\documentclass[amscd, leqno, amstex]{amsart}
\usepackage{amscd}
\usepackage{amsmath}
\usepackage{amssymb}
\theoremstyle{plain}
\newtheorem{theorem}{Theorem}[section]
\newtheorem{proposition}[theorem]{Proposition}
\newtheorem{lemma}[theorem]{Lemma}

\newtheorem{corollary}[theorem]{Corollary}
\theoremstyle{definition}
\newtheorem{definition}[theorem]{Definition}
\newtheorem{conjecture}[theorem]{Conjecture}
\theoremstyle{remark}
\newtheorem{example}[theorem]{Example}
\newtheorem{remark}[theorem]{Remark}
\newtheorem{question}[theorem]{Question}
\newenvironment{pf}{\begin{proof}}{\end{proof}}
\makeatletter

\@addtoreset{equation}{section}
\makeatother

\begin{document}

\title
{Analytic torsion for Borcea-Voisin threefolds}
\author{Ken-Ichi Yoshikawa}
\address{Department of Mathematics, 
Faculty of Science,
Kyoto University,
Kyoto 606-8502, JAPAN}
\email{yosikawa@math.kyoto-u.ac.jp}

\dedicatory
{Dedicated to Professor Jean-Michel Bismut}

\begin{abstract}
In their study of genus-one string amplitude $F_{1}$, Bershadsky-Cecotti-Ooguri-Vafa discovered a remarkable identification between
holomorphic Ray-Singer torsion and instanton numbers for Calabi-Yau threefolds.
The holomorphic torsion invariant of Calabi-Yau threefolds corresponding to $F_{1}$ is called BCOV invariant.
In this paper, we establish an identification between the BCOV invariants of Borcea-Voisin threefolds
and another holomorphic torsion invariants for $K3$ surfaces with involution.
We also introduce BCOV invariants for abelian Calabi-Yau orbifolds.
Between Borcea-Voisin orbifold and its crepant resolution, we compare their BCOV invariants.
\end{abstract}

\maketitle

\tableofcontents

\section*{Introduction}\label{sect:0}
\par
In \cite{BCOV93}, \cite{BCOV94}, Bershadsky-Cecotti-Ooguri-Vafa discovered a remarkable identification between holomorphic Ray-Singer torsion and
instanton numbers for Calabi-Yau threefolds.
For a Calabi-Yau threefold $X$, they introduced the following combination of Ray-Singer analytic torsions \cite{RaySinger73}, \cite{BGS88}
$$
T_{\rm BCOV}(X,\gamma)=\exp\{-\sum_{p,q}(-1)^{p+q}pq\,\zeta'_{p,q}(0)\},
$$
where $\zeta_{p,q}(s)$ is the spectral zeta function of the Laplacian acting on the $(p,q)$-forms on $X$ with respect to a Ricci-flat K\"ahler metric $\gamma$.
Although $T_{\rm BCOV}(X,\gamma)$ itself is {\em not} an invariant of $X$, its correction \cite{FLY08} 
$$
\tau_{\rm BCOV}(X)
=
{\rm Vol}(X,\gamma)^{-3+\frac{\chi(X)}{12}}
{\rm Vol}_{L^{2}}(H^{2}(X,{\bf Z}),[\gamma])^{-1}
T_{\rm BCOV}(X,\gamma)
$$ 
is an invariant of $X$, where ${\rm Vol}_{L^{2}}(H^{2}(X,{\bf Z}),[\gamma])$ is the covolume of the lattice $H^{2}(X,{\bf Z})\subset H^{2}(X,{\bf R})$ with respect to
the $L^{2}$-metric on the cohomology group induced by $\gamma$ and $\chi(X)$ is the topological Euler number of $X$. 
We call $\tau_{\rm BCOV}(X)$ the {\em BCOV invariant} of $X$. 
Because of its invariance property, $\tau_{\rm BCOV}$ gives rise to a function on the moduli space of Calabi-Yau threefolds.
The function $-\log\tau_{\rm BCOV}$ is identified with the physical quantity $F_{1}$.
Then the remarkable conjecture of Bershadsky-Cecotti-Ooguri-Vafa \cite{BCOV94} can be formulated as follows.
\par
Let $\varDelta^{*}\subset{\bf C}$ be the unit punctured disc. 
Let $\pi\colon{\mathcal X}\to(\varDelta^{*})^{n}$ be a large complex structure limit of Calabi-Yau threefolds \cite{Morrison93b} with fiber $X_{s}=\pi^{-1}(s)$,
$s\in(\varDelta^{*})^{n}$ and let $Y$ be the mirror Calabi-Yau threefold corresponding to the large complex structure limit point $0\in\varDelta^{n}$.
Define the function $F_{1}^{\rm top}(t)$ on the complexified K\"ahler cone $H^{2}(Y,{\bf R})+i\,{\mathcal K}_{Y}$ as the infinite product \cite{BCOV93}
$$
F_{1}^{\rm top}(t)
=
e^{2\pi i\int_{Y}\frac{1}{24}t\wedge c_{2}(Y)}
\prod_{d\in H_{2}(Y,{\bf Z})\setminus\{0\}}
\left(1-e^{2\pi i\langle t,d\rangle}\right)^{\frac{n_{0}(d)}{12}}
\prod_{\nu\geq1}\left(1-e^{2\pi i\nu\langle t,d\rangle}\right)^{n_{1}(d)}
$$
where ${\mathcal K}_{Y}$ is the K\"ahler cone of $Y$ and the numbers $\{n_{g}(d)\}$ are certain curve counting invariants of $Y$.
Then the conjecture of Bershadsky-Cecotti-Ooguri-Vafa claims that
$F_{1}^{\rm top}(t)$ converges when $\Im t\gg0$ and that the following equality holds near the large complex structure limit point $0\in\varDelta^{n}$
\begin{equation}
\label{eqn:BCOV:conj:2}
\tau_{\rm BCOV}(X_{s})
=
C\,
\left\|
F_{1}^{\rm top}\left(t(s)\right)^{2}
\left(
\frac{\varXi_{s}}{\langle A_{0}^{\lor},\varXi_{s}\rangle}
\right)^{3+n+\frac{\chi(X_{s})}{12}}
\otimes
\left(
\frac{\partial}{\partial t_{1}}\wedge\cdots\wedge\frac{\partial}{\partial t_{n}}
\right)_{t(s)}
\right\|^{2},
\end{equation}
where $C$ is a constant, $n=h^{2,1}(X_{s})=h^{1,1}(Y)$, 
$A_{0}^{\lor}\in\Gamma((\varDelta^{*})^{n},R^{3}\pi_{*}{\bf Z})$ is $\pi_{1}((\varDelta^{*})^{n})$-invariant,
$\{\varXi_{s}\}_{s\in(\varDelta^{*})^{n}}\in\Gamma((\varDelta^{*})^{n},\pi_{*}K_{{\mathcal X}/(\varDelta^{*})^{n}})$ is a nowhere vanishing 
relative canonical form of the family $\pi\colon{\mathcal X}\to(\varDelta^{*})^{n}$,
$t(s)=(t_{1}(s),\ldots,t_{n}(s))$ is the system of {\em canonical coordinates} \cite{Morrison93b} on $\varDelta^{n}$ and
$\|\cdot\|$ is the Hermitian metric induced from the $L^{2}$-metric on $\pi_{*}K_{{\mathcal X}/(\varDelta^{*})^{n}}$ and the Weil-Petersson metric on 
the holomorphic tangent bundle on $(\varDelta^{*})^{n}$.
After the works of Klemm-Mari\~no \cite{KlemmMarino08}, Maulik-Pandharipande \cite{MaulikPandharipande08}, Zinger \cite{Zinger09}, 
the curve counting invariants $\{n_{g}(d)\}$ appearing in $F_{1}^{\rm top}(t)$ are expected to be the Gopakumar-Vafa invariants of $Y$.
\par
To our knowledge, this conjecture of Bershadsky-Cecotti-Ooguri-Vafa is still widely open.
Although the curvature equation characterizing $\tau_{\rm BCOV}$ on the moduli space of Calabi-Yau threefolds is known 
\cite{BGS88}, \cite{BCOV94}, \cite{FangLu03}, \cite{FLY08}, because of the global nature of the differential equation and 
because the global structure of the moduli space of Calabi-Yau threefolds is not well understood in general, 
it is difficult to derive an explicit expression of the function $\tau_{\rm BCOV}$ in the canonical coordinates near the large complex structure limit point.
(See \cite{FLY08}, \cite{Yoshikawa09} for some Calabi-Yau threefolds whose BCOV invariants are explicitly expressed on the moduli space.)
\par
In the present paper, generalizing \cite{FLY08}, \cite{Yoshikawa09}, we shall give more examples of Calabi-Yau threefolds 
whose BCOV invariant can be explicitly expressed on a certain {\em locus} of their moduli space. 
The Calabi-Yau threefolds which we mainly treat in the present paper are those studied by Borcea \cite{Borcea97} and Voisin \cite{Voisin93}.
In their study of mirror symmetry, Borcea and Voisin introduced a class of Calabi-Yau threefolds
which are defined as the natural resolution
$$
\widetilde{X}_{(S,\theta,T)}\to X_{(S,\theta,T)}=(S\times T)/\theta\times(-1)_{T}.
$$
Here, $(S,\theta)$ is a $2$-elementary $K3$ surface, i.e., a $K3$ surface equipped with an anti-symplectic holomorphic involution,
and $T$ is an elliptic curve. We call $\widetilde{X}_{(S,\theta,T)}$ (resp $X_{(S,\theta,T)}$) the Borcea-Voisin threefold (resp. orbifold)
associated with $(S,\theta,T)$. 
By the computation of the Hodge numbers of $\widetilde{X}_{(S,\theta,T)}$ and hence the dimension of its Kuranishi space \cite{Voisin93}, 
a generic deformation of $\widetilde{X}_{(S,\theta,T)}$ is no longer of Borcea-Voisin type 
unless $S^{\theta}$, the fixed-point-set of the $\theta$-action on $S$, is either empty or rational.
Hence, in most cases, the Borcea-Voisin locus forms a proper subvariety of the deformation space of $\widetilde{X}_{(S,\theta,T)}$.
\par
By construction, once their deformation type is fixed, Borcea-Voisin threefolds are parametrized by a certain Zariski open subset of the locally symmetric variety 
associated with the product of the period domain for $2$-elementary $K3$ surfaces and that of elliptic curves. 
We regard this locally symmetric variety as the moduli space of Borcea-Voisin threefolds of fixed deformation type.
By Nikulin \cite{Nikulin80}, \cite{Nikulin83}, the deformation type of $(S,\theta)$ and hence that of $\widetilde{X}_{(S,\theta,T)}$ is determined by
the isometry class of the invariant sublattice of $H^{2}(S,{\bf Z})$ with respect to the $\theta$-action
$$
H^{2}(S,{\bf Z})^{\theta}=\{l\in H^{2}(S,{\bf Z});\,\theta^{*}l=l\}.
$$ 
For this reason, the isometry class of $H^{2}(S,{\bf Z})^{\theta}$ is called the type of $(S,\theta)$ and also the type of $\widetilde{X}_{(S,\theta,T)}$.
There exist $75$ distinct types \cite{Nikulin80}, \cite{Nikulin83}. We often identify a type with its representative, 
hence a primitive $2$-elementary Lorentzian sublattice of the $K3$-lattice.
\par
In their study of duality in string theory \cite{HarveyMoore98}, Harvey-Moore studied $F_{1}$ for a special class of Borcea-Voisin threefolds, 
called Enriques Calabi-Yau threefolds, and gave an identification between $F_{1}$ for Enriques Calabi-Yau threefolds 
and the equivariant determinant of Laplacian for $K3$ surfaces with fixed-point-free involution. (See also \cite{FLY08}, \cite{Yoshikawa08}.)
After Harvey-Moore, in \cite{Yoshikawa04}, another holomorphic torsion invariant was constructed for $2$-elementary $K3$ surfaces of type $M$  
(cf. Section~\ref{sect:2.2}):
$$
\tau_{M}(S,\theta)
=
{\rm Vol}(S,\kappa)^{\frac{14-r(M)}{4}}
\tau_{{\bf Z}_{2}}(S,\kappa)(\theta)\,
{\rm Vol}(S^{\theta},\kappa|_{S^{\theta}})
\tau(S^{\theta},\kappa|_{S^{\theta}}),
$$
where $\kappa$ is a $\theta$-invariant Ricci-flat K\"ahler form on $S$, $\tau_{{\bf Z}_{2}}(S,\kappa)(\theta)$ is the equivariant analytic torsion of
$(S,\theta)$ with respect to $\kappa$, $r(M)$ is the rank of $M$, 
and $\tau(S^{\theta},\kappa|_{S^{\theta}})$ is the analytic torsion of $(S^{\theta},\kappa|_{S^{\theta}})$.
By its invariance property, $\tau_{M}$ is regarded as a function on the moduli space of $2$-elementary $K3$ surfaces of type $M$.
In \cite{Yoshikawa04}, \cite{Yoshikawa12}, the automorphy of the function $\tau_{M}$ was established:
$$
\tau_{M}
=
\left\|\Phi_{M}\right\|^{-1/2\nu},
$$
where $\nu>0$ is a certain integer and $\|\Phi_{M}\|$ is the Petersson norm of an automorphic form $\Phi_{M}$ on the moduli space 
vanishing exactly on the discriminant locus with multiplicity $\nu$.
Recently, in \cite{Yoshikawa13}, \cite{MaYoshikawa14}, the structure of $\Phi_{M}$ as an automorphic form on the moduli space of $2$-elementary $K3$ surfaces 
of type $M$ is determined. 
Namely, $\Phi_{M}$ is always the tensor product of an explicit Borcherds product and a classical Siegel modular form. 
(See Section~\ref{sect:2.4}.)
\par
In this way, we have two invariants $\tau_{\rm BCOV}(\widetilde{X}_{(S,\theta,T)})$ and $\tau_{M}(S,\theta)$, 
which can be identified by Harvey-Moore \cite{HarveyMoore98} when $S/\theta$ is an Enriques surface.
Our first result is an extension of this identification of Harvey-Moore to arbitrary types of Borcea-Voisin threefolds 
and the corresponding $2$-elementary $K3$ surfaces.

\begin{theorem}
\label{MainThm:BCOV:Borcea:Voisin}
There exists a constant $C_{M}$ depending only on the isometry class $M$ such that for every Borcea-Voisin threefold $\widetilde{X}_{(S,\theta,T)}$ of type $M$,
$$
\tau_{\rm BCOV}(\widetilde{X}_{(S,\theta,T)})
=
C_{M}\,
\tau_{M}\left(S,\theta\right)^{-4}
\left\|
\eta(T)^{24}
\right\|^{2},
$$ 
where $\|\eta(T)\|$ is the value of the Petersson norm of the Dedekind $\eta$-function evaluated at the period of $T$. 
\end{theorem}

For an explicit formula for the BCOV invariant in terms of automorphic forms, see Corollary~\ref{cor:BCOV:invariant:Borcea:Voisin:threefolds}.
We remark that Theorem~\ref{MainThm:BCOV:Borcea:Voisin} was conjectured in \cite[Conj.\,5.17]{Yoshikawa09}.
After the conjecture of Bershadsky-Cecotti-Ooguri-Vafa \eqref{eqn:BCOV:conj:2},
Theorem~\ref{MainThm:BCOV:Borcea:Voisin} clarifies the significance of the invariant $\tau_{M}$ in mirror symmetry.
\par
Since the BCOV torsion makes sense for Calabi-Yau orbifolds, it is natural to ask the possibility of extending the construction of BCOV invariants 
to Calabi-Yau orbifolds.
In this paper, we give a partial answer to this question by constructing BCOV invariants for abelian Calabi-Yau orbifolds.
\par
Let $Y$ be an abelian Calabi-Yau orbifold and let $\Sigma Y$ be its inertia orbifold (cf. \cite{Kawasaki78}, \cite{Ma05}). 
We define the BCOV invariant of $Y$ as 
$$
\begin{aligned}
\tau_{\rm BCOV}^{\rm orb}(Y)
&=
T_{\rm BCOV}(Y,\gamma)\,
{\rm Vol}(Y,\gamma)^{-3+\frac{{\chi}^{\rm orb}(Y)}{12}}
{\rm Vol}_{L^{2}}(H^{2}(Y,{\bf Z}),[\gamma])^{-1}
\\
&\quad\times
\tau(\Sigma Y,\gamma|_{\Sigma Y})^{-1}
{\rm Vol}(\Sigma Y,\gamma|_{\Sigma Y})^{-1},
\end{aligned}
$$
where $\gamma$ is a Ricci-flat K\"ahler metric on $Y$ in the sense of orbifolds, ${\chi}^{\rm orb}(Y)$ is a rational number defined as
the integral of certain characteristic forms on $Y$ and $\Sigma Y$ (cf. \eqref{eqn:Euler:X:G}), and
$\tau(\Sigma Y,\gamma|_{\Sigma Y})$ is the analytic torsion of $(\Sigma Y,\gamma|_{\Sigma Y})$.
(See Section~\ref{sect:7}.)
When $Y$ is a global Calabi-Yau orbifold $Y=X/G$, then ${\chi}^{\rm orb}(Y)=\chi(X,G)$ is the string-theoretic orbifold Euler characteristic of $Y$ 
and is equal to the Euler characteristic of an arbitrary crepant resolution of $Y$ (cf. \cite{Roan96}).
In Sections~\ref{sect:6} and \ref{sect:7}, we shall prove that $\tau_{\rm BCOV}^{\rm orb}(Y)$ is independent of the choice of $\gamma$. 
Hence $\tau_{\rm BCOV}^{\rm orb}(Y)$ is an invariant of $Y$.
Then a natural question is a comparison of the two BCOV invariants $\tau_{\rm BCOV}(\widetilde{Y})$ and $\tau_{\rm BCOV}^{\rm orb}(Y)$,
where $\widetilde{Y}$ is a crepant resolution of $Y$.
In Section~\ref{sect:8}, following Harvey-Moore \cite{HarveyMoore98}, we shall prove the following.

\begin{theorem}
\label{thm:factorization:orb:BCOV}
There exists a constant $C'_{M}$ depending only on $M$ such that for every Borcea-Voisin orbifold $X_{(S,\theta,T)}$ of type $M$,
$$
\tau_{\rm BCOV}^{\rm orb}\left(X_{(S,\theta,T)}\right)
=
C'_{M}\,
\tau_{M}\left(S,\theta\right)^{-4}
\|\eta(T)^{24}\|^{2}.
$$
\end{theorem}

\begin{corollary}
\label{cor:BCOV=orb:BCOV}
There exists a constant $C''_{M}$ depending only on $M$ such that for every Borcea-Voisin orbifold $X_{(S,\theta,T)}$ of type $M$,
$$
\tau_{\rm BCOV}(\widetilde{X}_{(S,\theta,T)})=C''_{M}\,\tau_{\rm BCOV}^{\rm orb}\left(X_{(S,\theta,T)}\right).
$$
\end{corollary}

Our proof of Corollary~\ref{cor:BCOV=orb:BCOV} is indirect in the sense that it is a consequence of 
Theorems~\ref{MainThm:BCOV:Borcea:Voisin} and \ref{thm:factorization:orb:BCOV}. 
A direct proof of Corollary~\ref{cor:BCOV=orb:BCOV} along the line of \cite[Question 5.18]{Yoshikawa09} is strongly desired. 
\par
To prove Theorem~\ref{MainThm:BCOV:Borcea:Voisin}, we compare the complex Hessians of $\log\tau_{\rm BCOV}$ and $\log(\tau_{M}^{-4}\|\eta^{24}\|^{2})$
to conclude that $\log[\tau_{\rm BCOV}/(\tau_{M}^{-4}\|\eta^{24}\|^{2})]$ is a pluriharmonic function on the moduli space of Borcea-Voisin threefolds of type $M$.
Then we must prove that $\log\tau_{\rm BCOV}$ and $\log(\tau_{M}^{-4}\|\eta^{24}\|^{2})$ have the same singularity
on the discriminant locus of the moduli space. By the locality of the singularity of $\log\tau_{\rm BCOV}$ established in \cite{Yoshikawa15} and 
the corresponding result for $\log\tau_{M}$ in \cite{Yoshikawa04}, we can verify that $\tau_{\rm BCOV}$ and $\tau_{M}^{-4}\|\eta^{24}\|^{2}$ have the same singularity
on the discriminant locus by computing their singularities for some particular examples of Borcea-Voisin threefolds.
\par
The construction of BCOV invariants for abelian Calabi-Yau orbifolds is an application of the theory of
Quillen metrics \cite{BGS88}, \cite{Bismut95}, \cite{Ma00}, \cite{Ma05}.
Proof of Theorem~\ref{thm:factorization:orb:BCOV} is parallel to the one in \cite{HarveyMoore98}.
\par
This paper is organized as follows.
In Section~\ref{sect:1}, we recall BCOV invariants and their basic properties.
In Section~\ref{sect:2}, we recall the holomorphic torsion invariant $\tau_{M}$ for $2$-elementary $K3$ surfaces and its explicit formula.
In Section~\ref{sect:3}, we prove Theorem~\ref{MainThm:BCOV:Borcea:Voisin}.
In Section~\ref{sect:4}, we study a certain orbifold characteristic form associated to orbifold submersions.
In Section~\ref{sect:5}, we study the variation of equivariant BCOV torsion.
In Section~\ref{sect:6}, we extend BCOV invariants to global abelian Calabi-Yau orbifolds and establish its curvature formula.
In Section~\ref{sect:7}, we extend BCOV invariants to general abelian Calabi-Yau orbifolds.
In Section~\ref{sect:8}, we prove Theorem~\ref{thm:factorization:orb:BCOV}.
\par
{\bf Acknowledgements }
The author thanks Professor Akira Fujiki for helpful discussions about Kuranishi family.
He also thanks the referee for the simplification of the proofs of Lemmas~\ref{lemma:equivariant:Todd:chern:character} and
\ref{lemma:equivariant:Todd:chern:character:relative:dim=0} and for very helpful comments.

The author is partially supported by JSPS KAKENHI Grant Numbers 23340017, 22244003, 22224001, 25220701.

\section
{BCOV invariants for Calabi-Yau threefolds}
\label{sect:1}
\par
A compact connected K\"ahler manifold $X$ of dimension $3$ is called a {\em Calabi-Yau threefold}
if the following conditions are satisfied
$$
(1)\quad
K_{X}\cong{\mathcal O}_{X},
\qquad\qquad
(2)\quad
H^{1}(X,{\mathcal O}_{X})=H^{2}(X,{\mathcal O}_{X})=0.
$$
If a compact K\"ahler {\em orbifold} of dimension $3$ satisfies (1), (2), then $X$ is called a Calabi-Yau orbifold.
It is classical that every Calabi-Yau threefold is projective algebraic
and that every K\"ahler class of a Calabi-Yau threefold contains a unique Ricci-flat K\"ahler form \cite{Yau78}.

\subsection
{Analytic torsion}
\label{sect:1.1}
\par
Let $V$ be a compact K\"ahler manifold of dimension $n$ with K\"ahler metric 
$g=\sum_{\alpha,\beta}g_{\alpha\bar{\beta}}dz_{\alpha}\otimes d\bar{z}_{\beta}$.
Then the K\"ahler form of $g$ is defined as $\gamma:=\frac{i}{2}\sum_{\alpha,\beta}g_{\alpha\bar{\beta}}dz_{\alpha}\wedge d\bar{z}_{\beta}$.
Let $(\xi,h)$ be a holomorphic Hermitian vector bundle on $V$.
Let $\square_{q}=(\bar{\partial}+\bar{\partial}^{*})^{2}$ be the Laplacian acting on the $(0,q)$-forms on $V$ with values in $\xi$.
Let $\sigma(\square_{q})$ be the spectrum of $\square_{q}$ and let $E(\square_{q};\lambda)$ be the eigenspace of $\square_{q}$
corresponding to the eigenvalue $\lambda\in\sigma(\square_{q})$.
Then the spectral zeta function of $\square_{q}$ is defined as
$$
\zeta_{q}(s)
:=
\sum_{\lambda\in\sigma(\square_{q})\setminus\{0\}}\lambda^{-s}\dim E(\square_{q};\lambda).
$$
It is classical that $\zeta_{q}(s)$ converges absolutely on the half-plane $\Re s>\dim V$, that $\zeta_{q}(s)$ extends to 
a meromorphic function on ${\bf C}$ and that $\zeta_{q}(s)$ is holomorphic at $s=0$.
Ray-Singer \cite{RaySinger73} introduced the notion of {\em analytic torsion}.

\begin{definition}
\label{def:analytic:torsion}
The analytic torsion of $(\xi,h)$ is defined as the real number
$$
\tau(V,\xi)
:=
\exp\{-\sum_{q\geq0}(-1)^{q}q\,\zeta'_{q}(0)\}.
$$
\end{definition}

In mirror symmetry, the following combination of analytic torsions introduced by Bershadsky-Cecotti-Ooguri-Vafa \cite{BCOV94} plays a crucial role.
Write $\square_{p,q}=(\bar{\partial}+\bar{\partial}^{*})^{2}$ for the Laplacian acting on the $(p,q)$-forms on $V$
and $\zeta_{p,q}(s)$ for its spectral zeta function.

\begin{definition}
\label{def:BCOV:torsion}
The {\em BCOV torsion} of $(V,\gamma)$ is the real number defined as
$$
T_{\rm BCOV}(V,\gamma)
:=
\exp\{-\sum_{p,q\geq0}(-1)^{p+q}pq\,\zeta'_{p,q}(0)\}
=
\prod_{p\geq0}\tau(X,\Omega^{p}_{X})^{(-1)^{p}p}.
$$
\end{definition}

The BCOV torsion itself is not an invariant of a Calabi-Yau threefold, even if the K\"ahler form $\gamma$ is Ricci-flat.
However, by adding a small correction term, it becomes an invariant of a Calabi-Yau threefold \cite{FLY08}. Let us recall its construction.

\subsection
{The BCOV invariant}
\label{sect:1.2}
\par
Let $X$ be a Calabi-Yau threefold and let $\gamma$ be a K\"ahler form on $X$, which is not necessarily Ricci-flat.
Let $\eta$ be a nowhere vanishing canonical form on $X$.
Let $c_{3}(X,\gamma)$ be top Chern form of $(TX,\gamma)$.
In this paper, we follow the convention in Arakeolv geometry. Hence
$$
{\rm Vol}(X,\gamma):=\|1\|_{L^{2},\gamma}^{2}=(2\pi)^{-3}\int_{X}\frac{\gamma^{3}}{3!},
\qquad
\|\eta\|_{L^{2}}^{2}=(2\pi)^{-3}\int_{X}i^{3^{2}}\eta\wedge\overline{\eta}.
$$
We define the covolume ${\rm Vol}_{L^{2}}(H^{2}(X,{\bf Z}),[\gamma])$ of the lattice $H^{2}(X,{\bf Z})\subset H^{2}(X,{\bf R})$ 
as the volume of the compact real torus $H^{2}(X,{\bf R})/H^{2}(X,{\bf Z})$ with respect to the $L^{2}$-metric on $H^{2}(X,{\bf R})$ induced from $\gamma$
$$
{\rm Vol}_{L^{2}}(H^{2}(X,{\bf Z}),[\gamma])
:=
{\rm Vol}(H^{2}(X,{\bf R})/H^{2}(X,{\bf Z}),[\gamma])
=
\det\left(
\langle{\bf e}_{\alpha},{\bf e}_{\beta}\rangle_{L^{2},\gamma}
\right),
$$
where $\{{\bf e}_{\alpha}\}_{1\leq\alpha\leq b_{2}(X)}$ is a basis of $H^{2}(X,{\bf Z})_{\rm free}$. 
Here $\langle{\bf e}_{\alpha},{\bf e}_{\beta}\rangle_{L^{2},\gamma}$ is defined as
$$
\langle{\bf e}_{\alpha},{\bf e}_{\beta}\rangle_{L^{2},\gamma}
:=
(2\pi)^{-3}\int_{X}{\mathcal H}{\bf e}_{\alpha}\wedge\overline{*}({\mathcal H}{\bf e}_{\beta}),
$$
where ${\mathcal H}{\bf e}_{\alpha}$ is the harmonic representative of ${\bf e}_{\alpha}$ with respect to $\gamma$ and 
$*$ denotes the Hodge star-operator with respect to $\gamma$.
Notice that the $L^{2}$-metric on $H^{2}(X,{\bf R})$ induced from $\gamma$ depends only on the cohomology class $[\gamma]\in H^{2}(X,{\bf R})$.

\begin{definition}
\label{def:BCOV:invariant}
The BCOV invariant of $X$ is the real number defined by
$$
\begin{aligned}
\tau_{\rm BCOV}(X):
&=
{\rm Vol}(X,\gamma)^{-3+\frac{\chi(X)}{12}}
{\rm Vol}_{L^{2}}(H^{2}(X,{\bf Z}),[\gamma])^{-1}
T_{\rm BCOV}(X,\gamma)
\\
&\qquad
\times
\exp\left[-\frac{1}{12}\int_{X}
\log\left(
\frac{i\,\eta\wedge\bar{\eta}}{\gamma^{3}/3!}
\frac{{\rm Vol}(X,\gamma)}{\|\eta\|_{L^{2}}^{2}}
\right)\,c_{3}(X,\gamma)\right],
\end{aligned}
$$
where $\chi(X)=\int_{X}c_{3}(X)$ is the topological Euler number of $X$.
\end{definition}

\begin{remark}
\label{remark:BCOV:invariant:Ricci:flat:case}
The K\"ahler form $\gamma$ on $X$ is Ricci-flat if and only if
$$
\frac{i\,\eta\wedge\bar{\eta}}{\gamma^{3}/3!}
=
\frac{\|\eta\|_{L^{2}}^{2}}{{\rm Vol}(X,\gamma)}.
$$
When $\gamma$ is Ricci-flat, we get a simpler expression
$$
\tau_{\rm BCOV}(X)
=
{\rm Vol}(X,\gamma)^{-3+\frac{\chi(X)}{12}}
{\rm Vol}_{L^{2}}(H^{2}(X,{\bf Z}),[\gamma])^{-1}
T_{\rm BCOV}(X,\gamma).
$$
We also remark that ${\rm Vol}(X,\gamma)$ and ${\rm Vol}_{L^{2}}(H^{2}(X,{\bf Z}),[\gamma])$ are constant
under polarized deformations of Calabi-Yau threefolds \cite[Lemma 4.12]{FLY08}. 
Hence, for each moduli space of {\em polarized Ricci-flat} Calabi-Yau threefolds ${\mathcal M}$
(cf. \cite{Fujiki84}, \cite{Schumacher84}, \cite{Viehweg95}), there is a constant $C_{\mathcal M}>0$ depending only on the polarization such that
$$
\tau_{\rm BCOV}(X)=C_{\mathcal M}\,T_{\rm BCOV}(X,\gamma)
$$
for all $(X,\gamma_{L},c_{1}(L))\in{\mathcal M}$. 
Here $L$ is an ample line bundle on $X$ giving the polarization and $\gamma_{L}$ is the unique Ricci-flat K\"ahler form on $X$ with $[\gamma_{L}]=c_{1}(L)$.
In this sense, the BCOV invariant can be identified with the BCOV torsion for polarized Ricci-flat Calabi-Yau threefolds.
\end{remark}

\begin{theorem}
\label{Thm:BCOV:invariant}
For a Calabi-Yau threefold $X$, $\tau_{\rm BCOV}(X)$ is independent of the choice of a K\"ahler form on $X$.
Namely, $\tau_{\rm BCOV}(X)$ is an invariant of $X$. 
\end{theorem}

\begin{pf}
See \cite[Th.\,4.16]{FLY08}.
\end{pf}

For an expression of $\tau_{\rm BCOV}(X)$ in terms of arithmetic characteristic classes, we refer to \cite[Eq.(4)]{MaillotRoessler12}.
After Bershadsky-Cecotti-Ooguri-Vafa and Theorem~\ref{Thm:BCOV:invariant},
$\tau_{\rm BCOV}$ is regarded as a function on the moduli space of Calabi-Yau threefolds.

\subsection
{Singularity of BCOV invariants}
\label{sect:1.3}
\par
In this subsection, we recall some results in \cite{FLY08}, \cite{Yoshikawa15} about the boundary behavior of BCOV invariants,
which are applications of the Bismut-Lebeau embedding theorem for Quillen metrics \cite{BismutLebeau91}.
\par
Let $\pi\colon{\mathcal X}\to B$ be a surjective morphism from an irreducible projective fourfold ${\mathcal X}$ to a compact Riemann surface $B$.
Assume that there exists a finite subset $\Delta_{\pi}\subset B$ such that 
$\pi|_{B\setminus\Delta_{\pi}}\colon{\mathcal X}|_{B\setminus\Delta_{\pi}}\to B\setminus\Delta_{\pi}$ is a smooth morphism 
and such that $X_{t}:=\pi^{-1}(t)$ is a Calabi-Yau threefold for all $t\in B\setminus\Delta_{\pi}$.

\begin{theorem}
\label{thm:log:divergence:BCOV}
For every $0\in\Delta_{\pi}$, there exists $\alpha=\alpha_{0}\in{\bf Q}$ such that
$$
\log\tau_{\rm BCOV}(X_{t})=\alpha\,\log|t|^{2}+O\left(\log(-\log|t|)\right)
\qquad
(t\to0),
$$
where $t$ is a local parameter of $B$ centered at $0\in\Delta_{\pi}$. 
\end{theorem}

\begin{pf}
See \cite[Ths.\,0.1 and 0.2]{Yoshikawa15}.
\end{pf}

\par
Next we recall the {\em locality} of the logarithmic singularity of BCOV invariants.
Let ${\mathcal X}$ and ${\mathcal X}'$ be normal irreducible projective fourfolds. 
Let $B$ and $B'$ be compact Riemann surfaces.
Let $\pi\colon{\mathcal X}\to B$ and $\pi'\colon{\mathcal X}'\to B'$ be surjective holomorphic maps.
Let $\overline{\varSigma}_{\pi|_{{\mathcal X}\setminus{\rm Sing}\,{\mathcal X}}}$ 
(resp. $\overline{\varSigma}_{\pi'|_{{\mathcal X}'\setminus{\rm Sing}\,{\mathcal X}'}}$)
be the closure of the critical locus of $\pi|_{{\mathcal X}\setminus{\rm Sing}\,{\mathcal X}}$
(resp. $\pi'|_{{\mathcal X}'\setminus{\rm Sing}\,{\mathcal X}'}$) in ${\mathcal X}$ (resp. ${\mathcal X}'$).
Define the critical loci of $\pi$ and $\pi'$ as
$$
\varSigma_{\pi}:={\rm Sing}\,{\mathcal X}\cup\overline{\varSigma}_{\pi|_{{\mathcal X}\setminus{\rm Sing}\,{\mathcal X}}},
\qquad
\varSigma_{\pi'}:={\rm Sing}\,{\mathcal X}'\cup\overline{\varSigma}_{\pi'|_{{\mathcal X}'\setminus{\rm Sing}\,{\mathcal X}'}}
$$
and the discriminant loci of $\pi$ and $\pi'$ as $\Delta_{\pi}:=\pi(\varSigma_{\pi})$ and $\Delta_{\pi'}:=\pi'(\varSigma_{\pi'})$, respectively.
Let $0\in\Delta_{\pi}$ and $0'\in\Delta_{\pi'}$. Let $\varDelta$ be the unit disc of ${\bf C}$.
Let $V$ (resp. $V'$) be a neighborhood of $0$ (resp. $0'$) in $B$ (resp. $B'$) such that
$V\cong\varDelta$ and $V\cap\Delta_{\pi}=\{0\}$ (resp. $V'\cong\varDelta$ and $V'\cap\Delta_{\pi'}=\{0\}$). 
We assume the following: 
\begin{itemize}
\item[(A1)]
$\Delta_{\pi}\not=B$, $\Delta_{\pi'}\not=B'$, $\dim\varSigma_{\pi}\leq2$, $\dim\varSigma_{\pi'}\leq2$, and $X_{0}$ and $X'_{0'}$ are irreducible.
\item[(A2)]
$X_{t}$ and $X'_{t'}$ are Calabi-Yau threefolds for all $t\in B\setminus\Delta_{\pi}$ and $t'\in B'\setminus\Delta_{\pi'}$.
\item[(A3)]
$\pi^{-1}(V)\setminus\varSigma_{\pi}$ carries a nowhere vanishing canonical form $\varXi$.
Similarly, $(\pi')^{-1}(V')\setminus\varSigma_{\pi'}$ carries a nowhere vanishing canonical form $\varXi'$.
\item[(A4)]
The function germ of $\pi$ near $\varSigma_{\pi}\cap\pi^{-1}(V)$ and the function germ of $\pi'$ near $\varSigma_{\pi'}\cap(\pi')^{-1}(V')$ are isomorphic.
Namely, there exist a neighborhood $O$ of $\varSigma_{\pi}\cap\pi^{-1}(V)$ in $\pi^{-1}(V)$,
a neighborhood $O'$ of $\varSigma_{\pi'}\cap(\pi')^{-1}(V')$ in $(\pi')^{-1}(V')$, 
and an isomorphism $\varphi\colon O\to O'$ such that $\pi|_{O}=\pi'\circ\varphi|_{O'}$.
\end{itemize}
\par
For $b\in B$ and $b'\in B'$, we set $X_{b}:=\pi^{-1}(b)$ and $X'_{b'}:=(\pi')^{-1}(b')$.
For $b\in B\setminus\Delta_{\pi}$ and $b'\in B'\setminus\Delta_{\pi'}$, 
the BCOV invariants $\tau_{\rm BCOV}(X_{b})$ and $\tau_{\rm BCOV}(X'_{b'})$ are well defined.
Let $0\in\Delta_{\pi}$ and $0'\in\Delta_{\pi'}$. A local parameter of $B$ (resp. $B'$) centered at $0$ (resp. $0'$) is denoted by $t$.
Hence $t$ is a generator of the maximal ideal of ${\mathcal O}_{B,0}$ and that of ${\mathcal O}_{B',0'}$.

\begin{theorem}
\label{conjecture:locality:BCOV:invariant}
If the assumptions {\rm (A1)--(A4)} are satisfied, 
then $\log\tau_{\rm BCOV}(X_{t})$ and $\log\tau_{\rm BCOV}(X'_{t})$ have the same logarithmic singularities at $t=0$:
$$
\lim_{t\to0}\frac{\log\tau_{\rm BCOV}(X_{t})}{\log|t|}
=
\lim_{t\to0}\frac{\log\tau_{\rm BCOV}(X'_{t})}{\log|t|}.
$$
In particular,
$$
\log\tau_{\rm BCOV}(X_{t})-\log\tau_{\rm BCOV}(X'_{t})=O\left(\log(-\log|t|)\right)
\qquad
(t\to0).
$$
\end{theorem}

\begin{pf}
See \cite[Th.\,4.1]{Yoshikawa15}.
\end{pf}

\subsection
{Algebraic section on the moduli space corresponding to $\tau_{\rm BCOV}$}
\label{sect:1.5}
\par
In this subsection, we study the equation of currents satisfied by $\tau_{\rm BCOV}$ on a compactified moduli space of polarized Calabi-Yau threefolds.
\par
Let ${\frak M}$ be a coarse moduli space of polarized Calabi-Yau threefolds \cite{Fujiki84, Schumacher84, Viehweg95}:
\begin{itemize}
\item[(i)]
Every point of ${\frak M}$ corresponds to an isomorphism class of a polarized Calabi-Yau threefold.
\item[(ii)]
For any polarized family of Calabi-Yau threefolds $(\pi\colon{\mathcal X}\to B,{\mathcal L})$, where ${\mathcal L}$ is a relatively ample line bundle 
on ${\mathcal X}$, the classifying map $B\ni b\to[(X_{b},L_{b})]\in{\frak M}$ is holomorphic.
Here $X_{b}:=\pi^{-1}(b)$, $L_{b}:={\mathcal L}|_{X_{b}}$ and $[(X_{b},L_{b})]$ is the isomorphism class of $(X_{b},L_{b})$.
 \end{itemize}
Since every holomorphic line bundle $L$ on a Calabi-Yau threefold $X$ extends to a holomorphic line bundle on the Kuranishi family of $X$,
any $[(X,L)]\in{\frak M}$ has a neighborhood isomorphic to ${\rm Def}(X)/{\rm Aut}(X,L)$, where ${\rm Def}(X)$ is the Kuranishi space of $X$. 
Let ${\mathcal M}$ be a component of ${\frak M}$.
Since ${\rm Def}(X)$ is smooth \cite{Tian87, Todorov89} and since ${\rm Aut}(X,L)={\rm Aut}(X,\gamma_{L})$ is a compact Lie group with
trivial Lie algebra and hence is a finite group, ${\mathcal M}$ is an orbifold.
Set ${\mathcal M}_{\rm reg}:={\mathcal M}\setminus{\rm Sing}\,{\mathcal M}$. Then ${\mathcal M}_{\rm reg}$ is a complex manifold.
Since ${\mathcal M}$ is quasi-projective by Viehweg \cite[Th.\,1.13]{Viehweg95}, there exists a projective manifold $\overline{\mathcal M}$ 
containing ${\mathcal M}_{\rm reg}$ as a dense Zariski open subset such that
${\mathcal D}:=\overline{\mathcal M}\setminus{\mathcal M}_{\rm reg}$ is a normal crossing divisor of $\overline{\mathcal M}$.
For $[(X,L)]\in{\mathcal M}$, we set $h^{1,2}:=h^{1,2}(X)=\dim{\mathcal M}$ and $\chi:=\chi(X)$, where $\chi(X)$ is the topological Euler number of $X$.
Then $h^{1,2}$ and $\chi$ are constant on ${\mathcal M}$.

\subsubsection
{Hodge bundles and the Kodaira-Spencer maps}
\label{sect:1.5.1}
\par
For $[(X,L)]\in{\mathcal M}$, let 
$$
{\frak f}\colon({\frak X},X)\to({\rm Def}(X),[X])
$$ 
be the Kuranishi family of $X$.
Let ${\mathcal H}^{3}$ be the holomorphic vector bundle on ${\mathcal M}_{\rm reg}$ such that ${\mathcal H}^{3}_{[(X,L)]}=H^{3}(X,{\bf C})$
and let 
$$
0={\mathcal F}^{4}\subset{\mathcal F}^{3}\subset{\mathcal F}^{2}\subset{\mathcal F}^{1}\subset{\mathcal F}^{0}={\mathcal H}^{3}
$$ 
be the Hodge filtration. We have
$$
{\mathcal H}^{3}|_{{\rm Def}(X)}=R^{3}{\frak f}_{*}{\bf C}\otimes{\mathcal O}_{{\rm Def}(X)},
\qquad
{\mathcal F}^{p}/{\mathcal F}^{p+1}|_{{\rm Def}(X)}\cong R^{3-p}{\frak f}_{*}\Omega^{p}_{{\frak X}/{\rm Def}(X)},
$$
where $\Omega^{p}_{{\frak X}/{\rm Def}(X)}:=\Lambda^{p}\Omega^{1}_{{\frak X}/{\rm Def}(X)}$ and 
$\Omega^{1}_{{\frak X}/{\rm Def}(X)}:=\Omega^{1}_{\frak X}/{\frak f}^{*}\Omega^{1}_{{\rm Def}(X)}$.
The line bundle $\lambda:={\mathcal F}^{3}$ on ${\mathcal M}_{\rm reg}$ is called the Hodge bundle.
Then $\lambda|_{{\rm Def}(X)}={\frak f}_{*}K_{{\frak X}/{\rm Def}(X)}$.
\par
By Kawamata \cite[Th.\,17]{Kawamata81}, there exists a finite covering 
$$
\phi\colon\widetilde{\mathcal M}\to\overline{\mathcal M}
$$ 
with branch locus ${\mathcal D}$ such that the monodromy operator on $\phi^{*}{\mathcal H}^{3}$ along $\phi^{-1}({\mathcal D})$ is unipotent. 
Under this assumption, by the nilpotent orbit theorem of Schmid \cite[Th.\,4.12]{Schmid73}, 
the vector bundles $\phi^{*}{\mathcal H}^{3}$ and $\phi^{*}{\mathcal F}^{p}$ extend to 
holomorphic vector bundles $\widetilde{\mathcal H}^{3}$ and $\widetilde{\mathcal F}^{p}$ on $\widetilde{\mathcal M}$, respectively.
We define line bundles $\widetilde{\lambda}$ and $\widetilde{\mu}$ on $\widetilde{\mathcal M}$ by
$$
\widetilde{\lambda}:=\widetilde{\mathcal F}^{3},
\qquad
\widetilde{\mu}:=\det(\widetilde{\mathcal F}^{2}/\widetilde{\mathcal F}^{3}).
$$
Then $\widetilde{\lambda}|_{\phi^{-1}({\mathcal M}_{\rm reg})}=\phi^{*}\lambda$.
\par
Define the holomorphic vector bundle $N$ on ${\frak X}$ by the exact sequence:
$$
0\longrightarrow\Theta_{{\frak X}/{\rm Def}(X)}:=\ker{\frak f}_{*}\longrightarrow\Theta_{\frak X}\longrightarrow N\longrightarrow 0.
$$
The projection ${\frak f}_{*}\colon N\to{\frak f}^{*}\Theta_{{\rm Def}(X)}$ induces an isomorphism $\iota\colon{\frak f}_{*}N\cong\Theta_{{\rm Def}(X)}$. 
Let $\delta\colon{\frak f}_{*}N\to R^{1}{\frak f}_{*}\Theta_{{\frak X}/{\rm Def}(X)}$ be the connecting homomorphism.
The Kodaira-Spencer map is defined as the composite
\begin{equation}
\label{eqn:Kodaira:Spencer:map}
\begin{aligned}
\rho_{{\rm Def}(X)}
:=
\delta\circ\iota^{-1}\colon \Theta_{{\rm Def}(X)}\to R^{1}{\frak f}_{*}\Theta_{{\frak X}/{\rm Def}(X)}
&=
R^{1}{\frak f}_{*}\Omega_{{\frak X}/{\rm Def}(X)}^{2}\otimes({\frak f}_{*}K_{{\frak X}/{\rm Def}(X)})^{\lor}
\\
&\cong
({\mathcal F}^{2}/{\mathcal F}^{3})\otimes({\mathcal F}^{3})^{\lor}|_{{\rm Def}(X)}.
\end{aligned}
\end{equation}
Then the Kodaira-Spencer map induces an isomorphism of holomorphic vector bundles on ${\mathcal M}_{\rm reg}$:
$$
\rho\colon\Theta_{{\mathcal M}_{\rm reg}}
\to
({\mathcal F}^{2}/{\mathcal F}^{3})\otimes({\mathcal F}^{3})^{\lor}.
$$
The isomorphism $\rho$ and its lift
$$
\phi^{*}\rho\colon\Theta_{\widetilde{\mathcal M}}|_{\phi^{-1}({\mathcal M}_{\rm reg})}
\to
(\widetilde{\mathcal F}^{2}/\widetilde{\mathcal F}^{3})\otimes(\widetilde{\mathcal F}^{3})^{\lor}|_{\phi^{-1}({\mathcal M}_{\rm reg})}
$$
are again called the Kodaira-Spencer map. 
Since 
$$
\det(\phi^{*}\rho)
\in 
H^{0}\left(
\phi^{-1}({\mathcal M}_{\rm reg}),
\det((\widetilde{\mathcal F}^{2}/\widetilde{\mathcal F}^{3})\otimes(\widetilde{\mathcal F}^{3})^{\lor})
\otimes(\det\Theta_{\widetilde{\mathcal M}})^{\lor}
\right)
$$ 
has at most algebraic singularity along $\phi^{-1}({\mathcal D})$,  
$\det(\phi^{*}\rho)$ is a meromorphic section of the line bundle $\widetilde{\mu}\otimes\widetilde{\lambda}^{-h^{1,2}}\otimes K_{\widetilde{\mathcal M}}$.
Hence we have the canonical isomorphism
\begin{equation}
\label{eqn:canonical:isomorphism}
K_{\widetilde{\mathcal M}}^{-1}\otimes{\mathcal O}_{\widetilde{\mathcal M}}\left({\rm div}(\det(\phi^{*}\rho))\right)
\cong
\widetilde{\mu}\otimes\widetilde{\lambda}^{-h^{1,2}}.
\end{equation}

\subsubsection
{Weil-Petersson metric and its boundary behavior}
\label{sect:1.5.2}
\par
Since the third cohomology group of Calabi-Yau threefolds consists of primitive cohomology classes, 
${\mathcal F}^{p}/{\mathcal F}^{p+1}$ is equipped with the $L^{2}$-metric, which is independent of the choice of a K\"ahler metric on each fiber.
In particular, so is the Hodge bundle $\widetilde{\lambda}$. This metric is denoted by $h_{L^{2}}$ or $(\cdot,\cdot)_{L^{2}}$.
The $L^{2}$-metric $(\cdot,\cdot)_{L^{2}}$ on $({\mathcal F}^{p}/{\mathcal F}^{p+1})|_{[(X,L)]}=H^{q}(X,\Omega^{p}_{X})$, $p+q=3$, is expressed by
\begin{equation}
\label{eqn:L2:metric:middle:degree:primitive}
(u,v)_{L^{2}}:=(-1)^{3}i^{p-q}\int_{X}u\wedge\overline{v}.
\end{equation}
The Hermitian metric on the line bundle $\det({\mathcal F}^{p}/{\mathcal F}^{p+1})$ induced from $h_{L^{2}}$ is denoted by $h_{L^{2}}$ or $\|\cdot\|_{L^{2}}$.
The Weil-Petersson forms on  ${\mathcal M}_{\rm reg}$ and $\widetilde{\mathcal M}$ are defined as
$$
\omega_{\rm WP}=c_{1}(\lambda,h_{L^{2}}),
\qquad
\widetilde{\omega}_{\rm WP}=c_{1}(\widetilde{\lambda},h_{L^{2}}).
$$
Then $\phi^{*}\omega_{\rm WP}=\widetilde{\omega}_{\rm WP}|_{\phi^{-1}({\mathcal M}_{\rm reg})}$.
Let $\eta_{{\frak X}/{\rm Def}(X)}\in H^{0}({\rm Def}(X),{\frak f}_{*}K_{{\frak X}/{\rm Def}(X)})$ be a nowhere vanishing holomorphic section
and define the function $\|\eta_{{\frak X}/{\rm Def}(X)}\|_{L^{2}}^{2}$ on ${\rm Def}(X)$ by
$$
\|\eta_{{\frak X}/{\rm Def}(X)}\|_{L^{2}}^{2}([X_{t}]):=\|\eta_{{\frak X}/{\rm Def}(X)}|_{X_{t}}\|_{L^{2}}^{2}.
$$
Since $\omega_{\rm WP}=c_{1}(\lambda,h_{L^{2}})$, we have
\begin{equation}
\label{eqn:Weil:Petersson:form:2}
\omega_{\rm WP}
=
-dd^{c}\log\|\eta_{{\frak X}/{\rm Def}(X)}\|_{L^{2}}^{2}
=
c_{1}({\frak f}_{*}K_{{\frak X}/{\rm Def}(X)},h_{L^{2}}).
\end{equation}
By e.g. \cite[Th.\,2]{Tian87} and \cite[Def.\,4.3]{FLY08}, we get 
\begin{equation}
\label{eqn:Weil:Petersson:metric}
\omega_{\rm WP}(u,v)
=
\frac{(\rho_{{\rm Def}(X)}(u)\otimes\eta,\rho_{{\rm Def}(X)}(v)\otimes\eta)_{L^{2}}}{\|\eta\|_{L^{2}}^{2}}
\end{equation}
for all $u,v\in\Theta_{{\rm Def}(X),[X]}=H^{1}(X,\Theta_{X})$,
where $\eta\in H^{0}(X,K_{X})\setminus\{0\}$ and the numerator is the cup-product pairing between 
$H^{1}(X,\Omega_{X}^{2})$ and $\overline{H^{1}(X,\Omega_{X}^{2})}=H^{2}(X,\Omega_{X}^{1})$.
By the expression \eqref{eqn:Weil:Petersson:metric}, $\omega_{\rm WP}$ is a real analytic $(1,1)$-form on ${\rm Def}(X)$.
Let
$$
{\rm Ric}\,\omega_{\rm WP}:=c_{1}(\Theta_{{\rm Def}(X)},\omega_{\rm WP})
$$
be the Ricci-form of the Weil-Petersson form on $\Theta_{{\rm Def}(X)}$.
\par
By \eqref{eqn:Weil:Petersson:metric} and \cite[Prop.\,5.22, Proof of Cor.\,5.23]{CattaniKaplanSchmid86}, 
$\widetilde{\omega}_{\rm WP}$ has Poincar\'e growth near $\phi^{-1}({\mathcal D})$.
Namely, for any $x\in\phi^{-1}({\mathcal D})$, there is a neighborhood ${\mathcal U}$ of $x$ in $\widetilde{\mathcal M}$
such that ${\mathcal U}\cong\varDelta^{h^{1,2}}$, ${\mathcal U}\setminus\phi^{-1}({\mathcal D})\cong(\varDelta^{*})^{k}\times\varDelta^{h^{1,2}-k}$
for some $1\leq k\leq h^{1,2}$ and such that
$$
0\leq\omega_{\rm WP}|_{{\mathcal U}\setminus\phi^{-1}({\mathcal D})}\leq C\,\omega_{P}.
$$
Here $\omega_{P}$ is the K\"ahler form of the Poincar\'e metric on $(\varDelta^{*})^{k}\times\varDelta^{h^{1,2}-k}$, i.e.,
$$
\omega_{P}
=
i\sum_{\alpha=1}^{k}\frac{dt_{\alpha}\wedge d\bar{t}_{\alpha}}{|t_{\alpha}|^{2}(\log|t_{\alpha}|)^{2}}
+
i\sum_{\beta=k+1}^{h^{1,2}}dt_{\beta}\wedge d\bar{t}_{\beta}
$$
\par
Let ${\rm Ric}\,\widetilde{\omega}_{\rm WP}$ be the Ricci form of $\widetilde{\omega}_{\rm WP}$, namely the first Chern form of
the holomorphic tangent bundle equipped with the Weil-Petersson metric $(\Theta_{\phi^{-1}({\mathcal M}_{\rm reg})},\widetilde{\omega}_{\rm WP})$.
Then ${\rm Ric}\,\widetilde{\omega}_{\rm WP}=\phi^{*}{\rm Ric}\,\omega_{\rm WP}$.
By \eqref{eqn:Weil:Petersson:metric}, the Kodaira-Spencer map induces an isometry of holomorphic Hermitian vector bundles on
$\phi^{-1}({\mathcal M}_{\rm reg})$:
$$
\phi^{*}\rho\colon(\Theta_{\widetilde{\mathcal M}}|_{\phi^{-1}({\mathcal M}_{\rm reg})},\widetilde{\omega}_{\rm WP})
\cong
((\widetilde{\mathcal F}^{2}/\widetilde{\mathcal F}^{3})\otimes(\widetilde{\mathcal F}^{3})^{\lor}|_{\phi^{-1}({\mathcal M}_{\rm reg})},
(\cdot,\cdot)_{L^{2}}\otimes\|\cdot\|_{L^{2}}^{-1}).
$$
From this isometry, we have the following equality of $(1,1)$-forms on $\phi^{-1}({\mathcal M}_{\rm reg})$
$$
{\rm Ric}\,\widetilde{\omega}_{\rm WP}+h^{1,2}\,\widetilde{\omega}_{\rm WP}
=
c_{1}(\widetilde{\mathcal F}^{2}/\widetilde{\mathcal F}^{3},(\cdot,\cdot)_{L^{2}})
=
c_{1}(\widetilde{\mu},\|\cdot\|_{L^{2}}).
$$
By \cite[Th.\,1.1]{Lu01}, ${\rm Ric}\,\widetilde{\omega}_{\rm WP}+(h^{1,2}+3)\,\widetilde{\omega}_{\rm WP}$ is a positive $(1,1)$-form on
$\phi^{-1}({\mathcal M}_{\rm reg})$.
By \cite[Prop.\,5.22, Proof of Cor.\,5.23]{CattaniKaplanSchmid86} again and the above expression, 
the $(1,1)$-form ${\rm Ric}\,\widetilde{\omega}_{\rm WP}+(h^{1,2}+3)\,\widetilde{\omega}_{\rm WP}$
has Poincar\'e growth near the normal crossing divisor $\phi^{-1}({\mathcal D})$. 
Hence, on ${\mathcal U}$ as above, we have the following estimate:
$$
0
\leq
\{
{\rm Ric}\,\widetilde{\omega}_{\rm WP}+(h^{1,2}+3)\,\widetilde{\omega}_{\rm WP}
\}|_{{\mathcal U}\setminus\phi^{-1}({\mathcal D})}
\leq 
C\,\omega_{P}.
$$
\par
We identify the $(1,1)$-forms $\widetilde{\omega}_{\rm WP}$ and ${\rm Ric}\,\widetilde{\omega}_{\rm WP}+(h^{1,2}+3)\,\widetilde{\omega}_{\rm WP}$ 
on $\phi^{-1}({\mathcal M}_{\rm reg})$ with 
the closed positive $(1,1)$-currents on $\widetilde{\mathcal M}$ defined as 
their trivial extensions from $\phi^{-1}({\mathcal M}_{\rm reg})$ to $\widetilde{\mathcal M}$.

\subsubsection
{The differential equation satisfied by $\log\tau_{\rm BCOV}$ on $\widetilde{\mathcal M}$}
\label{sect:1.5.3}
\par
As usual, we define $d^{c}=\frac{1}{4\pi i}(\partial-\bar{\partial})$ for complex manifolds.
Hence $dd^{c}=\frac{1}{2\pi i}\bar{\partial}\partial$.
By the curvature formula for Quillen metrics \cite{BGS88}, the complex Hessian of $\log\tau_{\rm BOOV}$ can be computed
on the Kuranishi space of $X$:

\begin{theorem}
\label{thm:curvature:BCOV}
The following equality of $(1,1)$-forms on $({\rm Def}(X),[X])$ holds:
\begin{equation}
\label{eqn:curvature:BCOV}
-dd^{c}\log\tau_{\rm BCOV}={\rm Ric}\,\omega_{\rm WP}+\left(h^{1,2}+3+\frac{\chi}{12}\right)\omega_{\rm WP}.
\end{equation}
In particular, $\tau_{\rm BCOV}\in C^{\omega}({\rm Def}(X))$.
\end{theorem}

\begin{pf}
See \cite[Th.\,4.14]{FLY08}.
\end{pf}

This differential equation can be globalized to the compactified moduli space $\overline{\mathcal M}$.
For this, firstly, we prove the following:

\begin{theorem}
\label{thm:curvature:current:BCOV:invariant}
Let $\phi^{-1}({\mathcal D})=\bigcup_{k\in K}\widetilde{\mathcal D}_{k}$ be the irreducible decomposition.
Then there exists $\alpha_{k}\in{\bf Q}$ such that the following equality of $(1,1)$-currents on $\widetilde{\mathcal M}$ holds
\begin{equation}
\label{eqn:curvature:pullback:BCOV:inv}
-12dd^{c}(\phi^{*}\log\tau_{\rm BCOV})
=
(36+12h^{1,2}+\chi)\,\widetilde{\omega}_{\rm WP}
+
12\,{\rm Ric}\,\widetilde{\omega}_{\rm WP}
-
\sum_{k\in K}\alpha_{k}\delta_{\widetilde{\mathcal D}_{k}}.
\end{equation}
\end{theorem}

\begin{pf}
For $k\in K$, set $\widetilde{\mathcal D}_{k}^{0}:=\widetilde{\mathcal D}_{k}\setminus\bigcup_{l\not=k}\widetilde{\mathcal D}_{l}$.
For every $p\in\widetilde{\mathcal D}_{k}^{0}$, there is a neighborhood $U$ of $p$ in $\widetilde{\mathcal M}$ and an isomorphism
$U\cong\varDelta^{h^{1,2}}$ such that $U\cap\widetilde{\mathcal D}_{k}^{0}\cong\{z_{1}=0\}\cap\varDelta^{h^{1,2}}$ and
$U\cap\widetilde{\mathcal D}_{l}^{0}=\emptyset$ for $l\not=k$. We write $z'=(z_{2},\ldots,z_{h^{1,2}})$.
We define $Z:=\bigcup_{k\not=l}\widetilde{\mathcal D}_{k}\cap\widetilde{\mathcal D}_{l}$.
\par
Since $\widetilde{\omega}_{\rm WP}$ and ${\rm Ric}\,\widetilde{\omega}_{\rm WP}+(h^{1,2}+3)\widetilde{\omega}_{\rm WP}$
are closed positive $(1,1)$-currents on $\varDelta^{h^{1,2}}$, there exist by Siu \cite[Proof of Lemma 5.4]{Siu74}
plurisubharmonic functions $\widetilde{\psi}_{1}$, $\widetilde{\psi}_{2}$ on $\varDelta^{h^{1,2}}$ such that
\begin{equation}
\label{eqn:estimate:potential:3}
dd^{c}\widetilde{\psi}_{1}=\widetilde{\omega}_{\rm WP},
\qquad
dd^{c}\widetilde{\psi}_{2}={\rm Ric}\,\widetilde{\omega}_{\rm WP}+(h^{1,2}+3)\,\widetilde{\omega}_{\rm WP}
\end{equation}
as currents on $\varDelta^{h^{1,2}}$. 
Since $\widetilde{\omega}_{\rm WP}$ and ${\rm Ric}\,\widetilde{\omega}_{\rm WP}+(h^{1,2}+3)\widetilde{\omega}_{\rm WP}$ have Poincar\'e growth,
there exist constants $C_{3},C_{4}>0$ such that
\begin{equation}
\label{eqn:estimate:potential:4}
0\leq dd^{c}\widetilde{\psi}_{1}\leq C_{3}\,\omega_{P},
\qquad
0\leq dd^{c}\widetilde{\psi}_{2}\leq C_{4}\,\omega_{P}.
\end{equation}
Set 
$$
Q(z_{1},z'):=-\log(-\log|z_{1}|^{2})+\|z'\|^{2}.
$$
Since $\omega_{P}=dd^{c}Q$, we deduce from \eqref{eqn:estimate:potential:3}, \eqref{eqn:estimate:potential:4} that the functions
$$
\widetilde{\psi}_{1},
\qquad
C_{3}Q-\widetilde{\psi}_{1},
\qquad
\widetilde{\psi}_{2},
\qquad
C_{4}Q-\widetilde{\psi}_{2}
$$
are plurisubharmonic on $\varDelta^{h^{1,2}}$.
Since these functions are bounded from above on a neighborhood of $0\in\varDelta^{h^{1,2}}$,
there exist constants $C_{5},C_{6}>0$ with
\begin{equation}
\label{eqn:estimate:potential:5}
C_{5}\left\{-\log(-\log|z_{1}|)-1\right\}
\leq
\widetilde{\psi}_{m}(z_{1},z')
\leq 
C_{6}
\qquad
(m=1,2).
\end{equation}
\par
By Theorems~\ref{thm:curvature:BCOV}, \ref{thm:log:divergence:BCOV} and \eqref{eqn:estimate:potential:5}, 
there exists $\alpha(z')\in{\bf Q}$ such that on $\varDelta^{*}\times\varDelta^{h^{1,2}-1}$
$$
-dd^{c}\left(\phi^{*}\log\tau_{\rm BCOV}+\frac{\chi}{12}\,\widetilde{\psi}_{1}+\widetilde{\psi}_{2}\right)=0,
$$
$$
\left|
\left(\phi^{*}\log\tau_{\rm BCOV}+\frac{\chi}{12}\,\widetilde{\psi}_{1}+\widetilde{\psi}_{2}\right)(z_{1},z')-\alpha(z')\log|z_{1}|^{2}
\right|
\leq
C(z')\,\log(-\log|z_{1}|),
$$
where $C(z')$ is a (possibly discontinuous, unbounded) positive function on $\varDelta^{h^{1,2}-1}$.
By \cite[Lemma 5.9]{Yoshikawa09}, there exists $\alpha_{k}\in{\bf Q}$ and a pluriharmonic function $h\in C^{\omega}(\varDelta^{h^{1,2}})$ such that
\begin{equation}
\label{eqn:estimate:potential:6}
\phi^{*}\log\tau_{\rm BCOV}(z_{1},z')+\frac{\chi}{12}\,\widetilde{\psi}_{1}(z_{1},z')+\widetilde{\psi}_{2}(z_{1},z')
=
\alpha_{k}\log|z_{1}|^{2}+h(z_{1},z').
\end{equation}
By \eqref{eqn:estimate:potential:6}, we get the following equation of currents on $U$
$$
-12dd^{c}(\phi^{*}\log\tau_{\rm BCOV})
=
(36+12h^{1,2}+\chi)\,\widetilde{\omega}_{\rm WP}
+
12\,{\rm Ric}\,\widetilde{\omega}_{\rm WP}
-
\alpha_{k}\delta_{\widetilde{\mathcal D}_{k}}.
$$
Since $p$ is an arbitrary point of $\widetilde{\mathcal D}_{k}^{0}$, we get the equation of currents on $\widetilde{\mathcal M}\setminus Z$:
\begin{equation}
\label{eqn:curvature:current:BCOV:moduli}
-12dd^{c}(\phi^{*}\log\tau_{\rm BCOV})
=
(36+12h^{1,2}+\chi)\,\widetilde{\omega}_{\rm WP}
+
12\,{\rm Ric}\,\widetilde{\omega}_{\rm WP}
-
\sum_{k\in K}\alpha_{k}\delta_{\widetilde{\mathcal D}_{k}}.
\end{equation}
\par
Let $p\in Z$. Since the right hand side is a linear combination of closed positive $(1,1)$-currents on $\widetilde{\mathcal M}$,
there is a neighborhood $V$ of $p$ in $\widetilde{\mathcal M}$ such that the right hand side of \eqref{eqn:curvature:current:BCOV:moduli}
has a potential function $R$ on $V$.
Namely $R$ is a plurisubharmonic function on $V$ satisfying the following equation of currents on $U$
\begin{equation}
\label{eqn:potential}
dd^{c}R
=
(36+12h^{1,2}+\chi)\,\widetilde{\omega}_{\rm WP}
+
12\,{\rm Ric}\,\widetilde{\omega}_{\rm WP}
-
\sum_{k\in K}\alpha_{k}\delta_{\widetilde{\mathcal D}_{k}}.
\end{equation}
By \eqref{eqn:curvature:current:BCOV:moduli}, \eqref{eqn:potential}, 
$12\phi^{*}\log\tau_{\rm BCOV}+R$ is a pluriharmonic function on $V\setminus Z$. 
Since $Z$ has codimension $\geq2$ in $V$,
$12\phi^{*}\log\tau_{\rm BCOV}+R$ extends to a pluriharmonic function on $V$. 
Thus \eqref{eqn:curvature:current:BCOV:moduli} holds on $V$.
Since $p\in Z$ is arbitrary, \eqref{eqn:curvature:current:BCOV:moduli} holds on $\widetilde{\mathcal M}$.
\end{pf}

By the definition of trivial extension of currents, $\frac{1}{\deg\phi}\phi_{*}\widetilde{\omega}_{\rm WP}$ coincides with 
the trivial extension of $\omega_{\rm WP}$ from ${\mathcal M}_{\rm reg}$ to $\overline{\mathcal M}$.
This trivial extension is again denoted by $\omega_{\rm WP}$. 
Similarly, $\frac{1}{\deg\phi}\phi_{*}\{{\rm Ric}\,\widetilde{\omega}_{\rm WP}+(h^{1,2}+3)\widetilde{\omega}_{\rm WP}\}$ coincides with 
the trivial extension of ${\rm Ric}\,\omega_{\rm WP}+(h^{1,2}+3)\omega_{\rm WP}$ from ${\mathcal M}_{\rm reg}$ to $\overline{\mathcal M}$.
Set ${\mathcal D}_{k}:={\rm Supp}\,\phi(\widetilde{\mathcal D}_{k})$. 
Then $\phi_{k}:=\phi|_{\widetilde{\mathcal D}_{k}}$ is a surjective map from $\widetilde{\mathcal D}_{k}$ to ${\mathcal D}_{k}$ and 
we have the equation of currents $\phi_{*}\delta_{\widetilde{\mathcal D}_{k}}=\deg\phi_{k}\cdot\delta_{{\mathcal D}_{k}}$.
Applying $\phi_{*}$ to the both sides of \eqref{eqn:curvature:pullback:BCOV:inv}, we get

\begin{corollary}
\label{cor:curvature:current:BCOV:invariant}
The following equation of currents on $\overline{\mathcal M}$ holds:
$$
-12dd^{c}\log\tau_{\rm BCOV}
=
(36+12h^{1,2}+\chi)\,{\omega}_{\rm WP}
+
12\,{\rm Ric}\,{\omega}_{\rm WP}
-
\sum_{k\in K}\frac{\deg\phi_{k}}{\deg\phi}\alpha_{k}\,\delta_{{\mathcal D}_{k}}.
$$
\end{corollary}

\subsubsection
{The section on $\widetilde{\mathcal M}$ corresponding to $\tau_{\rm BCOV}$}
\label{sect:1.5.4}
\par
As a consequence of Corollary~\ref{cor:curvature:current:BCOV:invariant}, we have the following algebraicity of $\tau_{\rm BCOV}$.
We remark that this result (Theorem~\ref{thm:consequence:BCOV:conjecture}) is deduced directly 
from Maillot-R\"ossler's formula \cite[Eq.\,(4)]{MaillotRoessler12} for $\tau_{\rm BCOV}$ in terms of certain arithmetic characteristic classes.

\begin{theorem}
\label{thm:consequence:BCOV:conjecture}
There exist $\ell\in{\bf Z}_{>0}$, a character $\chi\in{\rm Hom}(\pi_{1}(\widetilde{\mathcal M}),U(1))$ and a meromorphic section $\sigma$
of the line bundle
$\widetilde{\lambda}^{(36+\chi)\ell}\otimes\widetilde{\mu}^{-12\ell}\otimes[\chi]$
such that the following equality of functions on $\phi^{-1}({\mathcal M}_{\rm reg})$ holds:
$$
\phi^{*}\tau_{\rm BCOV}^{12\ell}=\|\sigma\|^{2}.
$$
Here $[\chi]$ is the holomorphic line bundle on $\widetilde{\mathcal M}$ corresponding to the unitary character $\chi$, 
$|\cdot|$ is the Hermitian metric on $[\chi]$ induced from the standard norm $|\cdot|$ on ${\bf C}$, 
and the line bundle $\widetilde{\lambda}^{(36+\chi)\ell}\otimes\widetilde{\mu}^{-12\ell}\otimes[\chi]$
is equipped with the Hermitian metric induced from the $L^{2}$-metrics on $\widetilde{\lambda}$ and $\widetilde{\mu}$ 
and the metric $|\cdot|$ on $[\chi]$.
\end{theorem}

\begin{pf}
Set $\xi:=\widetilde{\lambda}^{(36+\chi)}\otimes\widetilde{\mu}^{-12}$.
Let $\Psi$ be a non-zero meromorphic section of $\xi$ and let $E:={\rm div}(\Psi)$. 
By Theorem~\ref{thm:curvature:current:BCOV:invariant} and the Poincar\'e-Lelong formula,
$\phi^{*}\log\tau_{\rm BCOV}$ and $\|\Psi\|^{2}$ satisfy the following equations of currents on $\widetilde{\mathcal M}$:
$$
-12dd^{c}(\phi^{*}\tau_{\rm BCOV})
=
(36+12h^{1,2}+\chi)\,\widetilde{\omega}_{\rm WP}
+
12\,{\rm Ric}\,\widetilde{\omega}_{\rm WP}
-
\sum_{k\in K}\alpha_{k}\delta_{{\mathcal D}_{k}},
$$
$$
-dd^{c}\log\|\Psi\|^{2}
=
(36+12h^{1,2}+\chi)\,\widetilde{\omega}_{\rm WP}
+
12\,{\rm Ric}\,\widetilde{\omega}_{\rm WP}
-
\delta_{E}.
$$
Set $\Delta:=\sum_{k\in K}\alpha_{k}D_{k}-E$. On $\widetilde{\mathcal M}$, we get
$$
-dd^{c}\log\left[\tau_{\rm BCOV}^{12}/\|\Psi\|^{2}\right]
=
-\delta_{\Delta}.
$$
Let $\ell\in{\bf Z}_{>0}$ be an integer such that $\ell\Delta$ is an integral divisor.
By Lemma~\ref{lemma:character:ass:log:1:form} below applied to $F:=(\tau_{\rm BCOV}^{12}/\|\Psi\|^{2})^{\ell}$,
there is a character $\chi\in{\rm Hom}(\pi_{1}(\widetilde{\mathcal M}),U(1))$ and a meromorphic section $s$ of $[\chi]$
with ${\rm div}(s)=\ell\Delta$ such that
$$
|s|^{2}=F=(\tau_{\rm BCOV}^{12}/\|\Psi\|^{2})^{\ell}.
$$
Hence $\sigma:=\Psi^{\ell}\otimes s$ is a meromorphic section of $\xi^{\ell}\otimes[\chi]$ with divisor
$$
{\rm div}(\sigma)={\rm div}(\Psi^{\ell}\otimes s)=\ell E+\ell\Delta=\ell\sum_{k\in K}\alpha_{k}D_{k}
$$
such that $\tau_{\rm BCOV}^{12\ell}=\|\Psi^{\ell}\otimes s\|^{2}=\|\sigma\|^{2}$. This completes the proof.
\end{pf}

\begin{lemma}
\label{lemma:character:ass:log:1:form}
Let $D$ be a divisor on a complex manifold $M$. Let $F$ be a positive function on $M\setminus D$ satisfying $\log F\in L^{1}_{\rm loc}(M)$
and the equation $dd^{c}\log F=\delta_{D}$ of currents on $M$.
Then there exist $\chi\in{\rm Hom}(\pi_{1}(M),U(1))$ and a meromorphic section $s$ of $[\chi]$ with $F=|s|^{2}$ and ${\rm div}(s)=D$.
\end{lemma}

\begin{pf}
The proof is standard and is omitted.
\end{pf}

Since
${\mathcal O}_{\widetilde{\mathcal M}}({\rm div}(\det\phi^{*}\rho)))|_{\phi^{-1}(\widetilde{\mathcal M}_{\rm reg})}\cong{\mathcal O}_{\widetilde{\mathcal M}}$
and since
$$
\widetilde{\lambda}^{(36+\chi)\ell}\otimes\widetilde{\mu}^{-12\ell}\otimes{\mathcal O}_{\widetilde{\mathcal M}}(\chi)
\cong
\widetilde{\lambda}^{(36-12h^{1,2}+\chi)\ell}
\otimes 
K_{\widetilde{\mathcal M}}^{12}
\otimes
{\mathcal O}_{\widetilde{\mathcal M}}(-12{\rm div}(\det\phi^{*}\rho)))
\otimes
[\chi]
$$
by \eqref{eqn:canonical:isomorphism},
$\sigma|_{\phi^{-1}(\widetilde{\mathcal M}_{\rm reg})}$ is a nowhere vanishing holomorphic section of the line bundle
$$
\widetilde{\lambda}^{(36-12h^{1,2}+\chi)\ell}\otimes K_{\widetilde{\mathcal M}}^{12\ell}\otimes[\chi]
=\phi^{*}(\lambda^{(36-12h^{1,2}+\chi)\ell}\otimes K_{\mathcal M}^{12\ell})\otimes[\chi]
$$
on $\phi^{-1}(\widetilde{\mathcal M}_{\rm reg})$. A corresponding statement for $A$-model can be found in \cite[Eq.\,(14)]{BCOV93}.

\begin{question}
By choosing $\widetilde{\mathcal M}$ appropriately, 
is the line bundle $[\chi]$ in Theorem~\ref{thm:consequence:BCOV:conjecture} a torsion element of ${\rm Pic}(\widetilde{\mathcal M})$?
\end{question}

\section
{Holomorphic torsion invariant for $2$-elementary $K3$ surfaces}
\label{sect:2}
\par
\subsection
{Lattices, domains of type IV and orthogonal modular varieties}
\label{sect:2.1}
\par
A free ${\bf Z}$-module of finite rank equipped with a non-degenerate integral symmetric bilinear form is called a lattice.
For a lattice $L$, its rank is denoted by $r(L)$ and its automorphism group is denoted by $O(L)$.
The set of roots of $L$ is defined as $\Delta_{L}:=\{d\in L;\,\langle d,d\rangle=-2\}$.
For a non-zero integer $k\in{\bf Z}$, $L(k)$ denotes the ${\bf Z}$-module $L$ equipped with the rescaled bilinear form $k\langle\cdot,\cdot\rangle_{L}$.
The dual lattice of $L$ is defined as $L^{\lor}:={\rm Hom}_{\bf Z}(L,{\bf Z})\subset L\otimes{\bf Q}$. 
The finite abelian group $A_{L}:=L^{\lor}/L$ is called the {\em discriminant group} of $L$. 
If $A_{L}\cong({\bf Z}/2{\bf Z})^{l}$ for some $l\in{\bf Z}_{\geq0}$, then $L$ is said to be {\em $2$-elementary}.
For an even $2$-elementary lattice $L$, we define $l(L):={\rm rank}_{{\bf Z}_{2}}A_{L}$ and we denote by $\delta(L)\in\{0,1\}$
the parity of the discriminant form on $A_{L}$ (cf. \cite{Nikulin80}).
If $L$ is indefinite, the isometry class of $L$ is determined by the triplet $({\rm sign}(L),l(L),\delta(L))$ (cf. \cite{Nikulin80}).
We define ${\Bbb U}:=({\bf Z}^{2},\binom{0\,1}{1\,0})$. 
For root systems $A_{\ell}$, $D_{\ell}$, $E_{\ell}$, their root lattices are denoted by 
${\Bbb A}_{\ell}$, ${\Bbb D}_{\ell}$, ${\Bbb E}_{\ell}$ and are assumed to be {\em negative-definite}.
The {\em $K3$-lattice} is defined as the even unimodular lattice
$$
{\Bbb L}_{K3}:={\Bbb U}\oplus{\Bbb U}\oplus{\Bbb U}\oplus{\Bbb E}_{8}\oplus{\Bbb E}_{8}.
$$
\par
For a lattice $\Lambda$ of ${\rm sign}(\Lambda)=(2,r(\Lambda)-2)$, we define
$$
\Omega_{\Lambda}
:=
\left\{
[\eta]\in{\bf P}(\Lambda\otimes{\bf C});\,
\langle\eta,\eta\rangle_{\Lambda}=0,\quad
\langle\eta,\bar{\eta}\rangle_{\Lambda}>0
\right\}.
$$
Then $\Omega_{\Lambda}$ consists of two connected components $\Omega_{\Lambda}=\Omega_{\Lambda}^{+}\amalg\Omega_{\Lambda}^{-}$,
each of which is isomorphic to a bounded symmetric domain of type IV of dimension $r(\Lambda)-2$.
The projective $O(\Lambda)$-action on $\Omega_{\Lambda}$ is proper and discontinuous. The quotient
$$
{\mathcal M}_{\Lambda}:=\Omega_{\Lambda}/O(\Lambda)=\Omega_{\Lambda}^{+}/O^{+}(\Lambda)
$$
is an analytic space of dimension $r(\Lambda)-2$,
where $O^{+}(\Lambda)$ is the subgroup of index $2$ of $O(\Lambda)$ preserving $\Omega_{\Lambda}^{\pm}$.
The {\em discriminant locus} of ${\mathcal M}_{\Lambda}$ is the {\em reduced} $O(\Lambda)$-invariant divisor of $\Omega_{\Lambda}$ defined as
$$
{\mathcal D}_{\Lambda}:=\sum_{d\in\Delta_{\Lambda}/\pm1}H_{d},
\qquad
H_{d}:=\{[\eta]\in\Omega_{\Lambda};\,\langle d,\eta\rangle=0\}.
$$
We define
$$
\Omega_{\Lambda}^{0}:=\Omega_{\Lambda}\setminus{\mathcal D}_{\Lambda},
\qquad
{\mathcal M}_{\Lambda}^{0}:=\Omega_{\Lambda}^{0}/O(\Lambda).
$$
By Baily-Borel, ${\mathcal M}_{\Lambda}$ and ${\mathcal M}_{\Lambda}^{0}$ are irreducible, normal quasi-projective varieties.
Let ${\mathcal M}_{\Lambda}^{*}$ be the Baily-Borel compactification of ${\mathcal M}_{\Lambda}$ and define the boundary locus by
$$
{\mathcal B}_{\Lambda}:={\mathcal M}_{\Lambda}^{*}\setminus{\mathcal M}_{\Lambda}.
$$
Then $\dim{\mathcal B}_{\Lambda}=1$ if $r(\Lambda)\geq4$ and $\dim{\mathcal B}_{\Lambda}=0$ if $r(\Lambda)=3$.
When $M$ is a primitive $2$-elementary Lorentzian sublattice of ${\Bbb L}_{K3}$ with $r(M)\geq18$, 
then ${\mathcal B}_{M^{\perp}}$ is irreducible by \cite[Prop.\,11.7]{Yoshikawa13}.
\par
For simplicity, assume the splitting of lattices
$$
\Lambda={\Bbb U}(N)\oplus L,
$$
where $N\in{\bf Z}_{>0}$ and $L$ is an even Lorentzian lattice. 
Let 
$$
{\mathcal C}_{L}:=\{x\in L\otimes{\bf R};\,\langle x,x\rangle_{L}>0\}
$$ 
be the positive cone of $L$, which consists of connected components ${\mathcal C}_{L}={\mathcal C}_{L}^{+}\amalg{\mathcal C}_{L}^{-}$.
We identify the tube domain $L\otimes{\bf R}+i\,{\mathcal C}_{L}$ with $\Omega_{\Lambda}$ via the map
\begin{equation}
\label{eqn:tube:domain:realization}
L\otimes{\bf R}+i\,{\mathcal C}_{L}
\ni z\to 
\left[
\left(
\left(
-\langle z,z\rangle_{L}/2,1/N
\right)
,
z
\right)
\right]
\in\Omega_{\Lambda}.
\end{equation}
The K\"ahler form of the Bergman metric on $L\otimes{\bf R}+i\,{\mathcal C}_{L}$ is the positive $(1,1)$-form defined as
$$
\omega_{\Lambda}(z):=-dd^{c}\log\langle\Im z,\Im z\rangle_{L}.
$$
Via the identification \eqref{eqn:tube:domain:realization}, we regard $\omega_{\Lambda}$ as the $(1,1)$-form on $\Omega_{\Lambda}$.

\subsection
{$2$-elementary $K3$ surfaces and their holomorphic torsion invariants}
\label{sect:2.2}
\par

\subsubsection
{$2$-elementary $K3$ surfaces and their moduli space}
\label{sect:2.2.1}
\par
Let $S$ be a $K3$ surface. Then $H^{2}(S,{\bf Z})$ endowed with the cup-product pairing is isometric to the $K3$-lattice.
Namely, there exists an isometry of lattices $\alpha\colon H^{2}(S,{\bf Z})\cong{\Bbb L}_{K3}$.
\par
Let $\theta\colon S\to S$ be a holomorphic involution.
The pair $(S,\theta)$ is called a {\em $2$-elementary $K3$ surface} if $\theta$ is anti-symplectic, i.e., $\theta^{*}|_{H^{0}(S,K_{S})}=-1$.
The type of an anti-symplectic involution $\theta$ is defined as the isometry class of the invariant sublattice $H^{2}(S,{\bf Z})^{+}$, where 
$$
H^{2}(S,{\bf Z})^{\pm}:=\{v\in H^{2}(S,{\bf Z});\,\theta^{*}v=\pm v\}.
$$
By Nikulin \cite{Nikulin80}, $\alpha(H^{2}(S,{\bf Z})^{+})\subset{\Bbb L}_{K3}$ is a primitive $2$-elementary Lorentzian sublattice.
Since the embedding of a primitive $2$-elementary Lorentzian lattice into ${\Bbb L}_{K3}$ is unique up to an action of $O({\Bbb L}_{K3})$ 
by \cite{Nikulin80}, the type of an anti-symplectic holomorphic involution is independent of the choice of an isometry $\alpha$. 
By Nikulin \cite{Nikulin80}, \cite{Nikulin83}, the topological type of an anti-symplectic holomorphic involution on a $K3$ surface is determined by its type
in the sense that if $(S,\theta)$ and $(S',\theta')$ are two $2$-elementary $K3$ surfaces of the same type, then $(S',\theta')$ is deformation equivalent to $(S,\theta)$.
By \cite{Nikulin80}, \cite{Nikulin83}, there exists $75$ distinct types of $2$-elementary $K3$ surfaces.
The moduli space of $2$-elementary $K3$ surfaces of type $M$ is constructed as follows.
\par
Let $(S,\theta)$ be a $2$-elementary $K3$ surface of type $M$ and
let $\alpha\colon H^{2}(S,{\bf Z})\cong{\Bbb L}_{K3}$ be an isometry with $\alpha(H^{2}(S,{\bf Z})^{-})=M^{\perp}$.
We define the period of $(S,\theta)$ as
$$
\pi_{M}(S,\theta)
:=
\left[\alpha\left(H^{0}(S,K_{S})\right)\right].
$$ 
Then any point of the discriminant locus $\overline{\mathcal D}_{M^{\perp}}$ is never realized as 
the period of a $2$-elementary $K3$ surface of type $M$.
By \cite[Th.\,1.8]{Yoshikawa04}, ${\mathcal M}_{M^{\perp}}^{0}$ is a coarse moduli space of $2$-elementary $K3$ surfaces of type $M$
via the period map.

\subsubsection
{The fixed point set of a $2$-elementary $K3$ surface}
\label{sect:2.2.2}
\par
Let $(S,\theta)$ be a $2$-elementary $K3$ surface of type $M$. Set
$$
S^{\theta}:=\{x\in S;\,\theta(x)=x\}.
$$
The topology of $S^{\theta}$ was determined by Nikulin \cite{Nikulin83} as follows:

\begin{proposition}
\label{prop:topology:fixed:points}
With the same notation as above, the following hold:
\newline{\rm (1)}
If $M\cong{\Bbb U}(2)\oplus{\Bbb E}_{8}(2)$, then $S^{\theta}=\emptyset$ and the quotient $S/\theta$ is an Enriques surface. 
\newline{\rm (2)}
If $M\cong{\Bbb U}\oplus{\Bbb E}_{8}(2)$, then $S^{\theta}=C_{1}^{(1)}\amalg C_{2}^{(1)}$, where $C_{1}^{(1)}$, $C_{2}^{(1)}$ are elliptic curves.
\newline{\rm (3)}
If $M\not\cong{\Bbb U}(2)\oplus{\Bbb E}_{8}(2),{\Bbb U}\oplus{\Bbb E}_{8}(2)$, then
\begin{equation}
\label{eqn:fixed:point:set:2-elementary:K3}
S^{\theta}
=
C^{g(M)}\amalg E_{1}\amalg\ldots\amalg E_{k(M)},
\end{equation}
where $C^{(g)}$ is a curve of genus $g$ and $E_{i}$ is a smooth rational curve with 
\begin{equation}
\label{eqn:fixed:point:set:2-elementary:K3:2}
g(M):=11-\frac{r(M)+l(M)}{2},
\qquad
k(M)=\frac{r(M)-l(M)}{2}.
\end{equation}
\end{proposition}

When $M\cong{\Bbb U}(2)\oplus{\Bbb E}_{8}(2)$, $M$ and $M^{\perp}$ are called {\em exceptional} in this paper.
Notice that, when $M\cong{\Bbb U}\oplus{\Bbb E}_{8}(2)$, the total genus of $S^{\theta}$ is still given by $g(M)=2$, 
so that the first equality of \eqref{eqn:fixed:point:set:2-elementary:K3:2} remains valid.
\par
{\em Warning} : In \cite{Yoshikawa13}, the lattices ${\Bbb U}\oplus{\Bbb E}_{8}(2)$ and ${\Bbb U}\oplus{\Bbb U}\oplus{\Bbb E}_{8}(2)$ 
are also called exceptional. In the present paper, these lattices are {\em not} exceptional.
\par
Let ${\frak S}_{g}$ be the Siegel upper half-space of degree $g$ and let
$$
{\mathcal A}_{g}:={\frak S}_{g}/{\rm Sp}_{2g}({\bf Z})
$$
be the Siegel modular variety of degree $g$, the coarse moduli space of principally polarized abelian varieties of dimension $g$.
Let $\omega_{{\mathcal A}_{g}}$ be the K\"ahler form on ${\mathcal A}_{g}$ in the sense of orbifolds induced from the ${\rm Sp}_{2g}({\bf Z})$-invariant
K\"ahler form on ${\frak S}_{g}$:
$$
\omega_{{\frak S}_{g}}:=-dd^{c}\log\det\Im\tau
$$
\par
After \eqref{eqn:fixed:point:set:2-elementary:K3}, we define the {\em Torelli map} 
$\overline{J}_{M}^{0}\colon{\mathcal M}_{M^{\perp}}^{0}\to{\mathcal A}_{g(M)}$ by
$$
\overline{J}_{M}^{0}(\pi_{M}(S,\theta)):=[\varOmega(S^{\theta})],
$$
where $\varOmega(S^{\theta})\in{\frak S}_{g}$ is the period of $S^{\theta}$ and $[\varOmega(S^{\theta})]\in{\mathcal A}_{g(M)}$ is the corresponding point.
Let $\varPi_{M^{\perp}}\colon\Omega_{M^{\perp}}\to{\mathcal M}_{M^{\perp}}$ be the projection.
We define the holomorphic map $J_{M}^{0}\colon\Omega_{M^{\perp}}^{0}\to{\mathcal A}_{g(M)}$, again called the Torelli map, by
$$
J_{M}^{0}:=\overline{J}_{M}\circ\varPi_{M^{\perp}}.
$$
By Borel's extension theorem, $J_{M}^{0}$ extends to a holomorphic map from 
$\Omega_{M^{\perp}}^{0}\cup{\mathcal D}_{M^{\perp}}^{0}$ to ${\mathcal A}_{g(M)}^{*}$, the Satake compactification of ${\mathcal A}_{g}$,
where ${\mathcal D}_{M^{\perp}}^{0}$ is the dense Zariski open subset of ${\mathcal D}_{M^{\perp}}$ defined as
$$
{\mathcal D}_{M^{\perp}}^{0}:=\bigcup_{d\in\Delta_{M^{\perp}}}H_{d}^{0},
\qquad
H_{d}^{0}:=H_{d}\setminus\bigcup_{\delta\in\Delta_{M^{\perp}}\setminus\{\pm d\}}H_{\delta}.
$$
This extension of $J_{M}^{0}$ is denoted by $J_{M}$. Then the semi-positive $(1,1)$-form $(J_{M}^{0})^{*}\omega_{{\mathcal A}_{g(M)}}$
on $\Omega_{M^{\perp}}^{0}$ extends trivially to a closed positive $(1,1)$-current on $\Omega_{M^{\perp}}$.
The trivial extension of $(J_{M}^{0})^{*}\omega_{{\mathcal A}_{g(M)}}$ from $\Omega_{M^{\perp}}^{0}$ to $\Omega_{M^{\perp}}$
is denoted by $J_{M}^{*}\omega_{{\mathcal A}_{g(M)}}$.

\subsubsection
{A holomorphic torsion invariant for $2$-elementary $K3$ surfaces}
\label{sect:2.2.3}
\par

\begin{definition}
\label{def:analytic:torsion:invariant:2-elementary:K3}
Let $\eta$ be a nowhere vanishing holomorphic $2$-form on $S$.
Let $\gamma$ be an $\theta$-invariant K\"ahler form on $S$.
Define
$$
\begin{aligned}
\tau_{M}(S,\theta)
&:=
{\rm Vol}(S,\gamma)^{\frac{14-r(M)}{4}}
\tau_{{\bf Z}_{2}}(S,\gamma)(\theta)\,
{\rm Vol}(S^{\theta},\gamma|_{S^{\theta}})
\tau(S^{\theta},\gamma|_{S^{\theta}})
\\
&\quad
\times
\exp\left[
\frac{1}{8}\int_{S^{\theta}}
\log\left.\left(\frac{\eta\wedge\overline{\eta}}
{\gamma^{2}/2!}\cdot\frac{{\rm Vol}(S,\gamma)}{\|\eta\|_{L^{2}}^{2}}
\right)\right|_{S^{\theta}}
c_{1}(S^{\theta},\gamma|_{S^{\theta}})
\right].
\end{aligned}
$$
Here ${\rm Vol}(S^{\theta},\gamma|_{S^{\theta}})$ and $\tau(S^{\theta},\gamma|_{S^{\theta}})$ are multiplicative with respect to
the decomposition into connected components.
\end{definition}

By \cite{Yoshikawa04}, $\tau_{M}(S,\theta)$ is independent of the choices of $\eta$ and $\gamma$ and is determined by the period of $(S,\theta)$.
In particular, $\tau_{M}(S,\theta)$ is an invariant of $(S,\theta)$. 
If $\gamma$ is an $\iota$-invariant Ricci-flat K\"ahler form on $S$, then
\begin{equation}
\label{eqn:tau:M:Ricci:flat}
\tau_{M}(S,\theta)
=
{\rm Vol}(S,\gamma)^{\frac{14-r(M)}{4}}\tau_{{\bf Z}_{2}}(S,\gamma)(\theta)
\cdot
{\rm Vol}(S^{\theta},\gamma|_{S^{\theta}})\tau(S^{\theta},\gamma|_{S^{\theta}}).
\end{equation}
Define the $O(M^{\perp})$-invariant smooth function $\widetilde{\tau}_{M}$ on $\Omega_{M^{\perp}}^{0}$ by
$$
\widetilde{\tau}_{M}:=\varPi_{M^{\perp}}^{*}\tau_{M}.
$$
By \cite[Eq.\,(7.1)]{Yoshikawa04}, \cite[Th.\,5.2]{Yoshikawa12}, the following equation of $(1,1)$-currents on $\Omega_{M^{\perp}}$ holds 
\begin{equation}
\label{eqn:curvature:tau:M}
dd^{c}\log\widetilde{\tau}_{M}
=
\frac{r(M)-6}{4}\,\omega_{M^{\perp}}
+
J_{M}^{*}\omega_{{\mathcal A}_{g(M)}}
-
\frac{1}{4}\delta_{{\mathcal D}_{M^{\perp}}}.
\end{equation}

\subsection
{Borcherds products for $2$-elementary lattices}
\label{sect:2.3}
\par
Let ${\frak H}\subset{\bf C}$ be the complex upper half-plane.
Recall that the Dedekind $\eta$-function and the Jacobi theta series are holomorphic functions on ${\frak H}$ defined as 
$$
\eta(\tau)=e^{2\pi i\tau/24}\prod_{n=1}^{\infty}
\left(
1-e^{2\pi in\tau}
\right),
\qquad
\vartheta_{{\Bbb A}_{1}^{+}+\frac{k}{2}}(\tau)=\sum_{n\in{\bf Z}}e^{2\pi i(n+\frac{k}{2})^{2}\tau}
\quad
(k=0,1).
$$
\par
For a $2$-elementary lattice $\Lambda$, we set
$$
\phi_{\Lambda}(\tau):=\eta(\tau)^{-8}\eta(2\tau)^{8}\eta(4\tau)^{-8}\,\theta_{{\Bbb A}_{1}^{+}}(\tau)^{12-r(\Lambda)}. 
$$ 
Let $\{{\frak e}_{\gamma}\}_{\gamma\in A_{\Lambda}}$ be the standard basis of the group ring ${\bf C}[A_{\Lambda}]$.
Let ${\rm Mp}_{2}({\bf Z})$ be the metaplectic double cover of ${\rm SL}_{2}({\bf Z})$ and 
let $\rho_{\Lambda}\colon{\rm Mp}_{2}({\bf Z})\to{\rm GL}({\bf C}[A_{\Lambda}])$ be the Weil representation \cite[Sect.\,4]{Borcherds98}.
Define a ${\bf C}[A_{\Lambda}]$-valued holomorphic functions on ${\frak H}$ by
$$
F_{\Lambda}
:=
\sum_{\gamma\in\widetilde{\Gamma}_{0}(4)\backslash{\rm Mp}_{2}({\bf Z})}
\phi_{\Lambda}|_{\gamma}\,
\rho_{\Lambda}(\gamma^{-1})\,{\frak e}_{0},
$$
where $|_{\gamma}$ is the Petersson slash operator.
When ${\rm sign}(\Lambda)=(2,r(\Lambda)-2)$, we deduce from \cite[Th.\,7.7]{Yoshikawa13} that $F_{\Lambda}(\tau)$ is a modular form 
of type $\rho_{\Lambda}$ of weight $(4-r(\Lambda))/2$.
To get the integrality of the Fourier coefficients, we set
$$
g(\Lambda):=(r(\Lambda)-l(\Lambda))/2.
$$
If $\Lambda=M^{\perp}$ for a primitive $2$-elementary Lorentzian sublattice $M\subset{\Bbb L}_{K3}$,
then $g(\Lambda)=g(M)$. By \cite[Sect.\,7]{MaYoshikawa14}, $2^{g(\Lambda)-1}F_{\Lambda}(\tau)$ always has integral Fourier expansion.
\par
Let $\alpha\in{\bf Z}_{>0}$ be such that $\alpha F_{\Lambda}$ has integral Fourier expansion.
By Borcherds \cite[Th.\,13.3]{Borcherds98}, the Borcherds lift $\Psi_{\Lambda}(\cdot,\alpha F_{\Lambda})$
is an automorphic form on $\Omega_{\Lambda}^{+}$ for $O^{+}(\Lambda)$.
For the weight and the singularities of $\Psi_{\Lambda}(\cdot,\alpha F_{\Lambda})$, see \cite[Th.\,7.1]{MaYoshikawa14}.
The infinite product expansion of $\Psi_{\Lambda}(\cdot,\alpha F_{\Lambda})$ is given as follows.
For simplicity, assume 
$$
\Lambda={\Bbb U}\oplus L.
$$
Then $F_{\Lambda}(\tau)=F_{L}(\tau)$.
Let $F_{\Lambda}(\tau)=\sum_{\gamma\in A_{\Lambda}}{\frak e}_{\gamma}\sum_{k\in{\bf Z}+\gamma^{2}/2}c_{\gamma}(k)\,e^{2\pi ik\tau}$
be the Fourier series expansion.
Under the identification $\Omega_{\Lambda}^{+}\cong L\otimes{\bf R}+i\,{\mathcal C}_{L}^{+}$ via \eqref{eqn:tube:domain:realization},
the following equality holds for $z\in L\otimes{\bf R}+i\,{\mathcal W}$ with $\langle\Im z,\Im z\rangle\gg0$
\begin{equation}
\label{eqn:infinite:product}
\Psi_{\Lambda}(z,\alpha F_{\Lambda})
=
e^{2\pi i\langle\alpha\varrho,z\rangle}
\prod_{\lambda\in L^{\lor},\,\lambda\cdot{\mathcal W}>0}
\left(
1-e^{2\pi i\langle\lambda,z\rangle}
\right)^{\alpha c_{\overline{\lambda}}(\lambda^2/2)},
\end{equation}
where $\overline{\lambda}:=\lambda+L\in A_{L}=A_{\Lambda}$,
the cone ${\mathcal W}\subset C_{L}^{+}$ is a {\em Weyl chamber} of $\alpha F_{L}$ and
$\alpha\varrho\in L\otimes{\bf Q}$ is the {\em Weyl vector} of $\alpha F_{L}$. See \cite[Sect.\,10]{Borcherds98} for these notions. 
\par
By \cite[Th.\,13.3]{Borcherds98}, there exists $w(\Lambda)\in{\bf Q}$ such that 
$\Psi_{\Lambda}(\cdot,\alpha F_{\Lambda})$ has weight $\alpha w(\Lambda)$.
The Petersson norm of $\Psi_{\Lambda}(\cdot,\alpha F_{\Lambda})$ is the $C^{\infty}$ function on
$L\otimes{\bf R}+i\,{\mathcal C}_{L}^{+}$ defined by
$$
\|\Psi_{\Lambda}(z,\alpha F_{\Lambda})\|^{2}:=\langle\Im z,\Im z\rangle_{L}^{\alpha w(\Lambda)}\left|\Psi_{\Lambda}(z,\alpha F_{\Lambda})\right|^{2}.
$$
Through \eqref{eqn:tube:domain:realization}, 
$\|\Psi_{\Lambda}(\cdot,\alpha F_{\Lambda})\|^{2}$ is viewed as an $O(\Lambda)$-invariant function on $\Omega_{\Lambda}$. 
Set
$$
\|\Psi_{\Lambda}(\cdot,F_{\Lambda})\|:=\|\Psi_{\Lambda}(\cdot,\alpha F_{\Lambda})\|^{1/\alpha}.
$$
In what follow, $\|\Psi_{\Lambda}(\cdot,F_{\Lambda})\|^{2}$ is viewed as a function on ${\mathcal M}_{\Lambda}$.

\subsection
{An explicit formula for $\tau_{M}$}
\label{sect:2.4}
\par
We recall the main result of \cite{MaYoshikawa14}.
Let $a,b\in\{0,1/2\}^{g}$. The pair $(a,b)$ is said to be even if $4{}^{t}a\cdot b\in2{\bf Z}$.
For even $(a,b)$, the corresponding Riemann theta constant $\theta_{a,b}(\varOmega)$ is defined as the theta series
$$
\theta_{a,b}(\varOmega)
:=
\sum_{n\in{\bf Z}^{g}}\exp\{\pi i{}^{t}(n+a)\varOmega(n+a)+2\pi i{}^{t}(n+a)b\}.
$$
For $\varOmega\in{\frak S}_{g}$, we define $\chi_{g}(\varOmega)$ and $\Upsilon_{g}(\varOmega)$ by
$$
\chi_{g}(\varOmega)
:=
\prod_{(a,b)\,{\rm even}}
\theta_{a,b}(\varOmega),
\qquad
\Upsilon_{g}(\varOmega)=\chi_{g}(\varOmega)^{8}\sum_{(a,b)\,{\rm even}}\theta_{a,b}(\varOmega)^{-8}.
$$
Then $\chi_{g}^{8}$ and $\Upsilon_{g}$ are Siegel modular forms of weights $2^{g+1}(2^{g}+1)$ and $2(2^{g}-1)(2^{g}+2)$, respectively.
Their Petersson norms are ${\rm Sp}_{2g}({\bf Z})$-invariant $C^{\infty}$ functions
$$
\|\chi_{g}(\varOmega)^{8}\|^{2}
:=
\left(\det\Im\varOmega\right)^{2^{g+1}(2^{g}+1)}
\left|
\chi_{g}(\varOmega)^{8}
\right|^{2},
$$
$$
\|\Upsilon_{g}(\varOmega)\|^{2}
:=
\left(\det\Im\varOmega\right)^{2(2^{g}-1)(2^{g}+2)}
\left|
\Upsilon_{g}(\varOmega)^{8}
\right|^{2}
$$
on ${\frak S}_{g}$.
We regard $\|\chi_{g}^{8}\|^{2},\|\Upsilon_{g}\|^{2}\in C^{\infty}({\mathcal A}_{g})$ in what follows.

\begin{theorem}
\label{thm:torsion:invariant:2-elementary:K3:Ma:Yoshikawa}
Let $M$ be a primitive $2$-elementary Lorentzian sublattice of ${\Bbb L}_{K3}$.
Set $r:=r(M)$, $\delta:=\delta(M)$, $g:=g(M)$ and $\Lambda:=M^{\perp}$.
Then there exists a constant $C_{M}$ depending only on $M$ such that the following equality of functions on ${\mathcal M}_{\Lambda}^{0}$ holds:
\newline{\rm (1) }
If $(r,\delta)\not=(2,0),(10,0)$, then
$$
\tau_{M}^{-2^{g}(2^{g}+1)}
=
C_{M}\,
\left\|
\Psi_{\Lambda}(\cdot,2^{g-1}F_{\Lambda})
\right\|
\cdot 
J_{M}^{*}\left\|
\chi_{g}^{8}
\right\|.
$$
\newline{\rm (2) }
If $(r,\delta)=(10,0)$, then
$$
\tau_{M}^{-(2^{g}-1)(2^{g}+2)}
=
C_{M}\,
\left\|
\Psi_{\Lambda}(\cdot,(2^{g-1}+1)F_{\Lambda})
\right\|
\cdot
J_{M}^{*}
\left\|
\Upsilon_{g}
\right\|.
$$
\newline{\rm (3) }
If $(r,\delta)=(2,0)$, then $M\cong{\Bbb U}$ or ${\Bbb U}(2)$ and
$$
\tau_{M}^{-(2^{g}-1)(2^{g}+2)}
=
C_{M}\,
\left\|
\Psi_{\Lambda}(\cdot,2^{g-1}F_{\Lambda}+f_{\Lambda})
\right\|
\cdot 
J_{M}^{*}
\left\|
\Upsilon_{g}
\right\|,
$$
where $f_{\Lambda}$ is the elliptic modular form of type $\rho_{\Lambda}$ given as follows:
\newline
When $\Lambda={\Bbb U}^{\perp}={\Bbb U}^{\oplus2}\oplus{\Bbb E}_{8}^{\oplus2}$, 
$$
f_{\Lambda}(\tau):=\theta_{{\Bbb E}_{8}^{+}}(\tau)/\eta(\tau)^{24},
\qquad
\theta_{{\Bbb E}_{8}^{+}}(\tau):=\sum_{\lambda\in{\Bbb E}_{8}^{+}}\exp(\pi i\langle\lambda,\lambda\rangle\tau).
$$
When $\Lambda={\Bbb U}(2)^{\perp}={\Bbb U}(2)\oplus{\Bbb U}\oplus{\Bbb E}_{8}^{\oplus2}$, 
$$
f_{\Lambda}(\tau)
:=
8
\sum_{\gamma\in A_{\Lambda}}
\{
\frac{1}{\eta(\tau/2)^{8}\eta(\tau)^{8}}
+
\frac{(-1)^{\gamma^{2}}}{\eta(\frac{\tau+1}{2})^{8}\eta(\tau+1)^{8}}
\}\,
{\frak e}_{\gamma}
+
\frac{{\frak e}_{0}}{\eta(\tau)^{8}\eta(2\tau)^{8}}.
$$
\end{theorem}

\begin{pf}
See \cite[Th.\,0.1]{MaYoshikawa14}.
\end{pf}

\section
{BCOV invariants for Borcea-Voisin threefolds}
\label{sect:3}
\par
In Section~\ref{sect:3}, $M\subset{\Bbb L}_{K3}$ denotes a primitive $2$-elementary Lorentzian sublattice. 
Throughout this section, we keep the following notation
$$
\Lambda:=M^{\perp},
\qquad
g:=g(M),
\qquad
r:=r(M),
\qquad
\delta:=\delta(M).
$$
Hence ${\rm sign}(\Lambda)=(2,20-r)$ and $\dim\Omega_{\Lambda}=20-r$.

\begin{definition}
\label{def:Borcea:Voisin:orbifold}
Let $(S,\theta)$ be a $2$-elementary $K3$ surface and let $T$ be an elliptic curve. The orbifold
$$
X_{(S,\theta,T)}:=S\times T/\theta\times(-1_{T}).
$$
is called a {\em Borcea-Voisin orbifold} associated with $(S,\theta,T)$.
The type of a Borcea-Voisin orbifold $X_{(S,\theta,T)}$ is defined as that of $(S,\theta)$. 
Hence the type of $X_{(S,\theta,T)}$ is an isometry class of a primitive $2$-elementary Lorentzian sublattice of ${\Bbb L}_{K3}$.
\end{definition}

Let $T[2]=T^{(-1_{T})}$ be the set of points of order $2$ of $T$. Then
$$
{\rm Sing}\,X_{(S,\theta,T)}=S^{\theta}\times T[2]
$$
is the $4$ copies of $S^{\theta}$. Let 
$$
p\colon\widetilde{X}_{(S,\theta,T)}\to X_{(S,\theta,T)}
$$
be the blowing-up of ${\rm Sing}(X_{(S,\theta,T)})=S^{\theta}\times T[2]$. 
Then $\widetilde{X}_{(S,\theta,T)}$ is a Calabi-Yau threefold, called a {\em Borcea-Voisin threefold},
whose mirror symmetry was studied by Borcea \cite{Borcea97}, Voisin \cite{Voisin93} and Gross-Wilson \cite{GrossWilson97}.
One of the main results of this paper is the following:

\begin{theorem}
\label{conj:BCOV=orb:BCOV}
There exists a constant $C_{M}$ depending only on the lattice $M$ such that
for every Borcea-Voisin orbifold $X_{(S,\theta,T)}$ of type $M$,
$$
\tau_{\rm BCOV}(\widetilde{X}_{(S,\theta,T)})
=
C_{M}\,
\tau_{M}\left(S,\theta\right)^{-4}
\left\|
\eta(\varOmega(T))^{24}
\right\|^{2}.
$$
\end{theorem}

\begin{corollary}
\label{cor:BCOV:invariant:Borcea:Voisin:threefolds}
There exists a constant $C_{M}$ depending only on $M$ such that the following equality holds 
for every Borcea-Voisin threefold $\widetilde{X}_{(S,\theta,T)}$ of type $M$:
\begin{itemize}
\item[(1)]
If $(r,\delta)\not=(2,0),(10,0)$, then
\begin{equation}
\label{eqn:BCOV:explicit:formula:1}
\begin{aligned}
\tau_{\rm BCOV}(\widetilde{X}_{(S,\theta,T)})^{2^{g-1}(2^{g}+1)}
&=
C_{M}\,
\left\|
\Psi_{\Lambda}
\left(
\pi_{M}(S,\theta),2^{g-1}F_{\Lambda}
\right)
\right\|^{2}
\left\|
\chi_{g}
\left(
\varOmega(S^{\theta})
\right)^{8}
\right\|^{2}
\\
&\qquad\times
\left\|
\eta
\left(
\varOmega(T)
\right)^{24}
\right\|^{2^{g}(2^{g}+1)}.
\end{aligned}
\end{equation}
\item[(2)]
If $(r,\delta)=(10,0)$, then
\begin{equation}
\label{eqn:BCOV:explicit:formula:2}
\begin{aligned}
\tau_{\rm BCOV}(\widetilde{X}_{(S,\theta,T)})^{(2^{g-1}+1)(2^{g}-1)}
&=
C_{M}\,
\left\|
\Psi_{\Lambda}
\left(
\pi_{M}(S,\theta),(2^{g-1}+1)F_{\Lambda}
\right)
\right\|^{2}
\\
&\qquad\times
\left\|
\Upsilon_{g}
\left(
\varOmega(S^{\theta})
\right)
\right\|^{2}
\left\|
\eta
\left(
\varOmega(T)
\right)^{24}
\right\|^{2(2^{g-1}+1)(2^{g}-1)}.
\end{aligned}
\end{equation}
\item[(3)]
If $(r,\delta)=(2,0)$, then
\begin{equation}
\label{eqn:BCOV:explicit:formula:3}
\begin{aligned}
\tau_{\rm BCOV}(\widetilde{X}_{(S,\theta,T)})^{(2^{g-1}+1)(2^{g}-1)}
&=
C_{M}\,
\left\|
\Psi_{\Lambda}
\left(
\pi_{M}(S,\theta),2^{g-1}F_{\Lambda}+f_{\Lambda}
\right)
\right\|^{2}
\\
&\qquad\times
\left\|
\Upsilon_{g}
\left(
\varOmega(S^{\theta})
\right)
\right\|^{2}
\left\|
\eta
\left(
\varOmega(T)
\right)^{24}
\right\|^{2(2^{g-1}+1)(2^{g}-1)},
\end{aligned}
\end{equation}
where $f_{\Lambda}$ is the same elliptic modular form as in Theorem~\ref{thm:torsion:invariant:2-elementary:K3:Ma:Yoshikawa} (3).
\end{itemize}
\end{corollary}

\begin{pf}
The result follows from Theorems~\ref{thm:torsion:invariant:2-elementary:K3:Ma:Yoshikawa} and \ref{conj:BCOV=orb:BCOV}.
\end{pf}

\begin{remark}
By the same argument as in \cite[Sect.\,6.1.1]{KlemmMarino08},
the B-model canonical coordinates on the moduli space of Borcea-Voisin threefolds with $g=0$ coincide with the standard linear coordinates
$(z_{1},\ldots,z_{m},\tau)$ on $\Omega_{\Lambda}\times{\frak H}$ induced by the isomorphism \eqref{eqn:tube:domain:realization},
where $m=\dim\Omega_{\Lambda}$.
As a result, we get $\|\varXi_{s}/\langle A_{0}^{\lor},\varXi_{s}\rangle\|^{2}=2\Im\tau\cdot\langle\Im z,\Im z\rangle$ and
$$
\left\|
\frac{\partial}{\partial t_{1}}\wedge\cdots\wedge\frac{\partial}{\partial t_{n}}
\right\|^{2}
=
\left(
-\frac{\partial^{2}\log\Im\tau}{\partial\tau\partial\bar{\tau}}
\right)
\det
\left(
-\frac{\partial^{2}\log\langle\Im z,\Im z\rangle}{\partial z_{\alpha}\partial\bar{z}_{\beta}}
\right)
=
(\Im\tau)^{-2}\langle\Im z,\Im z\rangle^{-m}
$$
in \eqref{eqn:BCOV:conj:2}. In particular, $\|\varXi_{s}/\langle A_{0}^{\lor},\varXi_{s}\rangle\|^{2}$ and
$\|(\partial/\partial t_{1})\wedge\cdots\wedge(\partial/\partial t_{n})\|^{2}$ are expressed by the Bergman kernel functions of 
$\Omega_{\Lambda}$ and ${\frak H}$.
If the conjecture of Bershadsky-Cecotti-Ooguri-Vafa \eqref{eqn:BCOV:conj:2} holds true for Borcea-Voisin threefolds with $g=0$, this fact implies that
{\em the infinite product $F_{1}^{\rm top}(t)$ in \eqref{eqn:BCOV:conj:2} must extend to an automorphic form on $\Omega_{\Lambda}\times{\frak H}$.}
In this sense, Conjecture \eqref{eqn:BCOV:conj:2} may be viewed as a conjectural extension of the theory of Borcherds products to 
a certain specific section on the moduli space of Calabi-Yau threefolds.
Corollary~\ref{cor:BCOV:invariant:Borcea:Voisin:threefolds} verifies the requirement of infinite product for Borcea-Voisin threefolds with $g=0$.
In particular, Conjecture \eqref{eqn:BCOV:conj:2} for Borcea-Voisin threefolds with $g=0$ is reduced to the conjecture that
the exponents of the infinite products $\Psi_{\Lambda}(\cdot,2^{g-1}F_{\Lambda})$ and $\eta(\tau)$ are given by the instanton numbers of the mirror
of such Borcea-Voisin threefolds. 
However, the Borcea-Voisin mirror construction does not apply to the case $g=0$ because of the non-existence of a primitive 
${\Bbb U}\subset\Lambda$ (cf. \cite[Sect.\,2.6]{Voisin93}). Thus we are naturally led to the following question:
{\em What is the mirror of Borcea-Voisin threefolds with $g=0$?} 
To our knowledge, this basic question is still open.
\end{remark}

\subsection
{A variational formula for $\tau_{\rm BCOV}$}
\label{sect:3.1}
\par
We keep the notation in Section~\ref{sect:2}.
\par
The modular curve ${\bf X}(1)$ is the quotient of ${\frak H}$ defined as
$$
{\bf X}(1):={\frak H}/{\rm SL}_{2}({\bf Z}).
$$ 
Let $\omega_{\rm hyp}$ be the K\"ahler form of the Poincar\'e metric on ${\frak H}$
$$
\omega_{\rm hyp}=-dd^{c}\log\Im \tau.
$$
\par
By definition, Borcea-Voisin threefolds of type $M$ are parametrized by the product ${\mathcal M}_{\Lambda}^{0}\times{\bf X}(1)$.
Hence $\tau_{\rm BCOV}\in C^{\omega}({\mathcal M}_{\Lambda}^{0}\times{\bf X}(1))$.

\begin{theorem}
\label{thm:comparison:curvature:BCOV:orb:BCOV}
Regard $\tau_{\rm BCOV}$ as an $O(\Lambda)\times{\rm SL}_{2}({\bf Z})$-invariant $C^{\omega}$ function on $\Omega_{\Lambda}^{0}\times{\frak H}$.
Then the following equations of $(1,1)$-forms on $\Omega_{\Lambda}^{0}\times{\frak H}$ hold:
\begin{equation}
\label{eqn:curvature:BCOV:invariant:Borcea:Voisin}
-dd^{c}\log\tau_{\rm BCOV}
=
(r(M)-6){\rm pr}_{1}^{*}\omega_{\Lambda}
+
4{\rm pr}_{1}^{*}J_{M}^{*}\omega_{{\mathcal A}_{g}}
+
12{\rm pr}_{2}^{*}\omega_{\rm hyp}.
\end{equation}
In particular, $\log\left[\tau_{\rm BCOV}/(\tau_{M}^{-4}\|\eta^{24}\|^{2})\right]$ is a pluriharmonic function on ${\mathcal M}_{\Lambda}^{0}\times{\bf X}(1)$.
\end{theorem}

\begin{pf}
Take a $2$-elementary $K3$ surface $(S,\theta)$ of type $M$ and an elliptic curve $T$.
Let $f\colon({\mathcal S},\theta)\to{\rm Def}(S,\theta)$ be the Kuranishi family of $(S,\theta)$ and
let $g\colon{\mathcal T}\to{\rm Def}(T)$ be the Kuranishi family of $T$.
Set $\widetilde{X}:=\widetilde{X}_{(S,\theta,T)}$ and let ${\frak f}\colon({\frak X},\widetilde{X})\to({\rm Def}(\widetilde{X},[\widetilde{X}])$ 
be the Kuranishi family of $\widetilde{X}$. We have the embedding of germs 
$$
\mu\colon{\rm Def}(S,\theta)\times{\rm Def}(T)\to{\rm Def}(\widetilde{X}),
\qquad
\mu(s,t):=[\widetilde{X}_{(S_{s},\theta_{s},T_{t})}]
$$ 
for $s\in{\rm Def}(S,\theta)$, $t\in{\rm Def}(T)$, where $(S_{s},\theta_{s})=f^{-1}(s)$ and $T_{t}=g^{-1}(t)$.
\par{\em (Step 1) }
Let
$$
E:=p^{-1}({\rm Sing}(X_{(S,\theta,T)}))=p^{-1}(S^{\theta}\times T[2])
$$
be the exceptional divisor of $p\colon\widetilde{X}_{(S,\theta,T)}\to X_{(S,\theta,T)}$. 
Then $E$ is a ${\bf P}^{1}$-bundle over $S^{\theta}\times T[2]$, whose structure is given as follows.
Let ${\mathcal N}:=N_{(S^{\theta}\times T[2])/(S\times T)}$ be the normal bundle of $S^{\theta}\times T[2]$ in $S\times T$.
Then ${\mathcal O}_{{\bf P}({\mathcal N})}(-1)\subset{\bf P}({\mathcal N})\times{\mathcal N}$.
The projection ${\mathcal O}_{{\bf P}({\mathcal N})}(-1)\to{\mathcal N}$ is the blowing-up of the zero section of ${\mathcal N}$ and
$$
{\mathcal O}_{{\bf P}({\mathcal N})}(-2)
=
{\mathcal O}_{{\bf P}({\mathcal N})}(-1)/\pm1
\to 
{\mathcal N}/\pm1
$$ 
is a {\em crepant} resolution.
Let $N_{E/\widetilde{X}}$ be the normal bundle of $E$ in $\widetilde{X}$. Then
\begin{equation}
\label{eqn:exceptional:divisor}
E={\bf P}(N_{(S^{\theta}\times T[2])/(S\times T)}),
\qquad
N_{E/\widetilde{X}}|_{E}={\mathcal O}_{E}(-2).
\end{equation}
We set $p_{E}:=p|_{E}$. Then $p_{E}\colon E={\bf P}({\mathcal N})\to S^{\theta}\times T[2]$ is the projection of ${\bf P}^{1}$-bundle.
\par
Let $i\colon E\hookrightarrow\widetilde{X}$ be the inclusion.
By Voisin \cite[Lemme 1.7]{Voisin93}, we have the decomposition
\begin{equation}
\label{eqn:decomposition:cohomology:Borcea:Voisin:manifold}
\begin{aligned}
H^{1}(\widetilde{X},\Omega_{\widetilde{X}}^{2})
&=
[H^{1}(S,\Omega^{1}_{S})^{-}\otimes H^{0}(T,K_{T})]
\oplus
[H^{0}(S,K_{S})\otimes H^{1}(T,{\mathcal O}_{T})]
\\
&\qquad
\oplus
i_{*}p_{E}^{*}H^{0}(S^{\theta}\times T[2],\Omega_{S^{\theta}\times T[2]}^{1}),
\end{aligned}
\end{equation}
which is orthogonal with respect to the $L^{2}$-metric on $H^{1}(\widetilde{X},\Omega_{\widetilde{X}}^{2})$.
The Kodaira-Spencer map induces the following isomorphisms
\begin{equation}
\label{eqn:Kodaira:Spencer:map:Borcea:Voisin:1}
\begin{aligned}
\rho_{1}
\colon
\Theta_{{\rm Def}(S,\theta),[(S,\theta)]}
&\cong 
H^{1}(S,\Theta_{S})^{+}
\cong
H^{1}(S,\Omega^{1}_{S})^{-}\otimes H^{0}(S,K_{S})^{\lor}
\\
&\cong
[H^{1}(S,\Omega^{1}_{S})^{-}\otimes H^{0}(T,K_{T})]\otimes H^{0}(\widetilde{X},K_{\widetilde{X}})^{\lor},
\end{aligned}
\end{equation}
\begin{equation}
\label{eqn:Kodaira:Spencer:map:Borcea:Voisin:2}
\begin{aligned}
\rho_{2}
\colon
\Theta_{{\rm Def}(T),[T]}
&\cong 
H^{1}(T,\Theta_{T})
=
H^{1}(T,{\mathcal O}_{T})\otimes H^{0}(T,K_{T})^{\lor}
\\
&\cong
[H^{0}(S,K_{S})\otimes H^{1}(T,{\mathcal O}_{T})]\otimes H^{0}(\widetilde{X},K_{\widetilde{X}})^{\lor}.
\end{aligned}
\end{equation}
By \eqref{eqn:decomposition:cohomology:Borcea:Voisin:manifold}, \eqref{eqn:Kodaira:Spencer:map:Borcea:Voisin:1},
\eqref{eqn:Kodaira:Spencer:map:Borcea:Voisin:2},  we get the following canonical identification
\begin{equation}
\label{eqn:Kodaira:Spencer:map:Borcea:Voisin:3}
\begin{aligned}
\rho_{3}\colon
&\Theta_{{\rm Def}(\widetilde{X}),[\widetilde{X}]}
\cong 
H^{1}(\widetilde{X},\Theta_{\widetilde{X}})
\cong
H^{1}(\widetilde{X},\Omega_{\widetilde{X}}^{2})\otimes H^{0}(\widetilde{X},K_{\widetilde{X}})^{\lor}
\cong
\\
&
\Theta_{{\rm Def}(S,\theta),[(S,\theta)]}
\oplus
\Theta_{{\rm Def}(T),[T]}
\oplus
\left[
i_{*}p_{E}^{*}H^{0}(S^{\theta}\times T[2],\Omega_{S^{\theta}\times T[2]}^{1})\otimes H^{0}(\widetilde{X},K_{\widetilde{X}})^{\lor}
\right],
\end{aligned}
\end{equation}
where the first isomorphism is given by the Kodaira-Spencer map and the last decomposition is orthogonal with respect to $\omega_{\rm WP}$.
By \eqref{eqn:Weil:Petersson:form:2}, we have
\begin{equation}
\label{eqn:restriction:Weil:Petersson:form}
(\mu^{*}\omega_{\rm WP})|_{\Theta_{{\rm Def}(S,\theta)}\oplus\Theta_{{\rm Def}(T)}}
=
{\rm pr}_{1}^{*}\omega_{\Lambda}+{\rm pr}_{2}^{*}\omega_{\rm hyp}.
\end{equation}
\par{\em (Step 2) }
Let $({\frak F},h_{\frak F})$ be the automorphic vector bundle of rank $g$ on ${\frak S}_{g}$ equipped with the Hermitian structure induced from
the polarization such that 
\begin{equation}
\label{eqn:Hodge:vector:bundle}
(f_{*}\Omega_{{\mathcal S}^{\theta}/{\rm Def}(S^{\theta})}^{1},h_{L^{2}})=J_{M}^{*}({\frak F},h_{\frak F}).
\end{equation}
(${\frak F}$ is the relative cotangent bundle of the universal family of principally polarized abelian varieties over ${\frak S}_{g}$.)
Let us see the isometry of Hermitian vector spaces
\begin{equation}
\label{eqn:isometry:cohomology:Borcea:Voisin:1}
\begin{aligned}
\,&
\left(
i_{*}p_{E}^{*}H^{0}(S^{\theta}\times T[2],\Omega_{S^{\theta}\times T[2]}^{1})\otimes H^{0}(\widetilde{X},K_{\widetilde{X}})^{\lor},
\omega_{\rm WP}
\right)
\\
&\cong
\left(
H^{0}(S^{\theta}\times T[2],\Omega_{S^{\theta}\times T[2]}^{1}),
2h_{L^{2}}
\right)
\otimes
\left(
H^{0}(\widetilde{X},K_{\widetilde{X}}),h_{L^{2}}
\right)^{\lor}.
\end{aligned}
\end{equation}
Recall that $i_{*}\colon H^{*}(E)\to H^{*+2}(\widetilde{X})$ is defined as the dual of $i^{*}\colon H^{*}(\widetilde{X})\to H^{*}(E)$ 
with respect to the Poincar\'e duality pairing:
For all $\varphi\in H^{*}(\widetilde{X})$, $\psi\in H^{4-*}(E)$,
$$
\int_{\widetilde{X}}i_{*}\psi\wedge\varphi=2\pi\int_{E}\psi\wedge i^{*}\varphi.
$$
Let $\Phi\in H^{2}_{c}(N_{E/\widetilde{X}})$ be the Thom form of Mathai-Quillen (cf. \cite[(1.37)]{BerlineGetzlerVergne92}).
Identify $N_{E/\widetilde{X}}$ with a tubular neighborhood of $E$ in $\widetilde{X}$
and regard $\Phi$ as a $C^{\infty}$ closed $2$-form on $\widetilde{X}$ supported on the closure of $N_{E/\widetilde{X}}$ by this identification.
Since $i_{*}\psi=[\Phi]\wedge\psi$ and $i^{*}\Phi=2\pi\,c_{1}(N_{E/\widetilde{X}})$ by \cite[(1.38)]{BerlineGetzlerVergne92}, 
we deduce from \eqref{eqn:exceptional:divisor} that for all $\omega,\omega'\in H^{0}(S^{\theta}\times T[2],\Omega_{S^{\theta}\times T[2]}^{1})$
$$
\begin{aligned}
\,&
\left\langle 
i_{*}p_{E}^{*}\omega,i_{*}p_{E}^{*}\overline{\omega'}
\right\rangle_{L^{2}}
=
-\frac{\sqrt{-1}}{(2\pi)^{3}}\int_{\widetilde{X}}
(i_{*}p_{E}^{*}\omega)\wedge\Phi\wedge p_{E}^{*}\overline{\omega'}
=
-\frac{\sqrt{-1}}{(2\pi)^{2}}\int_{E}p_{E}^{*}(\omega\wedge\overline{\omega'})\wedge i^{*}\Phi
\\
&=
-\frac{\sqrt{-1}}{2\pi}\int_{E}p_{E}^{*}(\omega\wedge\overline{\omega'})\wedge c_{1}({\mathcal O}_{E}(-2))
=
\frac{\sqrt{-1}}{\pi}\int_{S^{\theta}\times T[2]}\omega\wedge\overline{\omega'}
=
2\langle\omega,\overline{\omega'}\rangle_{L^{2}},
\end{aligned}
$$
where we used the projection formula to get the $4$-th equality.
This verifies \eqref{eqn:isometry:cohomology:Borcea:Voisin:1}.
By \eqref{eqn:restriction:Weil:Petersson:form}, \eqref{eqn:isometry:cohomology:Borcea:Voisin:1},
we have an isometry of holomorphic Hermitian vector bundles on ${\rm Def}(S,\theta)\times{\rm Def}(T)$
\begin{equation}
\label{eqn:isometry:cohomology:Borcea:Voisin:2}
\begin{aligned}
\mu^{*}\left(
\Theta_{{\rm Def}(\widetilde{X})},
\omega_{\rm WP}
\right)
&\cong
\left(
\Theta_{{\rm Def}(S,\theta)},\omega_{M^{\perp}}
\right)
\oplus
\left(
\Theta_{{\rm Def}(T)},\omega_{\rm hyp}
\right)
\\
&\qquad
\oplus
\left[
J_{M}^{*}({\frak F}^{\oplus4},2h_{{\frak F}^{\oplus4}})\otimes({\frak f}_{*}K_{{\frak X}/{\rm Def}(\widetilde{X})},h_{L^{2}})^{\lor}
\right].
\end{aligned}
\end{equation}
\par{\em (Step 3) }
Since 
${\rm Ric}\,\omega_{\rm WP}=c_{1}(\Theta_{{\rm Def}(\widetilde{X})},\omega_{\rm WP})$,
$c_{1}({\frak F},2h_{{\frak F}})=\omega_{{\mathcal A}_{g}}$
and
$\omega_{\rm WP}=c_{1}({\frak f}_{*}K_{{\frak X}/{\rm Def}(\widetilde{X})},h_{L^{2}})$,
we get by \eqref{eqn:isometry:cohomology:Borcea:Voisin:2}
\begin{equation}
\label{eqn:Ricci:form:Weil:Petersson:Borcea:Voisin}
\begin{aligned}
\mu^{*}{\rm Ric}\,\omega_{\rm WP}
&=
c_{1}
\left(
\Theta_{{\rm Def}(S,\theta)},\omega_{\Lambda}
\right)
+
c_{1}
\left(
\Theta_{{\rm Def}(T)},\omega_{\rm hyp}
\right)
\\
&\quad
+
J_{M}^{*}c_{1}({\frak F}^{\oplus4},2h_{{\frak F}^{\oplus4}})
-
{\rm rk}({\frak F}^{\oplus4})c_{1}({\frak f}_{*}K_{{\frak X}/{\rm Def}(\widetilde{X})},h_{L^{2}})
\\
&=
c_{1}(\Omega_{\Lambda},\omega_{\Lambda})
+
c_{1}({\frak H},\omega_{\frak H})
+
4J_{M}^{*}\omega_{{\frak S}_{g}}
-
4g\,\mu^{*}\omega_{\rm WP}
\\
&=
-(\dim\Omega_{\Lambda})\cdot\omega_{\Lambda}
-
2\omega_{\frak H}
+
4J_{M}^{*}\omega_{{\frak S}_{g}}
-4g(\omega_{\Lambda}+\omega_{\frak H})
\\
&=
-(\dim\Omega_{\Lambda}+4g)\,\omega_{\Lambda}
-
(4g+2)\omega_{\frak H}
+
4J_{M}^{*}\omega_{{\frak S}_{g}}.
\end{aligned}
\end{equation}
To get the third equality, we used the Einstein property of the bounded symmetric domains equipped with the Bergman metric
$$
c_{1}(\Omega_{\Lambda},\omega_{\Lambda})=-(\dim\Omega_{\Lambda})\cdot\omega_{\Lambda}=-(20-r)\omega_{\Lambda},
\qquad
c_{1}({\frak H},\omega_{\frak H})=-2\omega_{\frak H}.
$$
Since
$$
h^{2,1}(\widetilde{X})-4g=21-r,
\qquad
\frac{\chi(\widetilde{X})}{12}=r-10
$$
by \cite[Cor.\,1.8]{Voisin93}, we deduce from \eqref{eqn:Ricci:form:Weil:Petersson:Borcea:Voisin} and 
Theorem~\ref{thm:curvature:BCOV} that
\begin{equation}
\label{eqn:curvature:BCOV:invariant:Borcea:Voisin:2}
-dd^{c}\log\widetilde{\tau}_{\rm BCOV}
=
(r-6)\,\omega_{\Lambda}+12\,\omega_{\frak H}+4\,J_{M}^{*}\omega_{{\frak S}_{g}}.
\end{equation}
This completes the proof of \eqref{eqn:curvature:BCOV:invariant:Borcea:Voisin}.
\end{pf}

Set 
$$
F^{\Lambda}
:=
\log
\left[
\tau_{\rm BCOV}/(\tau_{M}^{-4}\|\eta^{24}\|^{2})
\right]
\in
C^{\omega}\left({\mathcal M}_{\Lambda}^{0}\times{\bf X}(1)\right).
$$
Then Theorem~\ref{conj:BCOV=orb:BCOV} is equivalent to the assertion that $F^{\Lambda}$ is a constant function on 
${\mathcal M}_{\Lambda}^{0}\times{\bf X}(1)$. In the rest of Section~\ref{sect:3.1}, we study the possible singularities of $F^{\Lambda}$.
\par
Let $\varpi\colon{\frak H}\to{\bf X}(1)$ be the projection and let ${\frak P}_{\Lambda}:=\varPi_{\Lambda}\times\varpi$ be the projection 
from $\Omega_{\Lambda}\times{\frak H}$ to ${\mathcal M}_{\Lambda}\times{\bf X}(1)$.

\begin{proposition}
\label{prop:pluriharmonicity:F:M:perp}
For any $d\in\Delta_{\Lambda}$, there exists $\alpha(d)\in{\bf Q}$ such that
\begin{equation}
\label{eqn:curvature:BCOV:invariant:Borcea:Voisin:2}
-dd^{c}
\left[
{\frak P}_{\Lambda}^{*}F^{\Lambda}
\right]
=
\sum_{d\in\Delta_{\Lambda}/\pm1}\alpha(d)\,\delta_{H_{d}\times{\frak H}}.
\end{equation}
Here $\alpha(g\cdot d)=\alpha(d)$ for all $g\in O(\Lambda)$.
In particular, $\partial F^{\Lambda}$ is a logarithmic $1$-form on ${\mathcal M}_{\Lambda}\times{\bf X}(1)$ with possible pole along
$\overline{\mathcal D}_{\Lambda}\times{\bf X}(1)$.
\end{proposition}

\begin{pf}
Let $d\in\Delta_{\Lambda}$ and let $z\in{\frak H}$. 
Let $\gamma\colon\varDelta\to\Omega_{\Lambda}$ be a holomorphic curve intersecting $H_{d}$ transversally at $\gamma(0)\in H_{d}^{0}$.
By Theorem~\ref{thm:log:divergence:BCOV} and \cite[Prop.\,5.5]{Yoshikawa09}, there exists $a_{\gamma,d,z}\in{\bf Q}$ such that
\begin{equation}
\label{eqn:singularity:BCOV:inv:discr}
\log\tau_{\rm BCOV}\left({\frak P}_{\Lambda}(\gamma(t),z)\right)=a_{\gamma,d,z}\log|t|^{2}+O\left(\log(-\log|t|)\right)
\qquad
(t\to0).
\end{equation}
By \cite[Th.\,6.5]{Yoshikawa04}, we get
\begin{equation}
\label{eqn:singularity:orb:BCOV:inv:discr}
\log\left[\tau_{M}(\gamma(t))^{-4}\|\eta(z)^{24}\|^{2}\right]=\frac{1}{2}\log|t|^{2}+O(1)
\qquad
(t\to0).
\end{equation}
By \eqref{eqn:singularity:BCOV:inv:discr}, \eqref{eqn:singularity:orb:BCOV:inv:discr}, we get
\begin{equation}
\label{eqn:singularity:orb:F:discr}
F^{\Lambda}\left({\frak P}_{\Lambda}(\gamma(t),z)\right)=(a_{\gamma,d,z}-\frac{1}{2})\,\log|t|^{2}+O\left(\log(-\log|t|)\right)
\qquad
(t\to0).
\end{equation}
Since ${\frak P}_{\Lambda}^{*}F^{\Lambda}$ is pluriharmonic on $\Omega_{\Lambda}^{0}\times{\frak H}$, 
the constant $a_{\gamma,d,z}-\frac{1}{2}$ depends only on $d\in\Delta_{\Lambda}$ by \cite[Lemma 5.9]{Yoshikawa09}. 
Hence we can define $\alpha(d):=-(a_{\gamma,d,z}-\frac{1}{2})$. 
Then \eqref{eqn:curvature:BCOV:invariant:Borcea:Voisin:2} follows from \eqref{eqn:singularity:orb:F:discr} and \cite[Lemma 5.9]{Yoshikawa09}.
The property $\alpha(g\cdot d)=\alpha(d)$ for all $g\in O(\Lambda)$ is a consequence of the $O(\Lambda)$-invariance of 
${\frak P}_{\Lambda}^{*}F^{\Lambda}$.
\end{pf}

\par
Recall that ${\mathcal M}_{\Lambda}^{*}$ is the Baily-Borel compactification of ${\mathcal M}_{\Lambda}$ and 
${\mathcal B}_{\Lambda}={\mathcal M}_{\Lambda}^{*}\setminus{\mathcal M}_{\Lambda}$ is its boundary locus.

\begin{lemma}
\label{thm:BCOV=orb:BCOV:U(2)+U(2)}
Let $C\subset{\mathcal M}_{\Lambda}^{*}$ be an irreducible curve such that $C\not\subset{\mathcal B}_{\Lambda}$.
For any $z\in{\bf X}(1)$, $F^{\Lambda}|_{C\times\{z\}}$ has at most logarithmic singularities at $(C\cap{\mathcal B}_{\Lambda})\times\{z\}$.
Namely, for any $x\in(C\cap{\mathcal B}_{\Lambda})\times\{z\}$, there exists $\alpha\in{\bf Q}$ such that
$$
F^{\Lambda}(t,z)=\alpha\,\log|t|^{2}+O\left(\log(-\log|t|)\right)
\qquad
(t\to0),
$$
where $t$ is a local parameter of $C$ centered at $x$.
\end{lemma}

\begin{pf}
By \cite[Th.\,3.1]{Yoshikawa12}, there is an irreducible curve $B$, a finite surjective map $\varphi\colon B\to C$,
a smooth projective threefold ${\mathcal S}$, a surjective holomorphic map $f\colon{\mathcal S}\to B$,
and a holomorphic involution $\theta\colon{\mathcal S}\to{\mathcal S}$ preserving the fibers of $f$ with the following properties:
\begin{itemize}
\item[(i)]
There is a non-empty Zariski open subset $B^{0}\subset B$ such that
$(S_{b},\theta_{b})$ is a $2$-elementary $K3$ surface of type $M$ for all $b\in B^{0}$.
Here, $S_{b}:=f^{-1}(b)$ and $\theta_{b}:=\theta|_{S_{b}}$.
\item[(ii)]
$\varphi\colon B^{0}\to C\cap{\mathcal M}_{\Lambda}^{0}$ is the period map for the family of $2$-elementary $K3$ surfaces
$f\colon({\mathcal S},\theta)|_{B^{0}}\to B^{0}$ of type $M$.
\end{itemize}
\par
Let $T$ be an elliptic curve.
Let ${\mathcal X}\to({\mathcal S}\times T)/(\theta\times-1_{T})$ be a resolution and let $g\colon{\mathcal X}\to B$ be the map induced from 
the map $f\circ{\rm pr}_{1}\colon{\mathcal S}\times T\to B$. 
We may assume that $X_{b}:=g^{-1}(b)$ is the Borcea-Voisin threefold $\widetilde{X}_{(S_{b},\theta_{b},T)}$ for all $b\in B^{0}$. 
\par
Let $p\in\varphi^{-1}(C\cap{\mathcal B}_{\Lambda})$
and let $(U,s)$ be a coordinate neighborhood of $B$ centered at $p$.
By Theorem~\ref{thm:log:divergence:BCOV} applied to the family of Calabi-Yau threefolds $g\colon{\mathcal X}\to B$, 
there exists $\beta\in{\bf Q}$ such that at $s\to0$,
\begin{equation}
\label{eqn:behavior:BCOV:booundary:locus}
\varphi^{*}(\log\tau_{\rm BCOV})(s)
=
\log\tau_{\rm BCOV}(X_{s})
=
\beta\,\log|s|^{2}+O(\log(-\log|s|)).
\end{equation}
By \eqref{eqn:curvature:tau:M} and \eqref{eqn:behavior:BCOV:booundary:locus}, 
$F^{\Lambda}|_{C\times\{z\}}$ has at most logarithmic singularities at $(C\cap{\mathcal B}_{\Lambda})\times\{z\}$.
This completes the proof.
\end{pf}

Set $+i\infty:={\bf X}^{*}(1)\setminus{\bf X}(1)$.

\begin{lemma}
\label{lemma:singularity:BCOV:invariant:boundary:X(1)}
For any $[\eta]\in{\mathcal M}_{\Lambda}^{0}$,
$F^{\Lambda}|_{[\eta]\times{\bf X}(1)}$ has at most a logarithmic singularity at $([\eta],+i\infty)$.
\end{lemma}

\begin{pf}
The result follows from \cite[Prop.\,5.6]{Yoshikawa09} and Theorems~\ref{thm:Kronecker:limit:formula} and \ref{thm:Harvey:Moore} below.
\end{pf}

\subsection
{Examples verifying Theorem~\ref{conj:BCOV=orb:BCOV}}
\label{sect:3.2}
\par

\begin{lemma}
\label{lemma:Grauert:Remmert:extension}
Let $F$ be a real-valued pluriharmonic function on ${\mathcal M}_{\Lambda}\times{\bf X}(1)$.
Assume that for any $[\eta]\in{\mathcal M}_{\Lambda}^{0}$, $F|_{[\eta]\times{\bf X}(1)}$ has at most logarithmic singularity at $([\eta],+i\infty)$.
\newline{\rm (1) }
If $r\leq17$, then $F$ is a constant function.
\newline{\rm (2) }
If $r\geq18$, assume moreover that for every complete irreducible curve $C\subset{\mathcal M}_{\Lambda}^{*}$ with $C\not\subset{\mathcal B}_{\Lambda}$ 
and for every $z\in{\bf X}(1)$, $F|_{C\times\{z\}}$ has at most logarithmic singularity at any point of $(C\cap{\mathcal B}_{\Lambda})\times\{z\}$. 
Namely, there exists $\alpha\in{\bf Q}$ such that
$$
F([\eta],j)=\alpha\log|j|^{2}+O\left(\log(-\log|j|)\right)
\qquad
(j\to+i\infty),
$$
where $j\colon{\bf X}^{*}(1)\cong{\bf P}^{1}$ is the isomorphism given by the $j$-invariant.
Then $F$ is a constant.
\end{lemma}

\begin{pf}
{\bf (1) }
Since $r\leq17$, ${\mathcal B}_{\Lambda}$ has codimension $\geq2$ in ${\mathcal M}_{\Lambda}^{*}$.
By Grauert-Remmert \cite[Satz 4]{GrauertRemmert65} and the normality of ${\mathcal M}_{\Lambda}^{*}$, 
$F$ extends to a pluriharmonic function on ${\mathcal M}_{\Lambda}^{*}\times{\bf X}(1)$.
By the compactness of ${\mathcal M}_{\Lambda}^{*}$, there exists a function $\psi$ on ${\bf X}(1)$ such that $F={\rm pr}_{2}^{*}\psi$.
Since $F$ is pluriharmonic, $\psi$ must be a harmonic function on ${\bf X}(1)$ because $\psi=F|_{[\eta]\times{\bf X}(1)}$, where $[\eta]\in{\mathcal M}_{\Lambda}^{0}$.
Since $F|_{[\eta]\times{\bf X}(1)}$ has at most logarithmic singularity at $([\eta],+i\infty)$,
$\psi$ has at most logarithmic singularity at $+i\infty$.
Namely, there exists $\alpha\in{\bf R}$ such that
\begin{equation}
\label{eqn:asymptotic:psi}
\psi(j)=\alpha\,\log|j|^{2}+O(\log\log|j|)
\qquad
(j\to\infty).
\end{equation}
Since $\psi$ is harmonic on ${\bf X}(1)$, it follows from \eqref{eqn:asymptotic:psi} that $\partial\psi$ is a logarithmic $1$-form on ${\bf X}^{*}(1)\cong{\bf P}^{1}$
with possible pole at $j=\infty$ and with residue $\alpha$. By the residue theorem applied to $\partial\psi$, we get $\alpha=0$.
Hence $\partial\psi$ is a holomorphic $1$-form on $X^{*}(1)$. As a result, $\partial\psi=0$ on ${\bf X}^{*}(1)$, so that $\psi$ is a constant function on ${\bf X}(1)$.
This proves that $F={\rm pr}_{2}^{*}\psi$ is a constant function on ${\mathcal M}_{\Lambda}\times{\bf X}(1)$.
\par{\bf (2) }
Let $z\in{\bf X}(1)$. By assumption and \cite[Lemma 5.9]{Yoshikawa09}, 
$F|_{{\mathcal M}_{\Lambda}^{*}\times\{z\}}$ has at most logarithmic singularity along ${\mathcal B}_{\Lambda}$.
Hence $\partial F|_{{\mathcal M}_{\Lambda}^{*}\times\{z\}}$ is a logarithmic $1$-form on ${\mathcal M}_{\Lambda}^{*}\times\{z\}$
with possible pole along the {\em irreducible} divisor ${\mathcal B}_{\Lambda}$.
We set 
$$
\alpha:={\rm Res}_{{\mathcal B}_{\Lambda}\times\{z\}}(\partial F|_{{\mathcal M}_{\Lambda}^{*}\times\{z\}}).
$$
Let $C\subset{\mathcal M}_{\Lambda}^{*}$ be a complete irreducible curve with $C\not={\mathcal B}_{\Lambda}$ and 
$C\cap{\mathcal B}_{\Lambda}\not=\emptyset$. 
By the residue theorem applied to the logarithmic $1$-form $\partial F|_{C\times\{z\}}$ on $C\times\{z\}$, we get $\alpha\cdot\#(C\cap{\mathcal B}_{\Lambda})=0$.
Hence $\alpha=0$, so that $F|_{{\mathcal M}_{\Lambda}^{*}\times\{z\}}$ is a pluriharmonic function 
on ${\mathcal M}_{\Lambda}^{*}\setminus{\rm Sing}\,{\mathcal M}_{\Lambda}^{*}$.
By \cite[Satz 4]{GrauertRemmert65} and the normality of ${\mathcal M}_{\Lambda}^{*}$, 
$F|_{{\mathcal M}_{\Lambda}^{*}\times\{z\}}$ extends to a pluriharmonic function on ${\mathcal M}_{\Lambda}^{*}$.
By the compactness of ${\mathcal M}_{\Lambda}^{*}$, $F|_{{\mathcal M}_{\Lambda}^{*}\times\{z\}}$ is a constant function on
${\mathcal M}_{\Lambda}^{*}$.
This implies the existence of a harmonic function $\psi$ on ${\bf X}(1)$ such that $F=({\rm pr}_{2})^{*}\psi$.
By the same argument as in (1), $\psi$ is a constant function on ${\bf X}(1)$.
This proves that $F={\rm pr}_{2}^{*}\psi$ is a constant function on ${\mathcal M}_{\Lambda}\times{\bf X}(1)$.
This completes the proof.
\end{pf}

\begin{theorem}
\label{thm:BCOV=orb:BCOV:irreducible:discriminant}
If $r\leq17$ and $\overline{\mathcal D}_{\Lambda}$ is irreducible, then Theorem~\ref{conj:BCOV=orb:BCOV} holds.
\end{theorem}

\begin{pf}
Let $\alpha\in{\bf Q}$ be the residue of the logarithmic $1$-form $\partial F^{\Lambda}$ along $\overline{\mathcal D}_{\Lambda}\times{\bf X}(1)$.
Since $r\leq17$ and hence ${\mathcal B}_{\Lambda}$ has codimension $\geq2$ in ${\mathcal M}_{\Lambda}^{*}$, 
there is an irreducible complete curve $C\subset {\mathcal M}_{\Lambda}^{*}$ such that $C\cap{\mathcal B}_{\Lambda}=\emptyset$ and 
$C\cap\overline{\mathcal D}_{\Lambda}\not=\emptyset$. 
Let $z\in{\bf X}(1)$ be an arbitrary point.
Since $\partial F^{\Lambda}|_{C\times\{z\}}$ is a logarithmic $1$-form on $C\times\{z\}$ with residue $\alpha$
at any pole of $\partial F^{\Lambda}|_{C\times\{z\}}$, the total residue of $\partial F^{\Lambda}|_{C\times\{z\}}$
is a non-zero multiple of $\alpha$. By the residue theorem, we get $\alpha=0$.
Thus $\partial F^{\Lambda}$ is a holomorphic $1$-form on ${\mathcal M}_{\Lambda}\times{\bf X}(1)$,
so that $F^{\Lambda}$ is a pluriharmonic function on ${\mathcal M}_{\Lambda}\times{\bf X}(1)$. 
Now, the result follows from Lemma~\ref{lemma:Grauert:Remmert:extension} (1).
\end{pf}

\subsection
{The behavior of BCOV invariants near the discriminant locus}
\label{sect:3.3}

\subsubsection
{Ordinary singular families of $2$-elementary $K3$ surfaces}
\label{sect:3.3.1}
\par
Let ${\mathcal S}$ be a smooth complex threefold and let $p\colon{\mathcal S}\to\varDelta$ be a proper surjective holomorphic function 
without critical points on ${\mathcal S}\setminus p^{-1}(0)$.
Let $\theta\colon{\mathcal S}\to{\mathcal S}$ be a holomorphic involution preserving the fibers of $p$.
Set $S_{t}=p^{-1}(t)$ and $\theta_{t}=\theta|_{S_{t}}$ for $t\in\varDelta$.
Then $p\colon({\mathcal S},\theta)\to\varDelta$ is called an {\it ordinary singular family}
of $2$-elementary $K3$ surfaces if $p$ has a unique, non-degenerate critical point on $S_{0}$ 
and if $(S_{t},\theta_{t})$ is a $2$-elementary $K3$ surface for all $t\in\varDelta^{*}$. 
Let $o\in{\mathcal S}$ be the unique critical point of $p$.
By \cite[Sect.\,2.2]{Yoshikawa04}, there exists a system of coordinates $({\mathcal U},(w_{1},w_{2},w_{3}))$ centered at $o$ such that
$$
\iota(w)=(-w_{1},-w_{2},-w_{3})
\quad\hbox{ or }\quad
(w_{1},w_{2},-w_{3}),
\qquad
z\in{\mathcal U}.
$$
If $\iota(w)=(-w_{1},-w_{2},-w_{3})$ on ${\mathcal U}$, $\iota$ is said to be of {\it type $(0,3)$}. 
If $\iota(w)=(w_{1},w_{2},-w_{3})$ on ${\mathcal U}$, $\iota$ is said to be of {\it type $(2,1)$}.

\begin{lemma}
\label{lemma:normalized:coodinates}
There exists a system of local coordinates $(z_{1},z_{2},z_{3})$ of ${\mathcal Z}$ centered at $o\in{\mathcal Z}$ and a coordinate $t$ of $\varDelta$
centered at $0\in\varDelta$ with the following properties.
\begin{itemize}
\item[(1)]
If $\iota$ is of type $(0,3)$, then
\begin{equation}
\label{eqn:normal:form:type(2,1)}
\iota(z_{1},z_{2},z_{3})=(-z_{1},-z_{2},-z_{3}),
\qquad
p(z_{1},z_{2},z_{3})=(z_{1})^{2}+(z_{2})^{2}+(z_{3})^{2}.
\end{equation}
\item[(2)]
If $\iota$ is of type $(2,1)$, then
\begin{equation}
\label{eqn:normal:form:type(0,3)}
\iota(z_{1},z_{2},z_{3})=(z_{1},z_{2},-z_{3}),
\qquad
p(z_{1},z_{2},z_{3})=(z_{1})^{2}+(z_{2})^{2}+(z_{3})^{2}.
\end{equation}
\end{itemize}
\end{lemma}

\begin{pf}
The proof is standard and is omitted.
\end{pf}

\subsubsection
{Two local models of critical points}
\label{sect:3.3.2}
\par
We introduce two local models of critical points appearing in certain degenerations of Borcea-Voisin threefolds.
\par
Let ${\bf B}\subset{\bf C}^{3}$ be the unit ball of radius $1$. Let $T$ be an elliptic curve.
Define involutions  $\iota^{(2,2)}$, $\iota^{(0,4)}$ on ${\bf B}\times T$ by
$$
\iota^{(2,2)}(z,w)=(z_{1},z_{2},-z_{3},-w),
\qquad
\iota^{(0,4)}(z,w)=(-z_{1},-z_{2},-z_{3},-w),
$$
where $z=(z_{1},z_{2},z_{3})\in{\bf B}$ and $w\in T$.
Set 
$$
{\mathcal V}^{(2,2)}:=({\bf B}\times T)/\iota^{(2,2)},
\qquad
{\mathcal V}^{(0,4)}:=({\bf B}\times T)/\iota^{(0,4)}.
$$
Then ${\mathcal V}^{(2,2)}$ and ${\mathcal V}^{(0,4)}$ are orbifolds. 
Since the nowhere vanishing canonical form $dz_{1}\wedge dz_{2}\wedge dz_{3}\wedge dw$ on ${\bf B}\times T$ is invariant under the $\iota^{(2,2)}$
and $\iota^{(0,4)}$-actions, it descends to a nowhere vanishing canonical form in the sense of orbifolds on ${\mathcal V}^{(2,2)}$ and ${\mathcal V}^{(0,4)}$,
respectively.
We write $\varXi^{(2,2)}$ (resp. $\varXi^{(0,4)}$) for the nowhere vanishing canonical form on ${\mathcal V}^{(2,2)}$ (resp. ${\mathcal V}^{(0,4)}$)
induced by $dz_{1}\wedge dz_{2}\wedge dz_{3}\wedge dw$.
\par
For $(z,w)\in{\bf B}\times T$, write $[(z,w)]^{(2,2)}\in{\mathcal V}^{(2,2)}$ and $[(z,w)]^{(0,4)}\in{\mathcal V}^{(0,4)}$ for the images of $(z,w)$
by the projections ${\bf B}\times T\to{\mathcal V}^{(2,2)}$ and ${\bf B}\times T\to{\mathcal V}^{(0,4)}$, respectively.
\newline
\par{\em (Case 1) }
Set $\Sigma:=\{[(z,w)]^{(2,2)}\in{\mathcal V}^{(2,2)};\,z_{1}=z_{2}=0\}$. Then ${\rm Sing}\,{\mathcal V}^{(2,2)}=\Sigma$.
Let $\sigma^{(2,2)}\colon\widetilde{\mathcal V}^{(2,2)}\to{\mathcal V}^{(2,2)}$ be the blowing-up along $\Sigma$, 
which is a resolution of the singularities of ${\mathcal V}^{(2,2)}$.
Define $F^{(2,2)}\in{\mathcal O}({\mathcal V}^{(2,2)})$ by
$$
F^{(2,2)}([(z,w)]^{(2,2)}):=(z_{1})^{2}+(z_{2})^{2}+(z_{3})^{2}
$$
and set 
$$
\widetilde{F}^{(2,2)}:=F^{(2,2)}\circ\sigma^{(2,2)}\in{\mathcal O}(\widetilde{\mathcal V}^{(2,2)}).
$$
\par
Set $W:=\{(u,v,r)\in{\bf C}^{3};\,uv-r^{2}=0\}$. 
Since ${\bf C}^{2}/\pm1\cong W$ via the map $\pm(z_{3},w)\mapsto((z_{3})^{2},w^{2},z_{3}w)$, 
we have an isomorphism of germs $({\mathcal V}^{(2,2)},x)\cong({\bf C}^{2}\times W,(0,0))$ for any $x\in{\rm Sing}\,{\mathcal V}^{(2,2)}$. 
Under this identification of germs, the function germ $F^{(2,2)}\colon({\bf C}^{4}/\iota^{(2,2)},0)=({\bf C}^{2}\times W,(0,0))\to({\bf C},0)$ is expressed as
\begin{equation}
\label{eqn:normal:form:projection}
F^{(2,2)}(z_{1},z_{2},u,v,r)=(z_{1})^{2}+(z_{2})^{2}+u,
\end{equation}
where $(z_{1},z_{2})\in{\bf C}^{2}$, $(u,v,r)\in W$. 
\par
Let $\sigma\colon(\widetilde{W},E)\to(W,0)$ be the blowing-up at the origin, where $E=\sigma^{-1}(0)\cong{\bf P}^{1}$.  
The isomorphism $({\mathcal V}^{(2,2)},x)\cong({\bf C}^{2}\times W,(0,0))$
induces an isomorphism of germs $(\widetilde{\mathcal V}^{(2,2)},\widetilde{x})\cong({\bf C}^{2}\times\widetilde{W},(0,\zeta))$ for any
$\widetilde{x}\in(\sigma^{(2,2)})^{-1}(x)$, where $\zeta\in\sigma^{-1}(0)$ is the point corresponding to $\widetilde{x}$. 
Let $\{{\mathcal U}_{0},{\mathcal U}_{1},{\mathcal U}_{2}\}$ be the open covering of ${\bf C}^{2}\times\widetilde{W}$ defined as
$$
{\mathcal U}_{0}:={\bf C}^{2}\times\sigma^{-1}(\{u\not=0\}),
\quad
{\mathcal U}_{1}:={\bf C}^{2}\times\sigma^{-1}(\{v\not=0\}),
\quad
{\mathcal U}_{2}:={\bf C}^{2}\times\sigma^{-1}(\{r\not=0\}).
$$
By \eqref{eqn:normal:form:projection}, $\widetilde{F}^{(2,2)}$ has no critical points on ${\mathcal U}_{0}\cup{\mathcal U}_{2}$.
On ${\mathcal U}_{1}$, we have the system of coordinates $(z_{1},z_{2},v,s:=r/v)$. Since
$$
\widetilde{F}^{(2,2)}(z_{1},z_{2},v,s)=(z_{1})^{2}+(z_{2})^{2}+u=(z_{1})^{2}+(z_{2})^{2}+\frac{r^{2}}{v}=(z_{1})^{2}+(z_{2})^{2}+vs^{2},
$$
we get 
$$
\varSigma_{\widetilde{F}^{(2,2)}}=\{z_{1}=z_{2}=s=0\}\cap{\mathcal U}_{1}.
$$
Hence $\dim\varSigma_{\widetilde{F}^{(2,2)}}=1$. In particular, the divisor $(\widetilde{F}^{(2,2)})^{-1}(0)$ is irreducible.
\par
Since the dualizing sheaf of $\widetilde{W}$ is trivial, so is the dualizing sheaf of ${\bf C}^{2}\times\widetilde{W}$.
Since $\varXi^{(2,2)}$ is a nowhere vanishing section of the dualizing sheaf of ${\bf C}^{2}\times\widetilde{W}$, 
$(\sigma^{(2,2)})^{*}\varXi^{(2,2)}$ extends to a nowhere vanishing canonical form on $\widetilde{\mathcal V}^{(2,2)}$.
If $\Upsilon$ is a nowhere vanishing canonical form on ${\mathcal V}^{(2,2)}\setminus{\rm Sing}\,{\mathcal V}^{(2,2)}$,
then $\Upsilon/\varXi^{(2,2)}$ is a nowhere vanishing holomorphic function on ${\mathcal V}^{(2,2)}\setminus{\rm Sing}\,{\mathcal V}^{(2,2)}$,
which extends to a nowhere vanishing holomorphic function on ${\mathcal V}^{(2,2)}$. Hence $(\sigma^{(2,2)})^{*}\Upsilon$ is also a nowhere vanishing
canonical form on $\widetilde{\mathcal V}^{(2,2)}$.
\newline
\par{\em (Case 2) }
Let $\omega_{1},\omega_{2},\omega_{3}\in T$ be non-zero points of order $2$. Then
$$
{\rm Sing}\,{\mathcal V}^{(0,4)}=\{[(0,0)]^{(0,4)},\,[(0,\omega_{1})]^{(0,4)},\,[(0,\omega_{2})]^{(0,4)},\,[(0,\omega_{3})]^{(0,4)}\}
$$
consists of $4$ isolated quotient singularities isomorphic to $({\bf C}^{4}/\pm1,0)$.
\par
Define $F^{(0,4)}\in{\mathcal O}({\mathcal V}^{(0,4)})$ by
$$
F^{(0,4)}([(z,w)]^{(0,4)}):=(z_{1})^{2}+(z_{2})^{2}+(z_{3})^{2}
$$
By this expression, we have
$$
\varSigma_{F^{(0,4)}}=\{[(0,w)]^{(0,4)}\in{\mathcal V}^{(0,4)};\,w\in T\}.
$$
In particular, $\dim\varSigma_{F^{(0,4)}}=1$.
Since the inverse image of $(F^{(0,4)})^{-1}(0)$ in ${\bf B}\times T$ is irreducible, $(F^{(0,4)})^{-1}(0)$ is an irreducible divisor of ${\mathcal V}^{(0,4)}$.

\subsubsection
{A degenerating family of Borcea-Voisin threefolds and BCOV invariants, I}
\label{sect:3.3.3}
\par
Let ${\mathcal S}$ and ${\mathcal S}'$ be smooth irreducible projective threefolds. 
Let $\theta\colon{\mathcal S}\to{\mathcal S}$ and $\theta'\colon{\mathcal S}'\to{\mathcal S}'$ be holomorphic involutions on ${\mathcal S}$ and ${\mathcal S}'$, 
respectively.
Let $B$ and $B'$ be compact Riemann surfaces and let $p\colon{\mathcal S}\to B$ and $p'\colon{\mathcal S}'\to B'$ be surjective holomorphic maps. 
Let $\Delta\subset B$ and $\Delta'\subset B'$ be the discriminant loci of $p\colon{\mathcal S}\to B$ and $p'\colon{\mathcal S}'\to B'$, respectively.
Let ${\frak p}\in\Delta$ and ${\frak p}'\in\Delta'$. For $b\in B$ and $b'\in B'$, set $S_{b}:=p^{-1}(b)$ and $S'_{b'}:=(p')^{-1}(b')$.
We assume the following:
\begin{itemize}
\item[(1)]
$\theta$ and $\theta'$ preserve the fibers of $p$ and $p'$, respectively. Set $\theta_{b}:=\theta|_{S_{b}}$ and $\theta'_{b'}:=\theta'|_{S'_{b'}}$ for 
$b\in B$ and $b'\in B'$.
\item[(2)]
There exist primitive $2$-elementary Lorentzian sublattices $M,M'\subset{\Bbb L}_{K3}$ such that
$(S_{b},\theta_{b})$ is a $2$-elementary $K3$ surface of type $M$ for all $B\setminus\Delta$
and $(S'_{b'},\theta'_{b'})$ is a $2$-elementary $K3$ surface of type $M'$ for all $B'\setminus\Delta'$.
\item[(3)]
There is a neighborhood $U$ of ${\frak p}$ in $B$ such that 
$$
p\colon(p^{-1}(U),\theta|_{p^{-1}(U)})\to U
$$ 
is an ordinary singular family of $2$-elementary $K3$ surfaces of type $M$.
Similarly, there is a neighborhood $U'$ of ${\frak p}'$ in $B'$ such that 
$$
p'\colon((p')^{-1}(U'),\theta'|_{(p')^{-1}(U')})\to U'
$$ 
is an ordinary singular family of $2$-elementary $K3$ surfaces of type $M'$.
\end{itemize}
\par
Let $T$ be an elliptic curve. We set 
$$
\iota:=\theta\times(-1)_{T},
\qquad
\iota':=\theta'\times(-1)_{T}
$$
and
$$
{\mathcal X}:=({\mathcal S}\times T)/\iota,
\qquad
{\mathcal X}':=({\mathcal S}'\times T)/\iota'.
$$
Let $\pi\colon{\mathcal X}\to B$ and $\pi'\colon{\mathcal X}'\to B'$ be the projections induced from the projections $p\colon{\mathcal S}\to B$
and $p'\colon{\mathcal S}'\to B'$, respectively. 
Since ${\mathcal S}\times T$ (resp. ${\mathcal S}'\times T$) is a complex manifold, the set of fixed points of the $\iota$-action 
(resp. $\iota'$-action) on ${\mathcal S}\times T$ (resp. ${\mathcal S}'\times T$), 
i.e., ${\mathcal S}^{\theta}\times T[2]$ (resp. $({\mathcal S}')^{\theta'}\times T[2]$) is the disjoint union of complex submanifolds.
\par
Let $Z:=({\mathcal S}^{\theta}\times T[2])^{\rm hol}$ and $Z':=(({\mathcal S}')^{\theta'}\times T[2])^{\rm hol}$ be the horizontal components.
Namely, $Z$ is the union of those connected components of ${\mathcal S}^{\theta}\times T[2]$ which are flat over $B$.
Similarly, $Z'$ is the union of those connected components of $({\mathcal S}')^{\theta'}\times T[2]$ which are flat over $B'$. 
Then $Z$ and $Z'$ are complex submanifolds of ${\mathcal S}\times T$ and ${\mathcal S}'\times T$ of codimension $2$, respectively.
\par
Let $\sigma\colon\widetilde{\mathcal X}\to{\mathcal X}$ be the blowing-up of ${\mathcal X}$ along $Z$ and
let $\sigma'\colon\widetilde{\mathcal X}'\to{\mathcal X}'$ be the blowing-up of ${\mathcal X}'$ along $Z'$. 
We set 
$$
\widetilde{\pi}:=\pi\circ\sigma\colon\widetilde{\mathcal X}\to B,
\qquad
\widetilde{\pi}':=\pi'\circ\sigma'\colon\widetilde{\mathcal X}'\to B'.
$$
By construction, $\widetilde{\pi}^{-1}(b)=\widetilde{X}_{(S_{b},\theta_{b},T)}$ for $b\in U\setminus\{\frak p\}$ and 
$(\widetilde{\pi}')^{-1}(b')=\widetilde{X}_{(S'_{b'},\theta'_{b'},T)}$ for $b'\in U'\setminus\{{\frak p}'\}$.
Fix isomorphisms of germs $(U,{\frak p})\cong({\bf C},0)$ and $(U',{\frak p}')\cong({\bf C},0)$.
Then, for $t\in{\bf C}$ with $0<|t|\ll1$, $\widetilde{X}_{(S_{t},\theta_{t},T)}$ is a Borcea-Voisin threefold of type $M$ and 
$\widetilde{X}_{(S'_{t},\theta'_{t},T)}$ is a Borcea-Voisin threefold of type $M'$.

\begin{theorem}
\label{thm:dependence:log:divergence}
If $\theta|_{p^{-1}(U)}$ and $\theta'|_{(p')^{-1}(U')}$ have the same type, then 
$$
\lim_{t\to0}\frac{\log\tau_{\rm BCOV}(\widetilde{X}_{(S_{t},\theta_{t},T)})}{\log|t|^{2}}
=
\lim_{t\to0}\frac{\log\tau_{\rm BCOV}(\widetilde{X}_{(S'_{t},\theta'_{t},T)})}{\log|t|^{2}}.
$$
In particular, as $t\to0$,
$$
\log\tau_{\rm BCOV}(\widetilde{X}_{(S_{t},\theta_{t},T)})
-
\log\tau_{\rm BCOV}(\widetilde{X}_{(S'_{t},\theta'_{t},T)})
=
O\left(\log(-\log|t|)\right).
$$
\end{theorem}

\begin{pf}
By Theorem~\ref{conjecture:locality:BCOV:invariant}, it suffices to verify conditions (A1), (A2), (A3), (A4) in Section~\ref{sect:1.3} 
for the families $\widetilde{\pi}\colon\widetilde{\mathcal X}\to B$ and $\widetilde{\pi}'\colon\widetilde{\mathcal X}'\to B'$.
Since $\widetilde{\pi}^{-1}(b)=\widetilde{X}_{(S_{b},\theta_{b},T)}$ for $b\in U\setminus\{\frak p\}$ and 
$(\widetilde{\pi}')^{-1}(b')=\widetilde{X}_{(S'_{b'},\theta'_{b'},T)}$ for $b'\in U'\setminus\{{\frak p}'\}$, condition (A2) holds.
\par
Set $\theta_{U}:=\theta|_{p^{-1}(U)}$ and $\theta'_{U'}:=\theta'|_{(p')^{-1}(U')}$.
Let $o$ be the unique critical point of $p|_{p^{-1}(U)}$ and let $o'$ be the unique critical point of $p'|_{(p')^{-1}(U')}$.
\par{\em (Step 1) }
By Lemma~\ref{lemma:normalized:coodinates}, there exist a neighborhood ${\mathcal U}$ of $o$ in ${\mathcal S}$
and a system of coordinates $(z_{1},z_{2},z_{3})$ on ${\mathcal U}$ centered at $o$ such that 
$({\mathcal U},(z_{1},z_{2},z_{3}))=({\bf B},(z_{1},z_{2},z_{3}))$ and 
such that either \eqref{eqn:normal:form:type(2,1)} or \eqref{eqn:normal:form:type(0,3)} holds on ${\mathcal U}$.
Similarly, there exist a neighborhood ${\mathcal U}'$ of $o'$ in ${\mathcal S}'$
and a system of coordinates $(z_{1},z_{2},z_{3})$ on ${\mathcal U}'$ centered at $o'$ such that 
$({\mathcal U}',(z_{1},z_{2},z_{3}))=({\bf B},(z_{1},z_{2},z_{3}))$ and 
such that either \eqref{eqn:normal:form:type(2,1)} or \eqref{eqn:normal:form:type(0,3)} holds on ${\mathcal U}'$.
Since $\theta_{U}$ and $\theta'_{U'}$ have the same type, we get 
\begin{equation}
\label{eqn:identification:nbhd:critical:locus}
({\mathcal U}\times T)/\iota
\cong
({\mathcal U}'\times T)/\iota'
\cong
\begin{cases}
\begin{array}{ll}
{\mathcal V}^{(2,2)}
&
\hbox{ if }\theta_{U},\,\theta'_{U'}\hbox{ are of type }(2,1),
\\
{\mathcal V}^{(0,4)}
&
\hbox{ if }\theta_{U},\,\theta'_{U'}\hbox{ are of type }(0,3).
\end{array}
\end{cases}
\end{equation}
Define open subsets $O\subset\widetilde{\mathcal X}$ and $O'\subset\widetilde{\mathcal X}'$ by
$$
O:=\sigma^{-1}\left(({\mathcal U}\times T)/\iota\right),
\qquad
O':=(\sigma')^{-1}\left(({\mathcal U}'\times T)/\iota'\right).
$$
By \eqref{eqn:identification:nbhd:critical:locus}, we get the following isomorphism
\begin{equation}
\label{eqn:identification:nbhd:critical:locus:2}
O
\cong
O'
\cong
\begin{cases}
\begin{array}{ll}
\widetilde{\mathcal V}^{(2,2)}
&
\hbox{ if }\theta_{U},\,\theta'_{U'}\hbox{ are of type }(2,1),
\\
{\mathcal V}^{(0,4)}
&
\hbox{ if }\theta_{U},\,\theta'_{U'}\hbox{ are of type }(0,3).
\end{array}
\end{cases}
\end{equation}
By \eqref{eqn:normal:form:type(2,1)}, \eqref{eqn:normal:form:type(0,3)}, we get under the isomorphism \eqref{eqn:identification:nbhd:critical:locus}
\begin{equation}
\label{eqn:identification:projection:near:critical:locus}
\pi|_{({\mathcal U}\times T)/\iota}
=
\pi'|_{({\mathcal U}'\times T)/\iota'}
=
\begin{cases}
\begin{array}{ll}
F^{(2,2)}
&
\hbox{ if }\theta_{U},\,\theta'_{U'}\hbox{ are of type }(2,1),
\\
F^{(0,4)}
&
\hbox{ if }\theta_{U},\,\theta'_{U'}\hbox{ are of type }(0,3).
\end{array}
\end{cases}
\end{equation}
Under the isomorphism \eqref{eqn:identification:nbhd:critical:locus:2}, we get by \eqref{eqn:identification:projection:near:critical:locus} 
the following isomorphism of pairs:
\begin{equation}
\label{eqn:identification:projection:near:critical:locus:2}
\left(
O,\widetilde{\pi}
\right)
\cong
\left(
O',\widetilde{\pi}'
\right)
\cong
\begin{cases}
\begin{array}{ll}
\left(
\widetilde{\mathcal V}^{(2,2)},
\widetilde{F}^{(2,2)}
\right)
&
\hbox{ if }\theta_{U},\,\theta'_{U'}\hbox{ are of type }(2,1),
\\
\left(
{\mathcal V}^{(0,4)},
F^{(0,4)}
\right)
&
\hbox{ if }\theta_{U},\,\theta'_{U'}\hbox{ are of type }(0,3).
\end{array}
\end{cases}
\end{equation}
This verifies condition (A4) in Section~\ref{sect:1.3}.  
\par{\em (Step 2) }
Since $S_{0}$ and $S'_{0}$ are singular $K3$ surfaces with a unique ordinary double point as its singular set,
the divisors $\widetilde{\pi}^{-1}(0)$ and $(\widetilde{\pi}')^{-1}(0)$ are irreducible by the descriptions of  
$(\widetilde{\mathcal V}^{(2,2)},\widetilde{F}^{(2,2)})$ and $(\widetilde{\mathcal V}^{(0,4)},\widetilde{F}^{(0,4)})$ in Section~\ref{sect:3.3.2}.
Similarly, by the descriptions of $(\widetilde{\mathcal V}^{(2,2)},\widetilde{F}^{(2,2)})$ and $(\widetilde{\mathcal V}^{(0,4)},\widetilde{F}^{(0,4)})$ 
in Section~\ref{sect:3.3.2}, we get 
\begin{equation}
\label{eqn:dimension:critical:loci}
\dim\varSigma_{\widetilde{\pi}}=\dim\varSigma_{\widetilde{\pi}'}=1.
\end{equation}
This verifies condition (A1) in Section~\ref{sect:1.3}.  
\par{\em (Step 3) }
By choosing $U$ and $U'$ sufficiently small, we may assume by \cite[Lemma 2.3]{Yoshikawa04} that $p^{-1}(U)$ and $(p')^{-1}(U')$ carry
nowhere vanishing canonical forms $\xi$ and $\xi'$, respectively, such that $\theta^{*}\xi=-\xi$ and $(\theta')^{*}\xi'=-\xi'$. 
Then $\xi\wedge dw$ and $\xi'\wedge dw$ are nowhere vanishing canonical forms on $p^{-1}(U)\times T$ and $(p')^{-1}(U')\times T$, respectively,
such that $\iota^{*}(\xi\wedge dw)=\xi\wedge dw$ and $(\iota')^{*}(\xi'\wedge dw)=\xi'\wedge dw$.
Hence $\xi\wedge dw$ (resp. $\xi'\wedge dw$) induces a nowhere vanishing canonical form $\Xi$ (resp. $\Xi'$) in the sense of orbifolds
on  $({\mathcal U}\times T)/\iota$ (resp. $({\mathcal U}'\times T)/\iota'$). 
\par
Let $\theta_{U}$ and $\theta'_{U'}$ be of type $(2,1)$.
Since every nowhere vanishing canonical form on $({\mathcal U}\times T)/\iota$ (resp. $({\mathcal U}'\times T)/\iota'$) lifts to a nowhere vanishing canonical form
on $O$ (resp. $O'$)  via $\sigma$ (resp. $\sigma'$) by Section~\ref{sect:3.3.2} Case 1, 
$\sigma^{*}\Xi$ (resp. $(\sigma')^{*}\Xi'$) is a nowhere vanishing canonical form on $O$ (resp. $O'$). 
\par
Let $\theta_{U}$ and $\theta'_{U'}$ be of type $(0,3)$. Then $O\cong O'\cong{\mathcal V}^{(0,4)}$ by \eqref{eqn:identification:nbhd:critical:locus:2}.
Hence $\Xi$ (resp. $\Xi'$) is a nowhere vanishing canonical form on $O\setminus\varSigma_{\widetilde{\pi}}$ 
(resp. $O'\setminus\varSigma_{\widetilde{\pi}'}$).
This verifies condition (A3) in Section~\ref{sect:1.3}.  
Since conditions (A1), (A2), (A3), (A4) are verified for $\widetilde{\pi}\colon\widetilde{\mathcal X}\to B$ and $\widetilde{\pi}'\colon\widetilde{\mathcal X}'\to B'$,
the result follows from Theorem~\ref{conjecture:locality:BCOV:invariant}.
\end{pf}

\subsubsection
{A degenerating family of Borcea-Voisin threefolds and BCOV invariants, II}
\label{sect:3.3.4}
\par
Let ${\mathcal S}$ be a smooth projective threefold equipped with a holomorphic involution $\theta\colon{\mathcal S}\to{\mathcal S}$.
Let $B$ be a compact Riemann surface and let $p\colon{\mathcal S}\to B$ be a surjective holomorphic map. 
Let $\Delta\subset B$ be the discriminant locus of $p\colon{\mathcal S}\to B$ and let ${\frak p}\in\Delta$.
We assume the following:
\begin{itemize}
\item[(1)]
$\theta$ preserves the fibers of $p$ and the pair $(S_{b},\theta_{b})$ is a $2$-elementary $K3$ surface of type $M$ for all $B\setminus\Delta$,
where $S_{b}:=p^{-1}(b)$ and $\theta_{b}:=\theta|_{S_{b}}$.
\item[(2)]
There is a neighborhood $U$ of ${\frak q}$ in $B$ such that 
$p\colon(p^{-1}(U),\theta|_{p^{-1}(U)})\to U$ is an ordinary singular family of $2$-elementary $K3$ surfaces of type $M$.
\end{itemize}

\begin{theorem}
\label{thm:singularity:BCOV:type(2,2)}
Let $T$ be an elliptic curve. If $\theta|_{p^{-1}(U)}$ is of type $(2,1)$, then
$$
\log\tau_{\rm BCOV}(\widetilde{X}_{(S_{t},\theta_{t},T)})=\frac{1}{2}\log|t|^{2}+O\left(\log(-\log|t|)\right)
\qquad
(t\to0).
$$
\end{theorem}

\begin{pf}
Set $M:={\Bbb U}\oplus{\Bbb E}_{8}(2)$. Then $\overline{\mathcal D}_{\Lambda}$ is irreducible by \cite[Prop.\,11.6]{Yoshikawa13}.
\par{\em (Step 1) }
By \cite[Th.\,2.8]{Yoshikawa04}, there exist a smooth projective threefold ${\mathcal S}'$ equipped with a holomorphic involution 
$\theta'\colon{\mathcal S}'\to{\mathcal S}'$, a pointed compact Riemann surface $(B',{\frak p}')$ equipped with a neighborhood $U'$ of ${\frak p}'$, 
and a surjective holomorphic map $p'\colon{\mathcal S}'\to B'$ with the following properties:
\begin{itemize}
\item[(1)]
$\theta'$ preserves the fibers of $p'$.
\item[(2)]
$p'\colon({\mathcal S}'_{U'},\theta'|_{{\mathcal S}'_{U'}})\to U'$ is an ordinary singular family of $2$-elementary $K3$ surfaces of type $M$.
\end{itemize}
Since $M={\Bbb U}\oplus{\Bbb E}_{8}(2)$, the fixed-point-set $(S'_{b'})^{\theta'_{b'}}$ consists of two disjoint elliptic curves for all $b'\in U'\setminus\{{\frak p}'\}$.
Assume that $\theta'_{U'}$ is of type $(0,3)$. Then the set of fixed points of the $\theta'_{{\frak p}'}$-action on $S'_{{\frak p}'}=(p')^{-1}({\frak p}')$ 
consists of two disjoint elliptic curves and the isolated point ${\rm Sing}\,S'_{{\frak p}'}$.
Let $\mu\colon\widetilde{S}'_{{\frak p}'}\to S'_{{\frak p}'}$ be the minimal resolution.
By \cite[Th.\,2.3 (1)]{Yoshikawa13}, the involution $\theta'_{{\frak p}'}$ on $S'_{{\frak p}'}$ lifts to an involution $\widetilde{\theta}'_{{\frak p}'}$ 
on $\widetilde{S}'_{{\frak p}'}$ and the pair $(\widetilde{S}'_{{\frak p}'},\widetilde{\theta}'_{{\frak p}'})$ is a $2$-elementary $K3$ surface.
Since $\theta'_{U'}$ is of type $(0,3)$, $\mu^{-1}({\rm Sing}\,S'_{{\frak p}'})\cong{\bf P}^{1}$ must be a component of the fixed-point-set
$(\widetilde{S}'_{{\frak p}'})^{\widetilde{\theta}'_{{\frak p}'}}$. Thus $(\widetilde{S}'_{{\frak p}'})^{\widetilde{\theta}'_{{\frak p}'}}$ consists of
two elliptic curves and a $(-2)$-curve. By Proposition~\ref{prop:topology:fixed:points}, this is impossible.
Hence $\theta'_{U'}$ is of type $(2,1)$.
(Even though ${\Bbb U}\oplus{\Bbb E}_{8}(2)$ is one of the exceptional lattices in the sense of \cite{Yoshikawa13}, 
the proof of \cite[Th.\,2.3 (2)]{Yoshikawa13} remains valid.
The fact that $\theta'_{U'}$ is of type $(2,1)$ also follows from \cite[Th.\,2.3 (2)]{Yoshikawa13}.)
\par{\em (Step 2) }
Since the elliptic curve $T$ is fixed, there exists by Theorem~\ref{thm:BCOV=orb:BCOV:irreducible:discriminant} a constant $C$ such that
\begin{equation}
\label{eqn:asymptotics:BCOV:family:associated:to:OSF}
\log\tau_{\rm BCOV}(\widetilde{X}_{(S'_{t},\theta'_{t},T)})
=
-4\log\tau_{M}(S'_{t},\theta'_{t})+C
\qquad
(\forall\,t\in U'\setminus\{{\frak p}'\}).
\end{equation}
Since $p'|_{U'}\colon({\mathcal S}'|_{U'},\theta')\to U'$ is an ordinary singular family of $2$-elementary $K3$ surfaces, we get by \cite[Th.\,6.5]{Yoshikawa04}
\begin{equation}
\label{eqn:asymptotics:equivariant:torsion:OSF:type(2,1)}
\log\tau_{M}(S'_{t},\theta'_{t})=-\frac{1}{8}\log|t|^{2}+O\left(\log(-\log|t|)\right)
\qquad
(t\to0).
\end{equation}
By \eqref{eqn:asymptotics:BCOV:family:associated:to:OSF}, \eqref{eqn:asymptotics:equivariant:torsion:OSF:type(2,1)}, we get
\begin{equation}
\label{eqn:asymptotics:BCOV:family:associated:to:OSF:2}
\log\tau_{\rm BCOV}(\widetilde{X}_{(S'_{t},\theta'_{t},T)})
=
\frac{1}{2}\log|t|^{2}+O\left(\log(-\log|t|)\right)
\qquad
(t\to0).
\end{equation}
Since $\theta'_{U'}$ is of type $(2,1)$, it follows from \eqref{eqn:asymptotics:BCOV:family:associated:to:OSF:2} and
Theorem~\ref{thm:dependence:log:divergence} that
$$
\begin{aligned}
\log\tau_{\rm BCOV}(\widetilde{X}_{(S_{t},\theta_{t},T)})
&=
\log\tau_{\rm BCOV}(\widetilde{X}_{(S'_{t},\theta'_{t},T)})+O\left(\log(-\log|t|)\right)
\\
&=
\frac{1}{2}\log|t|^{2}+O\left(\log(-\log|t|)\right)
\end{aligned}
$$
as $t\to0$. This completes the proof.
\end{pf}

\begin{theorem}
\label{thm:singularity:BCOV:type(0,4)}
Let $T$ be an elliptic curve. If $\theta|_{p^{-1}(U)}$ is of type $(0,3)$, then
$$
\log\tau_{\rm BCOV}(\widetilde{X}_{(S_{t},\theta_{t},T)})=\frac{1}{2}\log|t|^{2}+O\left(\log(-\log|t|)\right)
\qquad
(t\to0).
$$
\end{theorem}

\begin{pf}
Set $M:={\Bbb U}(2)\oplus{\Bbb E}_{8}(2)$. Then $\overline{\mathcal D}_{\Lambda}$ is irreducible by e.g. \cite[Prop.\,11.6]{Yoshikawa13}.
\par{\em (Step 1) }
By \cite[Th.\,2.8]{Yoshikawa04}, there exist a smooth projective threefold ${\mathcal S}'$ equipped with a holomorphic involution 
$\theta'\colon{\mathcal S}'\to{\mathcal S}'$, a pointed compact Riemann surface $(B',{\frak p}')$ equipped with a neighborhood $U'$ of ${\frak p}'$, 
and a surjective holomorphic map $p'\colon{\mathcal S}'\to B'$ with the following properties:
\begin{itemize}
\item[(1)]
$\theta'$ preserves the fibers of $p'$.
\item[(2)]
$p'\colon({\mathcal S}'_{U'},\theta'|_{{\mathcal S}'_{U'}})\to U'$ is an ordinary singular family of $2$-elementary $K3$ surfaces of type $M$.
\end{itemize}
Let $\Delta'$ be the discriminant locus of $p'\colon{\mathcal S}'\to B'$.
Since $M={\Bbb U}(2)\oplus{\Bbb E}_{8}(2)$, $\theta'_{b'}$ has no fixed points on $S'_{b'}$ for all $b'\in B'\setminus\Delta'$
by Proposition~\ref{prop:topology:fixed:points}.
Hence $({\mathcal S}')^{\theta'}$ has no horizontal components, which implies that $\theta'_{U'}:=\theta'|_{{\mathcal S}'_{U'}}$ is of type $(0,3)$.
(Since the lattice ${\Bbb U}(2)\oplus{\Bbb E}_{8}(2)$ is exceptional, \cite[Th.\,2.3 (2)]{Yoshikawa13} does not apply in this case.)
\par{\em (Step 2) }
Let $t$ be a local parameter of $B'$ centered at ${\frak p}'$ and set $S'_{t}:=(p')^{-1}(t)$ and $\theta'_{t}:=\theta'|_{S'_{t}}$.
Since the elliptic curve $T$ is fixed, there exist by Theorem~\ref{thm:BCOV=orb:BCOV:irreducible:discriminant} a constant $C$ such that 
\begin{equation}
\label{eqn:asymptotics:BCOV:family:associated:to:OSF:type(0,3)}
\log\tau_{\rm BCOV}(\widetilde{X}_{(S'_{t},\theta'_{t},T)})
=
-4\log\tau_{M}(S'_{t},\theta'_{t})+C
\qquad
(\forall\,t\in U'\setminus\{{\frak p}'\}).
\end{equation}
Since $p'|_{U'}\colon({\mathcal S}'|_{U'},\theta')\to U'$ is an ordinary singular family of $2$-elementary $K3$ surfaces, we get by \cite[Th.\,6.5]{Yoshikawa04}
\begin{equation}
\label{eqn:asymptotics:equivariant:torsion:OSF:type(0,3)}
\log\tau_{M}(S'_{t},\theta'_{t})=-\frac{1}{8}\log|t|^{2}+O\left(\log(-\log|t|)\right)
\qquad
(t\to0).
\end{equation}
By \eqref{eqn:asymptotics:BCOV:family:associated:to:OSF:type(0,3)}, \eqref{eqn:asymptotics:equivariant:torsion:OSF:type(0,3)}, we get
\begin{equation}
\label{eqn:asymptotics:BCOV:family:associated:to:OSF:3}
\log\tau_{\rm BCOV}(\widetilde{X}_{(S'_{t},\theta'_{t},T)})
=
\frac{1}{2}\log|t|^{2}+O\left(\log(-\log|t|)\right)
\qquad
(t\to0).
\end{equation}
Since $\theta'_{U'}$ is of type $(0,3)$, it follows from \eqref{eqn:asymptotics:BCOV:family:associated:to:OSF:3} and
Theorem~\ref{thm:dependence:log:divergence} that
$$
\begin{aligned}
\log\tau_{\rm BCOV}(\widetilde{X}_{(S_{t},\theta_{t},T)})
&=
\log\tau_{\rm BCOV}(\widetilde{X}_{(S'_{t},\theta'_{t},T)})+O\left(\log(-\log|t|)\right)
\\
&=
\frac{1}{2}\log|t|^{2}+O\left(\log(-\log|t|)\right)
\end{aligned}
$$
as $t\to0$. This completes the proof.
\end{pf}

\subsubsection
{The singularity of BCOV invariants near the discriminant locus}
\label{sect:3.3.5}

\begin{theorem}
\label{thm:conj:formula:singularity:BCOV:inv}
Let $M$ be a primitive $2$-elementary Lorentzian sublattice of ${\Bbb L}_{K3}$. 
Let $C\subset{\mathcal M}_{\Lambda}$ be a compact Riemann surface intersecting $\overline{\mathcal D}_{\Lambda}$ transversally
at sufficiently general point ${\frak p}\in C\cap\overline{\mathcal D}_{\Lambda}$. 
Let $\gamma\colon({\bf C},0)\to(C,{\frak p})$ be an isomorphism of germs.
Let $c\in{\bf X}(1)$.
Write $\tau_{\rm BCOV}(\gamma(t),c)$ for the BCOV invariant of the Borcea-Voisin threefold associated with the $2$-elementary $K3$ surface 
with period $\gamma(z)$ and the elliptic curve with period $c$. 
Then the following holds
\begin{equation}
\label{eqn:singularity:BCOV:torsion:Borcea:Voisin}
\log\tau_{\rm BCOV}(\gamma(t),c)=\frac{1}{2}\log|t|^{2}+O\left(\log(-\log|t|)\right)
\qquad
(t\to0).
\end{equation}
\end{theorem}

\begin{pf}
If \eqref{eqn:singularity:BCOV:torsion:Borcea:Voisin} holds for one isomorphism $\varpi\colon({\bf C},0)\to(C,{\frak p})$,
then it holds for every isomorphism of germs $\gamma\colon({\bf C},0)\to(C,{\frak p})$.
Hence it suffices to prove \eqref{eqn:singularity:BCOV:torsion:Borcea:Voisin} for one particular isomorphism $\varpi\colon({\bf C},0)\to(C,{\frak p})$.
\par
By \cite[Th.\,2.8]{Yoshikawa04}, there is a family of $2$-elementary $K3$ surfaces 
$p\colon({\mathcal S},\theta)\to B$ of type $M$ 
(with degenerate fibers) over a pointed compact Riemann surface $(B,{\frak q})$ with the following properties:
\begin{itemize}
\item[(1)]
Let $\Delta\subset B$ be the discriminant locus of $p\colon{\mathcal S}\to B$ and
let $\varpi\colon B\setminus\Delta\ni b\to\pi_{M}(S_{b},\theta_{b})\in{\mathcal M}_{\Lambda}^{0}$ be the period map for $p\colon({\mathcal S},\theta)\to B$,
where $S_{b}:=p^{-1}(b)$ and $\theta_{b}:=\theta|_{S_{b}}$.
Then $\varpi$ extends to a surjective holomorphic map from $B$ to $C$ with ${\frak p}=\varpi({\frak q})$ such that $\varpi$ is non-degenerate at ${\frak q}$.
In particular, $\varpi\colon B\to{\mathcal M}_{\Lambda}$ intersects $\overline{\mathcal D}_{\Lambda}$ transversally at ${\frak p}=\varpi({\frak q})$.
\item[(2)]
There is a neighborhood $U$ of ${\frak q}$ in B such that $p\colon(p^{-1}(U),\theta|_{p^{-1}(U)})\to U$ is an ordinary singular family of
$2$-elementary $K3$ surfaces of type $M$.
\end{itemize}
\par
Let $T$ be an elliptic curve with period $c$. Let $(U,t)$ be a coordinate neighborhood of ${\frak q}$ in $B$. By choosing $U$ sufficiently small, 
it follows from (1) that $\varpi$ induces an isomorphism of germs $(U,{\frak q})$ and $(C,{\frak p})$. By construction,
\begin{equation}
\label{eqn:identification}
\tau_{\rm BCOV}(\varpi(t),c)=\tau_{\rm BCOV}(\widetilde{X}_{(S_{t},\theta_{t},T)}).
\end{equation}
Since $p\colon(p^{-1}(U),\theta|_{p^{-1}(U)})\to U$ is an ordinary singular family of $2$-elementary $K3$ surfaces by (2), 
the result for the isomorphism of germs $\varpi\colon(U,{\frak q})\to(C,{\frak p})$ follows from \eqref{eqn:identification} and 
Theorems~\ref{thm:singularity:BCOV:type(2,2)} and \ref{thm:singularity:BCOV:type(0,4)}.
\end{pf}

\subsection
{Proof of Theorem~\ref{conj:BCOV=orb:BCOV}}
\label{sect:3.4}
\par
By Theorems~\ref{thm:comparison:curvature:BCOV:orb:BCOV} and \ref{thm:conj:formula:singularity:BCOV:inv},
we have the following equation of currents on ${\mathcal M}_{\Lambda}\times{\bf X}(1)$
$$
-dd^{c}\log\tau_{\rm BCOV}
=
(r(M)-6){\rm pr}_{1}^{*}\omega_{\Lambda}
+
4{\rm pr}_{1}^{*}J_{M}^{*}\omega_{{\mathcal A}_{g(M)}}
+
12{\rm pr}_{2}^{*}\omega_{\rm hyp}
-
\frac{1}{2}\delta_{\overline{\mathcal D}_{\Lambda}\times{\bf X}(1)}.
$$
Since $\log(\tau_{M}^{-4}\|\eta^{24}\|^{2})$ satisfies the same equation of currents on ${\mathcal M}_{\Lambda}\times{\bf X}(1)$
by \cite[Th.\,6.5]{Yoshikawa04}, we get the following equation of currents on ${\mathcal M}_{\Lambda}\times{\bf X}(1)$
$$
-dd^{c}\log[\tau_{\rm BCOV}/(\tau_{M}^{-4}\|\eta^{24}\|^{2})]=0.
$$
Hence 
$F^{\Lambda}=\log[\tau_{\rm BCOV}/(\tau_{M}^{-4}\|\eta^{24}\|^{2})]$ is a pluriharmonic function on ${\mathcal M}_{\Lambda}\times{\bf X}(1)$.
By Lemma~\ref{lemma:singularity:BCOV:invariant:boundary:X(1)}, $F^{\Lambda}$ verifies the assumption of Lemma~\ref{lemma:Grauert:Remmert:extension}.
\par
If $r\leq17$, then $F^{\Lambda}$ is a constant function on ${\mathcal M}_{\Lambda}\times{\bf X}(1)$ by Lemma~\ref{lemma:Grauert:Remmert:extension} (1). 
The assertion is proved when $r\leq17$.
Let $r\geq18$.
Then ${\mathcal B}_{\Lambda}={\mathcal M}_{\Lambda}^{*}\setminus{\mathcal M}_{\Lambda}$ is an irreducible divisor of ${\mathcal M}_{\Lambda}^{*}$
by \cite[Prop.\,11.7]{Yoshikawa13}. By Lemma~\ref{thm:BCOV=orb:BCOV:U(2)+U(2)},
$F^{\Lambda}$ is a pluriharmonic function on ${\mathcal M}_{\Lambda}\times{\bf X}(1)$ 
with at most logarithmic singularity along ${\mathcal B}_{\Lambda}\times{\bf X}(1)$. 
By Lemma~\ref{lemma:Grauert:Remmert:extension} (2), $F^{\Lambda}$ is a constant function on ${\mathcal M}_{\Lambda}\times{\bf X}(1)$.
This proves the assertion when $r\geq18$.
This completes the proof.
\qed

\section
{Orbifold submersions and characteristic forms}
\label{sect:4}
\par
In this section, we study a local model of holomorphic orbifold submersions (cf. \cite[Def.\,1.7]{Ma05}) of relative dimension $3$
and give explicit expressions of some orbifold characteristic forms associated to it (Theorem~\ref{thm:orbifold:Tod:ch:family}).
This result shall play crucial roles in Sections~\ref{sect:5}, \ref{sect:6}, \ref{sect:7}.
In what follows, for a group $G$, we set 
$$
G^{*}:=G\setminus\{1\}
$$

\subsection
{Characteristic forms}
\label{sect:4.1}
\par
For a holomorphic Hermitian vector bundle $(E,h)$ over a complex manifold $Z$, we denote by $R(E,h)\in A^{1,1}(Z,{\rm End}(E))$
the curvature form of $(E,h)$ with respect to the holomorphic Hermitian connection.
The $p$-th Chern form of $(E,h)$ is denoted by $c_{p}(E,h)\in A^{p,p}(Z)$.
The Todd and Chern character forms of $(E,h)$ are the differential forms on $Z$ defined as
$$
{\rm Td}(E,h):=\det\left(\frac{\frac{i}{2\pi}R(E,h)}{I_{E}-\exp(-\frac{i}{2\pi}R(E,h))}\right),
\quad
{\rm ch}(E,h):={\rm Tr}\left[\exp\left(\frac{i}{2\pi}R(E,h)\right)\right].
$$

\subsubsection
{Equivariant characteristic forms}
\label{sect:4.1.1}
\par
Let $G\subset{\rm Aut}(Z)$ be a finite group of automorphisms of $Z$. 
Assume that $(E,h)$ is a $G$-equivariant holomorphic Hermitian vector bundle. 
For $g\in G$, let $Z^{g}$ be the set of its fixed points:
$$
Z^{g}:=\{z\in Z;\,g(z)=z\}.
$$
If ${\rm ord}(g)=n$, we have the splitting of holomorphic vector bundles on $Z^{g}$
\begin{equation}
\label{eqn:isotypical:splitting}
E|_{Z^{g}}=\bigoplus_{\alpha}E(\theta_{\alpha}),
\qquad
\theta_{\alpha}\in\{0,\,\frac{2\pi}{n},\ldots,\frac{2\pi(n-1)}{n}\},
\end{equation}
where $E(\theta):=\{v\in E;\,g(v)=e^{i\theta}\,v\}$ is the eigenbundle of $E|_{Z^{g}}$.
The splitting \eqref{eqn:isotypical:splitting} is orthogonal with respect to the Hermitian metric $h$.
We set $h_{E(\theta)}:=h|_{E(\theta)}$. 
\par
Define the equivariant Todd and Chern character forms of $(E,h)$ as
$$
{\rm Td}_{g}(E,h):={\rm Td}(E(0),h_{E(0)})
\prod_{\theta_{\alpha}\not=0}\frac{\rm Td}{\rm e}\left(\frac{i}{2\pi}R(E(\theta_{\alpha}),h_{E(\theta_{\alpha})})+i\theta_{\alpha}\right),
$$
$$
{\rm ch}_{g}(E,h)
:=
{\rm Tr}\left[g\cdot\exp\left(\frac{i}{2\pi}R(E,h)\right)\right]
=
\sum_{\alpha}e^{i\theta_{\alpha}}{\rm Tr}\left[\exp\left(\frac{i}{2\pi}R(E(\theta_{\alpha}),h_{E(\theta_{\alpha})})\right)\right],
$$
where $({\rm Td}/{\rm e})(A+i\theta):=\det[I_{n}/(I_{n}-e^{-i\theta}\exp(-A))]$ for an $(n,n)$-matrix $A$.
By definition, ${\rm Td}_{g}(E,h)$ and ${\rm ch}_{g}(E,h)$ are differential forms in $\bigoplus_{p\geq0}A^{p,p}(Z^{g})$.

\subsubsection
{Orbifold characteristic forms}
\label{sect:4.1.2}
\par
We consider a local model of abelian orbifold. Assume that {\em $G$ is abelian}. 
For $z\in Z$, let $G_{z}$ be the stabilizer of $z$ in $G$.
Set $\overline{Z}:=Z/G$. To define its inertia orbifold $\Sigma\overline{Z}$ (cf. \cite[Sect.\,1]{Kawasaki78}, \cite[Sect.\,1.1]{Ma05}), set 
\begin{equation}
\label{eqn:twisted:sector:1}
\Sigma Z:=\{(z,g)\in Z\times G^{*};\,z\in Z^{g}\}=\amalg_{g\in G^{*}}Z^{g}\times\{g\}.
\end{equation}
Since $Z^{g}$ is a disjoint union of manifolds of various dimensions, so is $\Sigma Z$. 
On $\Sigma Z$, $G$ acts by $h\cdot(z,g):=(h\cdot z,hgh^{-1})$ for $(z,g)\in\Sigma Z$, $h\in G$. 
Since $G$ is abelian, the $G$-action on $\Sigma Z$ preserves every $Z^{g}=Z^{g}\times\{g\}$.
We define the orbifold $\Sigma\overline{Z}$ as
\begin{equation}
\label{eqn:twisted:sector:2}
\Sigma\overline{Z}:=(\Sigma Z)/G=(\amalg_{g\in G^{*}}Z^{g}\times\{g\})/G.
\end{equation}
\par
Define the map $\nu\colon\Sigma Z\to Z$ by $\nu(z,g)=z$. Then $\nu|_{Z^{g}}$ is the inclusion $Z^{g}\hookrightarrow Z$.
By the $G$-equivariance of $\nu$, it descends to a map $\overline{\nu}\colon\Sigma\overline{Z}\to\overline{Z}$.
Since $\Sigma Z=\{(z,g)\in Z\times G^{*};\,g\in G_{z}\}$, we have $\nu^{-1}(z)=\{(z,g);\,g\in G_{z}^{*}\}$.
Let $U_{z}$ be a small $G_{z}$-invariant neighborhood of $z\in\nu(\Sigma Z)$ in $Z$ and set $\overline{U}_{z}:=U_{z}/G_{z}$. 
Then 
\begin{equation}
\label{eqn:structure:inertia:orbifold}
\nu^{-1}(U_{z})=\amalg_{g\in G_{z}^{*}}U_{z}^{g}\times\{g\},
\qquad
\overline{\nu}^{-1}(\overline{U}_{z})=\amalg_{g\in G_{z}^{*}}(U_{z}^{g}/G_{z})\times\{g\}.
\end{equation}
\par
Let $\pi\colon Z\to\overline{Z}$ and $\pi^{\Sigma}\colon\Sigma Z\to\Sigma\overline{Z}$ be the projections.
For a $G$-invariant differential form $\omega$ on $Z$, $\pi_{*}\omega$ is the differential form on $\overline{Z}$ in the sense of orbifolds 
such that $\pi^{*}(\pi_{*}\omega)=|G|\,\omega$. Similarly, the map of differential forms $\pi^{\Sigma}_{*}$ is defined.
By definition, $\int_{\overline{Z}}\pi_{*}\omega=\int_{Z}\omega$ and $\int_{\Sigma\overline{Z}}\pi^{\Sigma}_{*}\omega'=\int_{\Sigma Z}\omega'$
for $G$-invariant differential forms with compact support $\omega\in A_{0}^{*}(Z)$ and $\omega'\in A_{0}^{*}(\Sigma Z)$.
\par
Let $(\overline{E},\overline{h})$ be the Hermitian orbifold vector bundle on $\overline{Z}$ associated to $(E,h)$. 
Let $\phi(\cdot)$ be a ${\rm GL}({\bf C}^{r})$-invariant polynomial on $\frak{gl}({\bf C}^{r})$, $r={\rm rk}\,E$. 
By the $G$-equivariance of $(E,h)$, the characteristic forms $\phi(E,h)$ and $\sum_{g\in G^{*}}\phi_{g}(E,h)$ are $G$-invariant 
differential forms on $Z$ and $\Sigma Z$, respectively. We define
\begin{equation}
\label{eqn:orb:char:form:1}
\phi^{\Sigma}(\overline{E},\overline{h})
:=
\frac{1}{|G|}\pi_{*}\phi(E,h)
+
\frac{1}{|G|}(\pi^{\Sigma})_{*}\{\sum_{g\in G^{*}}\phi_{g}(E,h)\}
\in
\bigoplus_{p\geq0}A^{p,p}(\overline{Z}\amalg\Sigma\overline{Z}).
\end{equation}
Set $\pi_{g}:=\pi$ for $g=1$ and $\pi_{g}:=\pi^{\Sigma}|_{Z^{g}\times\{g\}}$ for $g\in G^{*}$. Then 
\begin{equation}
\label{eqn:orb:char:form:2}
\phi^{\Sigma}(\overline{E},\overline{h})
=
\frac{1}{|G|}\sum_{g\in G}(\pi_{g})_{*}\phi_{g}(E,h).
\end{equation}
For more about orbifolds and characteristic forms on orbifolds, see e.g. \cite{Kawasaki78}, \cite{Ma05}.

\subsection
{Group action on threefold}
\label{sect:4.2}
\par
Let $X$ be a smooth threefold and let $G\subset{\rm Aut}(X)$ be a finite {\em abelian} subgroup.
Assume that for all $g\in G$ and $x\in X^{g}$,
\begin{equation}
\label{eqn:assumption:group:action}
g_{*}\in{\rm SL}(T_{x}X)\cong{\rm SL}({\bf C}^{3}).
\end{equation}
If $1$ is an eigenvalue of  $g_{*}\in{\rm SL}(T_{x}X)$ and $g\not=1$, its possible multiplicity is one by \eqref{eqn:assumption:group:action}.
Hence $X^{g}$ consists of at most finitely many, disjoint curves and isolated points.
Let $X^{g,(k)}$ be the union of all components of $X^{g}$ of dimension $k$. Then
$$
\Sigma X={\Sigma}^{(0)}X\amalg{\Sigma}^{(1)}X,
\qquad
\Sigma^{(k)}X:=\amalg_{g\in G^{*}}X^{g,(k)}.
$$
We denoted by
$$
\nu({\Sigma}^{(1)}X)=\bigcup_{\lambda\in\Lambda}C_{\lambda}
$$
the irreducible decomposition. Then $C_{\lambda}\subset M^{g}$ for some $g\in G^{*}$.
In Section~\ref{sect:4.2}, we study the local structure of $\nu({\Sigma}X)$ in $X$.
\par
By \eqref{eqn:assumption:group:action}, we have an inclusion $G_{x}\subset{\rm SL}(T_{x}X)\cong{\rm SL}({\bf C}^{3})$
and hence an isomorphism of germs with group action
\begin{equation}
\label{eqn:local:model}
(X,x,G_{x})\cong({\bf C}^{3},0,\Gamma),
\end{equation}
where $\Gamma\subset{\rm SL}({\bf C}^{3})$ is a finite abelian subgroup by our assumption. 
Since $\Gamma$ is abelian, there exist characters $\chi_{k}\in{\rm Hom}(\Gamma,{\bf C}^{*})$, $k=1,2,3$, such that
\begin{equation}
\label{eqn:abelian:subgroup:SL(3,C)}
\Gamma=\{{\rm diag}(\chi_{1}(g),\chi_{2}(g),\chi_{3}(g))\in{\rm SL}({\bf C}^{3});\,g\in\Gamma\},
\qquad
\chi_{1}\chi_{2}\chi_{3}=1.
\end{equation}
Let $({\frak C}_{k},0)\subset({\bf C}^{3},0)$ be the germ of the $x_{k}$-axis of $({\bf C}^{3},0)$.
By \eqref{eqn:abelian:subgroup:SL(3,C)}, if $x\in C_{\lambda}$, $\lambda\in\Lambda$, the germ $(C_{\lambda},x)$ must be of the form 
$({\frak C}_{k},0)$ in the expression \eqref{eqn:local:model}. Moreover, 
\begin{equation}
\label{eqn:axis:fixed:curve}
({\frak C}_{k},0)=(C_{\lambda},x)\hbox{ for some }\lambda\in\Lambda
\quad\Longleftrightarrow\quad
\ker\chi_{k}\not=\{1\}.
\end{equation}
By \eqref{eqn:axis:fixed:curve}, the number of the components of the germ $(\nu({\Sigma}^{(1)}X),x)$ is at most $3$ for all $x\in X$. 
For a given number $1\leq n\leq 3$, there exists $\Gamma\subset{\rm SL}({\bf C}^{3})$ such that the corresponding $(\nu({\Sigma}^{(1)}X),x)$ 
consists of exactly $n$ components.

\begin{example}
\label{example:subgroup:SL(3,C):1}
Let $\Gamma=\langle{\rm diag}(i,i,-1)\rangle$. Then $(\nu({\Sigma}^{(1)}X),x)={\frak C}_{3}$
and $x\in\nu({\Sigma}^{(0)}X)\cap\nu({\Sigma}^{(1)}X)$ is a smooth point of $\nu({\Sigma}^{(1)}X)$.
\end{example}

\begin{example}
\label{example:subgroup:SL(3,C):2}
Let $\Gamma=\langle{\rm diag}(-\omega,-1,\omega^{2})\rangle$, where $\omega$ is a primitive cube root of $1$.
Then $(\nu({\Sigma}^{(1)}X),x)={\frak C}_{2}\cup{\frak C}_{3}$ and $x\in\nu(\Sigma^{(0)}X)\cap\nu(\Sigma^{(1)}X)$.
\end{example}

\begin{example}
\label{example:subgroup:SL(3,C):3}
Let $\Gamma=\langle{\rm diag}(1,-1,-1),{\rm diag}(-1,1,-1),{\rm diag}(-1,-1,1)\rangle$. 
Then $(\nu({\Sigma}^{(1)}X),x)={\frak C}_{1}\cup{\frak C}_{2}\cup{\frak C}_{3}$ and 
$x\in\nu(\Sigma^{(1)}X)\setminus\nu(\Sigma^{(0)}X)$.
\end{example}

\par
For later uses, we introduce $\Gamma^{0}\subset\Gamma$ and $\delta_{k}(\Gamma)\in{\bf Q}$ $(k=1,2,3)$ as follows.

\begin{definition}
\label{def:subset:Gamma:0}
For a subgroup $\Gamma$ of ${\rm SL}({\bf C}^{3})$, define 
\begin{equation}
\label{def:Gamma:0}
\Gamma^{0}:=\{g\in\Gamma^{*};\,\det(g-1_{{\bf C}^{3}})\not=0\}.
\end{equation}
\end{definition}

\begin{definition}
\label{def:admissibility:subgroup:SL}
Let $\Gamma$ be a finite abelian subgroup of ${\rm SL}({\bf C}^{3})$ expressed as in \eqref{eqn:abelian:subgroup:SL(3,C)}.
If $\Gamma^{0}=\emptyset$, set $\delta_{1}(\Gamma)=\delta_{2}(\Gamma)=\delta_{3}(\Gamma)=0$.
If $\Gamma^{0}\not=\emptyset$, set
\begin{equation}
\label{eqn:admissibility}
\delta_{k}(\Gamma)
:=
\sum_{g\in\Gamma^{0}}\frac{\chi_{k}(g)}{(1-\chi_{k}(g))^{2}}
\qquad
(k=1,2,3).
\end{equation}
\end{definition}

For a finite abelian subgroup $\Gamma\subset{\rm SL}({\bf C}^{3})$, 
set $({\bf C}^{3})^{\Gamma^{*}}:=\bigcap_{\gamma\in\Gamma^{*}}({\bf C}^{3})^{\gamma}$.
To understand the structure of $\Sigma C_{\lambda}$ in Section~\ref{sect:4.5} below, we need the following:

\begin{lemma}
\label{lemma:invariant:subgrp:SL:C3:Gamma0:empty}
Let $\Gamma\subset{\rm SL}({\bf C}^{3})$ be a finite abelian subgroup expressed as in \eqref{eqn:abelian:subgroup:SL(3,C)}.
If $({\bf C}^{3})^{\Gamma^{*}}=\{0\}$ and $\Gamma^{0}=\emptyset$, then
$$
\Gamma=\{1,\,{\rm diag}(-1,-1,1),\,{\rm diag}(-1,1,-1),\,{\rm diag}(1,-1,-1)\}.
$$
\end{lemma}

\begin{pf}
Since $\Gamma^{0}=\emptyset$, we get $\Gamma=\ker\chi_{1}\cup\ker\chi_{2}\cup\ker\chi_{3}$.
Thus we have $\Gamma^{*}=(\ker\chi_{1})^{*}\cup(\ker\chi_{2})^{*}\cup(\ker\chi_{3})^{*}$.
Since $\Gamma^{*}\not=(\ker\chi_{k})^{*}$ for all $k\in\{1,2,3\}$ by the condition $({\bf C}^{3})^{\Gamma^{*}}=\{0\}$,
we get $(\ker\chi_{l})^{*}\cup(\ker\chi_{m})^{*}\not=\emptyset$ for any $l\not=m$. 
This implies $(\ker\chi_{l'})^{*}\not=\emptyset$ and $(\ker\chi_{m'})^{*}\not=\emptyset$ for some $l'\not=m'$.
We may assume $(\ker\chi_{1})^{*}\not=\emptyset$, $(\ker\chi_{2})^{*}\not=\emptyset$.
Take arbitrary elements
$g_{1}={\rm diag}(1,\zeta,\zeta^{-1})\in(\ker\chi_{1})^{*}$ and $g_{2}={\rm diag}(\xi,1,\xi^{-1})\in(\ker\chi_{2})^{*}$, 
where $\zeta,\xi\in{\bf C}^{*}\setminus\{1\}$.
Since $g_{1}g_{2}={\rm diag}(\xi,\zeta,\zeta^{-1}\xi^{-1})\not\in\ker\chi_{1}\cup\ker\chi_{2}$, we get $\xi=\zeta^{-1}$.
Since $g_{1}g_{2}^{-1}={\rm diag}(\zeta,\zeta,\zeta^{-2})\not\in\ker\chi_{1}\cup\ker\chi_{2}$, we get $\zeta=-1$. 
Since $g_{1},g_{2}$ are arbitrary, we get $(\ker\chi_{1})^{*}=\{{\rm diag}(1,-1,-1)\}$ and $(\ker\chi_{2})^{*}=\{{\rm diag}(-1,1,-1)\}$.
Since $g_{1}g_{2}\in(\ker\chi_{3})^{*}$, we get $(\ker\chi_{3})^{*}\not=\emptyset$.
By the same argument, we get $(\ker\chi_{3})^{*}=\{{\rm diag}(-1,-1,1)\}$.
Since $\Gamma^{*}=(\ker\chi_{1})^{*}\cup(\ker\chi_{2})^{*}\cup(\ker\chi_{3})^{*}$, we get the result.
\end{pf}

After Lemma~\ref{lemma:invariant:subgrp:SL:C3:Gamma0:empty}, we set
$$
{\Bbb G}
:=
\{1,\,{\rm diag}(-1,-1,1),\,{\rm diag}(-1,1,-1),\,{\rm diag}(1,-1,-1)\}
\cong
({\bf Z}/2{\bf Z})^{\oplus2}.
$$
In Section~\ref{sect:4.5.1}, we shall show that 
$G_{x}={\Bbb G}$ iff $x\in\bigcup_{\lambda\in\Lambda}\Sigma C_{\lambda}\setminus\nu(\Sigma^{(0)}X)$.
Set
\begin{equation}
\label{eqn:ext:twisted:sector:0}
\widetilde{\Sigma}^{(0)}X
:=
\Sigma^{(0)}X\amalg\{(x,g)\in\Sigma X;\,G_{x}={\Bbb G}\}.
\end{equation}
For $x\in X$, we set $X^{G_{x}^{*}}:=\bigcap_{g\in G_{x}^{*}}X^{g}$.
Since $\Gamma^{0}\not=\emptyset$ implies $({\bf C}^{3})^{\Gamma^{*}}=\{0\}$ and since the set germ $(X^{G_{x}^{*}},x)$ is 
either the point $\{x\}$ or the union of some coordinate axes ${\frak C}_{k}$ by \eqref{eqn:local:model},
\eqref{eqn:abelian:subgroup:SL(3,C)} when $G_{x}^{*}\not=\emptyset$, we deduce from \eqref{eqn:ext:twisted:sector:0} and
Lemma~\ref{lemma:invariant:subgrp:SL:C3:Gamma0:empty} that
\begin{equation}
\label{eqn:ext:twisted:sector:1}
\nu(\widetilde{\Sigma}^{(0)}X)=\{x\in X;\,\dim_{x}X^{G_{x}^{*}}=0\}
\end{equation}
under the convention $\dim\emptyset=-\infty$. We set
$$
X^{(1)}:=\amalg_{\lambda\in\Lambda}C_{\lambda},
\qquad
X^{(0)}:=\nu(\widetilde{\Sigma}^{(0)}X),
$$
where we consider their reduced structure.
Then $X^{(1)}$ is the normalization of $\nu(\Sigma^{(1)}X)$. 
By \eqref{eqn:axis:fixed:curve}, $C_{\lambda}\cap C_{\lambda'}\subset X^{(0)}$ for any $\lambda\not=\lambda'$.
We set $C_{\lambda}^{0}:=C_{\lambda}\setminus X^{(0)}$.
Since $\dim_{x}X^{G_{x}^{*}}\leq1$ by \eqref{eqn:assumption:group:action} and since $X^{G_{x}^{*}}\not=\emptyset$ iff
$x\in\nu(\Sigma X)$, we get by \eqref{eqn:ext:twisted:sector:1}
$$
\nu(\Sigma^{(1)}X)\setminus\nu(\widetilde{\Sigma}^{(0)}X)=\amalg_{\lambda\in\Lambda}C_{\lambda}^{0}=
\{x\in X;\,\dim_{x}X^{G_{x}^{*}}=1\}.
$$ 
In particular, $\nu({\Sigma}^{(1)}X)=\bigcup_{\lambda\in\Lambda}C_{\lambda}=\overline{\{x\in X;\,\dim_{x}X^{G_{x}^{*}}=1\}}$.
\par
For $C_{\lambda}$, $\lambda\in\Lambda$, set 
$$
G_{C_{\lambda}}:=\{g\in G;\,g|_{C_{\lambda}}={\rm id}_{C_{\lambda}}\}.
$$
Since $G_{C_{\lambda}}$ is identified with its image in ${\rm SL}(N_{C_{\lambda}/X})$, $G_{C_{\lambda}}$ is isomorphic to
a finite abelian subgroup of ${\rm SL}({\bf C}^{2})$. Hence $G_{C_{\lambda}}$ is cyclic.

\begin{lemma}
\label{lemma:generic:stabilizer}
If $x\in C_{\lambda}^{0}$, then $G_{x}=G_{C_{\lambda}}$.
\end{lemma}

\begin{pf}
Let $\Gamma\cong G_{x}$ be as in \eqref{eqn:local:model}, \eqref{eqn:abelian:subgroup:SL(3,C)}. 
Since $\dim({\bf C}^{3})^{\Gamma^{*}}=1$, $({\bf C}^{3})^{\Gamma^{*}}$ is one of the coordinate axis, say ${\frak C}_{1}$. 
Then $(X^{G_{x}^{*}},x)\cong({\frak C}_{1},0)$ and $\chi_{1}\equiv1$, $\chi_{3}=\chi_{2}^{-1}$, so that 
$\Gamma\cong{\rm Im}\,\chi_{2}=\{{\rm diag}(1,\zeta,\zeta^{-1});\,\zeta\in{\rm Im}\,\chi_{2}\}$. 
By this expression, $G_{x}=\langle g\rangle$ for some $g\in G$ and $g|_{C_{\lambda}}={\rm id}_{C_{\lambda}}$,
where $(C_{\lambda},x)=({\frak C}_{1},0)$ via \eqref{eqn:local:model}.
Since $g\in G_{C_{\lambda}}$ and $G_{x}=\langle g\rangle$, we get $G_{x}\subset G_{C_{\lambda}}$.
Since $x\in C_{\lambda}$, the inclusion $G_{x}\supset G_{C_{\lambda}}$ is trivial.
\end{pf}

By \eqref{eqn:structure:inertia:orbifold} and Lemma~\ref{lemma:generic:stabilizer}, 
the fiber of the natural projection $\Sigma^{(1)}X\to X^{(1)}$ over $C_{\lambda}$ is given by $G_{C_{\lambda}}^{*}$.

\subsection
{Equivariant submersions and characteristic forms}
\label{sect:4.3}
\par

\subsubsection
{Set up}
\label{sect:4.3.1}
\par
Let ${\mathcal X}$ and $B\cong\varDelta^{\dim B}$ be complex manifolds and let
$$
f\colon{\mathcal X}\to B
$$ 
be a holomorphic submersion such that $X_{b}:=f^{-1}(b)$ is a connected {\em threefold}.
We do not assume the properness of $f$. Hence both ${\mathcal X}$ and $X_{b}$ can be non-compact.
\par
Let $G$ be a finite abelian group of automorphisms of ${\mathcal X}$ such that $f\colon{\mathcal X}\to B$ is $G$-equivariant with respect to
the trivial $G$-action on $B$. Hence $G$ preserves the fibers of
$f\colon{\mathcal X}\to B$. 
The order of $g\in G$ is denoted by $n_{g}={\rm ord}(g)$.
\par
Set ${\mathcal Y}:={\mathcal X}/G$,  which is equipped with the projection induced from $f$
$$
\overline{f}\colon{\mathcal Y}\to B.
$$ 
Then the fibers of $\overline{f}\colon{\mathcal Y}\to B$ are the orbifolds $Y_{b}:=X_{b}/G$, $b\in B$ of dimension $3$.
When ${\mathcal X}$ is an open subset of ${\bf C}^{\dim B+3}$, 
$\overline{f}\colon{\mathcal Y}\to B$ is viewed as a local model of holomorphic submersion of orbifolds of relative dimension $3$.
As in Section~\ref{sect:4.2}, we assume that for all $g\in G$ and $x\in X_{b}^{g}$,
\begin{equation}
\label{eqn:assumption:group}
g_{*}\in{\rm SL}(T_{x}X_{b})\cong{\rm SL}({\bf C}^{3}).
\end{equation}

\subsubsection
{Some loci of ${\mathcal X}$ and ${\mathcal Y}$ and some subgroups of $G$}
\label{sect:4.3.2}
\par
Since $f\colon{\mathcal X}\to B$ 
is a $G$-equivariant submersion, ${\mathcal X}^{g}$ is a complex submanifold of ${\mathcal X}$ flat over $B$. 
For a component ${\mathcal Z}$ of ${\mathcal X}^{g}$,
$f\colon{\mathcal Z}\to B$ is a holomorphic submersion with connected fiber.
Let ${\mathcal X}^{g,(k)}$ be the union of all components of ${\mathcal X}^{g}$ of relative dimension $k$. 
Then 
$$
\Sigma{\mathcal X}=\Sigma^{(0)}{\mathcal X}\amalg\Sigma^{(1)}{\mathcal X},
\qquad
{\Sigma}^{(k)}{\mathcal X}
=
\amalg_{g\in G^{*}}{\mathcal X}^{g,(k)}\times\{g\}.
$$
Since $\Sigma{\mathcal X}$ is a complex submanifold with $G$-action, we get the orbifolds
$$
\Sigma{\mathcal Y}:=\Sigma{\mathcal X}/G,
\qquad
\Sigma^{(k)}{\mathcal Y}=\Sigma^{(k)}{\mathcal X}/G
$$
such that $\Sigma{\mathcal Y}=\Sigma^{(0)}{\mathcal Y}\amalg\Sigma^{(1)}{\mathcal Y}$. 
After \eqref{eqn:ext:twisted:sector:0}, \eqref{eqn:ext:twisted:sector:1}, we also introduce
$$
\widetilde{\Sigma}^{(0)}{\mathcal X}
:=
\{(x,g)\in{\mathcal X}\times G^{*};\,g\in G_{x}^{*},\quad\dim_{x}X_{f(x)}^{G_{x}^{*}}=0\},
\qquad
\widetilde{\Sigma}^{(0)}{\mathcal Y}
:=
\widetilde{\Sigma}^{(0)}{\mathcal X}/G.
$$
\par
Let $\{{\mathcal C}_{\lambda}\}_{\lambda\in\Lambda}$ be the set of irreducible components of $\nu({\Sigma}^{(1)}{\mathcal X})$.
Hence $\nu({\Sigma}^{(1)}{\mathcal X})=\bigcup_{\lambda\in\Lambda}{\mathcal C}_{\lambda}$. 
For ${\mathcal C}_{\lambda}\subset\nu({\Sigma}^{(1)}{\mathcal X})$, set
$$
G_{{\mathcal C}_{\lambda}}
:=
\{g\in G;\,g|_{{\mathcal C}_{\lambda}}={\rm id}_{{\mathcal C}_{\lambda}}\},
\qquad
\Gamma_{{\mathcal C}_{\lambda}}
:=
\{g\in G;\,g({\mathcal C}_{\lambda})={\mathcal C}_{\lambda}\}
$$
Then $G_{{\mathcal C}_{\lambda}}$ and $\Gamma_{{\mathcal C}_{\lambda}}$ are subgroups of $G$ with 
$G_{{\mathcal C}_{\lambda}}\subset\Gamma_{{\mathcal C}_{\lambda}}$. 
Set
$$
\overline{\Gamma}_{{\mathcal C}_{\lambda}}
:=
\Gamma_{{\mathcal C}_{\lambda}}/G_{{\mathcal C}_{\lambda}}.
$$
Since $\overline{\Gamma}_{{\mathcal C}_{\lambda}}\subset{\rm Aut}({\mathcal C}_{\lambda})$, 
we define an orbifold of dimension $\dim B+1$ by
$$
\overline{\mathcal C}_{\lambda}:={\mathcal C}_{\lambda}/\overline{\Gamma}_{{\mathcal C}_{\lambda}}.
$$
Since $G_{{\mathcal C}_{\lambda}}\subset{\rm SL}(N_{{\mathcal C}_{\lambda}/{\mathcal X}})$ and since $G_{{\mathcal C}_{\lambda}}$ is abelian,
$G_{{\mathcal C}_{\lambda}}$ is a cyclic group.
We set
$$
n_{\lambda}:=|G_{{\mathcal C}_{\lambda}}|.
$$
\par
On the set $\Lambda$, $G$ acts by ${\mathcal C}_{g\cdot\lambda}:=g\cdot{\mathcal C}_{\lambda}$. 
Let $\overline{\Lambda}:=\Lambda/G$ be the $G$-orbits of the $G$-actions on $\Lambda$. 
Write $\bar{\lambda}$ for the $G$-orbit $G\cdot\lambda$.
The preimage of $\bar{\lambda}$ in $\Lambda$ consists of $|G/\Gamma_{{\mathcal C}_{\lambda}}|$ distinct points.
We have
$$
\nu(\Sigma^{(1)}{\mathcal X})
=
\bigcup_{\bar{\lambda}\in\overline{\Lambda}}\bigcup_{g\in G/\Gamma_{{\mathcal C}_{\lambda}}}{\mathcal C}_{g\cdot\lambda}.
$$
Since ${\mathcal C}_{g\cdot\lambda}\cong{\mathcal C}_{\lambda}$ via the action of $g\in G$ and hence 
$\overline{\mathcal C}_{g\cdot\lambda}=\overline{\mathcal C}_{\lambda}$
for all $g\in G/\Gamma_{{\mathcal C}_{\lambda}}$,
it makes sense to define $\overline{\mathcal C}_{\bar{\lambda}}:=\overline{\mathcal C}_{\lambda}$.
Since
$$
\nu(\Sigma^{(1)}{\mathcal X})/G
=
(\bigcup_{\bar{\lambda}\in\overline{\Lambda}}\bigcup_{g\in G/\Gamma_{{\mathcal C}_{\lambda}}}{\mathcal C}_{g\cdot\lambda})/G
=
\bigcup_{\bar{\lambda}\in\overline{\Lambda}}{\mathcal C}_{\lambda}/\Gamma_{{\mathcal C}_{\lambda}}
=
\bigcup_{\bar{\lambda}\in\overline{\Lambda}}\overline{\mathcal C}_{\bar\lambda},
$$
we have 
$\overline{\nu}(\Sigma^{(1)}{\mathcal Y})=\bigcup_{\bar{\lambda}\in\overline{\Lambda}}\overline{\mathcal C}_{\bar\lambda}$.
Hence there is a one-to-one correspondence between $\overline{\Lambda}$ and the components of 
${\rm Sing}\,{\mathcal Y}$ of dimension $\dim B+1$, i.e., $\overline{\nu}(\Sigma^{(1)}{\mathcal Y})$. 
We define 
$$
{\mathcal X}^{(1)}:=\amalg_{\lambda\in\Lambda}{\mathcal C}_{\lambda},
\qquad
{\mathcal Y}^{(1)}
:=
\amalg_{\bar{\lambda}\in\overline{\Lambda}}\overline{\mathcal C}_{\bar\lambda}.
$$
Then ${\mathcal X}^{(1)}$ and ${\mathcal Y}^{(1)}$ are the normalizations of ${\nu}(\Sigma^{(1)}{\mathcal X})$ and
$\overline{\nu}(\Sigma^{(1)}{\mathcal Y})$, respectively. In the same way as above, ${\mathcal Y}^{(1)}={\mathcal X}^{(1)}/G$.
Let $p\colon\Sigma^{(1)}{\mathcal X}\to{\mathcal X}^{(1)}$ be the projection. 
By \eqref{eqn:structure:inertia:orbifold} and Lemma~\ref{lemma:generic:stabilizer}, 
$p^{-1}(x)=\{g\in G^{*};\,{\mathcal C}_{\lambda}\subset{\mathcal X}^{g}\}=G_{{\mathcal C}_{\lambda}}^{*}$
for any $x\in{\mathcal C}_{\lambda}$.
\par
Set
$$
{\mathcal X}^{(0)}:=\nu(\widetilde{\Sigma}^{(0)}{\mathcal X}),
\qquad
{\mathcal Y}^{(0)}:=\overline{\nu}(\widetilde{\Sigma}^{(0)}{\mathcal Y})={\mathcal X}^{(0)}/G.
$$
For a component ${\frak p}$ of ${\mathcal X}^{(0)}$, 
let $\overline{\frak p}\subset{\mathcal Y}^{(0)}$ be its image by the projection ${\mathcal X}\to{\mathcal Y}$.
Then $\overline{\frak p}$ is identified with the $G$-orbit of ${\frak p}$.
Since $B$ is contractible,
${\frak p}$ (resp. $\overline{\frak p}$) is a section of $f\colon{\mathcal X}\to B$ (resp. $\overline{f}\colon{\mathcal Y}\to B$)
by the $G$-equivariance of $f$.
We set 
$$
G_{\frak p}
:=
\{g\in G;\,{\frak p}\subset{\mathcal X}^{g}\}
=
\{g\in G;\,g|_{\frak p}={\rm id}_{\frak p}\}.
$$
Then the preimage of $\overline{\frak p}$ in ${\mathcal X}^{(0)}$ consists of $|G/G_{\frak p}|$ distinct points.

\subsubsection
{Some vector bundles and their characteristic forms}
\label{sect:4.3.3}
\par
Let $g\in G$.
Let $\Theta_{{\mathcal X}/B}$ (resp. $\Theta_{{\mathcal X}^{g}/B}$) be the relative holomorphic tangent bundle of the family $f\colon{\mathcal X}\to B$ 
(resp. $f\colon{\mathcal X}^{g}\to B$).
Let $\Omega^{1}_{{\mathcal X}/B}$ (resp. $\Omega^{1}_{{\mathcal X}^{g}/B}$) be the vector bundle of relative K\"ahler differentials of the family 
$f\colon{\mathcal X}\to B$ (resp. $f\colon{\mathcal X}^{g}\to B$).
Let $N_{{\mathcal X}^{g}/{\mathcal X}}$ (resp. $N^{*}_{{\mathcal X}^{g}/{\mathcal X}}$) be the normal (resp. conormal) bundle of ${\mathcal X}^{g}$
in ${\mathcal X}$. 
We have 
$$
{\rm rk}(\Theta_{{\mathcal X}^{g,(k)}/B})={\rm rk}(\Omega^{1}_{{\mathcal X}^{g,(k)}/B})=k,
\qquad
{\rm rk}(N_{{\mathcal X}^{g,(k)}/{\mathcal X}})={\rm rk}(N^{*}_{{\mathcal X}^{g,(k)}/{\mathcal X}})=3-k.
$$
Since $\Theta_{{\mathcal X}^{g}/B}$ (resp. $N_{{\mathcal X}^{g}/{\mathcal X}}$) is the (resp. union of) eigenbundle(s) of  
$\Theta_{{\mathcal X}/B}|_{{\mathcal X}^{g}}$ with respect to the $g$-action corresponding to the eigenvalue $1$ (resp. $\not=1$), we get
\begin{equation}
\label{eqn:splitting:relative:tangent:bundle}
\Theta_{{\mathcal X}/B}|_{{\mathcal X}^{g}}=\Theta_{{\mathcal X}^{g}/B}\oplus N_{{\mathcal X}^{g}/{\mathcal X}}.
\end{equation}
Similarly, we have the splitting
\begin{equation}
\label{eqn:splitting:relative:Kaehler:differentials}
\Omega^{1}_{{\mathcal X}/B}|_{{\mathcal X}^{g}} 
=
\Omega^{1}_{{\mathcal X}^{g}/B}\oplus N^{*}_{{\mathcal X}^{g}/{\mathcal X}}.
\end{equation}

\begin{lemma}
\label{lemma:eigenvalue:normal:bundle:rel:dim=1}
The eigenvalues of the $g$-action on $N_{{\mathcal X}^{g}/{\mathcal X}}|_{{\mathcal X}^{g,(1)}}$ are of the form
$$
\{\exp(2\pi ik/n_{g}),\exp(-2\pi ik/n_{g})\},
\qquad
(n_{g},k)=1,
$$
where $k\in\{1,\ldots,n_{g}-1\}$ may depend on the component of ${\mathcal X}^{g,(1)}$.
\end{lemma}

\begin{pf}
For simplicity, write $n$ for $n_{g}$.
Set $\zeta_{n}:=\exp(2\pi i/n)$.
Since $g$ has order $n$ and $g\in{\rm SL}(N_{{\mathcal X}^{g}/{\mathcal X}}|_{{\mathcal X}^{g,(1)}})$, its eigenvalues are of the form
$\{\zeta_{n}^{k},\zeta_{n}^{-k}\}$, $0\leq k<n$. 
Assume $g^{m}=1$ on $N_{{\mathcal X}^{g}/{\mathcal X}}|_{{\mathcal X}^{g,(1)}}$ for some $1\leq m<n$. 
Since $g^{m}=1$ on $\Theta_{\mathcal X}|_{{\mathcal X}^{g,(1)}}$, we get $g^{m}=1$ on ${\mathcal X}$, so that $m=nl$ for some $l\in{\bf Z}$.
This contradicts the choice of $m$. Hence $g^{m}\not=1$ on $N_{{\mathcal X}^{g}/{\mathcal X}}|_{{\mathcal X}^{g,(1)}}$ for any $1\leq m<n$,
which implies $(n,k)=1$.
\end{pf}

Let $h_{{\mathcal X}/B}$ be a $G$-invariant Hermitian metric on $\Theta_{{\mathcal X}/B}$, which is fiberwise K\"ahler.
The vector bundles 
$\Theta_{{\mathcal X}^{g}/B}$, $\Omega^{p}_{{\mathcal X}/B}$, $\Omega^{p}_{{\mathcal X}^{g}/B}$, 
$N_{{\mathcal X}^{g}/{\mathcal X}}$, $N^{*}_{{\mathcal X}^{g}/{\mathcal X}}$ are equipped with the Hermitian metrics induced from $h_{{\mathcal X}/B}$,
which are denoted by $h_{{\mathcal X}^{g}/B}$, $h_{\Omega^{p}_{{\mathcal X}/B}}$, $h_{\Omega^{p}_{{\mathcal X}^{g}/B}}$,
$h_{N_{{\mathcal X}^{g}/{\mathcal X}}}$, $h_{N^{*}_{{\mathcal X}^{g}/{\mathcal X}}}$, respectively.
Then the splittings \eqref{eqn:splitting:relative:tangent:bundle}, \eqref{eqn:splitting:relative:Kaehler:differentials}
are orthogonal with respect to the metrics $h_{{\mathcal X}/B}$ and $h_{\Omega^{1}_{{\mathcal X}/B}}$, respectively.
\par
In what follows, we write 
$c_{k}({\mathcal X}/B)$, $c_{k}({\mathcal X}^{g}/B)$, $c_{k}(\Omega^{p}_{{\mathcal X}/B})$, $c_{k}(\Omega^{p}_{{\mathcal X}^{g}/B})$,
$c_{k}(N_{{\mathcal X}^{g}/{\mathcal X}})$, $c_{k}(N^{*}_{{\mathcal X}^{g}/{\mathcal X}})$
for the $k$-th Chern forms of 
$(\Theta_{{\mathcal X}/B}, h_{{\mathcal X}/B})$,
$(\Theta_{{\mathcal X}^{g}/B}, h_{{\mathcal X}^{g}/B})$,
$(\Omega^{p}_{{\mathcal X}/B},h_{\Omega^{p}_{{\mathcal X}/B}})$,
$(\Omega^{p}_{{\mathcal X}^{\iota}/B},h_{\Omega^{p}_{{\mathcal X}^{\iota}/B}})$,
$(N_{{\mathcal X}^{g}/{\mathcal X}},h_{N_{{\mathcal X}^{g}/{\mathcal X}}})$,
$(N^{*}_{{\mathcal X}^{g}/{\mathcal X}},h_{N^{*}_{{\mathcal X}^{g}/{\mathcal X}}})$,
respectively. We have the following standard relations of Chern forms:
\begin{equation}
\label{eqn:relation:Chern:forms:1}
c_{1}({\mathcal X}/B)|_{{\mathcal X}^{g}}=c_{1}({\mathcal X}^{g}/B)+c_{1}(N_{{\mathcal X}^{g}/{\mathcal X}}),
\end{equation}
\begin{equation}
\label{eqn:relation:Chern:forms:2}
c_{1}(\Omega^{1}_{{\mathcal X}/B})=-c_{1}({\mathcal X}/B),
\qquad
c_{1}(\Omega^{1}_{{\mathcal X}^{g}/B})=-c_{1}({\mathcal X}^{g}/B),
\end{equation}
\begin{equation}
\label{eqn:relation:Chern:forms:3}
c_{1}(N^{*}_{{\mathcal X}^{g}/{\mathcal X}})=-c_{1}(N_{{\mathcal X}^{g}/{\mathcal X}}),
\qquad
c_{2}(N^{*}_{{\mathcal X}^{g}/{\mathcal X}})=c_{2}(N_{{\mathcal X}^{g}/{\mathcal X}}).
\end{equation}

\begin{lemma}
\label{lemma:Todd:chern:character:BCOV}
The following equality of differential forms on ${\mathcal X}$ holds:
$$
[{\rm Td}(\Theta_{{\mathcal X}/B},h_{{\mathcal X}/B})
\sum_{p\geq0}(-1)^{p}p\,{\rm ch}(\Omega^{p}_{{\mathcal X}/B},h_{\Omega_{{\mathcal X}/B}^{p}})]^{(8)}
=
-\frac{1}{12}c_{1}({\mathcal X}/B)c_{3}({\mathcal X}/B).
$$
\end{lemma}

\begin{pf}
See \cite[p.374]{BCOV94}.
\end{pf}

\subsubsection
{An equivariant characteristic form: case of relative dimension $1$}
\label{sect:4.3.4}
\par
Let ${\mathcal C}$ be a component of ${\mathcal X}^{(1)}$.
Let $g\in G_{\mathcal C}^{*}$. Then $n_{g}|n_{\mathcal C}$.
By Lemma~\ref{lemma:eigenvalue:normal:bundle:rel:dim=1},
there exists $k_{g}$ with $(k_{g},n_{g})=1$ such that the $g$-action on $N_{{\mathcal C}/{\mathcal X}}$ induces the splitting
$$
N_{{\mathcal C}/{\mathcal X}}
=
N_{{\mathcal C}/{\mathcal X}}(\theta)\oplus N_{{\mathcal C}/{\mathcal X}}(-\theta),
\qquad
\theta=\frac{2k_{g}\pi}{n_{g}}.
$$
We have the corresponding splitting
$$
N^{*}_{{\mathcal C}/{\mathcal X}}
=
N^{*}_{{\mathcal C}/{\mathcal X}}(\theta)\oplus N^{*}_{{\mathcal C}/{\mathcal X}}(-\theta)
$$
such that
\begin{equation}
\label{eqn:relation:dual:g:action}
N^{*}_{{\mathcal C}/{\mathcal X}}(\theta)=(N_{{\mathcal C}/{\mathcal X}}(-\theta))^{*},
\qquad
N^{*}_{{\mathcal C}/{\mathcal X}}(-\theta)=(N_{{\mathcal C}/{\mathcal X}}(\theta))^{*}.
\end{equation}
Set
$$
{\frak l}_{1}:=c_{1}(N_{{\mathcal C}/{\mathcal X}}(\theta),h_{N_{{\mathcal C}/{\mathcal X}}(\theta)}),
\qquad
{\frak l}_{2}:=c_{1}(N_{{\mathcal C}/{\mathcal X}}(-\theta),h_{N_{{\mathcal C}/{\mathcal X}}(-\theta)}).
$$
By \eqref{eqn:relation:Chern:forms:1}, we get
\begin{equation}
\label{eqn:relation:sum:Chern:forms}
{\frak l}_{1}+{\frak l}_{2}
=
c_{1}(N_{{\mathcal C}/{\mathcal X}},h_{N_{{\mathcal C}/{\mathcal X}}})
=
c_{1}({\mathcal X}/B)|_{{\mathcal C}}-c_{1}({\mathcal C}/B).
\end{equation}
By \eqref{eqn:relation:dual:g:action}, we get
\begin{equation}
\label{eqn:Chern:form:eigen:bundle}
c_{1}(N^{*}_{{\mathcal C}/{\mathcal X}}(\theta),h_{N_{{\mathcal C}/{\mathcal X}}(\theta)})=-{\frak l}_{2},
\qquad
c_{1}(N^{*}_{{\mathcal C}/{\mathcal X}}(-\theta),h_{N_{{\mathcal C}/{\mathcal X}}(-\theta)})=-{\frak l}_{1}.
\end{equation}

\begin{lemma}
\label{lemma:equivariant:Todd:chern:character}
For $g\in G_{\mathcal C}^{*}$, set
$$
\zeta_{g}:=e^{i\theta}.
$$
Then the following equality of $(2,2)$-forms on ${\mathcal C}$ holds:
$$
\begin{aligned}
\,&
[
{\rm Td}_{g}(\Theta_{{\mathcal X}/B},h_{{\mathcal X}/B})
\sum_{p\geq0}(-1)^{p}p\,{\rm ch}_{g}(\Omega^{p}_{{\mathcal X}/B},h_{\Omega^{p}_{{\mathcal X}/B}})|_{\mathcal C}
]^{(2,2)}
\\
&=
\frac{\zeta_{g}}{(1-\zeta_{g})^{2}}c_{1}({\mathcal X}/B)c_{1}({\mathcal C}/B)
-
\left\{
\frac{1}{12}+\frac{\zeta_{g}}{(1-\zeta_{g})^{2}}
\right\}
c_{1}({\mathcal C}/B)^{2}.
\end{aligned}
$$
\end{lemma}

\begin{pf}
By the definition of equivariant Todd form and \eqref{eqn:Chern:form:eigen:bundle}, we get
\begin{equation}
\label{eqn:equivariant:Todd}
\begin{aligned}
{\rm Td}_{g}(\Theta_{{\mathcal X}/B}, h_{{\mathcal X}/B})|_{\mathcal C}
&=
{\rm Td}(\Theta_{{\mathcal X}^{g}/B},h_{{\mathcal X}^{g}/B})|_{\mathcal C}
\cdot
\frac{1}{1-\zeta_{g}^{-1}e^{-{\frak l}_{1}}}
\cdot
\frac{1}{1-\zeta_{g} e^{-{\frak l}_{2}}}
\\
&=
\frac{{\rm Td}(\Theta_{{\mathcal X}^{g}/B},h_{{\mathcal X}^{g}/B})|_{\mathcal C}}
{\sum_{p\geq0}(-1)^{p}{\rm ch}_{g}(\Lambda^{p}N^{*}_{{\mathcal C}/{\mathcal X}},h_{\Lambda^{p}N^{*}_{\mathcal C}/{\mathcal X}})}.
\end{aligned}
\end{equation}
Here $h_{\Lambda^{p}N^{*}_{{\mathcal C}/{\mathcal X}}}$ is the Hermitian metric on $\Lambda^{p}N^{*}_{{\mathcal C}/{\mathcal X}}$
induced by $h_{N_{{\mathcal C}/{\mathcal X}}}$.
By \eqref{eqn:equivariant:Todd} and Maillot-R\"ossler \cite[Lemma 3.1]{MaillotRoessler04}, we get
\begin{equation}
\label{eqn:prod:equiv:Todd:sumChern}
\begin{aligned}
\,&
[
{\rm Td}_{g}(\Theta_{{\mathcal X}/B},h_{{\mathcal X}/B})
\sum_{q\geq0}(-1)^{q}q\,{\rm ch}_{g}(\Omega^{q}_{{\mathcal X}/B},h_{\Omega^{q}_{{\mathcal X}/B}})|_{\mathcal C}
]^{(2,2)}
\\
&=
\left[
{\rm Td}(\Theta_{{\mathcal X}^{g}/B},h_{{\mathcal X}^{g}/B})|_{\mathcal C}
\frac{\sum_{q\geq0}(-1)^{q}q\,{\rm ch}_{g}(\Omega^{q}_{{\mathcal X}/B},h_{\Omega^{q}_{{\mathcal X}/B}})|_{\mathcal C}}
{\sum_{p\geq0}(-1)^{p}{\rm ch}_{g}(\Lambda^{p}N^{*}_{{\mathcal C}/{\mathcal X}},h_{\Lambda^{p}N^{*}_{{\mathcal C}/{\mathcal X}}})}
\right]^{(2,2)}
\\
&=
-c_{1}({\mathcal C}/B)
\wedge
[
\zeta_{L}(\zeta_{g},-1)\,e^{-{\frak l}_{2}}
+
\zeta_{L}(\zeta_{g}^{-1},-1)\,e^{-{\frak l}_{1}}
+
\zeta_{\bf Q}(-1)\,e^{-c_{1}({\mathcal C}/B)}
]^{(1,1)}
\\
&=
c_{1}({\mathcal C}/B)
\wedge
\left\{
\zeta_{L}(\zeta_{g},-1){\frak l}_{2}
+
\zeta_{L}(\zeta_{g}^{-1},-1){\frak l}_{1}
+
\zeta_{\bf Q}(-1)\,c_{1}({\mathcal C}/B)
\right\},
\end{aligned}
\end{equation}
where $\zeta_{L}(z,s)=\sum_{n\geq1}z^{n}/n^{s}$ is the Lerch $\zeta$-function and $\zeta_{\bf Q}(s)$ is the Riemann $\zeta$-function.
Substituting
 $\zeta_{L}(\zeta_{g},-1)=\zeta_{L}(\zeta_{g}^{-1},-1)=\zeta_{g}/(1-\zeta_{g})^{2}$ (cf. \cite[p.744]{MaillotRoessler04})
and $\zeta_{\bf Q}(-1)=-1/12$ into \eqref{eqn:prod:equiv:Todd:sumChern} and using \eqref{eqn:relation:sum:Chern:forms}, 
we get the result.
\end{pf}

\subsubsection
{An equivariant characteristic form: case of relative dimension $0$}
\label{sect:4.3.5}
\par
Let ${\frak p}\subset{\mathcal X}^{(0)}$ be a component. Since $B\cong\varDelta^{\dim B}$ is contractible,
$$
{\frak p}\colon B\ni b\to{\frak p}\cap X_{b}\in X_{b}
$$ 
is viewed as a section of the map $f\colon{\mathcal X}\to B$.
As in \eqref{eqn:abelian:subgroup:SL(3,C)}, we can express
$$
G_{\frak p}=\{{\rm diag}(\chi_{{\frak p},1}(g),\chi_{{\frak p},2}(g),\chi_{{\frak p},3}(g));\,g\in G_{\frak p}\},
\qquad
\chi_{{\frak p},k}\in{\rm Hom}(G_{\frak p},{\bf C}^{*}).
$$
Then the normal bundle $N_{{\frak p}/{\mathcal X}}$ splits into the eigenbundles of the $G_{\frak p}$-action
$$
N_{{\frak p}/{\mathcal X}}
=
N_{{\frak p}/{\mathcal X}}(\chi_{{\frak p},1})
\oplus 
N_{{\frak p}/{\mathcal X}}(\chi_{{\frak p},2})
\oplus 
N_{{\frak p}/{\mathcal X}}(\chi_{{\frak p},3}).
$$
We get the corresponding splitting of the conormal bundle $N^{*}_{{\frak p}/{\mathcal X}}$
$$
N^{*}_{{\frak p}/{\mathcal X}}
=
N^{*}_{{\frak p}/{\mathcal X}}(\chi_{{\frak p},1}^{-1})
\oplus 
N^{*}_{{\frak p}/{\mathcal X}}(\chi_{{\frak p},2}^{-1})
\oplus 
N^{*}_{{\frak p}/{\mathcal X}}(\chi_{{\frak p},3}^{-1})
$$
such that $N^{*}_{{\frak p}/{\mathcal X}}(\chi_{{\frak p},k}^{-1})=(N_{{\frak p}/{\mathcal X}}(\chi_{{\frak p},k}))^{*}$.
For $k=1,2,3$, set 
\begin{equation}
\label{eqn:***}
\rho_{{\frak p},k}:=c_{1}(N_{{\frak p}/{\mathcal X}}(\chi_{{\frak p},k}),h_{N_{{\frak p}/{\mathcal X}}(\chi_{{\frak p},k})}),
\end{equation}
where $h_{N_{{\frak p}/{\mathcal X}}(\chi_{{\frak p},k})}:=h_{N_{{\frak p}/{\mathcal X}}}|_{N_{{\frak p}/{\mathcal X}}(\chi_{{\frak p},k})}$.
Then we have the following relations
$$
c_{1}(N^{*}_{{\frak p}/{\mathcal X}}(\chi_{{\frak p},k}^{-1}),h_{N^{*}_{{\frak p}/{\mathcal X}}(\chi_{{\frak p},k}^{-1})})=-\rho_{{\frak p},k},
$$
\begin{equation}
\label{eqn:relation:sum:Chern:forms:relative:dim:0}
\rho_{{\frak p},1}+\rho_{{\frak p},2}+\rho_{{\frak p},3}
=
c_{1}(N_{{\frak p}/{\mathcal X}},h_{N_{{\frak p}/{\mathcal X}}})
=
c_{1}({\mathcal X}/B)|_{\frak p}.
\end{equation}

\begin{lemma}
\label{lemma:equivariant:Todd:chern:character:relative:dim=0}
{\rm (1) }
For $g\in G_{{\frak p}}^{0}$, the following equality of $(1,1)$-forms on ${\frak p}$ holds:
$$
[
{\rm Td}_{g}(\Theta_{{\mathcal X}/B},h_{{\mathcal X}/B})
\sum_{p\geq0}(-1)^{p}p\,{\rm ch}_{g}(\Omega^{p}_{{\mathcal X}/B},h_{\Omega^{p}_{{\mathcal X}/B}})|_{\frak p}
]^{(1,1)}
=
\sum_{k=1}^{3}
\frac{\chi_{{\frak p},k}(g)}{(1-\chi_{{\frak p},k}(g))^{2}}\rho_{{\frak p},k}.
$$
{\rm (2) }
For any component ${\frak p}\subset{\mathcal X}^{(0)}$, one has
$$
[\sum_{g\in G_{\frak p}^{0}}
{\rm Td}_{g}(\Theta_{{\mathcal X}/B},h_{{\mathcal X}/B})
\sum_{p\geq0}(-1)^{p}p\,{\rm ch}_{g}(\Omega^{p}_{{\mathcal X}/B},h_{\Omega^{p}_{{\mathcal X}/B}})|_{\frak p}]^{(1,1)}
=
\sum_{k=1}^{3}\delta_{k}(G_{\frak p})\,\rho_{{\frak p},k}.
$$
\end{lemma}

\begin{pf}
{\bf (1) }
Set $\zeta_{k}:=\chi_{{\frak p},k}(g)$. By the definition of equivariant Todd form, 
$$
{\rm Td}_{g}(\Theta_{{\mathcal X}/B}, h_{{\mathcal X}/B})|_{\frak p}
=
\prod_{k=1}^{3}
\frac{1}{1-\zeta_{k}^{-1}e^{-\rho_{{\frak p},k}}}
=
\frac{1}{\sum_{p\geq0}(-1)^{p}{\rm ch}_{g}(\Lambda^{p}N^{*}_{{\frak p}/{\mathcal X}},h_{\Lambda^{p}N^{*}_{{\frak p}/{\mathcal X}}})}.
$$
By \cite[Lemma 3.1]{MaillotRoessler04}, we get
$$
\begin{aligned}
\,&
[
{\rm Td}_{g}(\Theta_{{\mathcal X}/B},h_{{\mathcal X}/B})
\sum_{p\geq0}(-1)^{p}p\,{\rm ch}_{g}(\Omega^{p}_{{\mathcal X}/B},h_{\Omega^{p}_{{\mathcal X}/B}})|_{\frak p}
]^{(1,1)}
\\
&=
\left[
\frac{\sum_{p\geq0}(-1)^{p}p\,{\rm ch}_{g}(\Omega^{p}_{{\mathcal X}/B},h_{\Omega^{p}_{{\mathcal X}/B}})|_{\frak p}}
{\sum_{p\geq0}(-1)^{p}{\rm ch}_{g}(\Lambda^{p}N^{*}_{{\frak p}/{\mathcal X}},h_{\Lambda^{p}N^{*}_{{\frak p}/{\mathcal X}}})}
\right]^{(1,1)}
=
-\sum_{k=1}^{3}
\zeta_{L}(\zeta_{k}^{-1},-1)\,e^{-\rho_{{\frak p},k}}.
\end{aligned}
$$
Since $\zeta_{L}(\zeta_{k}^{-1},-1)=\zeta_{k}/(1-\zeta_{k})^{2}$, we get the desired equality.
\par{\bf (2) }
When $G_{{\frak p}}^{0}\not=\emptyset$, the assertion follows from (1). 
When $G_{{\frak p}}^{0}=\emptyset$, the assertion also holds because we defined $\delta_{k}(G_{\frak p})=0$ $(k=1,2,3)$ in this case.
\end{pf}

\subsection
{An orbifold characteristic forms on ${\mathcal Y}\amalg{\mathcal Y}^{(1)}\amalg{\mathcal Y}^{(0)}$}
\label{sect:4.4}
\par 
Let 
\begin{center}
$\pi\colon{\mathcal X}\to{\mathcal Y}$,
\qquad
$\pi_{\lambda}\colon{\mathcal C}_{\lambda}\to\overline{\mathcal C}_{\bar{\lambda}}$
\quad
$(\lambda\in\Lambda)$,
\qquad
$\pi_{\frak p}\colon{\frak p}\to\overline{\frak p}$
\quad
(${\frak p}\subset{\mathcal X}^{(0)}$)
\end{center}
be the projections. 
Let $(E,h)$ be a $G$-equivariant holomorphic vector bundle on ${\mathcal X}$ with $G$-invariant Hermitian metric.
Then $\phi^{\Sigma}(\overline{E},\overline{h})$ is a form on ${\mathcal Y}\amalg\Sigma{\mathcal Y}$.
Let
$$
\varPi\colon{\mathcal Y}\amalg\Sigma{\mathcal Y}\to{\mathcal Y}\amalg{\mathcal Y}^{(1)}\amalg{\mathcal Y}^{(0)}
$$
be the projection. Then 
$\varPi_{*}\phi^{\Sigma}(\overline{E},\overline{h})\in A^{*}({\mathcal Y}\amalg{\mathcal Y}^{(1)}\amalg{\mathcal Y}^{(0)})$.
Since 
\begin{equation}
\label{eqn:decomposition:fixed:point}
{\mathcal X}^{g}
=
(\amalg_{\{\lambda\in\Lambda;\,g\in G_{{\mathcal C}_{\lambda}}^{*}\}}{\mathcal C}_{\lambda})
\amalg
(\amalg_{\{{\frak p}\subset\nu(\Sigma^{(0)}{\mathcal X});g\in G_{\frak p}^{0}\}}\{{\frak p}\})
\end{equation}
for any $g\in G^{*}$, it follows from \eqref{eqn:orb:char:form:2}, \eqref{eqn:decomposition:fixed:point} that
$$
\begin{aligned}
\varPi_{*}\phi^{\Sigma}(\overline{E},\overline{h})
&=
\frac{1}{|G|}\pi_{*}\phi(E,h)
+
\frac{1}{|G|}\sum_{g\in G^{*}}\sum_{\{\lambda\in\Lambda;\,g\in G_{{\mathcal C}_{\lambda}}^{*}\}}(\pi_{\lambda})_{*}(\phi_{g}(E,h)|_{{\mathcal C}_{\lambda}})
\\
&\quad+
\frac{1}{|G|}\sum_{g\in G^{*}}\sum_{\{{\frak p}\subset\nu(\Sigma^{(0)}{\mathcal X});\,g\in G_{\frak p}^{0}\}}(\pi_{\frak p})_{*}(\phi_{g}(E,h)|_{\frak p}).
\end{aligned}
$$
Since $G_{\frak p}^{0}=\emptyset$ when $G_{\frak p}={\Bbb G}$, we get $(\pi_{\frak p})_{*}(\phi_{g}(E,h)|_{\frak p})=0$ in this case and hence
\begin{equation}
\label{eqn:orb:char:form:3}
\begin{aligned}
\varPi_{*}\phi^{\Sigma}(\overline{E},\overline{h})
&=
\frac{1}{|G|}\pi_{*}\phi(E,h)
+
\frac{1}{|G|}\sum_{\lambda\in\Lambda}\sum_{g\in G_{{\mathcal C}_{\lambda}}^{*}}(\pi_{\lambda})_{*}(\phi_{g}(E,h)|_{{\mathcal C}_{\lambda}})
\\
&\quad+
\frac{1}{|G|}\sum_{{\frak p}\subset{\mathcal X}^{(0)}}\sum_{g\in G_{\frak p}^{0}}(\pi_{\frak p})_{*}(\phi_{g}(E,h)|_{\frak p}).
\end{aligned}
\end{equation}
By \eqref{eqn:orb:char:form:3}, we get
\begin{equation}
\label{eqn:Bismut:Gillet:Soule:Ma:2}
\begin{aligned}
\,&
\varPi_{*}{\rm Td}^{\Sigma}(\Theta_{{\mathcal Y}/B},h_{{\mathcal Y}/B})
\sum_{p\geq0}(-1)^{p}p\,{\rm ch}^{\Sigma}(\Omega_{{\mathcal Y}/B}^{p},h_{\Omega^{p}_{{\mathcal Y}/B}})
\\
&=
\frac{1}{|G|}
\pi_{*}\{
{\rm Td}(\Theta_{{\mathcal X}/B},h_{{\mathcal X}/B})
\sum_{p\geq0}(-1)^{p}p\,{\rm ch}(\Omega^{p}_{{\mathcal X}/B},h_{\Omega^{p}_{{\mathcal X}/B}})
\}
\\
&\quad
+\frac{1}{|G|}\sum_{\lambda\in\Lambda}\sum_{g\in G_{{\mathcal C}_{\lambda}}^{*}}
(\pi_{\lambda})_{*}\{
{\rm Td}_{g}(\Theta_{{\mathcal X}/B},h_{{\mathcal X}/B})
\sum_{p\geq0}(-1)^{p}p\,{\rm ch}_{g}(\Omega^{p}_{{\mathcal X}/B},h_{\Omega^{p}_{{\mathcal X}/B}})|_{{\mathcal C}_{\lambda}}
\}
\\
&\quad
+\frac{1}{|G|}\sum_{{\frak p}\subset{\mathcal X}^{(0)}}\sum_{g\in G_{\frak p}^{0}}
(\pi_{\frak p})_{*}\{
{\rm Td}_{g}(\Theta_{{\mathcal X}/B},h_{{\mathcal X}/B})
\sum_{p\geq0}(-1)^{p}p\,{\rm ch}_{g}(\Omega^{p}_{{\mathcal X}/B},h_{\Omega^{p}_{{\mathcal X}/B}})|_{\frak p}
\}.
\end{aligned}
\end{equation}

Since the projection $\pi_{\frak p}\colon{\frak p}\to\overline{\frak p}$ is an isomorphism, 
we identify the $(1,1)$-form $\rho_{{\frak p},k}$ with $(\pi_{\frak p})_{*}(\rho_{{\frak p},k})$.
For $\overline{\frak p}\subset{\mathcal Y}^{(0)}$ and $k\in\{1,2,3\}$, we set
$$
\overline{\rho}_{\overline{\frak p},k}:=(\pi_{\frak p})_{*}(\rho_{{\frak p},k}).
$$
Then $\overline{\rho}_{\overline{\frak p},k}$ is a smooth $(1,1)$-form on $\overline{\frak p}$.
Here $\overline{\frak p}$ is not viewed as an orbifold (of multiplicity $|G_{\frak p}|$)
but the complex manifold underlying it isomorphic to $B$.
If ${\frak p},{\frak p}'\in{\mathcal X}^{(0)}$ lie on the same $G$-orbit, we have
$\overline{\rho}_{\overline{\frak p},k}=\overline{\rho}_{\overline{\frak p}',k}$.

\begin{definition}
\label{def:notation:relative:degree}
For a differential form $\omega$ on ${\mathcal Y}\amalg{\mathcal Y}^{(1)}\amalg{\mathcal Y}^{(0)}$, 
$[\omega]^{2\dim\overline{f}_{*}+(p,p)}$ is the differential form on ${\mathcal Y}\amalg{\mathcal Y}^{(1)}\amalg{\mathcal Y}^{(0)}$ defined as
$$
[\omega]^{2\dim\overline{f}_{*}+(p,p)}
:=
[\omega|_{\mathcal Y}]^{(p+3,p+3)}
+
\sum_{\bar{\lambda}\in\overline{\Lambda}}
[\omega|_{\overline{\mathcal C}_{\bar{\lambda}}}]^{(p+1,p+1)}
+
\sum_{\overline{\frak p}\subset{\mathcal Y}^{(0)}}
[\omega|_{\overline{\frak p}}]^{(p,p)}.
$$
\end{definition}

\begin{proposition}
\label{prop:orb:Tod:ch}
The following equality of forms on ${\mathcal Y}\amalg{\mathcal Y}^{(1)}\amalg{\mathcal Y}^{(0)}$ holds
$$
\begin{aligned}
\,&
[
\varPi_{*}
\{
{\rm Td}^{\Sigma}(\Theta_{{\mathcal Y}/B},h_{{\mathcal Y}/B})
\sum_{p\geq0}(-1)^{p}p\,{\rm ch}^{\Sigma}(\Omega_{{\mathcal Y}/B}^{p},h_{\Omega^{p}_{{\mathcal Y}/B}})
\}
]^{2\dim\overline{f}_{*}+(1,1)}
\\
&=
-\frac{1}{12}c_{1}({\mathcal Y}/B)c_{3}({\mathcal Y}/B)
-
\frac{1}{12}
\sum_{\bar{\lambda}\in\overline{\Lambda}}\frac{n_{\lambda}^{2}-1}{n_{\lambda}}
c_{1}({\mathcal Y}/B)c_{1}(\overline{\mathcal C}_{\lambda}/B)
|_{\overline{\mathcal C}_{\bar\lambda}}
\\
&\quad
+
\frac{1}{12}\sum_{\bar{\lambda}\in\overline{\Lambda}}
(n_{\lambda}-1)\,c_{1}(\overline{\mathcal C}_{\lambda}/B)^{2}
|_{\overline{\mathcal C}_{\bar\lambda}}
+
\sum_{\overline{\frak p}\subset{\mathcal Y}^{(0)}}\sum_{k=1}^{3}
\frac{\delta_{k}(G_{{\frak p}})}{|G_{{\frak p}}|}\,\rho_{\overline{\frak p},k}.
\end{aligned}
$$
\end{proposition}

\begin{pf}
Substituting the equalities in Lemmas~\ref{lemma:Todd:chern:character:BCOV}, \ref{lemma:equivariant:Todd:chern:character}, 
\ref{lemma:equivariant:Todd:chern:character:relative:dim=0} into the right hand side of \eqref{eqn:Bismut:Gillet:Soule:Ma:2}, 
we get
\begin{equation}
\label{eqn:Bismut:Gillet:Soule:Ma:3}
\begin{aligned}
\,&
[
\varPi_{*}\{
{\rm Td}^{\Sigma}(\Theta_{{\mathcal Y}/B},h_{{\mathcal Y}/B})
\sum_{p\geq0}(-1)^{p}p\,{\rm ch}^{\Sigma}(\Omega_{{\mathcal Y}/B}^{p},h_{\Omega^{p}_{{\mathcal Y}/B}})
\}
]^{2\dim\overline{f}_{*}+(1,1)}
=
\\
&
-\frac{1}{12|G|}
\pi_{*}\{
c_{1}({\mathcal X}/B)c_{3}({\mathcal X}/B)
\}
+
\frac{1}{|G|}\sum_{\lambda\in\Lambda}
\sum_{g\in G_{{\mathcal C}_{\lambda}}^{*}}\frac{\zeta_{g}}{(1-\zeta_{g})^{2}}
(\pi_{\lambda})_{*}
\{
c_{1}({\mathcal X}/B)c_{1}({\mathcal C}_{\lambda}/B)
\}
\\
&
-
\frac{1}{|G|}\sum_{\lambda\in\Lambda}
(\sum_{g\in G_{{\mathcal C}_{\lambda}}^{*}}
\left\{
\frac{1}{12}+\frac{\zeta_{g}}{(1-\zeta_{g})^{2}}
\right\})
(\pi_{\lambda})_{*}
\{
c_{1}({\mathcal C}_{\lambda}/B)^{2}
\}
+
\sum_{{\frak p}\subset{\mathcal X}^{(0)}}\sum_{k=1}^{3}
\frac{\delta_{k}(G_{\frak p})}{|G|}
(\pi_{\frak p})_{*}{\rho}_{{\frak p},k}.
\end{aligned}
\end{equation}
Recall the following classical formula for the Dedekind sum
\begin{equation}
\label{eqn:Dedekind:sum}
\sum_{k=1}^{n-1}\frac{\zeta^{k}}{(1-\zeta^{k})^{2}}=-\frac{n^{2}-1}{12},
\qquad
\zeta=\exp(2\pi i\,m/n),
\quad
(m,n)=1.
\end{equation}
Since $G_{{\mathcal C}_{\lambda}}$ is a cyclic group of order $n_{\lambda}$, 
we get by \eqref{eqn:Bismut:Gillet:Soule:Ma:3}, \eqref{eqn:Dedekind:sum}
\begin{equation}
\label{eqn:Bismut:Gillet:Soule:Ma:4}
\begin{aligned}
\,&
[
{\rm Td}^{\Sigma}(\Theta_{{\mathcal Y}/B},h_{{\mathcal Y}/B})
\sum_{p\geq0}(-1)^{p}p\,{\rm ch}^{\Sigma}(\Omega_{{\mathcal Y}/B}^{p},h_{\Omega^{p}_{{\mathcal Y}/B}})
]^{2\dim\overline{f}_{*}+(1,1)}
\\
&=
-\frac{1}{12|G|}\pi_{*}\{c_{1}({\mathcal X}/B)c_{3}({\mathcal X}/B)\}
-
\frac{1}{12|G|}
\sum_{\lambda\in\Lambda}(n_{\lambda}^{2}-1)
(\pi_{\lambda})_{*}
\{
c_{1}({\mathcal X}/B)c_{1}({\mathcal C}_{\lambda}/B)
\}
\\
&\quad
+
\frac{1}{12|G|}
\sum_{\lambda\in\Lambda}
\left(n_{\lambda}^{2}-n_{\lambda}
\right)
(\pi_{\lambda})_{*}
\{
c_{1}({\mathcal C}_{\lambda}/B)^{2}
\}
+
\sum_{{\frak p}\subset{\mathcal X}^{(0)}}\sum_{k=1}^{3}
\frac{\delta_{k}(G_{\frak p})}{|G|}\,\overline{\rho}_{\bar{\frak p},k}.
\end{aligned}
\end{equation}
Some terms in the right hand side of \eqref{eqn:Bismut:Gillet:Soule:Ma:4} are calculated as follows:
\begin{equation}
\label{eqn:relation:char:forms:X:Y:1}
\begin{aligned}
\,&
\sum_{\lambda\in\Lambda}(n_{\lambda}^{2}-1)
(\pi_{\lambda})_{*}
\{
c_{1}({\mathcal X}/B)c_{1}({\mathcal C}_{\lambda}/B)
\}
\\
&=
\sum_{\bar{\lambda}\in\overline{\Lambda}}(n_{\lambda}^{2}-1)
\sum_{g\in G/\Gamma_{{\mathcal C}_{\lambda}}}
(\pi_{\lambda})_{*}
\{
c_{1}({\mathcal X}/B)c_{1}({\mathcal C}_{g\cdot\lambda}/B)
\}
\\
&=
\sum_{\bar{\lambda}\in\overline{\Lambda}}
(n_{\lambda}^{2}-1)\cdot|G/\Gamma_{{\mathcal C}_{\lambda}}|\cdot\deg\pi_{\lambda}
\cdot
c_{1}({\mathcal Y}/B)c_{1}(\overline{\mathcal C}_{\bar{\lambda}}/B)
\\
&=
|G|\sum_{\bar{\lambda}\in\overline{\Lambda}}\frac{n_{\lambda}^{2}-1}{n_{\lambda}}
c_{1}({\mathcal Y}/B)c_{1}(\overline{\mathcal C}_{\bar{\lambda}}/B),
\end{aligned}
\end{equation}
\begin{equation}
\label{eqn:relation:char:forms:X:Y:2}
\begin{aligned}
\,&
\sum_{\lambda\in\Lambda}(n_{\lambda}^{2}-n_{\lambda})
(\pi_{\lambda})_{*}
\{
c_{1}({\mathcal C}_{\lambda}/B)^{2}
\}
=
\sum_{\bar{\lambda}\in\overline{\Lambda}}(n_{\lambda}^{2}-n_{\lambda})
\sum_{g\in G/\Gamma_{{\mathcal C}_{\lambda}}}
(\pi_{\lambda})_{*}
\{
c_{1}({\mathcal C}_{g\cdot\lambda}/B)^{2}
\}
\\
&=
\sum_{\bar{\lambda}\in\overline{\Lambda}}
(n_{\lambda}^{2}-n_{\lambda})\cdot|G/\Gamma_{{\mathcal C}_{\lambda}}|\cdot\deg\pi_{\lambda}
\cdot
c_{1}(\overline{\mathcal C}_{\bar\lambda}/B)^{2}
=
|G|\sum_{\bar{\lambda}\in\overline{\Lambda}}(n_{\lambda}-1)
\,
c_{1}(\overline{\mathcal C}_{\bar\lambda}/B)^{2},
\end{aligned}
\end{equation}
\begin{equation}
\label{eqn:relation:char:forms:X:Y:3}
\sum_{{\frak p}\subset{\mathcal X}^{(0)}}\sum_{k=1}^{3}\delta_{k}(G_{\frak p})\overline{\rho}_{\bar{\frak p},k}
=
\sum_{\overline{\frak p}\subset{\mathcal Y}^{(0)}}|G/G_{\frak p}|
\sum_{k=1}^{3}\delta_{k}(G_{\frak p})\overline{\rho}_{\bar{\frak p},k}
=
|G|\sum_{\overline{\frak p}\subset{\mathcal Y}^{(0)}}
\sum_{k=1}^{3}\frac{\delta_{k}(G_{\frak p})}{|G_{\frak p}|}\overline{\rho}_{\overline{\frak p},k}
\end{equation}
where we used 
$\deg\pi_{\lambda}=|\overline{\Gamma}_{{\mathcal C}_{\lambda}}|=|\Gamma_{{\mathcal C}_{\lambda}}/G_{{\mathcal C}_{\lambda}}|$
to get \eqref{eqn:relation:char:forms:X:Y:1}, \eqref{eqn:relation:char:forms:X:Y:2}.
Substituting the equality $\pi_{*}\{c_{1}({\mathcal X}/B)c_{3}({\mathcal X}/B)\}=|G|\cdot c_{1}({\mathcal Y}/B)c_{3}({\mathcal Y}/B)$ 
and \eqref{eqn:relation:char:forms:X:Y:1}, \eqref{eqn:relation:char:forms:X:Y:2}, \eqref{eqn:relation:char:forms:X:Y:3}
into \eqref{eqn:Bismut:Gillet:Soule:Ma:4}, we get the result.
\end{pf}

\subsection
{An orbifold characteristic form on $\overline{\mathcal C}_{\bar{\lambda}}\amalg\Sigma\overline{\mathcal C}_{\bar{\lambda}}$}
\label{sect:4.5}
\par
The relative tangent bundle $\Theta_{\overline{\mathcal C}_{\bar{\lambda}}/B}$
of the map $\overline{f}\colon\overline{\mathcal C}_{\bar{\lambda}}\to B$
is equipped with the Hermitian metric $h_{\overline{\mathcal C}_{\bar{\lambda}}/B}$ induced from $h_{{\mathcal X}/B}$.
As before, the Chern form $c_{1}(\overline{\mathcal C}_{\bar{\lambda}}/B,h_{\overline{\mathcal C}_{\bar{\lambda}}/B})$ is denoted by 
$c_{1}(\overline{\mathcal C}_{\bar{\lambda}}/B)$.

\subsubsection
{The structure of $\Sigma{\mathcal C}_{\lambda}$}
\label{sect:4.5.1}
\par
Let $b\in B$. 
In Section~\ref{sect:4.5.1}, for the sake of simplicity, we write $C_{\lambda}$ and $X$ for $C_{b,\lambda}$ and $X_{b}$, respectively.
Let ${\frak q}\in C_{\lambda}$. Then $G_{{\mathcal C}_{\lambda}}=G_{C_{\lambda}}\subset G_{\frak q}$.
Write $(\overline{\Gamma}_{C_{\lambda}})_{\frak q}$ for the stabilizer of ${\frak q}$ in $\overline{\Gamma}_{C_{\lambda}}$. 
By definition,
\begin{equation}
\label{eqn:stabilizer:p:alpha:1}
(\overline{\Gamma}_{C_{\lambda}})_{\frak q}
=
(\Gamma_{C_{\lambda}}\cap G_{\frak q})/G_{C_{\lambda}}.
\end{equation}
By \eqref{eqn:abelian:subgroup:SL(3,C)}, there is a coordinate neighborhood $(U,(z_{1},z_{2},z_{3}))$ of $X$ 
centered at ${\frak q}$ such that $G_{\frak q}=\{{\rm diag}(\chi_{1}(g),\chi_{2}(g),\chi_{3}(g));\,g\in G_{\frak q}\}\subset{\rm SL}({\bf C}^{3})$, 
where $\chi_{i}\in{\rm Hom}(G_{\frak q},{\bf C}^{*})$ and $\chi_{1}\chi_{2}\chi_{3}=1$.
A system of coordinates with this property is said to be {\em normal}.
Let $g_{0}\in G_{C_{\lambda}}$ be a generator. 
Since $G_{C_{\lambda}}\subset G_{\frak q}$ and $\dim C_{\lambda}=1$, 
we have the following in the normal coordinates $(U,(z_{1},z_{2},z_{3}))$ centered at ${\frak q}$:
\begin{itemize}
\item[(a)]
$C_{\lambda}$ is given by the eigenspace of ${\rm diag}(\chi_{1}(g_{0}),\chi_{2}(g_{0}),\chi_{3}(g_{0}))$ 
corresponding to the eigenvalue $1$. 
\item[(b)]
There exists $k_{\lambda}\in\{1,2,3\}$ such that 
$\chi_{k_{\lambda}}(g_{0})=1$, $\chi_{l_{\lambda}}(g_{0})\not=1$, $\chi_{m_{\lambda}}(g_{0})\not=1$, 
where $\{k_{\lambda},l_{\lambda},m_{\lambda}\}=\{1,2,3\}$.
\item[(c)]
The germ $(C_{\lambda})_{\frak q}$ is given by the coordinate axis ${\frak C}_{k_{\lambda}}=\{z_{l_{\lambda}}=z_{m_{\lambda}}=0\}$.
\end{itemize}
Since every element of $G_{\frak q}$ is of the form ${\rm diag}(\chi_{1}(g),\chi_{2}(g),\chi_{3}(g))$ and hence
$G_{\frak q}$ preserves all the coordinate axes of the normal coordinates, we get
\begin{equation}
\label{eqn:stabilizer:p:alpha:2}
G_{\frak q}\subset\Gamma_{C_{\lambda}}.
\end{equation} 
By \eqref{eqn:stabilizer:p:alpha:1}, \eqref{eqn:stabilizer:p:alpha:2}, we get
\begin{equation}
\label{eqn:stabilizer:p:alpha:3}
(\overline{\Gamma}_{C_{\lambda}})_{\frak q}
=
G_{\frak q}/G_{C_{\lambda}}.
\end{equation}
With respect to the $\overline{\Gamma}_{C_{\lambda}}$-action on $C_{\lambda}$, we consider the set $\Sigma C_{\lambda}$ 
as in Section~\ref{sect:4.1.2}.
By definition, 
$\Sigma C_{\lambda}=\{({\frak q},g)\in C_{\lambda}\times(\overline{\Gamma}_{C_{\lambda}})^{*};\,
g\in(\overline{\Gamma}_{C_{\lambda}})_{\frak q}^{*}\}$.

\begin{proposition}
\label{prop:stabilizer:curve}
One has $\nu(\Sigma C_{\lambda})=C_{\lambda}\cap X^{(0)}$.
\end{proposition}

\begin{pf}
Let ${\frak q}\in\Sigma C_{\lambda}$. Let $(z_{1},z_{2},z_{3})$ be normal coordinates of $X$ centered at ${\frak q}$
satisfying (a), (b), (c). We may assume $k_{\lambda}=3$. By (a), (c), $g_{0}={\rm diag}(\chi_{1}(g_{0}),\chi_{1}(g_{0})^{-1},1)$ and
$(C_{\lambda})_{\frak q}=\{z_{1}=z_{2}=0\}$.
Since $(\overline{\Gamma}_{C_{\lambda}})_{\frak q}\not=\{1\}$, we can take $g\in G_{\frak q}\setminus G_{C_{\lambda}}$.
If $\chi_{3}(g)=1$, then $g$ acts trivially on $(C_{\lambda})_{\frak q}$, so that $g\in G_{C_{\lambda}}$. 
This contradicts the choice of $g$. Thus $\chi_{3}(g)\not=1$.
Since $G_{\frak q}$ contains $g={\rm diag}(\chi_{1}(g),\chi_{2}(g),\chi_{3}(g))$ and 
$g_{0}={\rm diag}(\chi_{1}(g_{0}),\chi_{1}(g_{0})^{-1},1)\in G_{C_{\lambda}}^{*}$
with $\chi_{3}(g)\not=1$ and $\chi_{1}(g_{0})\not=1$, we get $X^{G_{\frak q}^{*}}\subset X^{g}\cap X^{g_{0}}=\{{\frak q}\}$ near ${\frak q}$. 
Hence $X^{G_{\frak q}^{*}}=\{{\frak q}\}$. By \eqref{eqn:ext:twisted:sector:1}, we get ${\frak q}\in X^{(0)}$.
\par
Conversely, let ${\frak q}\in C_{\lambda}\cap X^{(0)}$.
Let $(z_{1},z_{2},z_{3})$ be normal coordinates of $X$ centered at ${\frak q}$ satisfying (a), (b), (c) and assume $k_{\lambda}=3$.
Hence $C_{\lambda}$ is the $x_{3}$-axis.
Since $g={\rm diag}(\chi_{1}(g),\chi_{2}(g),\chi_{3}(g))$ for all $g\in G_{\frak q}$, we have $G_{C_{\lambda}}=\ker\chi_{3}$.
Then we get $(\overline{\Gamma}_{C_{\lambda}})_{\frak q}\cong{\rm Im}\,\chi_{3}$ by \eqref{eqn:stabilizer:p:alpha:3}.
If ${\rm Im}\,\chi_{3}=\{1\}$, then $G_{\frak q}=\ker\chi_{3}=G_{C_{\lambda}}$, so that $X^{G_{\frak q}^{*}}\supset C_{\lambda}$.
By comparing \eqref{eqn:ext:twisted:sector:1}, this contradicts ${\frak q}\in X^{(0)}$.
Thus $(\overline{\Gamma}_{C_{\lambda}})_{\frak q}\not=\{1\}$ and hence ${\frak q}\in\Sigma C_{\lambda}$.
\end{pf}

\begin{corollary}
\label{cor:twisted:sector:twisted:sector}
Let $x\in X$. Then $G_{x}={\Bbb G}$ iff $x\in\bigcup_{\lambda\in\Lambda}\nu(\Sigma C_{\lambda})\setminus\nu(\Sigma^{(0)}X)$.
\end{corollary}

\begin{pf}
By Proposition~\ref{prop:stabilizer:curve}, 
$\bigcup_{\lambda\in\Lambda}\nu(\Sigma C_{\lambda})\setminus\nu(\Sigma^{(0)}X)
=
\bigcup_{\lambda\in\Lambda}C_{\lambda}\cap\nu(\widetilde{\Sigma}^{(0)}X)\setminus\nu(\Sigma^{(0)}X)
=
\bigcup_{\lambda\in\Lambda}\{x\in C_{\lambda};\,G_{x}={\Bbb G}\}
=
\{x\in X;\,G_{x}={\Bbb G}\}$,
because the condition $G_{x}={\Bbb G}$ implies $x\in C_{\lambda}$ for some $\lambda\in\Lambda$.
\end{pf}

By \eqref{eqn:ext:twisted:sector:1} and Corollary~\ref{cor:twisted:sector:twisted:sector}, we have
$$
{\mathcal X}^{(0)}=\nu(\Sigma^{(0)}{\mathcal X})\cup\bigcup_{\lambda\in\Lambda}\nu(\Sigma{\mathcal C}_{\lambda}),
\qquad
{\mathcal Y}^{(0)}
=
\overline{\nu}(\Sigma^{(0)}{\mathcal Y})
\cup
\bigcup_{\bar{\lambda}\in\overline{\Lambda}}\nu(\Sigma\overline{\mathcal C}_{\bar\lambda}).
$$

\subsubsection
{Eigenbundle of $N_{{\frak p}/{\mathcal X}}$ and the locus $\Sigma{\mathcal X}^{(1)}$}
\label{sect:4.5.2}
\par
Let ${\frak p}\subset{\mathcal X}^{(0)}$ be a component. 
Choosing $B$ small, we identify a small tubular neighborhood of 
${\frak p}$ in ${\mathcal X}$ with a neighborhood of the zero section of $N_{{\frak p}/{\mathcal X}}$.
Let us fix a system of fiber coordinates of $N_{{\frak p}/{\mathcal X}}$, which is fiberwise normal.
Hence $G_{\frak p}\subset{\rm SL}(N_{{\frak p}/{\mathcal X}})$ is expressed as
\begin{equation}
\label{eqn:expression:orbi:group}
G_{\frak p}=\{{\rm diag}(\chi_{{\frak p},1}(g),\chi_{{\frak p},2}(g),\chi_{{\frak p},3}(g));\,g\in G_{\frak p}\},
\qquad
\chi_{{\frak p},1}\chi_{{\frak p},2}\chi_{{\frak p},3}=1
\end{equation}
in this system of fiber coordinates, where $\chi_{{\frak p},k}\in{\rm Hom}(G_{\frak p},{\bf C}^{*})$. We set
\begin{equation}
\label{eqn:order:ker:Im:chi}
n_{{\frak p},k}:=|\ker\chi_{{\frak p},k}|,
\qquad
\nu_{{\frak p},k}:=|{\rm Im}\,\chi_{{\frak p},k}|.
\end{equation}
Let ${\frak C}_{{\frak p},1}$, ${\frak C}_{{\frak p},2}$, ${\frak C}_{{\frak p},3}$ be the germs at ${\frak p}$
of the coordinate axes in the normal coordinates. 
Since ${\frak C}_{{\frak p},k}$ is linear, we identify it with $N_{{\frak p}/{\mathcal X}}(\chi_{{\frak p},k})$,
the eigenbundle of $N_{{\frak p}/{\mathcal X}}$ corresponding to the character $\chi_{{\frak p},k}$.
We set
$G_{{\frak C}_{{\frak p},k}}:=\{g\in G;\,g|_{{\frak C}_{{\frak p},k}}={\rm id}_{{\frak C}_{{\frak p},k}}\}$.
By \eqref{eqn:expression:orbi:group}, we get
\begin{equation}
\label{eqn:stabilizer:coord:axis}
G_{{\frak C}_{{\frak p},k}}=\ker\chi_{{\frak p},k}.
\end{equation}
Since $T{\frak C}_{{\frak p},k}/B=N_{{\frak p}/{\mathcal X}}(\chi_{{\frak p},k})$, we get by \eqref{eqn:***}
\begin{equation}
\label{eqn:Chern:form:coord:axis}
c_{1}({\frak C}_{{\frak p},k}/B)=c_{1}(N_{{\frak p}/{\mathcal X}}(\chi_{{\frak p},k}))={\rho}_{{\frak p},k}
\qquad
(k=1,2,3).
\end{equation}
By \eqref{eqn:axis:fixed:curve},
for any ${\frak p}\in{\mathcal X}^{(0)}$ and $k\in\{1,2,3\}$, we have the following:
\begin{equation}
\label{eqn:criterion:coordinate:axis}
{\frak C}_{{\frak p},k}\subset{\mathcal X}^{(1)}
\Longleftrightarrow
n_{{\frak p},k}>1,
\qquad\qquad
{\frak C}_{{\frak p},k}\not\subset{\mathcal X}^{(1)}
\Longleftrightarrow
n_{{\frak p},k}=1.
\end{equation}
\par
Let us consider the case ${\frak p}\subset{\mathcal C}_{\lambda}\subset{\mathcal X}^{(1)}$,
where ${\mathcal C}_{\lambda}$ is a component of ${\mathcal X}^{(1)}$.
By (c) in Section~\ref{sect:4.5.1}, there exists $k_{\lambda}\in\{1,2,3\}$ such that ${\mathcal C}_{\lambda}={\frak C}_{{\frak p},k_{\lambda}}$ 
as germs at ${\frak p}$. 
Then we have $G_{{\mathcal C}_{\lambda}}=G_{{\frak C}_{{\frak p},k_{\lambda}}}$ and hence
$G_{{\mathcal C}_{\lambda}}=\ker\chi_{{\frak p},k_{\lambda}}$ by \eqref{eqn:stabilizer:coord:axis}.
If ${\mathcal C}_{\lambda}={\frak C}_{{\frak p},k_{\lambda}}$ as germs at ${\frak p}$, we write
$$
\chi_{{\frak p},\lambda}:=\chi_{{\frak p},k_{\lambda}},
\qquad
n_{{\frak p},\lambda}:=n_{{\frak p},k_{\lambda}},
\qquad
\nu_{{\frak p},\lambda}:=\nu_{{\frak p},k_{\lambda}}.
$$
Since $n_{\lambda}=|G_{{\mathcal C}_{\lambda}}|=|\ker\chi_{{\frak p},k_{\lambda}}|=n_{{\frak p},k_{\lambda}}$ in this case, 
we get by \eqref{eqn:Chern:form:coord:axis}, \eqref{eqn:criterion:coordinate:axis} the following key identity:
For any ${\frak p}\subset{\mathcal X}^{(0)}$, one has
\begin{equation}
\label{eqn:key:identity}
\sum_{\{\lambda\in\Lambda;\,{\mathcal C}_{\lambda}\supset{\frak p}\}}
\frac{(n_{\lambda}-1)(\nu_{{\frak p},\lambda}^{2}-1)}{\nu_{{\frak p},\lambda}}\,
c_{1}({\mathcal C}_{\lambda}/B)|_{\frak p}
=
\sum_{k=1}^{3}
\frac{(n_{{\frak p},k}-1)(\nu_{{\frak p},k}^{2}-1)}{\nu_{{\frak p},k}}\,
\rho_{{\frak p},k}.
\end{equation}

\subsubsection
{The orbifold Todd form of the family $\overline{f}\colon\overline{\mathcal C}_{\bar{\lambda}}\to B$}
\label{sect:4.5.3}
\par
As in Section~\ref{sect:4.1.2}, we set
$$
\Sigma\overline{\mathcal C}_{\bar{\lambda}}
:=
\Sigma{\mathcal C}_{\lambda}/\overline{\Gamma}_{{\mathcal C}_{\lambda}},
\qquad
{\mathcal C}_{\lambda}^{(0)}:=\nu(\Sigma{\mathcal C}_{\lambda}),
\qquad
\overline{\mathcal C}_{\lambda}^{(0)}:=\overline{\nu}(\Sigma\overline{\mathcal C}_{\bar\lambda}).
$$
By Proposition~\ref{prop:stabilizer:curve}, we have
\begin{equation}
\label{eqn:twisted:sector:curve}
{\mathcal C}_{\lambda}^{(0)}
=
{\mathcal C}_{\lambda}
\cap
{\mathcal X}^{(0)},
\qquad
\overline{\mathcal C}_{\bar{\lambda}}^{(0)}
=
\overline{\mathcal C}_{\bar{\lambda}}
\cap
{\mathcal Y}^{(0)}.
\end{equation}
By \eqref{eqn:orb:char:form:1},  
${\rm Td}^{\Sigma}(\overline{\mathcal C}_{\bar{\lambda}}/B,h_{\overline{\mathcal C}_{\bar{\lambda}}/B})$ 
is a form on 
$\overline{\mathcal C}_{\bar{\lambda}}\amalg\Sigma\overline{\mathcal C}_{\bar{\lambda}}$.
As in Section~\ref{sect:4.4}, let
$$
\varPi\colon
\overline{\mathcal C}_{\bar{\lambda}}\amalg\Sigma\overline{\mathcal C}_{\bar{\lambda}}
\to
\overline{\mathcal C}_{\bar{\lambda}}\amalg\overline{\mathcal C}_{\bar{\lambda}}^{(0)}
$$ 
be the projection. Then
\begin{equation}
\label{eqn:orbifold:Todd:form:fixed:curve}
\begin{aligned}
\varPi_{*}{\rm Td}^{\Sigma}(\overline{\mathcal C}_{\bar{\lambda}}/B,h_{\overline{\mathcal C}_{\bar{\lambda}}/B})
&=
\frac{1}{|\overline{\Gamma}_{{\mathcal C}_{\lambda}}|}
(\pi_{\lambda})_{*}{\rm Td}(T{\mathcal C}_{\lambda}/B,h_{{\mathcal C}_{{\lambda}}/B})
\\
&\quad
+
\frac{1}{|\overline{\Gamma}_{{\mathcal C}_{\lambda}}|}
\sum_{{\frak p}\subset{\mathcal C}_{\lambda}^{(0)}}\sum_{g\in(\overline{\Gamma}_{{\mathcal C}_{\lambda}})_{\frak p}^{*}}
(\pi_{\frak p})_{*}
{\rm Td}_{g}(T{\mathcal C}_{\lambda}/B,h_{{\mathcal C}_{{\lambda}}/B}).
\end{aligned}
\end{equation}
We define the differential form
$[\varPi_{*}{\rm Td}^{\Sigma}(\overline{\mathcal C}_{\bar{\lambda}}/B,h_{\overline{\mathcal C}_{\bar{\lambda}}/B})]^{2\dim\overline{f}_{*}+(1,1)}$
on $\overline{\mathcal C}_{\bar{\lambda}}\amalg\overline{\mathcal C}_{\bar{\lambda}}^{(0)}$ in the same way as in
Definition~\ref{def:notation:relative:degree}.

\begin{proposition}
\label{prop:orb:Td:ch:fixed:curve}
The following equality of forms on $\overline{\mathcal C}_{\bar{\lambda}}\amalg\overline{\mathcal C}_{\bar{\lambda}}^{(0)}$ 
holds:
$$
[
\varPi_{*}
{\rm Td}^{\Sigma}(\overline{\mathcal C}_{\bar{\lambda}}/B,h_{\overline{\mathcal C}_{\bar{\lambda}}/B})
]^{2\dim\overline{f}_{*}+(1,1)}
=
\frac{1}{12}c_{1}(\overline{\mathcal C}_{\bar{\lambda}}/B)^{2}
+
\frac{1}{12}\sum_{\overline{\frak p}\subset\overline{\mathcal C}_{\bar{\lambda}}^{(0)}}
\frac{\nu_{{\frak p},\lambda}^{2}-1}{\nu_{{\frak p},\lambda}}
\,c_{1}(\overline{\mathcal C}_{\bar{\lambda}}/B)|_{\overline{\frak p}}.
$$
\end{proposition}

\begin{pf}
Let ${\frak p}\subset{\mathcal C}_{\lambda}^{(0)}$ and express $G_{\frak p}$ as in \eqref{eqn:expression:orbi:group}. 
As in Section~\ref{sect:4.5.2}, we have 
${\mathcal C}_{\lambda}={\frak C}_{{\frak p},k_{\lambda}}$ for some $k_{\lambda}\in\{1,2,3\}$ as germs at ${\frak p}$.
Since $G_{{\mathcal C}_{\lambda}}=G_{{\frak C}_{{\frak p},k_{\lambda}}}$ in this case, it follows from
\eqref{eqn:stabilizer:p:alpha:3}, \eqref{eqn:stabilizer:coord:axis} that
the $(\overline{\Gamma}_{{\mathcal C}_{\lambda}})_{\frak p}$-action on the germ ${\mathcal C}_{\lambda}$ 
at ${\frak p}$ is given by the homomorphism $\chi_{{\frak p},\lambda}=\chi_{{\frak p},k_{\lambda}}$.
Hence we have
\begin{equation}
\label{eqn:stabilizer:p:alpha:4}
(\overline{\Gamma}_{{\mathcal C}_{\lambda}})_{\frak p}\cong{\rm Im}\,\chi_{{\frak p},\lambda}\subset{\bf C}^{*},
\qquad
G_{{\mathcal C}_{\lambda}}=\ker\chi_{{\frak p},\lambda}.
\end{equation} 
By the definition of equivariant Todd form and \eqref{eqn:Dedekind:sum}, \eqref{eqn:stabilizer:p:alpha:4}, we get
\begin{equation}
\label{eqn:equiv:Todd:fixed:curve}
\begin{aligned}
\sum_{g\in(\overline{\Gamma}_{{\mathcal C}_{\lambda}})_{\frak p}^{*}}
{\rm Td}_{g}(T{\mathcal C}_{\lambda}/B,h_{{\mathcal C}_{{\lambda}}/B})|_{\frak p}
&=
-\sum_{\zeta\in({\rm Im}\chi_{{\frak p},\lambda})^{*}}\frac{\zeta}{(\zeta-1)^{2}}\,
c_{1}({\mathcal C}_{\lambda}/B)|_{\frak p}
\\
&=
\frac{\nu_{{\frak p},\lambda}^{2}-1}{12}\,
c_{1}({\mathcal C}_{\lambda}/B)|_{\frak p},
\end{aligned}
\end{equation}
where we used $|{\rm Im}\chi_{{\frak p},\lambda}|=|{\rm Im}\chi_{{\frak p},k_{\lambda}}|=\nu_{{\frak p},\lambda}$.
By \eqref{eqn:orbifold:Todd:form:fixed:curve} and \eqref{eqn:equiv:Todd:fixed:curve}, we get
\begin{equation}
\label{eqn:curvature:orb:torsion:fixed:curve}
\begin{aligned}
\,&
[
\varPi_{*}{\rm Td}^{\Sigma}(\overline{\mathcal C}_{\bar{\lambda}}/B,h_{\overline{\mathcal C}_{\bar{\lambda}}/B})
]^{2\dim\overline{f}_{*}+(1,1)}
\\
&=
\frac{1}{12|\overline{\Gamma}_{{\mathcal C}_{\lambda}}|}(\pi_{\lambda})_{*}c_{1}({\mathcal C}_{\lambda}/B)^{2}
+
\sum_{{\frak p}\subset{\mathcal C}_{\lambda}^{(0)}}
\frac{\nu_{{\frak p},\lambda}^{2}-1}{12|\overline{\Gamma}_{{\mathcal C}_{\lambda}}|}\,
(\pi_{\frak p})_{*}c_{1}({\mathcal C}_{\lambda}/B)|_{\frak p}
\\
&=
\frac{1}{12}c_{1}(\overline{\mathcal C}_{\bar{\lambda}}/B)^{2}
+
\sum_{\overline{\frak p}\subset\overline{\mathcal C}_{\bar{\lambda}}^{(0)}}
|\overline{\Gamma}_{{\mathcal C}_{\lambda}}/(\overline{\Gamma}_{{\mathcal C}_{\lambda}})_{\frak p}|
\frac{\nu_{{\frak p},\lambda}^{2}-1}{12|\overline{\Gamma}_{{\mathcal C}_{\lambda}}|}\,
c_{1}(\overline{\mathcal C}_{\bar{\lambda}}/B)|_{\overline{\frak p}}.
\end{aligned}
\end{equation}
Since $\nu_{{\frak p},\lambda}=|(\overline{\Gamma}_{{\mathcal C}_{\lambda}})_{\frak p}|$, the result follows from
\eqref{eqn:curvature:orb:torsion:fixed:curve}.
\end{pf}

\begin{proposition}
\label{prop:orb:Td:ch:fixed:curve:2}
The following equality of forms on 
$\overline{\mathcal C}_{\bar{\lambda}}\amalg\overline{\mathcal C}_{\bar{\lambda}}^{(0)}$ 
holds:
$$
\begin{aligned}
\,&
\sum_{\bar{\lambda}\in\overline{\Lambda}}(n_{\lambda}-1)\,
[
\varPi_{*}
{\rm Td}^{\Sigma}(\overline{\mathcal C}_{\bar{\lambda}}/B,h_{\overline{\mathcal C}_{\bar{\lambda}}/B})
]^{2\dim\overline{f}_{*}+(1,1)}
\\
&=
\frac{1}{12}\sum_{\bar{\lambda}\in\overline{\Lambda}}(n_{\lambda}-1)\,c_{1}(\overline{\mathcal C}_{\bar{\lambda}}/B)^{2}
+
\frac{1}{12}\sum_{\overline{\frak p}\subset{\mathcal Y}^{(0)}}\sum_{k=1}^{3}
\frac{(n_{{\frak p},k}-1)(\nu_{{\frak p},k}^{2}-1)}{\nu_{{\frak p},k}}
\,\overline{\rho}_{\overline{\frak p},k}.
\end{aligned}
$$
\end{proposition}

\begin{pf}
By Proposition~\ref{prop:orb:Td:ch:fixed:curve}, the result follows from the following equality
$$
\begin{aligned}
\,&
\sum_{\bar{\lambda}\in\overline{\Lambda}}(n_{\lambda}-1)
\sum_{\overline{\frak p}\subset\overline{\mathcal C}_{\bar{\lambda}}^{(0)}}
\frac{\nu_{{\frak p},\lambda}^{2}-1}{\nu_{{\frak p},\lambda}}\,
c_{1}(\overline{\mathcal C}_{\bar{\lambda}}/B)|_{\overline{\frak p}}
=
\sum_{\bar{\lambda}\in\overline{\Lambda}}
\sum_{\overline{\frak p}\subset\overline{\mathcal C}_{\bar{\lambda}}\cap{\mathcal Y}^{(0)}}
\frac{(n_{\lambda}-1)(\nu_{{\frak p},\lambda}^{2}-1)}{\nu_{{\frak p},\lambda}}\,
c_{1}(\overline{\mathcal C}_{\bar{\lambda}}/B)|_{\overline{\frak p}}
\\
&=
\sum_{\overline{\frak p}\subset{\mathcal Y}^{(0)}}
\sum_{\overline{\mathcal C}_{\bar{\lambda}}\supset\overline{\frak p}}
\frac{(n_{\lambda}-1)(\nu_{{\frak p},\lambda}^{2}-1)}{\nu_{{\frak p},\lambda}}\,
c_{1}(\overline{\mathcal C}_{\bar{\lambda}}/B)|_{\overline{\frak p}}
=
\sum_{\overline{\frak p}\subset{\mathcal Y}^{(0)}}\sum_{k=1}^{3}
\frac{(n_{{\frak p},k}-1)(\nu_{{\frak p},k}^{2}-1)}{\nu_{{\frak p},k}}\,
\overline{\rho}_{\overline{\frak p},k},
\end{aligned}
$$
where the first equality follows from \eqref{eqn:twisted:sector:curve}
and the third equality follows from \eqref{eqn:key:identity}.
This completes the proof.
\end{pf}

\subsection
{A number associated to a finite abelian subgroup of ${\rm SL}({\bf C}^{3})$ }
\label{sect:4.6}
\par
Let $\Gamma\subset{\rm SL}({\bf C}^{3})$ be a finite abelian subgroup and consider its expression
\eqref{eqn:abelian:subgroup:SL(3,C)}.
Recall that $\Gamma^{0}$ and $\delta_{k}(\Gamma)$ were introduced in Definitions~\ref{def:subset:Gamma:0}
and \ref{def:admissibility:subgroup:SL}. We define
\begin{equation}
\label{eqn:invariant:epsilon:0}
\epsilon_{k}(\Gamma)
:=
\frac{\delta_{k}(\Gamma)}{|\Gamma|}
-
\frac{(|\ker\chi_{k}|-1)(|{\rm Im}\,\chi_{k}|^{2}-1)}{12|{\rm Im}\,\chi_{k}|}
\qquad
(k=1,2,3).
\end{equation}

\begin{proposition}
\label{prop:invariant:subgrp:SL:C3}
If $({\bf C}^{3})^{\Gamma^{*}}=\{0\}$, then
\begin{equation}
\label{eqn:invariant:epsilon:1}
\epsilon_{1}(\Gamma)=\epsilon_{2}(\Gamma)=\epsilon_{3}(\Gamma)
=
-\frac{1}{12|\Gamma|}
\left\{
|\Gamma|^{2}+2-\sum_{k=1}^{3}|\ker\chi_{k}|^{2}
\right\}.
\end{equation}
\end{proposition}

\begin{pf}
Assume $\Gamma^{0}\not=\emptyset$.
For simplicity, we put $\nu_{k}=|{\rm Im}\,\chi_{k}|$ and $n_{k}=|\ker\chi_{k}|$.
\par{\em (Step 1) }
Firstly, we prove $\ker\chi_{k}\cap\ker\chi_{l}=\{1\}$ for $k\not=l$. Let $g\in\ker\chi_{k}\cap\ker\chi_{l}$ $(k\not=l)$ and 
let $\{m\}=\{1,2,3\}\setminus\{k,l\}$.
Since $\chi_{m}=\chi_{k}^{-1}\chi_{l}^{-1}$, we get $g\in\ker\chi_{m}$. Hence $g\in\ker\chi_{1}\cap\ker\chi_{2}\cap\ker\chi_{3}$.
Since $\ker\chi_{1}\cap\ker\chi_{2}\cap\ker\chi_{3}=\{1\}$ by \eqref{eqn:abelian:subgroup:SL(3,C)}, 
we get $\ker\chi_{k}\cap\ker\chi_{l}=\{1\}$ for $k\not=l$.
As a result, $(\ker\chi_{k})^{*}\cap(\ker\chi_{l})^{*}=\emptyset$ for $k\not=l$.
Since $\Gamma^{0}=\Gamma\setminus(\ker\chi_{1}\cup\ker\chi_{2}\cup\ker\chi_{3})$, we get 
$$
\Gamma^{0}=\Gamma^{*}\setminus\{(\ker\chi_{1})^{*}\amalg(\ker\chi_{2})^{*}\amalg(\ker\chi_{3})^{*}\},
$$
which yields the following equality for $k=1,2,3$
\begin{equation}
\label{eqn:Dedekind:sum:subgrp:SL:C3:0}
\delta_{k}(\Gamma)
=
\sum_{g\in\Gamma\setminus\ker\chi_{k}}
\frac{\chi_{k}(g)}{(\chi_{k}(g)-1)^{2}}
-
\sum_{l\not=k}
\sum_{g\in(\ker\chi_{l})^{*}}
\frac{\chi_{k}(g)}{(\chi_{k}(g)-1)^{2}}.
\end{equation}
\par
{\em (Step 2) }
Since $\chi_{k}\colon\Gamma\to{\bf C}^{*}$ is a homomorphism, we get by \eqref{eqn:Dedekind:sum}
\begin{equation}
\label{eqn:Dedekind:sum:subgrp:SL:C3:1}
\sum_{g\in\Gamma\setminus\ker\chi_{k}}\frac{\chi_{k}(g)}{(\chi_{k}(g)-1)^{2}}
=
|\ker\chi_{k}|\cdot\sum_{\zeta\in({\rm Im}\,\chi_{k})^{*}}\frac{\zeta}{(\zeta-1)^{2}}
=
-\frac{1}{12}n_{k}(\nu_{k}^{2}-1).
\end{equation}
Since $\ker\chi_{k}\cap\ker\chi_{l}=\{1\}$ for $l\not=k$, the homomorphism $\chi_{k}\colon\ker\chi_{l}\to{\bf C}^{*}$ is injective in this case.
Hence, when $k\not=l$, $\chi_{k}(\ker\chi_{l})\subset{\bf C}^{*}$ is a subgroup of order $n_{l}$, which implies 
$\chi_{k}(\ker\chi_{l})=\{\zeta\in{\bf C}^{*};\,\zeta^{n_{l}}=1\}$. By \eqref{eqn:Dedekind:sum}, we get for $l\not=k$
\begin{equation}
\label{eqn:Dedekind:sum:subgrp:SL:C3:2}
\sum_{g\in(\ker\chi_{l})^{*}}\frac{\chi_{k}(g)}{(\chi_{k}(g)-1)^{2}}
=
\sum_{\zeta\in\chi_{k}(\ker\chi_{l})\setminus\{1\}}\frac{\zeta}{(\zeta-1)^{2}}
=
\sum_{\zeta^{n_{l}}=1,\,\zeta\not=1}\frac{\zeta}{(\zeta-1)^{2}}
=
-\frac{n_{l}^{2}-1}{12}.
\end{equation}
Substituting \eqref{eqn:Dedekind:sum:subgrp:SL:C3:1} and \eqref{eqn:Dedekind:sum:subgrp:SL:C3:2} into \eqref{eqn:Dedekind:sum:subgrp:SL:C3:0},
we get 
\begin{equation}
\label{eqn:Dedekind:sum:subgrp:SL:C3:3}
\delta_{k}(\Gamma)
=
-\frac{1}{12}
\{
n_{k}(\nu_{k}^{2}-1)+2-\sum_{l\not=k}n_{l}^{2}
\}.
\end{equation}
By \eqref{eqn:invariant:epsilon:0}, \eqref{eqn:Dedekind:sum:subgrp:SL:C3:3} and $n_{k}\nu_{k}=|\Gamma|$, we get
$\epsilon_{k}(\Gamma)=-(|\Gamma|^{2}+2-\sum_{k=1}^{3}n_{k}^{2})/12|\Gamma|$.
This proves the result when $\Gamma^{0}\not=\emptyset$.
\par
{\em (Step 3) }
Let us prove the remaining case. We assume $\Gamma^{0}=\emptyset$.
By Lemma~\ref{lemma:invariant:subgrp:SL:C3:Gamma0:empty}, we can verify \eqref{eqn:invariant:epsilon:1}
in the case $\Gamma^{0}=\emptyset$ and $({\bf C}^{3})^{\Gamma^{*}}=\{0\}$, where
$\epsilon_{1}(\Gamma)=\epsilon_{2}(\Gamma)=\epsilon_{3}(\Gamma)=-1/8$.
This completes the proof.
\end{pf}

\begin{definition}
\label{def:invariant:subgrp:SL:C3}
For a finite abelian subgroup $\Gamma\subset{\rm SL}({\bf C}^{3})$ with $({\bf C}^{3})^{\Gamma^{*}}=\{0\}$, define
$$
\epsilon(\Gamma)
:=
-\frac{1}{12|\Gamma|}
\left\{
|\Gamma|^{2}+2-\sum_{k=1}^{3}|\ker\chi_{k}|^{2}
\right\}.
$$
\end{definition}

\subsection
{A formula for an orbifold characteristic form on ${\mathcal Y}\amalg{\mathcal Y}^{(1)}\amalg{\mathcal Y}^{(0)}$ }
\label{sect:4.7}
\par

\begin{theorem}
\label{thm:orbifold:Tod:ch:family}
The following equality of forms on ${\mathcal Y}\amalg{\mathcal Y}^{(1)}\amalg{\mathcal Y}^{(0)}$ holds
$$
\begin{aligned}
\,&
[
\varPi_{*}
{\rm Td}^{\Sigma}(\Theta_{{\mathcal Y}/B},h_{{\mathcal Y}/B})
\sum_{p\geq0}(-1)^{p}p\,{\rm ch}^{\Sigma}(\Omega_{{\mathcal Y}/B}^{p},h_{\Omega^{p}_{{\mathcal Y}/B}})
]^{2\dim\overline{f}_{*}+(1,1)}
\\
&\quad
-
\sum_{\bar{\lambda}\in\overline{\Lambda}}(n_{\lambda}-1)\,
[
\varPi_{*}
{\rm Td}^{\Sigma}(\overline{\mathcal C}_{\bar{\lambda}}/B,h_{\overline{\mathcal C}_{\bar{\lambda}}/B})
]^{2\dim\overline{f}_{*}+(1,1)}
\\
&=
-\frac{1}{12}c_{1}({\mathcal Y}/B)c_{3}({\mathcal Y}/B)
-
\frac{1}{12}\sum_{\bar{\lambda}\in\overline{\Lambda}}\frac{n_{\lambda}^{2}-1}{n_{\lambda}}
\{c_{1}({\mathcal Y}/B)c_{1}(\overline{\mathcal C}_{\bar{\lambda}}/B)\}|_{\overline{\mathcal C}_{\bar{\lambda}}}
\\
&\quad
+
\sum_{\overline{\frak p}\subset{\mathcal Y}^{(0)}}
\epsilon(G_{\frak p})\,c_{1}({\mathcal Y}/B)|_{\overline{\frak p}}.
\end{aligned}
$$
\end{theorem}

\begin{pf}
By Propositions~\ref{prop:orb:Tod:ch} and \ref{prop:orb:Td:ch:fixed:curve:2}, we get
\begin{equation}
\label{eqn:orb:fixed:curve:Td:ch}
\begin{aligned}
\,&
[
\varPi_{*}
{\rm Td}^{\Sigma}(\Theta_{{\mathcal Y}/B},h_{{\mathcal Y}/B})
\sum_{p\geq0}(-1)^{p}p\,{\rm ch}^{\Sigma}(\Omega_{{\mathcal Y}/B}^{p},h_{\Omega^{p}_{{\mathcal Y}/B}})
]^{2\dim\overline{f}_{*}+(1,1)}
\\
&\quad
-
\sum_{\bar{\lambda}\in\overline{\Lambda}}(n_{\lambda}-1)\,
[
\varPi_{*}
{\rm Td}^{\Sigma}(\overline{\mathcal C}_{\bar{\lambda}}/B,h_{\overline{\mathcal C}_{\bar{\lambda}}/B})
]^{2\dim\overline{f}_{*}+(1,1)}
\\
&=
-\frac{1}{12}c_{1}({\mathcal Y}/B)c_{3}({\mathcal Y}/B)
-
\frac{1}{12}\sum_{\bar{\lambda}\in\overline{\Lambda}}\frac{n_{\lambda}^{2}-1}{n_{\lambda}}
\{c_{1}({\mathcal Y}/B)c_{1}(\overline{\mathcal C}_{\bar{\lambda}}/B)\}|_{\overline{\mathcal C}_{\bar{\lambda}}}
\\
&\quad
+
\sum_{\overline{\frak p}\subset{\mathcal Y}^{(0)}}\sum_{k=1}^{3}
\left\{
\frac{\delta_{k}(G_{\frak p})}{|G_{\frak p}|}
-
\frac{(n_{{\frak p},k}-1)(\nu_{{\frak p},k}^{2}-1)}{12\nu_{{\frak p},k}}
\right\}
\overline{\rho}_{\overline{\frak p},k}.
\end{aligned}
\end{equation}
Substituting $\epsilon_{k}(G_{\frak p})=\epsilon(G_{\frak p})$ $(k=1,2,3)$ and 
$\sum_{k=1}^{3}\overline{\rho}_{\overline{\frak p},k}=c_{1}({\mathcal Y}/B)|_{\overline{\frak p}}$ 
into \eqref{eqn:orb:fixed:curve:Td:ch}, we get the result.
\end{pf}

\section
{Equivariant BCOV torsion and its curvature}
\label{sect:5}
\par

\subsection
{Equivariant analytic torsion}
\label{sect:5.1}
\par
Following Bismut \cite{Bismut95}, we recall equivariant analytic torsion.
Let $(X,g_{X})$ be a compact K\"ahler manifold. Let $G$ be a finite group of automorphisms of $(X,g_{X})$.
We write $\widehat{G}$ for the set of irreducible representations of $G$. The character of $W\in\widehat{G}$ is denoted by $\chi_{W}$. 
Since the $G$-action on $X$ is holomorphic, $G$ acts on $A^{p,q}(X)$. 
\par
Let $\square_{p,q}=(\bar{\partial}+\bar{\partial}^{*})^{2}$ be the Laplacian of $(X,g_{X})$ acting on $A^{p,q}(X)$. 
Let $\sigma(\square_{p,q})$ be the spectrum of $\square_{p,q}$ and let $E(\square_{p,q};\lambda)$ be the eigenspace of $\square_{p,q}$
corresponding to the eigenvalue $\lambda\in\sigma(\square_{p,q})$. 
Since $G$ preserves $g_{X}$ and hence $G$ acts on $E(\square_{p,q};\lambda)$, we get the orthogonal splitting
$$
E(\square_{p,q};\lambda)=\bigoplus_{W\in\widehat{G}}{\rm Hom}_{G}(W,E(\square_{p,q};\lambda))\otimes W.
$$
Define the equivariant $\zeta$-function of $\square_{p,q}$ as 
$$
\begin{aligned}
\zeta_{p,q;G}(s)(g)
&:=
\sum_{\lambda\in\sigma(\square_{p,q})\setminus\{0\}}\lambda^{-s}\,
{\rm Tr}\,
\left[
g|_{E(\square_{p,q};\lambda)}
\right]
\\
&=
\sum_{\lambda\in\sigma(\square_{p,q})\setminus\{0\}}\lambda^{-s}\,
\sum_{W\in\widehat{G}}\chi_{W}(g)\dim{\rm Hom}_{G}(W,E(\square_{p,q};\lambda)).
\end{aligned}
$$
It is classical that $\zeta_{p,q;G}(s)(g)$ converges absolutely on the half-plane $\Re s>\dim X$, extends to a meromorphic function on ${\bf C}$,
and is holomorphic at $s=0$.
The {\em equivariant analytic torsion} \cite{Bismut95} of $(X,g_{X},G)$ is the class function on $G$ defined as
$$
\tau_{G}(X,g_{X})(g)
:=
\exp\{-\sum_{q\geq0}(-1)^{q}q\,\zeta'_{0,q;G}(0)(g)\}
$$
and the {\em equivariant BCOV torsion} of $(X,g_{X})$ is the class function on $G$ defined as
$$
T_{{\rm BCOV},G}(X,g_{X})(g)
:=
\exp\{-\sum_{p,q\geq0}(-1)^{p+q}pq\,\zeta'_{p,q;G}(0)(g)\}.
$$
\par
Set 
$$
Y:=X/G
$$ 
and let $p\colon X\to Y$ be the projection. Let $g_{Y}$ be the K\"ahler metric on $Y$ in the sense of orbifolds induced from $g_{X}$.
Let $A^{p,q}(Y)$ be the space of $C^{\infty}$ $(p,q)$-forms on $Y$ in the sense of orbifolds
and let $\square_{p,q}^{\rm orb}=(\bar{\partial}+\bar{\partial}^{*})^{2}$ be the Laplacian of $(Y,g_{Y})$ acting on $A^{p,q}(Y)$.
We can define the spectral zeta function $\zeta_{p,q}^{\rm orb}(s)$ of $\square_{p,q}^{\rm orb}$ in the same way as above or as in Section~\ref{sect:1.1}.
Then $\zeta_{p,q}^{\rm orb}(s)$ extends to a meromorphic function on ${\bf C}$ and is holomorphic at $s=0$,
so that the BCOV torsion $T_{\rm BCOV}(Y,g_{Y})$ of $(Y,g_{Y})$ can be defined by the same formula as in Definition~\ref{def:BCOV:torsion}.
\par
We set
$$
A^{p,q}(X)^{G}:=\{\varphi\in A^{p,q}(X);\,g^{*}\varphi=\varphi\,(\forall\,g\in G)\}
$$
and $\square_{p,q}^{G}:=\square_{p,q}|_{A^{p,q}(X)^{G}}$.
Since $p^{*}\colon A^{p,q}(Y)\to A^{p,q}(X)^{G}$ is an isomorphism with $\square_{p,q}^{G}\circ p^{*}=p^{*}\circ\square_{p,q}^{\rm orb}$, 
the two operators $\square_{p,q}^{\rm orb}$ and $\square_{p,q}^{G}$ are isospectral. 
Hence 
$$
\zeta_{p,q}^{\rm orb}(s)
=
\sum_{\lambda\in\sigma(\square_{p,q})\setminus\{0\}}\lambda^{-s}\,\dim E(\square_{p,q};\lambda)^{G}
=
\frac{1}{|G|}\sum_{g\in G}\zeta_{p,q;G}(s)(g).
$$
By this expression, we have
\begin{equation}
\label{eqn:BCOV:torsion:global:orbifold}
\log T_{\rm BCOV}(Y,g_{Y})=\frac{1}{|G|}\sum_{g\in G}\log T_{{\rm BCOV},G}(X,g_{X})(g).
\end{equation}

\subsection
{A variational formula for equivariant BCOV torsion}
\label{sect:5.2}
\par
In the rest of this section, as an application of the curvature theorem of Bismut-Gillet-Soul\'e \cite{BGS88} for Quillen metrics 
and its equivariant extension by Ma \cite{Ma00}, we shall derive the curvature formula for the equivariant BCOV torsion.
\par
We keep the notation and assumption in Sections~\ref{sect:4.3}--\ref{sect:4.7}. In Section~\ref{sect:5.2}, we assume moreover that
$f\colon{\mathcal X}\to B$ is {\em locally-projective}. 
Hence $X_{b}$ is a smooth projective threefold for all $b\in B$.
Let $h_{{\mathcal X}/B}$ be a $G$-invariant Hermitian metric on $\Theta_{{\mathcal X}/B}$, which is fiberwise K\"ahler.
Let $h_{{\mathcal Y}/B}$ be the Hermitian metric on $\Theta_{{\mathcal Y}/B}$ induced by $h_{{\mathcal X}/B}$.
Define $T_{{\rm BCOV},G}({\mathcal X}/B)(g),T_{\rm BCOV}({\mathcal Y}/B)\in C^{\infty}(B)$ as
$$
T_{{\rm BCOV},G}({\mathcal X}/B)(g)(b):=T_{{\rm BCOV},G}(X_{b},h_{{\mathcal X}/B}|_{X_{b}})(g),
$$
$$
\begin{aligned}
\log T_{\rm BCOV}({\mathcal Y}/B)(b)
&:=
\log T_{{\rm BCOV}}(Y_{b}/G,h_{{\mathcal Y}/B}|_{Y_{b}})
\\
&
=
\frac{1}{|G|}\sum_{g\in G}\log T_{{\rm BCOV},G}(X_{b},h_{{\mathcal X}/B}|_{X_{b}})(g).
\end{aligned}
$$
By the $G$-equivariance of $f\colon{\mathcal X}\to B$, 
the direct image sheaf $R^{q}f_{*}\Omega^{p}_{{\mathcal X}/B}$ is a $G$-equivariant holomorphic vector bundle on $B$
for all $p,q\geq0$.
For $W\in\widehat{G}$, set
$$
(R^{q}f_{*}\Omega^{p}_{{\mathcal X}/B})_{W}:={\rm Hom}_{G}(W,R^{q}f_{*}\Omega^{p}_{{\mathcal X}/B})\otimes W
$$
and let
\begin{equation}
\label{eqn:isotypical:decomposition}
R^{q}f_{*}\Omega^{p}_{{\mathcal X}/B}=\bigoplus_{W\in\widehat{G}}(R^{q}f_{*}\Omega^{p}_{{\mathcal X}/B})_{W}
\end{equation}
be the isotypical decomposition of the vector bundle $R^{q}f_{*}\Omega^{p}_{{\mathcal X}/B}$.
We set 
$$
(R^{q}f_{*}\Omega^{p}_{{\mathcal X}/B})^{G}
:=
\{\varphi\in R^{q}f_{*}\Omega^{p}_{{\mathcal X}/B};\,g\cdot\varphi=\varphi\,(\forall\,g\in G)\}.
$$ 
The $L^{2}$-metrics on $R^{q}f_{*}\Omega^{p}_{{\mathcal X}/B}$, $(R^{q}f_{*}\Omega^{p}_{{\mathcal X}/B})_{W}$, 
$(R^{q}f_{*}\Omega^{p}_{{\mathcal X}/B})^{G}$ with respect to the metric $h_{{\mathcal X}/B}$ are denoted by $h_{L^{2}}$.

\begin{proposition}
\label{prop:2nd:variation:orb:BCOV:torsion}
The following equality of $(1,1)$-forms on $B$ holds
$$
\begin{aligned}
\,&
-dd^{c}\log T_{\rm BCOV}({\mathcal Y}/B)
+
\sum_{p,q\geq0}(-1)^{p+q}p\,c_{1}((R^{q}f_{*}\Omega^{p}_{{\mathcal X}/B})^{G},h_{L^{2}})
\\
&=
[\overline{f}_{*}\{{\rm Td}^{\Sigma}(\Theta_{{\mathcal Y}/B},h_{{\mathcal Y}/B})
\sum_{p\geq0}(-1)^{p}p\,{\rm ch}^{\Sigma}(\Omega^{p}_{{\mathcal Y}/B},h_{\Omega^{p}_{{\mathcal Y}/B}})\}]^{(1,1)}.
\end{aligned}
$$
\end{proposition}

\begin{pf}
By the curvature formulae for Quillen metrics \cite{BGS88} and equivariant Quillen metrics \cite{Ma00}
applied to the $G$-equivariant morphism $f\colon{\mathcal X}\to B$ and
the $G$-equivariant vector bundles $\Omega^{p}_{{\mathcal X}/B}$ $(p\geq0)$, 
we get for all $g\in G$
$$
\begin{aligned}
\,&
-dd^{c}\log T_{{\rm BCOV},G}({\mathcal X}/B)(g)
+
\sum_{p,q\geq0}(-1)^{p+q}p\sum_{W\in\widehat{G}}\frac{\chi_{W}(g)}{{\rm rk}\,W}
c_{1}((R^{q}f_{*}\Omega^{p}_{{\mathcal X}/B})_{W},h_{L^{2}})
\\
&=
[f_{*}\{{\rm Td}_{g}(\Theta_{{\mathcal X}/B},h_{{\mathcal X}/B})
\sum_{p\geq0}(-1)^{p}p\,{\rm ch}_{g}(\Omega^{p}_{{\mathcal X}/B},h_{\Omega^{p}_{{\mathcal X}/B}})\}]^{(1,1)}.
\end{aligned}
$$
This, together with \eqref{eqn:BCOV:torsion:global:orbifold}, \eqref{eqn:orb:char:form:2} and the equality
$$
\frac{1}{|G|}
\sum_{g\in G}\sum_{W\in\widehat{G}}\frac{\chi_{W}(g)}{{\rm rk}\,W}
c_{1}((R^{q}f_{*}\Omega^{p}_{{\mathcal X}/B})_{W},h_{L^{2}})
=
c_{1}((R^{q}f_{*}\Omega^{p}_{{\mathcal X}/B})^{G},h_{L^{2}}),
$$
yields the result.
\end{pf}

\par
For a component ${\mathcal C}_{\lambda}\subset{\mathcal X}^{(1)}$, we set 
$$
C_{b,\lambda}:={\mathcal C}_{\lambda}\cap X_{b},
\qquad
\overline{C}_{b,\bar{\lambda}}:=C_{b,\lambda}/\overline{\Gamma}_{{\mathcal C}_{\lambda}}.
$$
Since $X_{b}$ is a smooth projective threefold, $C_{b,\lambda}$ is a smooth projective curve and 
$\overline{C}_{b,\bar{\lambda}}$ is a compact orbifold curve.
We have the canonical isomorphism
$$
H^{0}(\overline{C}_{b,\bar{\lambda}},K_{\overline{C}_{b,\bar{\lambda}}})
\cong
H^{0}(C_{b,\lambda},K_{C_{b,\lambda}})^{\overline{\Gamma}_{{\mathcal C}_{\lambda}}}
$$
via the projection $C_{b,\lambda}\to\overline{C}_{b,\bar{\lambda}}$. 
Here $\overline{C}_{b,\bar{\lambda}}$ is viewed as the smooth curve underlying the orbifold curve $\overline{C}_{b,\bar{\lambda}}$.
For $\bar{\lambda}\in\overline{\Lambda}$, we set 
$$
\overline{h}_{\bar{\lambda}}
:=
\dim H^{0}(\overline{C}_{b,\bar{\lambda}},K_{\overline{C}_{b,\bar{\lambda}}})
=
\dim H^{0}(C_{b,\lambda},K_{C_{b,\lambda}})^{\overline{\Gamma}_{{\mathcal C}_{\lambda}}}.
$$
The Torelli map for the family $\overline{f}\colon\overline{\mathcal C}_{\bar{\lambda}}\to B$ is defined by
$$
J_{\overline{\mathcal C}_{\bar{\lambda}}/B}\colon B\ni b\to
[{\rm Jac}(\overline{C}_{b,\bar{\lambda}})]=[{\rm Jac}(C_{b,\lambda})^{\overline{\Gamma}_{{\mathcal C}_{\lambda}}}]
\in
{\mathcal A}_{\overline{h}_{\bar{\lambda}}},
$$
where $[{\rm Jac}(\overline{C}_{b,\bar{\lambda}})]$ denotes the isomorphism class of the Jacobian variety 
${\rm Jac}(\overline{C}_{b,\bar{\lambda}})$ of the smooth curve underlying $\overline{C}_{b,\bar{\lambda}}$.
\par
Let $h_{L^{2}}=h_{\overline{f}_{*}K_{\overline{\mathcal C}_{\bar{\lambda}}/B}}$ be the $L^{2}$-metric on 
$\overline{f}_{*}K_{\overline{\mathcal C}_{\bar{\lambda}}/B}$, which is independent of the choice of $h_{\overline{\mathcal C}_{\bar\lambda}/B}$,
the Hermitian metric on $T\overline{\mathcal C}_{\bar\lambda}/B$ induced by $h_{{\mathcal X}/B}$. 
Let $\varphi_{1}^{(\bar{\lambda})},\ldots,\varphi_{\overline{h}_{\bar{\lambda}}}^{(\bar{\lambda})}
\in H^{0}(B,\overline{f}_{*}K_{\overline{\mathcal C}_{\bar{\lambda}}/B})$ 
be a basis of $\overline{f}_{*}K_{\overline{\mathcal C}_{\bar{\lambda}}/B}$ as a free ${\mathcal O}_{B}$-module. 
For $b\in B$, set $\varphi_{k}^{(\bar{\lambda})}(b):=\varphi_{k}^{(\bar{\lambda})}|_{\overline{C}_{b,\bar{\lambda}}}$.
Then $\{\varphi_{1}^{(\bar{\lambda})}(b),\ldots,\varphi_{\overline{h}_{\lambda}}^{(\bar{\lambda})}(b)\}$ is a basis of 
$H^{0}(\overline{C}_{b,\bar{\lambda}},K_{\overline{C}_{b,\bar{\lambda}}})$ such that
$$
\left\|
\varphi_{1}^{(\bar{\lambda})}(b)\wedge\cdots\wedge\varphi_{\overline{h}_{\lambda}}^{(\bar{\lambda})}(b)
\right\|_{L^{2}}^{2}
=
\det
\left(
\frac{i}{2\pi}\int_{\overline{C}_{b,\bar{\lambda}}}\varphi_{k}^{(\bar{\lambda})}(b)\wedge\overline{\varphi_{l}^{(\bar{\lambda})}(b)}
\right)_{1\leq k,l\leq\overline{h}_{\bar{\lambda}}}.
$$
By this expression and the definition of $\omega_{{\mathcal A}_{\overline{h}_{\bar{\lambda}}}}$, we get
\begin{equation}
\label{eqn:curvature:direct:image:fixed:points}
c_{1}(\overline{f}_{*}K_{\overline{\mathcal C}_{\bar{\lambda}}/B},h_{L^{2}})
=
c_{1}((f_{*}K_{{\mathcal C}_{\lambda}/B})^{\overline{\Gamma}_{{\mathcal C}_{\lambda}}},h_{L^{2}})
=
J_{\overline{\mathcal C}_{\bar{\lambda}}/B}^{*}\omega_{{\mathcal A}_{\overline{h}_{\bar{\lambda}}}}.
\end{equation}

\par 
Let $\tau_{\overline{\Gamma}_{{\mathcal C}_{\lambda}}}({\mathcal C}_{\lambda}/B)(g)$ 
and ${\rm Vol}({\mathcal C}_{\lambda}/B)$ be the $C^{\infty}$ function on $B$ defined as
$$
\tau_{\overline{\Gamma}_{{\mathcal C}_{\lambda}}}({\mathcal C}_{\lambda}/B)(g)(b)
:=
\tau_{\overline{\Gamma}_{{\mathcal C}_{\lambda}}}(C_{b,\lambda},h_{{\mathcal X}/B}|_{C_{b,\lambda}})(g)
\qquad
(g\in\overline{\Gamma}_{{\mathcal C}_{\lambda}}),
$$
$$
{\rm Vol}({\mathcal C}_{\lambda}/B)(b)
:=
{\rm Vol}(C_{b,\lambda},h_{{\mathcal X}/B}|_{C_{b,\lambda}})
=
\frac{1}{2\pi}\int_{C_{b,\lambda}}\gamma_{{\mathcal X}/B}|_{C_{b,\lambda}},
$$
where $\gamma_{{\mathcal X}/B}$ is the K\"ahler form associated with $h_{{\mathcal X}/B}$.
Let $\tau(\overline{\mathcal C}_{\bar{\lambda}}/B)$ and ${\rm Vol}(\overline{\mathcal C}_{\bar{\lambda}}/B)$ 
be the $C^{\infty}$ function on $B$ defined as
$$
\tau(\overline{\mathcal C}_{\bar{\lambda}}/B)(b)
:=
\tau(\overline{C}_{b,\bar{\lambda}},h_{\overline{\mathcal C}_{\bar{\lambda}}/B}|_{\overline{C}_{b,\bar{\lambda}}}),
$$
$$
{\rm Vol}(\overline{\mathcal C}_{\bar{\lambda}}/B)(b)
:=
{\rm Vol}(\overline{C}_{b,\bar{\lambda}},h_{\overline{\mathcal C}_{\bar{\lambda}}/B}|_{\overline{C}_{b,\bar{\lambda}}})
=
\frac{1}{|\overline{\Gamma}_{{\mathcal C}_{\lambda}}|}{\rm Vol}({\mathcal C}_{\lambda}/B)(b).
$$

\begin{proposition}
\label{prop:variation:torsion:orbifold:curve:1}
The following equality of $(1,1)$-forms on $B$ holds:
$$
-dd^{c}\log
\{
\tau(\overline{\mathcal C}_{\bar{\lambda}}/B){\rm Vol}(\overline{\mathcal C}_{\bar{\lambda}}/B)
\}
+
J_{\overline{\mathcal C}_{\bar{\lambda}}/B}^{*}\omega_{{\mathcal A}_{\overline{h}_{\bar{\lambda}}}}
=
[
\overline{f}_{*}{\rm Td}^{\Sigma}(\overline{\mathcal C}_{\bar{\lambda}}/B,h_{\overline{\mathcal C}_{\bar{\lambda}}/B})
]^{(1,1)}.
$$
\end{proposition}

\begin{pf}
Write $R_{\lambda}$ for the set of irreducible representations of $\overline{\Gamma}_{{\mathcal C}_{\lambda}}$.
Applying the curvature formula for the equivariant Quillen metrics \cite{BGS88},\cite{Ma00} to 
the $\overline{\Gamma}_{{\mathcal C}_{\lambda}}$-equivariant holomorphic submersion $f\colon{\mathcal C}_{\lambda}\to B$, 
we get for $g\in\overline{\Gamma}_{{\mathcal C}_{\lambda}}$
\begin{equation}
\label{eqn:curvature:equivariant:torsion:fixed:curve}
\begin{aligned}
\,&
-dd^{c}\log\tau_{\overline{\Gamma}_{{\mathcal C}_{\lambda}}}({\mathcal C}_{\lambda}/B)(g)
+
\sum_{W\in R_{\lambda}}\frac{\chi_{W}(g)}{{\rm rk}\,W}\,c_{1}((f_{*}{\mathcal O}_{{\mathcal C}_{\lambda}})_{W},h_{L^{2}})
\\
&\quad
+
\sum_{W\in R_{\lambda}}\frac{\chi_{W}(g)}{{\rm rk}\,W}\,c_{1}((f_{*}K_{{\mathcal C}_{\lambda}/B})_{W},h_{L^{2}})
=
f_{*}{\rm Td}_{g}(T{\mathcal C}_{\lambda}/B).
\end{aligned}
\end{equation}
Since 
$|\overline{\Gamma}_{{\mathcal C}_{\lambda}}|\,\log\tau(\overline{\mathcal C}_{\bar\lambda}/B)
=
\sum_{g\in\overline{\Gamma}_{{\mathcal C}_{\lambda}}}\log\tau_{\overline{\Gamma}_{{\mathcal C}_{\lambda}}}({\mathcal C}_{\lambda}/B)(g)$
by \eqref{eqn:BCOV:torsion:global:orbifold}
and
$$
\frac{1}{|\overline{\Gamma}_{{\mathcal C}_{\lambda}}|}
\sum_{W\in R_{\lambda}}
\sum_{g\in\overline{\Gamma}_{{\mathcal C}_{\lambda}}}
\frac{\chi_{W}(g)}{{\rm rk}\,W}\,c_{1}((f_{*}{\mathcal O}_{{\mathcal C}_{\lambda}})_{W},h_{L^{2}})
=
-dd^{c}\log{\rm Vol}(\overline{\mathcal C}_{\bar{\lambda}}/B),
$$
$$
\frac{1}{|\overline{\Gamma}_{{\mathcal C}_{\lambda}}|}
\sum_{W\in R_{\lambda}}
\sum_{g\in\overline{\Gamma}_{{\mathcal C}_{\lambda}}}
\frac{\chi_{W}(g)}{{\rm rk}\,W}\,c_{1}((f_{*}K_{{\mathcal C}_{\lambda}/B})_{W},h_{L^{2}})
=
J_{\overline{\mathcal C}_{\bar{\lambda}}/B}^{*}\omega_{{\mathcal A}_{\overline{h}_{\bar{\lambda}}}},
$$
we get the result by \eqref{eqn:curvature:equivariant:torsion:fixed:curve}, \eqref{eqn:orb:char:form:2},
\eqref{eqn:orbifold:Todd:form:fixed:curve}.
\end{pf}

\begin{theorem}
\label{thm:curvature:orbifold:BCOV:torsion:fixed:curve}
The following equality of $(1,1)$-forms on $B$ holds
$$
\begin{aligned}
\,&
-dd^{c}\log T_{\rm BCOV}({\mathcal Y}/B)
+
\sum_{\bar{\lambda}\in\overline{\Lambda}}
(n_{\lambda}-1)\,dd^{c}\log\left\{\tau(\overline{\mathcal C}_{\bar{\lambda}}/B){\rm Vol}(\overline{\mathcal C}_{\bar{\lambda}}/B)\right\}
\\
&=
-\frac{1}{12}\overline{f}_{*}\{c_{1}({\mathcal Y}/B)c_{3}({\mathcal Y}/B)\}
-
\frac{1}{12}\sum_{\bar{\lambda}\in\overline{\Lambda}}\frac{n_{\lambda}^{2}-1}{n_{\lambda}}
(\overline{f}|_{\overline{\mathcal C}_{\bar{\lambda}}})_{*}\{c_{1}({\mathcal Y}/B)c_{1}(\overline{\mathcal C}_{\bar{\lambda}}/B)\}
\\
&\quad
+
\sum_{\overline{\frak p}\subset{\mathcal Y}^{(0)}}
\epsilon(G_{\frak p})\,(\overline{f}|_{\overline{\frak p}})_{*}c_{1}({\mathcal Y}/B)
-
\sum_{p,q\geq0}(-1)^{p+q}p\,c_{1}((R^{q}f_{*}\Omega^{p}_{{\mathcal X}/B})^{G},h_{L^{2}})
\\
&\quad
+
\sum_{\bar{\lambda}\in\overline{\Lambda}}(n_{\lambda}-1)\,J_{\overline{\mathcal C}_{\bar{\lambda}}/B}^{*}\omega_{{\mathcal A}_{\overline{h}_{\bar{\lambda}}}}.
\end{aligned}
$$
\end{theorem}

\begin{pf}
The result follows from Theorem~\ref{thm:orbifold:Tod:ch:family} and
Propositions~\ref{prop:2nd:variation:orb:BCOV:torsion} and \ref{prop:variation:torsion:orbifold:curve:1}.
\end{pf}

\section
{BCOV invariants for global abelian Calabi-Yau orbifolds}
\label{sect:6}
\par
In this section, we introduce BCOV invariants for global Calabi-Yau orbifolds and calculate their curvatures.

\subsection
{Global abelian Calabi-Yau orbifolds}
\label{sect:6.1}
\par
Let $X$ be a smooth projective threefold with trivial canonical line bundle. 
Let $G\subset{\rm Aut}(X)$ be a finite abelian group. We assume that 
\begin{equation}
\label{eqn:condition:cohomology:group:action}
H^{0}(X,\Omega^{p}_{X})^{G}=H^{0}(X,\Omega^{p}_{X})
\quad
(p=0,3),
\qquad
H^{0}(X,\Omega^{p}_{X})^{G}=0
\quad
(p=1,2)
\end{equation}
and we set
$$
Y:=X/G.
$$
Since $G$ preserves canonical forms on $X$ by \eqref{eqn:condition:cohomology:group:action}, 
$X^{g}$ is the disjoint union of finitely many curves and finite points for any $g\in G^{*}$.
We define ${\Sigma}^{(1)}X$, ${\Sigma}^{(0)}X$, $\widetilde{\Sigma}^{(1)}X$, ${\Sigma}^{(1)}Y$, ${\Sigma}^{(0)}Y$, $\widetilde{\Sigma}^{(1)}Y$
and $X^{(1)}$, $X^{(0)}$, $Y^{(1)}$, $Y^{(0)}$ in the same way as before in Sections~\ref{sect:4.1} and \ref{sect:4.2}. 
Let $\{C_{\lambda}\}_{\lambda\in\Lambda}$ be the components of $X^{(1)}$.
Then $X^{(1)}=\amalg_{\lambda\in\Lambda}C_{\lambda}$ and $C_{\lambda}\subset X^{g}$ for some $g\in G^{*}$.
Similarly, $X^{(0)}$ is the set of those points ${\frak p}\in X$ such that ${\frak p}$ is isolated in $X^{G_{\frak p}^{*}}$.
As in Section~\ref{sect:4}, let $\overline{\Lambda}=\Lambda/G$ be the $G$-orbits of the $G$-action on $\Lambda$.
Then $\overline{\Lambda}$ is identified with the set of one-dimensional irreducible components of ${\rm Sing}\,Y$.
\par
For $C_{\lambda}$, $\lambda\in\Lambda$,  we can associate the finite subgroups $G_{C_{\lambda}},\Gamma_{C_{\lambda}}\subset G$
and $\overline{\Gamma}_{C_{\lambda}}\subset{\rm Aut}(C_{\lambda})$ in the same way as in Section~\ref{sect:4.3.2}.
We set
$$
\overline{C}_{\lambda}:=C_{\lambda}/\overline{\Gamma}_{C_{\lambda}}.
$$
Then $\overline{C}_{\lambda}\cong\overline{C}_{\lambda'}$ via the $G$-action if and only if $\bar{\lambda}=\bar{\lambda'}$ in $\overline{\Lambda}$.
We write $\overline{C}_{\bar{\lambda}}$ for $\overline{C}_{\lambda}$.
We have $\Sigma^{(1)}Y=\amalg_{\bar{\lambda}\in\overline{\Lambda}}\overline{C}_{\bar{\lambda}}$.
Similarly, ${\mathcal O}_{Y,{\frak p}}\cong{\mathcal O}_{Y,{\frak p}'}$ via the $G$-action 
if and only if ${\frak p}$, ${\frak p}'$ lie on the same $G$-orbit. 
The point of $Y$ corresponding to ${\frak p}$ is denoted by $\overline{\frak p}$.
\par
Since $G_{C_{\lambda}}$ is a finite abelian subgroup of ${\rm SL}({\bf C}^{2})$, $G_{C_{\lambda}}$ is a cyclic group of order
$$
n_{\lambda}:=|G_{C_{\lambda}}|.
$$
By Proposition~\ref{prop:invariant:subgrp:SL:C3}, $G_{\frak p}$, ${\frak p}\in X^{(0)}$ is an abelian subgroup of 
${\rm SL}(T_{\frak p}X)$ such that
$\epsilon_{1}(G_{\frak p})=\epsilon_{2}(G_{\frak p})=\epsilon_{3}(G_{\frak p})=\epsilon(G_{\frak p})$.
For any $g\in G^{*}$, we get the decomposition
$$
X^{g}
=
(\amalg_{\{\lambda\in\Lambda;\,g\in G_{C_{\lambda}}^{*}\}}C_{\lambda})
\amalg
(\amalg_{\{{\frak p}\in X^{(0)};\,g\in G_{\frak p}^{0}\}}{\frak p}).
$$
\par
Let $\gamma$ be a $G$-invariant K\"ahler form on $X$ and let $\overline{\gamma}$ be the K\"ahler form on $Y$ in the sense of orbifolds
induced from $\gamma$.
We define 
\begin{equation}
\label{eqn:Euler:X:G}
\chi^{\rm orb}(Y)
:=
\int_{Y}c_{3}(Y,\overline{\gamma})
+
\sum_{\bar{\lambda}\in\overline{\Lambda}}\frac{n_{\lambda}^{2}-1}{n_{\lambda}}
\int_{\overline{C}_{\bar{\lambda}}}c_{1}(\overline{C}_{\lambda},\overline{\gamma}|_{\overline{C}_{\bar{\lambda}}})
-
12\sum_{\overline{\frak p}\in Y^{(0)}}\epsilon(G_{\frak p}),
\end{equation}
where $c_{3}(Y,\overline{\gamma})$ and $c_{1}(\overline{C}_{\bar{\lambda}},\overline{\gamma}|_{\overline{C}_{\bar{\lambda}}})$ 
denote the corresponding Chern forms. 
Hence their pullbacks to $X$ and $C_{\lambda}$ are the Chern forms $c_{3}(X,\gamma)$ and $c_{1}(C_{\lambda},\gamma|_{C_{\lambda}})$, respectively.
Let $\chi(\cdot)$ denote the ordinary topological Euler characteristic. By \eqref{eqn:Euler:X:G}, Definition~\ref{def:invariant:subgrp:SL:C3} 
and the Gauss-Bonnet-Chern formula, we easily get
\begin{equation}
\label{eqn:orb:Euler:1}
|G|\,\chi^{\rm orb}(Y)
=
\chi(X)
+
\sum_{\lambda\in\Lambda}(n_{\lambda}^{2}-1)\,\chi(C_{\lambda})
+
\sum_{{\frak p}\in X^{(0)}}\{(|G_{\frak p}|^{2}-1)-\sum_{k=1}^{3}(n_{{\frak p},k}^{2}-1)\},
\end{equation}
where we used the notation in Section~\ref{sect:4.5.2}, in particular \eqref{eqn:order:ker:Im:chi}.

\begin{proposition}
\label{prop:comparison:Euler:char}
Let $\widetilde{Y}$ be any crepant resolution of $Y=X/G$. Then
$$
\chi^{\rm orb}(Y)=\chi(\widetilde{Y}).
$$
\end{proposition}

\begin{pf}
By Roan \cite{Roan96}, we have 
\begin{equation}
\label{eqn:orb:Euler:2}
|G|\,\chi(\widetilde{Y})
=
\sum_{g,h\in G}\chi(X^{g}\cap X^{h})
=
\chi(X)+\sum_{(g,h)\in(G\times G)^{*}}\chi(X^{g}\cap X^{h}).
\end{equation}
Let $(g,h)\in(G\times G)^{*}$. Let $C\subset X$ be an irreducible curve.
Since $C\subset X^{g}$ if and only if $g\in G_{C}$, we get
\begin{equation}
\label{eqn:inclusion:fixed:curve}
C\subset X^{g}\cap X^{h}
\qquad\Longleftrightarrow\qquad
(g,h)\in(G_{C}\times G_{C})^{*}.
\end{equation}
Similarly, for any ${\frak p}\in X$, we have
\begin{equation}
\label{eqn:inclusion:fixed:point}
{\frak p}\in X^{g}\cap X^{h}
\qquad\Longleftrightarrow\qquad
(g,h)\in(G_{\frak p}\times G_{\frak p})^{*}.
\end{equation}
For ${\frak p}\in X$, let $\{C_{\lambda_{{\frak p},k}}\}_{k\in K_{\frak p}}\subset X^{(1)}$, $K_{\frak p}\subset\{1,2,3\}$ 
be the set of fixed curves passing through ${\frak p}$, where $K_{\frak p}$ can be empty.
By \eqref{eqn:inclusion:fixed:curve}, \eqref{eqn:inclusion:fixed:point},
we have 
\begin{equation}
\label{eqn:condition:group}
\begin{aligned}
\,&
{\frak p}\hbox{ is an isolated point of }X^{g}\cap X^{h}
\\
&\Longleftrightarrow
\qquad
(g,h)\in(G_{\frak p}\times G_{\frak p})^{*}
\setminus
\bigcup_{k\in K_{\frak p}}(G_{C_{\lambda_{{\frak p},k}}}\times G_{C_{\lambda_{{\frak p},k}}})^{*},
\end{aligned}
\end{equation}
in which case we have ${\frak p}\in X^{(0)}$. 
By \eqref{eqn:inclusion:fixed:curve}, \eqref{eqn:inclusion:fixed:point}, \eqref{eqn:condition:group}, we get
\begin{equation}
\label{eqn:orb:Euler:3}
\begin{aligned}
\,&
\sum_{(g,h)\in(G\times G)^{*}}\chi(X^{g}\cap X^{h})
=
\sum_{\lambda\in\Lambda}
\#
\left\{
(g,h)\in(G_{C_{\lambda}}\times G_{C_{\lambda}})^{*}
\right\}
\cdot
\chi(C_{\lambda})
\\
&
+
\sum_{{\frak p}\in X^{(0)}}
\#\left\{
(g,h)\in
(G_{\frak p}\times G_{\frak p})^{*}
\setminus
\bigcup_{k\in K_{\frak p}}(G_{C_{\lambda_{{\frak p},k}}}\times G_{C_{\lambda_{{\frak p},k}}})^{*}
\right\}
\cdot
\chi({\frak p})
\\
&\quad=
\sum_{\lambda\in\Lambda}(n_{\lambda}^{2}-1)\,\chi(C_{\lambda})
+
\sum_{{\frak p}\in X^{(0)}}\{(|G_{\frak p}|^{2}-1)-\sum_{k\in K_{\frak p}}(n_{{\frak p},k}^{2}-1)\}
\\
&\quad=
\sum_{\lambda\in\Lambda}(n_{\lambda}^{2}-1)\,\chi(C_{\lambda})
+
\sum_{{\frak p}\in X^{(0)}}\{(|G_{\frak p}|^{2}-1)-\sum_{k=1}^{3}(n_{{\frak p},k}^{2}-1)\},
\end{aligned}
\end{equation}
where the last equality follows from \eqref{eqn:criterion:coordinate:axis}.
Substituting \eqref{eqn:orb:Euler:3} into \eqref{eqn:orb:Euler:2} and comparing it with \eqref{eqn:orb:Euler:1}, we get the result.
\end{pf}

\subsection
{BCOV invariants}
\label{sect:6.2}
\par
Let $\eta\in H^{0}(Y,\Omega^{3}_{Y})$ be a nowhere vanishing canonical form on $Y$ in the sense of orbifolds.
We introduce 
\begin{equation}
\label{eqn:anomaly:1}
A(Y,\overline{\gamma})
:=
\exp
\left[
-\frac{1}{12}\int_{Y}
\log
\left(
\frac{i\,\eta\wedge\overline{\eta}}{\overline{\gamma}^{3}/3!}
\cdot
\frac{{\rm Vol}(Y,\overline{\gamma})}{\|\eta\|_{L^{2}}^{2}}
\right)
\,
c_{3}(Y,\overline{\gamma})
\right],
\end{equation}
\begin{equation}
\label{eqn:anomaly:2}
A(\overline{C}_{\bar{\lambda}},\overline{\gamma}|_{\overline{C}_{\bar{\lambda}}})
:=
\exp
\left[
-\frac{1}{12}\int_{\overline{C}_{\bar{\lambda}}}
\log
\left(
\frac{i\,\eta\wedge\overline{\eta}}{\overline{\gamma}^{3}/3!}
\cdot
\frac{{\rm Vol}(Y,\overline{\gamma})}{\|\eta\|_{L^{2}}^{2}}
\right)
\,
c_{1}(\overline{C}_{\bar{\lambda}},\overline{\gamma}|_{\overline{C}_{\bar{\lambda}}})
\right],
\end{equation}
\begin{equation}
\label{eqn:anomaly:3}
A(\overline{\frak p},\overline{\gamma}|_{\overline{\frak p}})
:=
\left.
\frac{i\,\eta\wedge\overline{\eta}}{\overline{\gamma}^{3}/3!}
\right|_{\overline{\frak p}}
\cdot
\frac{{\rm Vol}(Y,\overline{\gamma})}{\|\eta\|_{L^{2}}^{2}}
\qquad
(\overline{\frak p}\in Y^{(0)}).
\end{equation}
\par
For generic $y\in\overline{C}_{\bar{\lambda}}$, the germ $(Y,y)$ is isomorphic to $({\bf C},0)\times A_{\mu_{\bar{\lambda}}}$, where
$$
\mu_{\bar{\lambda}}:=n_{\lambda}-1=|G_{C_{\lambda}}^{*}|
$$
and $A_{\mu}:=({\bf C}^{2},0)/\langle{\rm diag}(e^{2\pi i/(\mu+1)},e^{-2\pi i/(\mu+1)})\rangle$ is the $2$-dimensional $A_{\mu}$-singularity.

\begin{definition}
\label{def:orbifold:BCOV:invariant:1}
For a K\"ahler form $\overline{\gamma}$ on $Y$ in the sense of orbifolds, define the {\em orbifold BCOV invariant} of $(Y,\overline{\gamma})$ by
$$
\begin{aligned}
\tau_{\rm BCOV}^{\rm orb}(Y,\overline{\gamma})
&:=
T_{\rm BCOV}(Y,\overline{\gamma})\,
{\rm Vol}(Y,\overline{\gamma})^{-3+\frac{\chi^{\rm orb}(Y)}{12}}
{\rm Vol}_{L^{2}}(H^{2}(Y,{\bf Z}),[\overline{\gamma}])^{-1}
A(Y,\overline{\gamma})
\\
&\quad\times
\prod_{\bar{\lambda}\in\overline{\Lambda}}
\left\{
\tau(\overline{C}_{\bar{\lambda}},\overline{\gamma}|_{\overline{C}_{\bar{\lambda}}})
{\rm Vol}(\overline{C}_{\lambda},\overline{\gamma}|_{\overline{C}_{\lambda}})
\right\}^{-\mu_{\lambda}}
A(\overline{C}_{\bar{\lambda}},\overline{\gamma}|_{\overline{C}_{\bar{\lambda}}})^{\frac{n_{\lambda}^{2}-1}{n_{\lambda}}}
\\
&\quad\times
\prod_{\overline{\frak p}\in Y^{(0)}}A(\overline{\frak p},\overline{\gamma}|_{\overline{\frak p}})^{\epsilon(G_{\frak p})}.
\end{aligned}
$$
\end{definition}

When $G=\{1\}$, $\tau_{\rm BCOV}^{\rm orb}(Y,\overline{\gamma})$ coincides with the BCOV invariant of $X$.
In this sense, the orbifold BCOV invariant is an extension of the ordinary BCOV invariant to global abelian Calabi-Yau orbifolds.
As in Remark~\ref{remark:BCOV:invariant:Ricci:flat:case}, if $\overline{\gamma}$ is Ricci-flat, then 
\begin{equation}
\label{eqn:BCOV:orb:Ricci:flat}
\begin{aligned}
\tau_{\rm BCOV}^{\rm orb}(Y,\overline{\gamma})
&=
T_{\rm BCOV}(Y,\overline{\gamma})\,
{\rm Vol}(Y,\overline{\gamma})^{-3+\frac{{\chi}^{\rm orb}(Y)}{12}}
{\rm Vol}_{L^{2}}(H^{2}(Y,{\bf Z}),[\overline{\gamma}])^{-1}
\\
&\quad\times
\prod_{\bar{\lambda}\in\overline{\Lambda}}
\left\{
\tau(\overline{C}_{\bar{\lambda}},\overline{\gamma}|_{\overline{C}_{\bar{\lambda}}})
{\rm Vol}(\overline{C}_{\bar{\lambda}},\overline{\gamma}|_{\overline{C}_{\bar{\lambda}}})
\right\}^{-\mu_{\bar{\lambda}}}.
\end{aligned}
\end{equation}
\par
In the rest of this section, we shall derive a variational formula for $\tau_{\rm BCOV}^{\rm orb}(Y,\overline{\gamma})$ 
and prove its independence of the choice of K\"ahler form. 
Namely, $\tau_{\rm BCOV}^{\rm orb}(Y,\overline{\gamma})$ is an invariant of $Y$.
In fact, the definition of $\tau_{\rm BCOV}^{\rm orb}(Y,\overline{\gamma})$ can be extended to arbitrary abelian Calabi-Yau orbifolds.
In Section~\ref{sect:7}, we shall prove that $\tau_{\rm BCOV}^{\rm orb}(Y,\overline{\gamma})$ is independent of the choice of $\overline{\gamma}$
even for general abelian Calabi-Yau orbifolds. Hence $\tau_{\rm BCOV}^{\rm orb}(Y,\overline{\gamma})$ is an invariant of abelian Calabi-Yau orbifolds.

\subsection
{Set up}
\label{sect:6.3}
\par
Let $f\colon{\mathcal X}\to B$ be a locally-projective smooth morphism from a complex manifold ${\mathcal X}$ to a complex manifold $B$
such that $X_{b}=f^{-1}(b)$ is a connected threefold with trivial canonical bundle 
$$
K_{X_{b}}\cong{\mathcal O}_{X_{b}},
\qquad
\forall\,b\in B.
$$
\par
Let $G$ be a finite abelian subgroup of ${\rm Aut}({\mathcal X})$ such that 
the projection $f\colon{\mathcal X}\to B$ is $G$-equivariant with respect to the trivial $G$-action on $B$. 
Hence $G$ preserves the fibers of $f$. Moreover, we assume that
\begin{equation}
\label{eqn:condition:involution}
g^{*}|_{H^{0}(X_{b},K_{X_{b}})}=1,
\qquad
H^{0}(X_{b},\Omega_{X_{b}}^{1})^{G}=H^{0}(X_{b},\Omega_{X_{b}}^{2})^{G}=\{0\}
\end{equation}
for all $g\in G$ and $b\in B$.
Then \eqref{eqn:assumption:group} is satisfied by the first condition of \eqref{eqn:condition:involution}.
\par
By the $G$-equivariance of the family $f\colon{\mathcal X}\to B$, 
the locally free sheaf $R^{q}f_{*}\Omega^{p}_{{\mathcal X}/B}$ is equipped with the $G$-action.
By the second condition of \eqref{eqn:condition:involution} and the Hodge symmetry, we have the following vanishing for 
$p\not=q$ and $p+q\not=3$
\begin{equation}
\label{eqn:vanishing:direct:image}
(R^{q}f_{*}\Omega^{p}_{{\mathcal X}/B})^{G}=0.
\end{equation}
In what follows, we keep the notation in Sections~\ref{sect:4.3}, \ref{sect:4.5}, \ref{sect:5.2}. Hence we set
$$
{\mathcal Y}={\mathcal X}/G
$$
and the projection from ${\mathcal Y}$ to $B$ induced from $f$ is denoted by $\overline{f}\colon{\mathcal Y}\to B$.
The relative holomorphic tangent bundle $\Theta_{{\mathcal X}/B}$ is equipped with a $G$-invariant Hermitian metric $h_{{\mathcal X}/B}$, 
which is fiberwise K\"ahler. Its fiberwise K\"ahler form is denoted by $\gamma_{{\mathcal X}/B}$. 
The corresponding Hermitian metric on $\Theta_{{\mathcal Y}/B}$ and its fiberwise K\"ahler form are denoted by
$\overline{h}_{{\mathcal Y}/B}$ and $\overline{\gamma}_{{\mathcal Y}/B}$, respectively.
Let $\Omega_{{\mathcal Y}/B}^{p}$ be the sheaf of relative holomorphic $p$-forms on ${\mathcal Y}$ in the sense of orbifolds
\cite[Def.\,1.7]{Steenbrink77}.
If $\pi\colon{\mathcal X}\to{\mathcal Y}$ denotes the projection, then 
we have the canonical identification $H^{q}(U,\Omega_{{\mathcal Y}/B}^{p})=H^{q}(\pi^{-1}(U),\Omega_{{\mathcal X}/B}^{p})^{G}$
for any open subset $U\subset{\mathcal Y}$ and for any $q\geq0$ (cf. \cite[Lemma\,1.8]{Steenbrink77}).
As a result, we get the canonical identification of direct images for all $p,q\geq0$:
\begin{equation}
\label{eqn:equality:direct:image}
R^{q}\overline{f}_{*}\Omega_{{\mathcal Y}/B}^{p}=(R^{q}f_{*}\Omega_{{\mathcal X}/B}^{p})^{G}.
\end{equation}
\par
Let ${\frak b}\in B$ be an arbitrary point and set $X:=X_{\frak b}$.
Let 
$$
{\frak f}\colon({\frak X},X)\to({\rm Def}(X),[X])
$$ 
be the Kuranishi family of $X$.
By \cite{Tian87}, \cite{Todorov89}, $({\rm Def}(X),[X])$ is smooth of dimension $h^{1}(\Theta_{X})=h^{1,2}(X)$ 
and is identified with the set germ at $0\in H^{1}(X,\Theta_{X})$. 
Since $G$ is finite, there is a $G$-invariant Hermitian metric on $X$.
By making use of this $G$-invariant Hermitian metric in the construction of Kuranishi family \cite[Chap.\,X]{Kuranishi71},
we see that the $G$-action on $X$ lifts to holomorphic $G$-actions on ${\frak X}$ and ${\rm Def}(X)$ 
so that ${\frak f}\colon{\frak X}\to{\rm Def}(X)$ is $G$-equivariant.
\par
Let 
$$
{\frak B}:={\rm Def}(X)^{G}=\{x\in{\rm Def}(X);\,g(x)=x\,\,(\forall\,g\in G)\}
$$ 
be the set of fixed points of the $G$-action. We define the map of germs 
$$
\varpi\colon(B,{\frak b})\to({\frak B},[X]),
\qquad
\varpi(b):=[X_{b}],
$$ 
where $[X_{b}]\in{\rm Def}(X)$ is the point corresponding to the complex structure on $X_{b}$.
Then the family of threefolds with $G$-action $f\colon({\mathcal X},G)\to B$ is induced from the Kuranishi family by the map $\varpi$.

\subsection
{The first Chern form of the relative tangent bundle}
\label{sect:6.4}
\par
Let $\omega_{\rm WP}^{G}$ be the positive $(1,1)$-form on ${\frak B}$ induced from the Weil-Petersson form $\omega_{\rm WP}$
(cf. Section~\ref{sect:1.5.2}), i.e.,
$$
\omega_{\rm WP}^{G}:=\omega_{\rm WP}|_{\frak B}.
$$
Let $\eta_{{\frak X}/{\rm Def}(X)}\in H^{0}({\rm Def}(X),{\frak f}_{*}K_{{\frak X}/{\rm Def}(X)})$ be a nowhere vanishing holomorphic section
and set $\eta_{{\mathcal X}/B}:=\varpi^{*}(\eta_{{\frak X}/{\rm Def}(X)})$. By its $G$-invariance, $\eta_{{\mathcal X}/B}$ descends to
a nowhere vanishing relative canonical form $\eta_{{\mathcal X}/B}$ on ${\mathcal Y}$ in the sense of orbifolds.
Since
$$
\eta_{{\mathcal X}/B}=\eta_{{\frak X}/{\rm Def}(X)}|_{\mu(B)},
$$ 
we deduce from  \eqref{eqn:Weil:Petersson:form:2} that
\begin{equation}
\label{eqn:Weil:Petersson:form:+}
\varpi^{*}\omega_{\rm WP}^{G}
=
-dd^{c}\log\|\eta_{{\mathcal X}/B}\|_{L^{2}}^{2}
=
c_{1}(f_{*}K_{{\mathcal X}/B},h_{L^{2}}).
\end{equation}
\par
By \cite[Eq.\,(4.2)]{FLY08} and \eqref{eqn:Weil:Petersson:form:+}, the following equality of $(1,1)$-forms on ${\mathcal X}$ holds:
\begin{equation}
\label{eqn:relative:first:Chern:form}
\begin{aligned}
c_{1}({\mathcal X}/B)
&=
-f^{*}
\left\{
\varpi^{*}\omega_{\rm WP}^{G}+dd^{c}\log{\rm Vol}({\mathcal X}/B)
\right\}
\\
&\quad
+dd^{c}\log
\left\{
\frac{i\,\eta_{{\mathcal X}/B}\wedge\overline{\eta}_{{\mathcal X}/B}}{\gamma_{{\mathcal X}/B}^{3}/3!}
f^{*}\left(\frac{{\rm Vol}({\mathcal X}/B)}{\|\eta_{{\mathcal X}/B}\|_{L^{2}}^{2}}\right)
\right\},
\end{aligned}
\end{equation}
which yields the equation
\begin{equation}
\label{eqn:relative:first:Chern:form:quotient}
\begin{aligned}
c_{1}({\mathcal Y}/B)
&=
-\overline{f}^{*}
\left\{
\varpi^{*}\omega_{\rm WP}^{G}+dd^{c}\log{\rm Vol}({\mathcal Y}/B)
\right\}
\\
&\quad
+dd^{c}\log
\left\{
\frac{i\,\eta_{{\mathcal Y}/B}\wedge\overline{\eta}_{{\mathcal Y}/B}}{\overline{\gamma}_{{\mathcal Y}/B}^{3}/3!}
\overline{f}^{*}\left(\frac{{\rm Vol}({\mathcal Y}/B)}{\|\eta_{{\mathcal Y}/B}\|_{L^{2}}^{2}}\right)
\right\}.
\end{aligned}
\end{equation}

\subsection
{The curvature of $R^{q}\overline{f}_{*}\Omega_{{\mathcal Y}/B}^{p}$: the case $p+q=3$}
\label{sect:6.5}
\par
Recall that the Kodaira-Spencer map 
$$
\rho\colon\Theta_{{\rm Def}(X)}\to R^{1}f_{*}\Omega_{{\frak X}/{\rm Def}(X)}^{2}\otimes(f_{*}K_{{\frak X}/{\rm Def}(X)})^{\lor}
$$
was defined in Section~\ref{sect:1.5.1}.
Since $\rho$ is $G$-equivariant by \eqref{eqn:Kodaira:Spencer:map}, we get an isomorphism of holomorphic vector bundles on ${\frak B}$
$$
\begin{aligned}
\rho
\colon 
\Theta_{\frak B}
=
(\Theta_{{\rm Def}(X)}|_{\frak B})^{G}
&\cong
(R^{1}{\frak f}_{*}\Theta_{{\frak X}/{\rm Def}(X)})^{G}
\\
&=
(R^{1}f_{*}\Omega_{{\frak X}/{\rm Def}(X)}^{2})^{G}\otimes(f_{*}K_{{\frak X}/{\rm Def}(X)})^{\lor}|_{\frak B}.
\end{aligned}
$$  
By \eqref{eqn:Weil:Petersson:metric}, we get an isometry of holomorphic Hermitian vector bundles on ${\frak B}$
\begin{equation}
\label{eqn:KS:isometry:+}
\begin{aligned}
\rho
\colon
(\Theta_{\frak B},\omega_{\rm WP}^{G})
&\cong
((R^{1}f_{*}\Theta_{{\frak X}/{\rm Def}(X)})^{G},h_{L^{2}})
\\
&\cong
((R^{1}f_{*}\Omega_{{\frak X}/{\rm Def}(X)}^{2})^{G},h_{L^{2}})\otimes((f_{*}K_{{\frak X}/{\rm Def}(X)})^{\lor},h_{L^{2}}^{-1}).
\end{aligned}
\end{equation}
Recall that the Ricci from of $({\frak B},\omega_{\rm WP}^{G})$ is defined as 
$$
{\rm Ric}\,\omega_{\rm WP}^{G}:=c_{1}(\Theta_{\frak B},\omega_{\rm WP}^{G}).
$$

\begin{lemma}
\label{lemma:curvature:direct:image:Kuranishi:family}
Set $h^{p,q}(Y):=\dim H^{q}(Y,\Omega_{Y}^{p})=\dim H^{q}(X,\Omega_{X}^{p})^{G}$.
When $p+q=3$, the following equality of $(1,1)$-forms on ${\frak B}$ holds:
$$
c_{1}((R^{q}f_{*}\Omega^{p}_{{\frak X}/{\rm Def}(X)})^{G},h_{L^{2}})
=
\begin{cases}
\begin{array}{ll}
\omega_{\rm WP}^{G}&(p,q)=(3,0),
\\
{\rm Ric}\,\omega_{\rm WP}^{G}+h^{1,2}(Y)\,\omega_{\rm WP}^{G}&(p,q)=(2,1),
\\
-{\rm Ric}\,\omega_{\rm WP}^{G}-h^{1,2}(Y)\,\omega_{\rm WP}^{G}&(p,q)=(1,2),
\\
-\omega_{\rm WP}^{G}&(p,q)=(0,3).
\end{array}
\end{cases}
$$
\end{lemma}

\begin{pf}
The result for $(p,q)=(3,0)$ follows from \eqref{eqn:Weil:Petersson:form:+}.
Since $h^{1,2}(Y)={\rm rk}\,(\Theta_{{\rm Def}(X)})^{G}$, we get by \eqref{eqn:KS:isometry:+} the result for $(p,q)=(2,1)$:
$$
\begin{aligned}
c_{1}((R^{1}f_{*}\Omega^{2}_{{\frak X}/{\rm Def}(X)})^{G},h_{L^{2}})
&=
c_{1}((\Theta_{{\rm Def}(X)})^{G},\omega_{\rm WP}^{G})
+
h^{1,2}(Y)c_{1}((f_{*}K_{{\frak X}/{\rm Def}(X)}),h_{L^{2}})|_{\frak B}
\\
&=
{\rm Ric}\,\omega_{\rm WP}^{G}+h^{1,2}(Y)\,\omega_{\rm WP}^{G}.
\end{aligned}
$$
Since the Serre duality $(R^{q}f_{*}\Omega^{p}_{{\frak X}/{\rm Def}(X)})^{\lor}\cong R^{3-q}f_{*}\Omega^{3-p}_{{\frak X}/{\rm Def}(X)}$
is a $G$-equivariant isometry with respect to the $L^{2}$-metrics on the direct image sheaves,
the results for $(p,q)=(1,2),(0,3)$ follow from those for $(p,q)=(2,1),(3,0)$.
\end{pf}

\begin{lemma}
\label{lemma:curvature:direct:image}
When $p+q=3$, the following equality of $(1,1)$-forms on $B$ holds:
$$
c_{1}((R^{q}f_{*}\Omega^{p}_{{\mathcal X}/B})^{G},h_{L^{2}})
=
\begin{cases}
\begin{array}{ll}
\varpi^{*}\omega_{\rm WP}^{G}&(p,q)=(3,0),
\\
\varpi^{*}{\rm Ric}\,\omega_{\rm WP}^{G}+h^{1,2}(Y)\,\varpi^{*}\omega_{\rm WP}^{G}&(p,q)=(2,1),
\\
-\varpi^{*}{\rm Ric}\,\omega_{\rm WP}^{G}-h^{1,2}(Y)\,\varpi^{*}\omega_{\rm WP}^{G}&(p,q)=(1,2),
\\
-\varpi^{*}\omega_{\rm WP}^{G}&(p,q)=(0,3).
\end{array}
\end{cases}
$$
\end{lemma}

\begin{pf}
Since $H^{q}(Y_{b},\Omega_{Y_{b}}^{p})=H^{q}(X_{b},\Omega_{X_{b}}^{p})^{G}$, $b\in B$, consists of primitive cohomology classes when $p+q=3$ 
and since the $L^{2}$-metric on the primitive cohomology coincides
with the cup-product, the map of germs $\varpi\colon(B,{\frak b})\to({\rm Def}(X),[X])$ induces an isometry of holomorphic Hermitian vector bundles on $(B,{\frak b})$
$$
((R^{q}f_{*}\Omega^{p}_{{\mathcal X}/B})^{G},h_{L^{2}})
=
\varpi^{*}((R^{q}f_{*}\Omega^{p}_{{\frak X}/{\rm Def}(X)})^{G},h_{L^{2}}).
$$
Hence the result follows from Lemma~\ref{lemma:curvature:direct:image:Kuranishi:family}. 
\end{pf}

\subsection
{The curvature of $R^{q}\overline{f}_{*}\Omega_{{\mathcal Y}/B}^{p}$: the case $p=q$}
\label{sect:6.6}
\par
For $b\in B$, the covolume of the lattice $H^{2}(Y_{b},{\bf Z})\subset H^{2}(Y_{b},{\bf R})$ relative to the $L^{2}$-metric
is defined in the same way as in Section~\ref{sect:1.2}.
If $\{{\bf e}_{1}(b),\ldots,{\bf e}_{r}(b)\}$ is a basis of $H^{2}(Y_{b},{\bf Z})/{\rm Torsion}$, then
$$
{\rm Vol}_{L^{2}}(H^{2}(Y_{b},{\bf Z}),\overline{\gamma}_{b})
:=
\det\left(\langle{\bf e}_{k}(b),{\bf e}_{l}(b)\rangle_{L^{2},\overline{\gamma}_{b}}\right)
=
\left\|
{\bf e}_{1}(b)\wedge\cdots\wedge{\bf e}_{r}(b)
\right\|_{L^{2},\overline{\gamma}_{b}}^{2},
$$
where $r:={\rm rk}_{\bf Z}\,H^{2}(Y_{b},{\bf Z})$ and $\langle\cdot,\cdot\rangle_{L^{2},\overline{\gamma}_{b}}$ denotes the $L^{2}$-metric
on $H^{2}(Y_{b},{\bf R})$ with respect to $\overline{\gamma}_{b}:=\overline{\gamma}_{{\mathcal Y}/B}|_{Y_{b}}$. 
Then $\langle\cdot,\cdot\rangle_{L^{2},\overline{\gamma}_{b}}$ depends only on the K\"ahler class $[\overline{\gamma}_{b}]$.
We define ${\rm Vol}_{L^{2}}(R^{2}\overline{f}_{*}{\bf Z})\in C^{\infty}(B)$ by
$$
{\rm Vol}_{L^{2}}(R^{2}\overline{f}_{*}{\bf Z})(b):={\rm Vol}_{L^{2}}(H^{2}(Y_{b},{\bf Z}),\overline{\gamma}_{b}),
\qquad
b\in B.
$$
By \cite[Lemma 4.12]{FLY08}, ${\rm Vol}_{L^{2}}(R^{2}\overline{f}_{*}{\bf Z})$ is constant on $B$ if the relative K\"ahler form 
$\overline{\gamma}_{{\mathcal Y}/B}$ is induced from a K\"ahler form on ${\mathcal X}$, i.e., if there is a K\"ahler form $\gamma_{\mathcal X}$ 
on ${\mathcal X}$ satisfying $\gamma_{{\mathcal X}/B}=\gamma_{\mathcal X}|_{\Theta_{{\mathcal X}/B}}$.

\begin{lemma}
\label{lemma:curvature:direct:image:(p,p)}
The following equality of $(1,1)$-forms on $B$ holds:
$$
c_{1}(R^{p}\overline{f}_{*}\Omega^{p}_{{\mathcal Y}/B},h_{L^{2}})
=
\begin{cases}
\begin{array}{ll}
-dd^{c}\log{\rm Vol}({\mathcal Y}/B)&(p=0),
\\
-dd^{c}\log{\rm Vol}_{L^{2}}R^{2}\overline{f}_{*}{\bf Z}&(p=1),
\\
dd^{c}\log{\rm Vol}_{L^{2}}R^{2}\overline{f}_{*}{\bf Z}&(p=2),
\\
dd^{c}\log{\rm Vol}({\mathcal Y}/B)&(p=3).
\end{array}
\end{cases}
$$
\end{lemma}

\begin{pf} 
Since $(f_{*}{\mathcal O}_{\mathcal X})^{G}=\overline{f}_{*}{\mathcal O}_{\mathcal Y}={\mathcal O}_{B}\cdot 1$
and since $\|1\|_{L^{2}}^{2}={\rm Vol}({\mathcal Y}/B)$, we get the result for $p=0$.
Similarly, since $(f_{*}\Omega_{{\mathcal X}/B}^{2})^{G}=(R^{2}f_{*}{\mathcal O}_{\mathcal X})^{G}=0$ by \eqref{eqn:condition:involution}, 
we get $\overline{f}_{*}\Omega_{{\mathcal Y}/B}^{2}=R^{2}\overline{f}_{*}{\mathcal O}_{\mathcal Y}=0$ by \eqref{eqn:equality:direct:image}.
Let $L_{1},\ldots,L_{r}$ be holomorphic line bundles on $Y_{b}$ such that
$\{c_{1}(L_{1}),\ldots,c_{1}(L_{r})\}$ is a basis of $H^{2}(Y_{b},{\bf Z})/{\rm Tors}(H^{2}(Y_{b},{\bf Z}))$.
Since $R^{2}\overline{f}_{*}{\mathcal O}_{\mathcal Y}=0$, it follows from the exactness of the sequence
$R^{1}\overline{f}_{*}{\mathcal O}_{\mathcal Y}^{*}\to R^{2}\overline{f}_{*}{\bf Z}\to R^{2}\overline{f}_{*}{\mathcal O}_{\mathcal Y}$
that every $L_{i}$ extends to a holomorphic line bundle ${\mathcal L}_{i}$ on ${\mathcal Y}$. 
Then $R^{1}\overline{f}_{*}\Omega^{1}_{{\mathcal Y}/B}=R^{2}\overline{f}_{*}{\bf Z}\otimes{\mathcal O}_{B}=
{\mathcal O}_{B}c_{1}({\mathcal L}_{1})+\cdots+{\mathcal O}_{B}c_{1}({\mathcal L}_{r})$ by \eqref{eqn:vanishing:direct:image},
\eqref{eqn:equality:direct:image}, so that we get
$\det R^{1}\overline{f}_{*}\Omega^{1}_{{\mathcal Y}/B}={\mathcal O}_{B}\cdot c_{1}({\mathcal L}_{1})\wedge\cdots\wedge c_{1}({\mathcal L}_{r})$
and
$$
c_{1}(R^{1}\overline{f}_{*}\Omega^{1}_{{\mathcal Y}/B},h_{L^{2}})
=
-dd^{c}\log\left\|c_{1}({\mathcal L}_{1})\wedge\cdots\wedge c_{1}({\mathcal L}_{r})\right\|_{L^{2}}^{2}
=
-dd^{c}\log{\rm Vol}_{L^{2}}(R^{2}\overline{f}_{*}{\bf Z}).
$$
This proves the assertion for $p=1$. The assertions for $p=2,3$ follow from those for $p=1,0$ by the Serre duality.
See \cite[p.200]{FLY08} for the details.
\end{pf}

\begin{corollary}
\label{cor:weighted:alternating:sum:chern:form:direct:image:+}
The following equality of $(1,1)$-forms on $B$ holds:
$$
\begin{aligned}
\sum_{p,q\geq0}(-1)^{p+q}p\,c_{1}(R^{q}\overline{f}_{*}\Omega^{p}_{{\mathcal Y}/B},h_{L^{2}})
&=
-\varpi^{*}{\rm Ric}\,\omega_{\rm WP}^{G}-\left(h^{1,2}(Y)+3\right)\,\varpi^{*}\omega_{\rm WP}^{G}
\\
&\quad
+
dd^{c}\log\left\{{\rm Vol}({\mathcal Y}/B)^{3}\,{\rm Vol}_{L^{2}}(R^{2}\overline{f}_{*}{\bf Z})\right\}.
\end{aligned}
$$
\end{corollary}

\begin{pf}
The result follows from Lemmas~\ref{lemma:curvature:direct:image} and \ref{lemma:curvature:direct:image:(p,p)}.
\end{pf}

\subsection
{The direct images of some Chern forms}
\label{sect:6.7}
\par

\begin{definition}
\label{def:Bott-Chern:anomaly}
Define $C^{\infty}$ functions $A({\mathcal Y}/B)$, $A(\overline{\mathcal C}_{\bar{\lambda}}/B)$, $A(\overline{\frak p}/B)$ on $B$ by
$$
A({\mathcal Y}/B)
:=
\exp
\left[
-\frac{1}{12}\overline{f}_{*}
\left\{
\log
\left(
\frac{i\,\eta_{{\mathcal Y}/B}\wedge\overline{\eta}_{{\mathcal Y}/B}}{\overline{\gamma}_{{\mathcal Y}/B}^{3}/3!}
\cdot
\overline{f}^{*}\frac{{\rm Vol}({\mathcal Y}/B)}{\|\eta_{{\mathcal Y}/B}\|_{L^{2}}^{2}}
\right)
\,
c_{3}({\mathcal Y}/B)
\right\}
\right],
$$
$$
A(\overline{\mathcal C}_{\bar{\lambda}}/B)
:=
\exp
\left[
-\frac{1}{12}\overline{f}_{*}
\left.
\left\{
\log
\left(
\frac{i\,\eta_{{\mathcal Y}/B}\wedge\overline{\eta}_{{\mathcal Y}/B}}{\overline{\gamma}_{{\mathcal Y}/B}^{3}/3!}
\cdot
\overline{f}^{*}\frac{{\rm Vol}({\mathcal Y}/B)}{\|\eta_{{\mathcal Y}/B}\|_{L^{2}}^{2}}
\right)
\,
c_{1}(\overline{\mathcal C}_{\bar{\lambda}}/B)
\right|_{\overline{\mathcal C}_{\bar{\lambda}}}
\right\}
\right],
$$
$$
A(\overline{\frak p}/B)
:=
\left.
\overline{f}_{*}
\left(
\frac{i\,\eta_{{\mathcal Y}/B}\wedge\overline{\eta}_{{\mathcal Y}/B}}{\gamma_{{\mathcal Y}/B}^{3}/3!}
\right|_{\overline{\frak p}}
\cdot
\overline{f}^{*}\frac{{\rm Vol}({\mathcal Y}/B)}{\|\eta_{{\mathcal Y}/B}\|_{L^{2}}^{2}}
\right)
\qquad
(\overline{\frak p}\subset{\mathcal Y}^{(0)}).
$$
\end{definition}

\begin{proposition}
\label{prop:direct:image:Chern:forms:1}
The following equality of $(1,1)$-forms on ${\mathcal Y}$ holds:
$$
\overline{f}_{*}\{c_{1}({\mathcal Y}/B)c_{3}({\mathcal Y}/B)\}
=
-\int_{Y}c_{3}(Y)
\left\{
\varpi^{*}\omega_{\rm WP}^{G}
+
dd^{c}\log{\rm Vol}({\mathcal Y}/B)
\right\}
-
12\,dd^{c}\log A({\mathcal Y}/B).
$$
\end{proposition}

\begin{pf}
Substituting \eqref{eqn:relative:first:Chern:form:quotient} into $\overline{f}_{*}\{c_{1}({\mathcal Y}/B)c_{3}({\mathcal Y}/B)\}$
and using the projection formula and the commutativity of $dd^{c}$ and $\overline{f}_{*}$, we get the result.
See also \cite[p.197 Eq.(4.3)]{FLY08}. 
\end{pf}

\begin{proposition}
\label{prop:direct:image:Chern:forms:2}
The following equality of $(1,1)$-forms on $B$ holds:
$$
\begin{aligned}
(\overline{f}|_{\overline{\mathcal C}_{\bar{\lambda}}})_{*}\{c_{1}({\mathcal Y}/B)c_{1}(\overline{\mathcal C}_{\bar{\lambda}}/B)\}
&=
-\int_{\overline{C}_{\bar{\lambda}}}c_{1}(\overline{C}_{\bar{\lambda}})
\left\{
\varpi^{*}\omega_{\rm WP}^{G}
+
dd^{c}\log{\rm Vol}({\mathcal Y}/B)
\right\}
\\
&\quad
-
12\,dd^{c}\log A(\overline{\mathcal C}_{\bar{\lambda}}/B).
\end{aligned}
$$
\end{proposition}

\begin{pf}
By \eqref{eqn:relative:first:Chern:form:quotient}, we get
\begin{equation}
\label{eqn:direct:image:equivariant:Todd:Chern:character}
\begin{aligned}
\,&
(\overline{f}|_{\overline{\mathcal C}_{\bar{\lambda}}})_{*}\{c_{1}({\mathcal Y}/B)c_{1}(\overline{\mathcal C}_{\bar{\lambda}}/B)\}
\\
&=
-(\overline{f}|_{\overline{\mathcal C}_{\bar{\lambda}}})_{*}
\left[
(\overline{f}|_{\overline{\mathcal C}_{\bar{\lambda}}})^{*}
\left\{
\varpi^{*}\omega_{\rm WP}^{G}+dd^{c}\log{\rm Vol}({\mathcal Y}/B)
\right\}
\wedge
c_{1}(\overline{\mathcal C}_{\bar{\lambda}}/B)
\right]
\\
&\quad
+
(\overline{f}|_{\overline{\mathcal C}_{\bar{\lambda}}})_{*}
\left[
dd^{c}\log
\left\{
\frac{i\,\eta_{{\mathcal Y}/B}\wedge\overline{\eta}_{{\mathcal Y}/B}}{\overline{\gamma}_{{\mathcal Y}/B}^{3}/3!}
\overline{f}^{*}\left(\frac{{\rm Vol}({\mathcal Y}/B)}{\|\eta_{{\mathcal Y}/B}\|_{L^{2}}^{2}}\right)
\right\}
\wedge
c_{1}(\overline{\mathcal C}_{\bar{\lambda}}/B)
\right]
\\
&=
-\int_{\overline{C}_{\bar{\lambda},b}}c_{1}(\overline{C}_{\bar{\lambda},b})
\cdot
\left\{
\varpi^{*}\omega_{\rm WP}^{G}
+
dd^{c}\log{\rm Vol}({\mathcal Y}/B)
\right\}
-
12\,dd^{c}\log A(\overline{\mathcal C}_{\bar{\lambda}}/B).
\end{aligned}
\end{equation}
To get the second equality, we used the projection formula and the commutativity of $dd^{c}$ and 
$(\overline{f}|_{\overline{\mathcal C}_{\bar{\lambda}}})_{*}$.
\end{pf}

\begin{proposition}
\label{prop:direct:image:Chern:forms:3}
The following equality of $(1,1)$-forms on $B$ holds:
$$
(\overline{f}|_{\overline{\frak p}})_{*}\{c_{1}({\mathcal Y}/B)\}
=
-\varpi^{*}\omega_{\rm WP}^{G}-dd^{c}\log{\rm Vol}({\mathcal Y}/B)+dd^{c}\log A(\overline{\frak p}/B).
$$
\end{proposition}

\begin{pf}
The result follows from \eqref{eqn:relative:first:Chern:form:quotient} and the projection formula.
\end{pf}

\subsection
{The curvature formula}
\label{sect:5.8}
\par
We define $\tau_{\rm BCOV}^{\rm orb}({\mathcal Y}/B)\in C^{\infty}(B)$ as
$$
\begin{aligned}
\tau_{\rm BCOV}^{\rm orb}({\mathcal Y}/B)
&:=
T_{\rm BCOV}({\mathcal Y}/B)\,
{\rm Vol}({\mathcal Y}/B)^{-3+\frac{{\chi}^{\rm orb}(Y)}{12}}
{\rm Vol}_{L^{2}}(R^{2}\overline{f}_{*}{\bf Z})^{-1}
A({\mathcal Y}/B)
\\
&\quad\times
\prod_{\bar{\lambda}\in\overline{\Lambda}}
\left\{
\tau(\overline{\mathcal C}_{\bar{\lambda}}/B){\rm Vol}(\overline{\mathcal C}_{\bar{\lambda}}/B)
\right\}^{-\mu_{\bar{\lambda}}}
A(\overline{\mathcal C}_{\bar{\lambda}}/B)^{\frac{n_{\lambda}^{2}-1}{n_{\lambda}}}
\prod_{\overline{\frak p}\subset{\mathcal Y}^{(0)}}A(\overline{\frak p}/B)^{\epsilon(G_{\frak p})}.
\end{aligned}
$$

\begin{theorem}
\label{thm:variational:formula:orbifold:BCOV:invariant}
The following equality of $(1,1)$-forms on $B$ holds
$$
\begin{aligned}
-dd^{c}\log\tau_{\rm BCOV}^{\rm orb}({\mathcal Y}/B)
&=
\left(
\frac{{\chi}^{\rm orb}(Y)}{12}+h^{1,2}(Y)+3
\right)
\varpi^{*}\omega_{\rm WP}^{G}
+
\varpi^{*}{\rm Ric}\,\omega_{\rm WP}^{G}
\\
&\quad
+
\sum_{\bar{\lambda}\in\overline{\Lambda}}
\mu_{\bar{\lambda}}\,J_{\overline{\mathcal C}_{\bar{\lambda}}/B}^{*}\omega_{{\mathcal A}_{\overline{h}_{\bar{\lambda}}}}.
\end{aligned}
$$
\end{theorem}

\begin{pf}
By Theorem~\ref{thm:curvature:orbifold:BCOV:torsion:fixed:curve}, we get
\begin{equation}
\label{eqn:curvature:orbifold:BCOV:1}
\begin{aligned}
\,&
-dd^{c}\log T_{\rm BCOV}({\mathcal Y}/B)
+
\sum_{\bar{\lambda}\in\overline{\Lambda}}
\mu_{\lambda}\,dd^{c}\log\left\{\tau(\overline{\mathcal C}_{\bar{\lambda}}/B){\rm Vol}(\overline{\mathcal C}_{\bar{\lambda}}/B)\right\}
\\
&=
-\frac{1}{12}\overline{f}_{*}\{c_{1}({\mathcal Y}/B)c_{3}({\mathcal Y}/B)\}
-
\frac{1}{12}\sum_{\bar{\lambda}\in\overline{\Lambda}}\frac{n_{\lambda}^{2}-1}{n_{\lambda}}
(\overline{f}|_{\overline{\mathcal C}_{\bar{\lambda}}})_{*}\{c_{1}({\mathcal Y}/B)c_{1}(\overline{\mathcal C}_{\bar{\lambda}}/B)\}
\\
&\quad
+
\sum_{\overline{\frak p}\subset{\mathcal Y}^{(0)}}\epsilon(G_{\frak p})\,(\overline{f}|_{\overline{\frak p}})_{*}\{c_{1}({\mathcal Y}/B)\}
-
\sum_{p,q\geq0}(-1)^{p+q}p\,c_{1}(\det R^{q}\overline{f}_{*}\Omega^{p}_{{\mathcal Y}/B},h_{L^{2}})
\\
&\quad
+
\sum_{\bar{\lambda}\in\overline{\Lambda}}
\mu_{\bar{\lambda}}\,J_{\overline{\mathcal C}_{\bar{\lambda}}/B}^{*}\omega_{{\mathcal A}_{\overline{h}_{\bar{\lambda}}}}.
\end{aligned}
\end{equation}
Substituting the formulae in Propositions~\ref{prop:direct:image:Chern:forms:1}, \ref{prop:direct:image:Chern:forms:2}, \ref{prop:direct:image:Chern:forms:3}
and Corollary~\ref{cor:weighted:alternating:sum:chern:form:direct:image:+} into \eqref{eqn:curvature:orbifold:BCOV:1},
we get
\begin{equation}
\label{eqn:curvature:orbifold:BCOV:3}
\begin{aligned}
\,&
-dd^{c}\log T_{\rm BCOV}({\mathcal Y}/B)
+
\sum_{\bar{\lambda}\in\overline{\Lambda}}
\mu_{\lambda}\,dd^{c}\log\left\{\tau(\overline{\mathcal C}_{\bar{\lambda}}/B){\rm Vol}(\overline{\mathcal C}_{\bar{\lambda}}/B)\right\}
\\
&=
\left(
\frac{\int_{Y}c_{3}(Y)}{12}
+
\frac{1}{12}\sum_{\bar{\lambda}\in\overline{\Lambda}}\frac{n_{\lambda}^{2}-1}{n_{\lambda}}
\int_{\overline{C}_{\bar{\lambda}}}c_{1}(\overline{C}_{\bar{\lambda}})
-
\sum_{\overline{\frak p}\subset{\mathcal Y}^{(0)}}\epsilon(G_{\frak p})
+
h^{1,2}(Y)+3
\right)
\varpi^{*}\omega_{\rm WP}^{G}
\\
&\quad
+
\varpi^{*}{\rm Ric}\,\omega_{\rm WP}^{G}
+
\sum_{\bar{\lambda}\in\overline{\Lambda}}\mu_{\lambda}\,J_{\overline{\mathcal C}_{\bar{\lambda}}/B}^{*}
\omega_{{\mathcal A}_{\overline{h}_{\bar{\lambda}}}}
\\
&\quad
+
\left(
\frac{\int_{Y}c_{3}(Y)}{12}
+
\frac{1}{12}\sum_{\bar{\lambda}\in\overline{\Lambda}}\frac{n_{\lambda}^{2}-1}{n_{\lambda}}\int_{\overline{C}_{\bar{\lambda}}}c_{1}(\overline{C}_{\bar{\lambda}})
-
\sum_{\overline{\frak p}\subset{\mathcal Y}^{(0)}}\epsilon(G_{\frak p})
-
3
\right)
dd^{c}\log{\rm Vol}({\mathcal Y}/B)
\\
&\quad
-dd^{c}\log{\rm Vol}_{L^{2}}(R^{2}\overline{f}_{*}{\bf Z})
+
dd^{c}\log A({\mathcal Y}/B)
+
\sum_{\bar{\lambda}\in\overline{\Lambda}}\frac{n_{\lambda}^{2}-1}{n_{\lambda}}\,dd^{c}\log A(\overline{\mathcal C}_{\bar{\lambda}}/B)
\\
&\quad
-12\sum_{\overline{\frak p}\in{\mathcal Y}^{(0)}}\epsilon(G_{\frak p})\,dd^{c}\log A(\overline{\frak p}/B).
\end{aligned}
\end{equation}
Substituting \eqref{eqn:Euler:X:G} into \eqref{eqn:curvature:orbifold:BCOV:3}, we get the result.
\end{pf}

\begin{theorem}
\label{Thm:invariance:property:orbifold:BCOV:invariant}
The number $\tau_{\rm BCOV}^{\rm orb}(Y,\overline{\gamma})$ is independent of the choice of a K\"ahler form $\overline{\gamma}$ on $Y$.
\end{theorem}

\begin{pf}
Let $\gamma_{0}$ and $\gamma_{\infty}$ be two $G$-invariant K\"ahler forms on $X$, which are identified with the corresponding Hermitian metrics.
Set ${\mathcal X}:=X\times{\bf P}^{1}$ and let $f:={\rm pr}_{2}\colon{\mathcal X}\to{\bf P}^{1}$ be the projection.
By the triviality of the family $f\colon{\mathcal X}\to{\bf P}^{1}$, the $G$-action on $X$ extends to a $G$-action on ${\mathcal X}$. 
Define  the Hermitian metric $\gamma_{{\mathcal X}/{\bf P}^{1}}:=\{\gamma_{t}\}_{t\in{\bf P}^{1}}$ on $\Theta_{{\mathcal X}/{\bf P}^{1}}$ by
$\gamma_{t}:=(\gamma_{0}+|t|^{2}\gamma_{\infty})/(1+|t|^{2})$.
Then $\gamma_{{\mathcal X}/{\bf P}^{1}}$ is fiberwise K\"ahler.
We set $Y:=X/G$, ${\mathcal Y}:={\mathcal X}/G$ and
$\tau_{\rm BCOV}^{\rm orb}({\mathcal Y}/{\bf P}^{1})(t):=\tau_{\rm BCOV}^{\rm orb}(Y,\overline{\gamma}_{t})$ for $t\in{\bf P}^{1}$.
By Theorem~\ref{thm:variational:formula:orbifold:BCOV:invariant} applied to the trivial family $f\colon{\mathcal Y}\to{\bf P}^{1}$,
we get on ${\bf P}^{1}$
$$
-dd^{c}\log\tau_{\rm BCOV}^{\rm orb}({\mathcal Y}/{\bf P}^{1})=0,
$$
because $\varpi$ and $J_{\overline{\mathcal C}_{\lambda}/{\bf P}^{1}}$ are constant maps. 
Since $\tau_{\rm BCOV}^{\rm orb}({\mathcal Y}/{\bf P}^{1})$ is a harmonic function on ${\bf P}^{1}$, it is a constant.
Hence $\tau_{\rm BCOV}^{\rm orb}(Y,\overline{\gamma}_{0})=\tau_{\rm BCOV}^{\rm orb}(Y,\overline{\gamma}_{\infty})$.
\end{pf}

\begin{definition}
\label{def:orbifold:BCOV:invariant:2}
Let $X$ be a smooth projective threefold with trivial canonical line bundle. 
Let $G$ be a finite abelian group of automorphisms of $X$ satisfying \eqref{eqn:condition:cohomology:group:action}.
Then the {\em (orbifold) BCOV invariant} of $X/G$ is defined as
$$
\tau_{\rm BCOV}^{\rm orb}(X/G):=\tau_{\rm BCOV}^{\rm orb}(X/G,\overline{\gamma}),
$$
where $\gamma$ is an arbitrary $G$-invariant K\"ahler form on $X$.
\end{definition}

\section
{BCOV invariants for general abelian Calabi-Yau orbifolds}
\label{sect:7}
\par
In this section, we extend BCOV invariants to general (i.e., possibly non-global) abelian Calabi-Yau orbifolds 
by using the anomaly formula for orbifolds \cite[Th.\,0.1]{Ma05}.
\par
Let $Y$ be a $3$-dimensional Calabi-Yau orbifold. We assume that $Y$ is {\em abelian}. 
Namely, for any $y\in Y$, there exist a neighborhood $\overline{U}_{y}$ of $y$ in $Y$, 
a finite abelian subgroup $G_{y}\subset{\rm SL}({\bf C}^{3})$ and a $G_{y}$-invariant open neighborhood $U_{y}$ of $0$ in ${\bf C}^{3}$ 
such that $(\overline{U}_{y},y)\cong(U_{y}/G_{y},0)$. 
We define $\Sigma U_{y}=\Sigma^{(1)}U_{y}\amalg\Sigma^{(0)}U_{y}$ and
$\widetilde{\Sigma}^{(0)}U_{y}$ as in Sections~\ref{sect:4.1} and \ref{sect:4.2}. 
We set $\widetilde{\Sigma}U_{y}:=\Sigma^{(1)}U_{y}\amalg\widetilde{\Sigma}^{(0)}U_{y}$.
Similarly, we define 
$\Sigma^{(k)}\overline{U}_{y}:=\Sigma^{(k)}U_{y}/G_{y}$ and $\widetilde{\Sigma}^{(0)}\overline{U}_{y}:=\widetilde{\Sigma}^{(0)}U_{y}/G_{y}$.
We set $\Sigma\overline{U}_{y}:=\Sigma^{(0)}\overline{U}_{y}\amalg{\Sigma}^{(1)}\overline{U}_{y}$ and 
$\widetilde{\Sigma}\overline{U}_{y}:=\widetilde{\Sigma}^{(0)}\overline{U}_{y}\amalg{\Sigma}^{(1)}\overline{U}_{y}$.
Gluing these $\Sigma\overline{U}_{y}$ and $\widetilde{\Sigma}\overline{U}_{y}$, we get complex orbifolds
$\Sigma Y:=\bigcup_{y\in Y}{\Sigma}\overline{U}_{y}$ and $\widetilde{\Sigma}Y:=\bigcup_{y\in Y}\widetilde{\Sigma}\overline{U}_{y}$.
Then $\Sigma Y={\Sigma}^{(1)}Y\amalg\Sigma^{(0)}Y$ and $\widetilde{\Sigma}Y=\Sigma^{(1)}Y\amalg\widetilde{\Sigma}^{(0)}Y$.
Let $\overline{\nu}_{y}\colon\Sigma\overline{U}_{y}\to\overline{U}_{y}$ be the map as in Section~\ref{sect:4.1.2}.
Gluing $\overline{\nu}_{y}$, we get a map $\overline{\nu}\colon\Sigma Y\to Y$. 
Let $\overline{\nu}(\Sigma^{(1)}Y)=\bigcup_{\bar{\lambda}\in\overline{\Lambda}}\overline{C}_{\bar\lambda}$ be the irreducible decomposition.
By the local description of $\overline{\nu}(\Sigma^{(1)}\overline{U}_{y})$ in Sections~\ref{sect:4.2} and \ref{sect:4.5.1},
$\overline{C}_{\bar\lambda}$ can have singular points isomorphic to the union of some coordinate axes of ${\bf C}^{3}$.
Since $Y$ is an abelian orbifold, there exists $n_{\bar\lambda}\in{\bf Z}_{>1}$ for each $\bar{\lambda}\in\overline{\Lambda}$ such that 
$G_{y}$ is a cyclic group of order $n_{\lambda}$ for generic $y\in\overline{C}_{\bar\lambda}$. In particular, $n_{\bar\lambda}=|G_{y}|$
for generic $y\in\overline{C}_{\bar\lambda}$.
Let $\widehat{C}_{\bar\lambda}$ be the normalization of $\overline{C}_{\bar\lambda}$. We set
$$
Y^{(1)}:=\amalg_{\bar{\lambda}\in\overline{\Lambda}}\widehat{C}_{\bar{\lambda}},
\qquad
Y^{(0)}:=\overline{\nu}(\widetilde{\Sigma}^{(0)}Y).
$$
Define $\iota\colon Y^{(1)}\to\overline{\nu}(\Sigma^{(1)}Y)$ as the normalization. 
By the definition of $\overline{\nu}_{y}$, $\overline{\nu}$ induces an unramified covering
$p\colon\Sigma^{(1)}Y\to Y^{(1)}$ with $\iota\circ p=\overline{\nu}$.
By \eqref{eqn:structure:inertia:orbifold}, the restriction
$p|_{p^{-1}(\widehat{C}_{\bar{\lambda}})}\colon p^{-1}(\widehat{C}_{\bar{\lambda}})\to\widehat{C}_{\bar{\lambda}}$ 
is an \'etale map of degree $n_{\bar\lambda}-1$.
If $Y$ is a global orbifold, $p^{-1}(\widehat{C}_{\bar{\lambda}})$ consists of $n_{\bar\lambda}-1$ components
and each component is isomorphic to $\widehat{C}_{\bar{\lambda}}$.
Since $Y$ is not necessarily a global orbifold, the number of the components of $p^{-1}(\widehat{C}_{\bar{\lambda}})$
can be smaller than $n_{\bar\lambda}-1$ in general. 
\par 
For a K\"ahler form $\overline{\gamma}$ on $Y$ in the sense of orbifolds, we define ${\chi}^{\rm orb}(Y)$ by \eqref{eqn:Euler:X:G},
$A(Y,\overline{\gamma})$ by \eqref{eqn:anomaly:1}, $A(\overline{C}_{\bar{\lambda}},\overline{\gamma}|_{\overline{C}_{\bar{\lambda}}})$
by \eqref{eqn:anomaly:2}, $A(\overline{\frak p},\overline{\gamma}|_{\overline{\frak p}})$ $(\overline{\frak p}\in Y^{(0)})$
by \eqref{eqn:anomaly:3}. Then ${\chi}^{\rm orb}(Y)$ is independent of the choice of $\overline{\gamma}$. 
Let 
$$
\Sigma^{(1)}Y=\amalg_{\alpha\in A}D_{\alpha}
$$ 
be the decomposition into connected components. When $p(D_{\alpha})=\widehat{C}_{\lambda}$, the restriction
$p|_{D_{\alpha}}\colon D_{\alpha}\to\widehat{C}_{\lambda}$ is an unramified covering.
We set
$$
\tau(\Sigma Y,\overline{\gamma}|_{\Sigma Y})
:=
\prod_{\alpha\in A}\tau(D_{\alpha},\overline{\gamma}|_{D_{\alpha}}),
\qquad
{\rm Vol}(\Sigma Y,\overline{\gamma}|_{\Sigma Y})
:=
\prod_{\alpha\in A}{\rm Vol}(D_{\alpha},\overline{\gamma}|_{D_{\alpha}})
$$

\begin{definition}
\label{def:orbifold:BCOV:invariant:general}
Define the {\em orbifold BCOV invariant} of $(Y,\overline{\gamma})$ by
$$
\begin{aligned}
\tau_{\rm BCOV}^{\rm orb}(Y,\overline{\gamma})
&:=
T_{\rm BCOV}(Y,\overline{\gamma})\,
{\rm Vol}(Y,\overline{\gamma})^{-3+\frac{\chi^{\rm orb}(Y)}{12}}
{\rm Vol}_{L^{2}}(H^{2}(Y,{\bf Z}),[\overline{\gamma}])^{-1}
A(Y,\overline{\gamma})
\\
&\quad\times
\tau(\Sigma Y,\overline{\gamma}|_{\Sigma Y})^{-1}{\rm Vol}(\Sigma Y,\overline{\gamma}|_{\Sigma Y})^{-1}
\\
&\quad\times
\prod_{\bar{\lambda}\in\overline{\Lambda}}
A(\overline{C}_{\bar{\lambda}},\overline{\gamma}|_{\overline{C}_{\bar{\lambda}}})^{\frac{n_{\lambda}^{2}-1}{n_{\lambda}}}
\prod_{\overline{\frak p}\in Y^{(0)}}A(\overline{\frak p},\overline{\gamma}|_{\overline{\frak p}})^{\epsilon(G_{\frak p})}.
\end{aligned}
$$
\end{definition}

For a holomorphic orbi-vector bundle $E$ on a complex orbifold $M$, Hermitian metrics $h$, $h'$ on $E$
and ${\rm GL}({\bf C}^{r})$-invariant polynomial $\phi(\cdot)$ on $\frak{gl}({\bf C}^{r})$, $r:={\rm rk}(E)$,
we denote by $\widetilde{\phi}^{\Sigma}(E;h,h')\in\oplus_{p\geq0}A^{p,p}(M\amalg\Sigma M)/{\rm Im}\partial+{\rm Im}\bar{\partial}$ 
the Bott-Chern secondary class \cite[I e), f)]{BGS88}, \cite[Sect.\,1.2]{Ma05} such that
$$
-dd^{c}\widetilde{\phi}^{\Sigma}(E;h,h')=\phi^{\Sigma}(E,h)-\phi^{\Sigma}(E,h').
$$
\par
Following \cite{FLY08}, we define the complex line $\lambda_{\rm BCOV}(Y)$ as
\begin{equation}
\label{eqn:BCOV:bundle}
\begin{aligned}
\,&
\lambda_{\rm BCOV}(Y)
:=
\bigotimes_{p\geq0}\lambda(\Omega_{Y}^{p})^{(-1)^{p}p}
=
\bigotimes_{p,q\geq0}(\det H^{q}(Y,\Omega_{Y}^{p}))^{(-1)^{p+q}p}
=
\\
&
\det H^{2}(Y,\Omega_{Y}^{1})^{-1}\otimes(\det H^{1}(Y,\Omega_{Y}^{2}))^{-2}\otimes(\det H^{0}(Y,\Omega_{Y}^{3}))^{-3}
\otimes\bigotimes_{r=1}^{3}(\det H^{2r}(Y,{\bf C}))^{r}.
\end{aligned}
\end{equation}
We denote by $\|\cdot\|_{\lambda_{\rm BCOV}(Y),L^{2},\overline{\gamma}}$ and $\|\cdot\|_{\lambda_{\rm BCOV}(Y),Q,\overline{\gamma}}$
the $L^{2}$ and Quillen metrics \cite{BGS88}, \cite{Ma05} on $\lambda_{\rm BCOV}(Y)$ with respect to $\overline{\gamma}$, respectively. 
By definition,
$$
\|\cdot\|_{\lambda_{\rm BCOV}(Y),Q,\overline{\gamma}}^{2}
=
T_{\rm BCOV}(Y,\overline{\gamma})\,
\|\cdot\|_{\lambda_{\rm BCOV}(Y),L^{2},\overline{\gamma}}^{2}.
$$
Similarly, the Quillen metric on 
$\lambda({\mathcal O}_{D_{\alpha}})=
\det H^{0}(D_{\alpha},{\mathcal O}_{D_{\alpha}})
\otimes
\det H^{1}(D_{\alpha},{\mathcal O}_{D_{\alpha}})^{-1}$
with respect to $\overline{\gamma}|_{D_{\alpha}}$ (resp. $\overline{\gamma}'|_{D_{\alpha}}$) is denoted by 
$\|\cdot\|_{\lambda({\mathcal O}_{D_{\alpha}}),Q,\overline{\gamma}}$
(resp. $\|\cdot\|_{\lambda({\mathcal O}_{D_{\alpha}}),Q,\overline{\gamma}'}$).
Set
$$
\lambda({\mathcal O}_{\Sigma Y})
:=
\bigotimes_{\alpha\in A}\lambda({\mathcal O}_{D_{\alpha}}),
\qquad
\|\cdot\|_{\lambda({\mathcal O}_{\Sigma Y}),Q,\overline{\gamma}}
:=
\bigotimes_{\alpha\in A}
\|\cdot\|_{\lambda({\mathcal O}_{D_{\alpha}}),Q,\overline{\gamma}}.
$$

\begin{lemma}
\label{lemma:anomaly:formula:Ma}
Let $\overline{\gamma}$ and $\overline{\gamma}'$ be K\"ahler forms on $Y$ in the sense of orbifolds.
Then 
$$
\begin{aligned}
\,&
\log
\left(
\frac{\|\cdot\|_{\lambda_{\rm BCOV}(Y),Q,\overline{\gamma}}}
{\|\cdot\|_{\lambda_{\rm BCOV}(Y),Q,\overline{\gamma}'}}
\right)^{2}
-
\log
\left(
\frac{\|\cdot\|_{\lambda({\mathcal O}_{\Sigma Y}),Q,\overline{\gamma}}}
{\|\cdot\|_{\lambda({\mathcal O}_{\Sigma Y}),Q,\overline{\gamma}'}}
\right)^{2}
\\
&=
-\frac{1}{12}\int_{Y}\widetilde{c_{1}c_{3}}(Y;\overline{\gamma},\overline{\gamma}')
-
\frac{1}{12}\sum_{\bar{\lambda}\in\overline{\Lambda}}\frac{n_{\overline{\lambda}}^{2}-1}{n_{\bar{\lambda}}}
\int_{\overline{C}_{\bar{\lambda}}}\widetilde{c_{1}(Y)c_{1}(\overline{C}_{\bar{\lambda}})}(\overline{\gamma},\overline{\gamma}')
\\
&\quad
+
\sum_{\overline{\frak p}\in Y^{(0)}}\epsilon(G_{\overline{\frak p}})\,
\widetilde{c_{1}(Y)}(\overline{\gamma},\overline{\gamma}')|_{\overline{\frak p}}.
\end{aligned}
$$
\end{lemma}

\begin{pf}
For $\overline{\frak p}\in Y^{(0)}$, we fix an isomorphism 
$(\overline{U}_{\overline{\frak p}},\overline{\frak p})\cong(U_{\frak p}/G_{\frak p},{\frak p})$,
where the origin of $U_{\overline{\frak p}}$ is denoted by ${\frak p}$ and we set $U_{\frak p}:=U_{\overline{\frak p}}$
and $G_{\frak p}:=G_{\overline{\frak p}}$.
As in Section~\ref{sect:4.3.5}, we get the splitting
$$
N_{{\frak p}/U_{\frak p}}
=
N_{{\frak p}/U_{\frak p}}(\chi_{{\frak p},1})
\oplus
N_{{\frak p}/U_{\frak p}}(\chi_{{\frak p},2})
\oplus
N_{{\frak p}/U_{\frak p}}(\chi_{{\frak p},3}).
$$
Since $\overline{\gamma}$ and $\overline{\gamma}'$ are K\"ahler forms on $Y$ in the sense of orbifolds, they are viewed as
positive $(1,1)$-forms $\gamma$ and $\gamma'$ on $U_{\frak p}$, respectively. We define
$$
\widetilde{\rho}_{{\frak p},k}(\overline{\gamma},\overline{\gamma}')({\frak p})
:=
\widetilde{c}_{1}(N_{{\frak p}/U_{\frak p}}(\chi_{{\frak p},k});\gamma,\gamma')({\frak p}),
$$
which is the Bott-Chern class associated to the characteristic form \eqref{eqn:***}. 
By \eqref{eqn:relation:sum:Chern:forms:relative:dim:0},
$$
\widetilde{\rho}_{{\frak p},1}({\gamma},{\gamma}')({\frak p})
+
\widetilde{\rho}_{{\frak p},2}({\gamma},{\gamma}')({\frak p})
+
\widetilde{\rho}_{{\frak p},3}({\gamma},{\gamma}')({\frak p})
=
\widetilde{c_{1}(Y)}({\gamma},{\gamma}')|_{\overline{\frak p}}.
$$
\par
We apply the anomaly formula of Ma \cite[Th.\,0.1]{Ma05} to the line $\lambda_{\rm BCOV}(Y)$.
Recall that $\varPi\colon Y\amalg\Sigma Y\to Y\amalg Y^{(1)}\amalg Y^{(0)}$ is the projection.
By Proposition~\ref{prop:orb:Tod:ch}, we get
\begin{equation}
\label{eqn:anomaly:formula:BCOV:bundle}
\begin{aligned}
\,&
\log
\left(
\frac{\|\cdot\|_{\lambda_{\rm BCOV}(Y),Q,\overline{\gamma}}}
{\|\cdot\|_{\lambda_{\rm BCOV}(Y),Q,\overline{\gamma}'}}
\right)^{2}
=
\int_{Y\amalg\Sigma Y}
\sum_{p\geq0}(-1)^{p}p\,
\widetilde{
{\rm Td}^{\Sigma}(\Theta_{Y}){\rm ch}^{\Sigma}(\Omega_{Y}^{p})
}
(\overline{\gamma},\overline{\gamma}')
\\
&=
\int_{Y\amalg Y^{(1)}\amalg Y^{(0)}}
\varPi_{*}\{
\sum_{p\geq0}(-1)^{p}p\,
\widetilde{
{\rm Td}^{\Sigma}(\Theta_{Y}){\rm ch}^{\Sigma}(\Omega_{Y}^{p})
}
(\overline{\gamma},\overline{\gamma}')
\}
\\
&=
-\frac{1}{12}\int_{Y}\widetilde{c_{1}c_{3}}(Y;\overline{\gamma},\overline{\gamma}')
-
\frac{1}{12}\sum_{\bar{\lambda}\in\overline{\Lambda}}\frac{n_{\bar{\lambda}}^{2}-1}{n_{\bar{\lambda}}}
\int_{\overline{C}_{\bar{\lambda}}}\widetilde{c_{1}(Y)c_{1}(\overline{C}_{\bar{\lambda}})}(\overline{\gamma},\overline{\gamma}')
\\
&\quad
+
\frac{1}{12}\sum_{\bar{\lambda}\in\overline{\Lambda}}\mu_{\bar{\lambda}}
\int_{\overline{C}_{\bar{\lambda}}}
\widetilde{c_{1}^{2}}(\overline{C}_{\bar{\lambda}};
\overline{\gamma}|_{\overline{C}_{\bar{\gamma}}},\overline{\gamma}'|_{\overline{C}_{\bar{\gamma}}})
+
\sum_{\overline{\frak p}\in Y^{(0)}}\sum_{k=1}^{3}\frac{\delta(G_{\frak p})}{|G_{\frak p}|}
\widetilde{\rho}_{{\frak p},k}({\gamma},{\gamma}')({\frak p}),
\end{aligned}
\end{equation}
where we used the fact that $\varPi\colon Y\amalg\Sigma Y\to Y\amalg Y^{(1)}\amalg Y^{(0)}$ is an unramified covering
to get the third equality.
\par
Write 
$G_{\frak p}=\{{\rm diag}(\chi_{{\frak p},1}(g),\chi_{{\frak p},2}(g),\chi_{{\frak p},3}(g));\,g\in G_{\frak p}\}$. 
As in Section~\ref{sect:4.5.2}, we set $n_{{\frak p},k}:=|\ker\chi_{{\frak p},k}|$ and $\nu_{{\frak p},k}=|{\rm Im}\,\chi_{{\frak p},k}|$.
Applying the anomaly formula \cite[Th.\,0.1]{Ma05} to $\lambda({\mathcal O}_{\Sigma Y})$
and using Proposition~\ref{prop:orb:Td:ch:fixed:curve:2}, we get
\begin{equation}
\label{eqn:anomaly:formula:fixed:curve}
\begin{aligned}
\,&
\log
\left(
\frac{\|\cdot\|_{\lambda({\mathcal O}_{\Sigma Y}),Q,\overline{\gamma}}}
{\|\cdot\|_{\lambda({\mathcal O}_{\Sigma Y}),Q,\overline{\gamma}'}}
\right)^{2}
=
\sum_{\alpha\in A}
\log
\left(
\frac{\|\cdot\|_{\lambda({\mathcal O}_{D_{\alpha}}),Q,\overline{\gamma}}}
{\|\cdot\|_{\lambda({\mathcal O}_{D_{\alpha}}),Q,\overline{\gamma}'}}
\right)^{2}
\\
&=
\sum_{\alpha\in A}
\int_{D_{\alpha}\amalg\Sigma D_{\alpha}}
\widetilde{
{\rm Td}^{\Sigma}(TD_{\alpha})
}
(\overline{\gamma},\overline{\gamma}')
=
\sum_{\bar{\lambda}\in\overline{\Lambda}}\mu_{\bar\lambda}\,
\int_{\overline{C}_{\bar{\lambda}}\amalg\Sigma\overline{C}_{\bar{\lambda}}}
\widetilde{
{\rm Td}^{\Sigma}(T\overline{C}_{\bar{\lambda}})
}
(\overline{\gamma},\overline{\gamma}')
=
\\
&
\frac{1}{12}\sum_{\bar{\lambda}\in\overline{\Lambda}}\mu_{\bar{\lambda}}
\int_{\overline{C}_{\bar{\lambda}}}
\widetilde{c_{1}^{2}}(\overline{C}_{\bar{\lambda}};
\overline{\gamma}|_{\overline{C}_{\bar{\lambda}}},\overline{\gamma}'|_{\overline{C}_{\bar{\lambda}}})
+
\frac{1}{12}\sum_{\overline{\frak p}\in Y^{(0)}}\sum_{k=1}^{3}\frac{(n_{{\frak p},k}-1)(\nu_{{\frak p},k}^{2}-1)}{\nu_{{\frak p},k}}
\widetilde{\rho}_{{\frak p},k}({\gamma},{\gamma}')({\frak p}),
\end{aligned}
\end{equation}
where we used the fact that $\overline{\nu}\colon\amalg_{\alpha\in\phi^{-1}(\bar{\lambda})}D_{\alpha}\to\widehat{C}_{\bar{\lambda}}$ 
is an unramified covering of degree $\mu_{\bar\lambda}$ to get the third equality.
The result follows from \eqref{eqn:relation:sum:Chern:forms:relative:dim:0}, 
\eqref{eqn:anomaly:formula:BCOV:bundle}, \eqref{eqn:anomaly:formula:fixed:curve} and 
Proposition~\ref{prop:invariant:subgrp:SL:C3}.
\end{pf}

\begin{lemma}
\label{lemma:anomaly:1}
Let $\overline{\gamma}$ and $\overline{\gamma}'$ be K\"ahler forms on $Y$ in the sense of orbifolds. Then 
\begin{equation}
\label{eqn:anomaly:CY:orbifold:Bott:Chern}
\begin{aligned}
\,&
\frac{\chi^{\rm orb}(Y)}{12}\log\frac{{\rm Vol}(Y,\overline{\gamma})}{{\rm Vol}(Y,\overline{\gamma}')}
+
\log\frac{A(Y,\overline{\gamma})}{A(Y,\overline{\gamma}')}
+
\sum_{\bar{\lambda}\in\overline{\Lambda}}
\frac{n_{\bar{\lambda}}^{2}-1}{n_{\bar{\lambda}}}
\log
\frac
{A(\overline{C}_{\bar{\lambda}},\overline{\gamma}|_{\overline{C}_{\bar{\lambda}}})}
{A(\overline{C}_{\bar{\lambda}},\overline{\gamma}'|_{\overline{C}_{\bar{\lambda}}})}
\\
&\quad+
\sum_{\overline{\frak p}\in Y^{(0)}}\epsilon(G_{\frak p})
\log
\frac
{A(\overline{\frak p},\overline{\gamma}|_{\overline{\frak p}})}
{A(\overline{\frak p},\overline{\gamma}'|_{\overline{\frak p}})}
\\
&=
\frac{1}{12}\int_{Y}\widetilde{c_{1}c_{3}}(Y;\overline{\gamma},\overline{\gamma}')
+
\frac{1}{12}\sum_{\bar{\lambda}\in\overline{\Lambda}}\frac{n_{\bar{\lambda}}^{2}-1}{n_{\bar{\lambda}}}
\int_{\overline{C}_{\bar{\lambda}}}
\widetilde{c_{1}(Y)c_{1}(\overline{C}_{\bar{\lambda}})}(\overline{\gamma},\overline{\gamma}')
\\
&\quad-
\sum_{\overline{\frak p}\in Y^{(0)}}\epsilon(G_{\frak p})\,
\widetilde{c_{1}(Y)}(\overline{\gamma},\overline{\gamma}')|_{\overline{\frak p}}.
\end{aligned}
\end{equation}
\end{lemma}

\begin{pf}
Write $I=I(\overline{\gamma},\overline{\gamma}')$ for the left hand side of \eqref{eqn:anomaly:CY:orbifold:Bott:Chern}.
Let $\eta\in H^{0}(Y,\Omega_{Y}^{3})\setminus\{0\}$ and let $\overline{\kappa}$ be a {\em Ricci-flat} K\"ahler form on $Y$ in the sense of orbifolds
such that
\begin{equation}
\label{eqn:Ricci:flatness}
\overline{\kappa}^{3}/3!=i\,\eta\wedge\overline{\eta}.
\end{equation}
By the definitions of $A(Y,\overline{\gamma})$, $A(\overline{C}_{\bar{\lambda}},\overline{\gamma}|_{\overline{C}_{\bar{\lambda}}})$,
$A(\overline{\frak p},\overline{\gamma}|_{\overline{\frak p}})$ and $\chi^{\rm orb}(Y)$, we get
\begin{equation}
\label{eqn:formula:I}
\begin{aligned}
I
&=
\frac{\chi^{\rm orb}(Y)}{12}\log\frac{{\rm Vol}(Y,\overline{\gamma})}{{\rm Vol}(Y,\overline{\gamma}')}
-\frac{\chi^{\rm orb}(Y)}{12}\log\frac{{\rm Vol}(Y,\overline{\gamma})}{\|\eta\|_{L^{2}}^{2}}
+\frac{\chi^{\rm orb}(Y)}{12}\log\frac{{\rm Vol}(Y,\overline{\gamma}')}{\|\eta\|_{L^{2}}^{2}}
\\
&
-\frac{1}{12}\int_{Y}
\left\{
\log\left(\frac{i\,\eta\wedge\overline{\eta}}{\overline{\gamma}^{3}/3!}\right)c_{3}(Y,\overline{\gamma})
-
\log\left(\frac{i\,\eta\wedge\overline{\eta}}{(\overline{\gamma}')^{3}/3!}\right)c_{3}(Y,\overline{\gamma}')
\right\}
\\
&
-\frac{1}{12}\sum_{\bar{\lambda}\in\overline{\Lambda}}\frac{n_{\bar{\lambda}}^{2}-1}{n_{\bar{\lambda}}}
\int_{\overline{C}_{\bar{\lambda}}}
\left\{
\log\left(\frac{i\,\eta\wedge\overline{\eta}}{\overline{\gamma}^{3}/3!}\right)
c_{1}(\overline{C}_{\bar{\lambda}},\overline{\gamma}|_{\overline{C}_{\bar{\lambda}}})
-
\log\left(\frac{i\,\eta\wedge\overline{\eta}}{(\overline{\gamma}')^{3}/3!}\right)
c_{1}(\overline{C}_{\bar{\lambda}},\overline{\gamma}'|_{\overline{C}_{\bar{\lambda}}})
\right\}
\\
&
+\sum_{\overline{\frak p}\in Y^{(0)}}
\epsilon(G_{\frak p})
\left\{
\log\left(\frac{i\,\eta\wedge\overline{\eta}}{\overline{\gamma}^{3}/3!}\right)(\overline{\frak p})
-
\log\left(\frac{i\,\eta\wedge\overline{\eta}}{(\overline{\gamma}')^{3}/3!}\right)(\overline{\frak p})
\right\}.
\end{aligned}
\end{equation}
\par
Since $c_{1}(Y,\overline{\kappa})=0$ by the definition of $\kappa$ and since
\begin{equation}
\label{eqn:Bott:Chern:c1}
\widetilde{c}_{1}(Y;\overline{\gamma},\overline{\kappa})
=
\log(\overline{\gamma}^{3}/\overline{\kappa}^{3})
=
\log[(\overline{\gamma}^{3}/3!)/(i\,\eta\wedge\overline{\eta})]
\end{equation}
and 
$\widetilde{c_{1}c_{3}}(Y;\overline{\gamma},\overline{\kappa})
=
\widetilde{c}_{1}(Y;\overline{\gamma},\overline{\kappa})c_{3}(Y,\overline{\gamma})
+
c_{1}(Y,\overline{\kappa})\widetilde{c}_{3}(Y;\overline{\gamma},\overline{\kappa})$,
we deduce from \eqref{eqn:Ricci:flatness}, \eqref{eqn:Bott:Chern:c1} that
\begin{equation}
\label{eqn:Bott:Chern:c1:c3}
\begin{aligned}
\,&
\widetilde{c_{1}c_{3}}(Y;\overline{\gamma},\overline{\gamma}')
=
\widetilde{c}_{1}(Y;\overline{\gamma},\overline{\kappa})c_{3}(Y,\overline{\gamma})
-
\widetilde{c}_{1}(Y;\overline{\gamma}',\overline{\kappa})c_{3}(Y,\overline{\gamma}')
\\
&=
\log[(\overline{\gamma}^{3}/3!)/(i\,\eta\wedge\overline{\eta})]\,c_{3}(Y,\overline{\gamma})
-
\log[((\overline{\gamma}')^{3}/3!)/(i\,\eta\wedge\overline{\eta})]\,c_{3}(Y,\overline{\gamma}').
\end{aligned}
\end{equation}
Similarly, since $c_{1}(Y,\overline{\kappa})=0$ and hence
$$
\begin{aligned}
\widetilde{c_{1}(Y)c_{1}(\overline{C}_{\bar{\lambda}})}(\overline{\gamma},\overline{\gamma}')
&=
\widetilde{c_{1}(Y)c_{1}(\overline{C}_{\bar{\lambda}})}(\overline{\gamma},\overline{\kappa})
-
\widetilde{c_{1}(Y)c_{1}(\overline{C}_{\bar{\lambda}})}(\overline{\gamma}',\overline{\kappa})
\\
&=
\widetilde{c}_{1}(Y;\overline{\gamma},\overline{\kappa})c_{1}(\overline{C}_{\bar{\lambda}},\overline{\gamma}|_{\overline{C}_{\bar{\lambda}}})
-
\widetilde{c}_{1}(Y;\overline{\gamma}',\overline{\kappa})c_{1}(\overline{C}_{\bar{\lambda}},\overline{\gamma}'|_{\overline{C}_{\bar{\lambda}}}),
\end{aligned}
$$
we get by \eqref{eqn:Bott:Chern:c1}
\begin{equation}
\label{eqn:Bott:Chern:c1:c1}
\widetilde{c_{1}(Y)c_{1}(\overline{C}_{\bar{\lambda}})}(\overline{\gamma},\overline{\gamma}')
=
\log\left(\frac{\overline{\gamma}^{3}/3!}{i\,\eta\wedge\overline{\eta}}\right)
c_{1}(\overline{C}_{\bar{\lambda}},\overline{\gamma}|_{\overline{C}_{\bar{\lambda}}})
-
\log\left(\frac{(\overline{\gamma}')^{3}/3!}{i\,\eta\wedge\overline{\eta}}\right)\,
c_{1}(\overline{C}_{\bar{\lambda}},\overline{\gamma}'|_{\overline{C}_{\bar{\lambda}}}).
\end{equation}
Finally, we get by \eqref{eqn:Bott:Chern:c1}
\begin{equation}
\label{eqn:Bott:Chern:c1:fixed:pt}
\begin{aligned}
\widetilde{c_{1}(Y)}(\overline{\gamma},\overline{\gamma}')|_{\overline{\frak p}}
&=
\widetilde{c_{1}(Y)}(\overline{\gamma},\overline{\kappa})|_{\overline{\frak p}}
-
\widetilde{c_{1}(Y)}(\overline{\gamma}',\overline{\kappa})|_{\overline{\frak p}}
\\
&=
\log[(\overline{\gamma}^{3}/3!)/(i\,\eta\wedge\overline{\eta})](\overline{\frak p})
-
\log[\{(\overline{\gamma}')^{3}/3!\}/(i\,\eta\wedge\overline{\eta})](\overline{\frak p}).
\end{aligned}
\end{equation}
Substituting \eqref{eqn:Bott:Chern:c1:c3}, \eqref{eqn:Bott:Chern:c1:c1}, \eqref{eqn:Bott:Chern:c1:fixed:pt} into \eqref{eqn:formula:I},
we get the result.
\end{pf}

\begin{theorem}
\label{thm:invariance:BCOV:inv:CY:orb}
The number $\tau_{\rm BCOV}^{\rm orb}(Y,\overline{\gamma})$ is independent of the choice of a K\"ahler form on $Y$
in the sense of orbifolds.
\end{theorem}

\begin{pf}
Let $\overline{\gamma}$ and $\overline{\gamma}'$ be arbitrary K\"ahler forms on $Y$ in the sense of orbifolds. 
Let $\{c_{1}(L_{1}),\ldots,c_{1}(L_{b_{2}})\}$ be a basis of $H^{2}(Y,{\bf Z})_{\rm fr}$ and let $\{c_{1}(L_{1})^{\lor},\ldots,c_{1}(L_{b_{2}})^{\lor}\}$
be its dual basis of $H^{4}(Y,{\bf Z})_{\rm fr}$ with respect to the cup-product. Then
$\det H^{0}(Y,{\bf C})$ has the canonical element $\sigma_{0}:=1_{Y}$,
$\det H^{2}(Y,{\bf C})$ has the canonical element $\sigma_{2}:=c_{1}(L_{1})\wedge\cdots\wedge c_{1}(L_{b_{2}})$,
$\det H^{4}(Y,{\bf C})$ has the canonical element $\sigma_{4}:=c_{1}(L_{1})^{\lor}\wedge\cdots\wedge c_{1}(L_{b_{2}})^{\lor}$,
and $\det H^{6}(Y,{\bf C})$ has the canonical element $\sigma_{6}:=[Y]$, the fundamental class.
If $\{e_{1},\ldots,e_{b_{2}}\}$ is an orthonormal basis of $H^{2}(Y,{\bf C})$ with respect to $\overline{\gamma}$, then its dual basis
$\{e_{1}^{\lor},\ldots,e_{b_{2}}^{\lor}\}$ is given by $\{*e_{1},\ldots,*e_{b_{2}}\}$, where $*$ denotes the Hodge $*$-operator
with respect to the K\"ahler form $\overline{\gamma}$.
From this, we get 
$\|\sigma_{4}\|_{L^{2},\overline{\gamma}}^{2}=\|\sigma_{2}\|_{L^{2},\overline{\gamma}}^{-2}={\rm Vol}_{L^{2}}(H^{2}(Y,{\bf Z}),[\overline{\gamma}])^{-1}$.
Similarly, we have
$\|\sigma_{6}\|_{L^{2},\overline{\gamma}}^{2}=\|\sigma_{0}\|_{L^{2},\overline{\gamma}}^{-2}={\rm Vol}(Y,\overline{\gamma})^{-1}$.
Set
$$
\ell_{Y}:=(\det H^{2}(Y,\Omega_{Y}^{1}))^{-1}\otimes(\det H^{1}(Y,\Omega_{Y}^{2}))^{-2}\otimes(\det H^{0}(Y,\Omega_{Y}^{3}))^{-3}.
$$
Since the $L^{2}$-metric on $H^{q}(Y,\Omega_{Y}^{p})$ is independent of the choice of a K\"ahler form on $Y$ for $p+q=3$ and hence
so is the $L^{2}$-metric $\|\cdot\|_{L^{2},\ell_{Y}}$ on the line $\ell_{Y}$, we get 
\begin{equation}
\label{eqn:ratio:Quillen:metric:1}
\begin{aligned}
\frac{\|\cdot\|_{\lambda_{\rm BCOV}(Y),Q,\overline{\gamma}}^{2}}
{\|\cdot\|_{\lambda_{\rm BCOV}(Y),Q,\overline{\gamma}'}^{2}}
&=
\frac
{T_{\rm BCOV}(Y,\overline{\gamma})\cdot\|\cdot\|_{L^{2},\ell_{Y}}^{2}\cdot\prod_{r=1}^{3}\|\sigma_{r}\|_{L^{2},\overline{\gamma}}^{2r}}
{T_{\rm BCOV}(Y,\overline{\gamma}')\cdot\|\cdot\|_{L^{2},\ell_{Y}}^{2}\cdot\prod_{r=1}^{3}\|\sigma_{r}\|_{L^{2},\overline{\gamma}'}^{2r}}
\\
&=
\frac
{T_{\rm BCOV}(Y,\overline{\gamma}){\rm Vol}(Y,\overline{\gamma})^{-3}{\rm Vol}_{L^{2}}(H^{2}(Y,{\bf Z}),[\overline{\gamma}])^{-1}}
{T_{\rm BCOV}(Y,\overline{\gamma}'){\rm Vol}(Y,\overline{\gamma}')^{-3}{\rm Vol}_{L^{2}}(H^{2}(Y,{\bf Z}),[\overline{\gamma}'])^{-1}}
\end{aligned}
\end{equation}
by \eqref{eqn:BCOV:bundle}.
Similarly, since the $L^{2}$-metric on the line
$$
\ell_{\Sigma Y}:=\bigotimes_{\alpha\in A}\ell_{D_{\alpha}},
\qquad
\ell_{D_{\alpha}}:=\det H^{1}(D_{\alpha},{\mathcal O}_{D_{\alpha}})^{-1}
$$
is independent of the choice of a K\"ahler form on $\Sigma Y$ and since 
$\det H^{0}(\Sigma Y,{\mathcal O}_{\Sigma Y})=\bigotimes_{\alpha\in A}\det H^{0}(D_{\alpha},{\mathcal O}_{D_{\alpha}})$ 
has the canonical element $1_{\Sigma Y}=\bigotimes_{\alpha\in A}1_{D_{\alpha}}$, we get
\begin{equation}
\label{eqn:ratio:Quillen:metric:2}
\frac{\|\cdot\|_{\lambda({\mathcal O}_{\Sigma Y}),Q,\overline{\gamma}}^{2}}
{\|\cdot\|_{\lambda({\mathcal O}_{\Sigma Y}),Q,\overline{\gamma}'}^{2}}
=
\frac
{\tau(\Sigma Y,\overline{\gamma}|_{\Sigma Y})
\|\cdot\|_{L^{2},\ell_{\Sigma Y}}^{2}
\|1_{\Sigma Y}\|_{L^{2},\gamma}^{2}
}
{\tau(\Sigma Y,\overline{\gamma}'|_{\Sigma Y})
\|\cdot\|_{L^{2},\ell_{\Sigma Y}}^{2}
\|1_{\Sigma Y}\|_{L^{2},\gamma'}^{2}
}
=
\frac
{\tau(\Sigma Y,\overline{\gamma}|_{\Sigma Y})
{\rm Vol}(\Sigma Y,\overline{\gamma}|_{\Sigma Y})}
{\tau(\Sigma Y,\overline{\gamma}'|_{\Sigma Y})
{\rm Vol}(\Sigma Y,\overline{\gamma}'|_{\Sigma Y})}.
\end{equation}
By \eqref{eqn:ratio:Quillen:metric:1}, \eqref{eqn:ratio:Quillen:metric:2} and Definition~\ref{def:orbifold:BCOV:invariant:1},
we get
\begin{equation}
\label{eqn:anomaly:formula:BCOV:inv}
\begin{aligned}
\log\frac{\tau_{\rm BCOV}^{\rm orb}(Y,\overline{\gamma})}{\tau_{\rm BCOV}^{\rm orb}(Y,\overline{\gamma}')}
&=
\log
\left(
\frac{\|\cdot\|_{\lambda_{\rm BCOV}(Y),Q,\overline{\gamma}}}
{\|\cdot\|_{\lambda_{\rm BCOV}(Y),Q,\overline{\gamma}'}}
\right)^{2}
-
\log
\left(
\frac{\|\cdot\|_{\lambda({\mathcal O}_{\Sigma Y}),Q,\overline{\gamma}}}
{\|\cdot\|_{\lambda({\mathcal O}_{\Sigma Y}),Q,\overline{\gamma}'}}
\right)^{2}
\\
&\quad
+
\frac{\chi^{\rm orb}(Y)}{12}\log\frac{{\rm Vol}(Y,\overline{\gamma})}{{\rm Vol}(Y,\overline{\gamma}')}
+
\log\frac{A(Y,\overline{\gamma})}{A(Y,\overline{\gamma}')}
\\
&\quad
+
\sum_{\bar{\lambda}\in\overline{\Lambda}}
\frac{n_{\bar{\lambda}}^{2}-1}{n_{\bar{\lambda}}}
\log
\frac
{A(\overline{C}_{\bar{\lambda}},\overline{\gamma}|_{\overline{C}_{\bar{\lambda}}})}
{A(\overline{C}_{\bar{\lambda}},\overline{\gamma}'|_{\overline{C}_{\bar{\lambda}}})}
+
\sum_{\overline{\frak p}\in Y^{(0)}}\epsilon(G_{\frak p})
\log
\frac
{A(\overline{\frak p},\overline{\gamma}|_{\overline{\frak p}})}
{A(\overline{\frak p},\overline{\gamma}'|_{\overline{\frak p}})}.
\end{aligned}
\end{equation}
Comparing \eqref{eqn:anomaly:formula:BCOV:inv} and Lemmas~\ref{lemma:anomaly:formula:Ma} and \ref{lemma:anomaly:1},
we get 
$$
\log[\tau_{\rm BCOV}^{\rm orb}(Y,\overline{\gamma})/\tau_{\rm BCOV}^{\rm orb}(Y,\overline{\gamma}')]=0.
$$
Since the K\"ahler forms $\overline{\gamma}$ and $\overline{\gamma}'$ are arbitrary, we get the result.
\end{pf}

After Theorem~\ref{thm:invariance:BCOV:inv:CY:orb}, the following definition makes sense.

\begin{definition}
\label{def:BCOV:inv:CY:orbifold:general}
Define the BCOV invariant of an abelian Calabi-Yau orbifold as
$$
\tau_{\rm BCOV}^{\rm orb}(Y):=\tau_{\rm BCOV}^{\rm orb}(Y,\overline{\gamma}).
$$
\end{definition}

To reformulate the invariance property of $\tau_{\rm BCOV}^{\rm orb}(Y)$, we make the following

\begin{definition}
\label{def:orb:BCOV:line}
Define the orbifold BCOV line of $Y$ as
$$
\lambda_{\rm BCOV}^{\rm orb}(Y)
:=
\lambda_{\rm BCOV}(Y)
\otimes
\lambda({\mathcal O}_{\Sigma Y})^{-1}
$$
and the Hermitian structure $\|\cdot\|_{\rm BCOV}$ on $\lambda_{\rm BCOV}^{\rm orb}(Y)$ as
$$
\begin{aligned}
\|\cdot\|_{\rm BCOV}^{2}
&:=
{\rm Vol}(Y,\overline{\gamma})^{\frac{\chi^{\rm orb}(Y)}{12}}A(Y,\overline{\gamma})
\prod_{\bar{\lambda}\in\overline{\Lambda}}
A(\overline{C}_{\lambda},\overline{\gamma}|_{\overline{C}_{\bar{\lambda}}})^{\frac{n_{\bar{\lambda}}^{2}-1}{n_{\bar{\lambda}}}}
\prod_{\overline{\frak p}\in Y^{(0)}}
A(\overline{\frak p},\overline{\gamma}|_{\overline{\frak p}})^{\epsilon(G_{\frak p})}
\\
&\qquad\times
\|\cdot\|_{\lambda_{\rm BCOV}(Y),Q,\overline{\gamma}}^{2}
\otimes
\|\cdot\|_{\lambda({\mathcal O}_{\Sigma Y}),Q,\overline{\gamma}}^{-2}.
\end{aligned}
$$
\end{definition}

\begin{theorem}
\label{thm:orb:BCOV:line}
The Hermitian structure $\|\cdot\|_{\rm BCOV}$ on $\lambda_{\rm BCOV}^{\rm orb}(Y)$ is independent of the choice of
a K\"ahler metric on $Y$ in the sense of orbifolds. In particular, the Hermitian line
$(\lambda_{\rm BCOV}^{\rm orb}(Y),\|\cdot\|_{\rm BCOV})$ is an invariant of $Y$.
\end{theorem}

\begin{pf}
The result follows from \eqref{eqn:anomaly:formula:BCOV:inv} and Theorem~\ref{thm:invariance:BCOV:inv:CY:orb}.
\end{pf}

By \eqref{eqn:ratio:Quillen:metric:1}, \eqref{eqn:ratio:Quillen:metric:2}, $\tau_{\rm BCOV}^{\rm orb}(Y)$ is obtained as the ratio of
two intrinsic Hermitian structures on the complex line $\ell_{Y}\otimes\ell_{\Sigma Y}^{-1}$ as follows:
\begin{equation}
\label{eqn:orb:BCOV:inv:ratio:metric}
\tau_{\rm BCOV}^{\rm orb}(Y)
=
\frac
{\|(\,\,\cdot\,\,)\otimes\bigotimes_{r=0}^{3}\sigma_{r}^{\otimes r}
\otimes
1_{\Sigma Y}^{\otimes-1}
\|_{\rm BCOV}^{2}}
{\|\cdot\|_{L^{2},\ell_{Y}\otimes\ell_{\Sigma Y}^{-1}}^{2}}.
\end{equation}
\par
A typical example of abelian Calabi-Yau orbifold is a general anti-canonical divisor of a toric variety of certain type.
See e.g. \cite[Prop.\,4.13 (i)]{CoxKatz99}. To extend BCOV invariants to non-abelian Calabi-Yau orbifolds, one has to deal with
three-dimensional quotient singularities obtained by finite non-abelian subgroups of ${\rm SL}({\bf C}^{3})$. For the classification of such groups,
see e.g. \cite[Sect.\,3]{Roan96}.

\section
{BCOV invariants for Borcea-Voisin orbifolds}
\label{sect:8}
\par

\subsection
{Analytic torsion for elliptic curves}
\label{sect:8.1}
\par
Let $T$ be an elliptic curve and let $\gamma$ be a K\"ahler form on $T$. 
Let $\xi$ be a non-zero holomorphic $1$-form on $T$. We set
$$
\tau_{\rm ell}(T)
:=
{\rm Vol}(T,\gamma)\,\tau(T,\gamma)\,
\exp\left[\frac{1}{12}\int_{T}
\log\left(
\frac{i\,\xi\wedge\overline{\xi}}{\gamma}
\right)\,c_{1}(T,\gamma)\right].
$$
Since $\chi(T)=\int_{T}c_{1}(T,\gamma)=0$, $\tau(T)_{\rm ell}$ is independent of the choice of $\xi$.
\par
When $T\cong{\bf C}/{\bf Z}+\tau{\bf Z}$, $\tau\in{\frak H}$, we define $g_{T}:=dz\otimes d\bar{z}/\Im\tau$.
Then its K\"ahler form is given by $\gamma_{T}:=\frac{i}{2}dz\wedge d\bar{z}/\Im\tau$.
Since ${\rm Vol}(T,\gamma_{T})=(2\pi)^{-1}\int_{T}\gamma_{T}=(2\pi)^{-1}$, we have
\begin{equation}
\label{eqn:torsion:flat:elliptic:curve}
\tau_{\rm ell}(T)=(2\pi)^{-1}\tau(T,\gamma_{T}).
\end{equation}
The Laplacian $\square^{T}_{0,0}$ of $(T,g_{T})$ acting on $C^{\infty}(T)$ is given by
$$
\square_{0,0}^{T}
=
-2\Im\tau\frac{\partial^{2}}{\partial z\partial\bar{z}}
=
-\frac{\Im\tau}{2}
\left(
\frac{\partial^{2}}{\partial x^{2}}+\frac{\partial^{2}}{\partial y^{2}}
\right).
$$
We have $\sigma(\square_{0,0}^{T})=\{\nu_{m,n}(T);\,(m,n)\in{\bf Z}^{2}\}$, where $\nu_{m,n}(T):=\frac{2\pi^{2}}{\Im\tau}|m\tau+n|^{2}$
for $(m,n)\in{\bf Z}^{2}$.

\begin{theorem}
\label{thm:Kronecker:limit:formula}
For every elliptic curve $T$,
$$
\tau_{\rm ell}(T)
=
\left(4\pi\,\|\eta\left(\varOmega(T)\right)^{4}\|\right)^{-1}.
$$
\end{theorem}

\begin{pf}
Let $\zeta_{0,0}^{T}(s)$ be the spectral zeta function of $\square_{0,0}^{T}$.
By the Kronecker limit formula \cite[p.75 Eq.\,(17)]{Weil99}, we have
$$
\left.\frac{d}{ds}\right|_{s=0}\zeta_{0,0}^{T}(s)
=
\left.\frac{d}{ds}\right|_{s=0}\sum_{(m,n)\in{\bf Z}^{2}\setminus\{0\}}\frac{(\Im\tau)^{s}/(2\pi^{2})^{s}}{|m\tau+n|^{2s}}
=
-\log\left(
2\|\eta(\tau)\|^{4}
\right).
$$
This, together with \eqref{eqn:torsion:flat:elliptic:curve} and the anomaly formula \cite[Th.\,0.2]{BGS88}, yields the result.
\end{pf}

\subsection
{An observation of Harvey-Moore for Borcea-Voisin orbifolds}
\label{sect:8.2}
\par
In the rest of this paper, we prove the following:

\begin{theorem}
\label{thm:Harvey:Moore}
Let $(S,\theta)$ be a $2$-elementary $K3$ surface of type $M$ and let $T$ be an elliptic curve. 
Let $\rho=\rho(X_{(S,\theta,T)})=r(M)+1$ be the Picard number of $X_{(S,\theta,T)}$.
Then the following equality holds:
$$
\tau^{\rm orb}_{\rm BCOV}\left(X_{(S,\theta,T)}\right)
=
2^{14}(2\pi)^{2\rho}|A_{M}|^{-1}\,\tau_{M}(S,\theta)^{-4}\tau_{\rm ell}(T)^{-12}.
$$
\end{theorem}

Theorem~\ref{thm:Harvey:Moore} was proved by Harvey-Moore \cite[Sect.\,V]{HarveyMoore98} when $\theta\colon S\to S$ is free from fixed points, 
equivalently when the type of $X_{(S,\theta,T)}$ is ${\Bbb U}(2)\oplus{\Bbb E}_{8}(2)$.

\subsubsection
{Symmetries of the spectrum of Laplacians}
\label{sect:8.2.1}
\par
Let $\gamma_{T}$ and $\gamma_{S}$ be {\em Ricci-flat} K\"ahler forms on $T={\bf C}/{\bf Z}+\tau{\bf Z}$ and $S$, respectively, such that 
$$
\gamma_{T}:=\frac{i\,dz\wedge d\bar{z}}{2\Im\tau},
\qquad
\theta^{*}\gamma_{S}=\gamma_{S}.
$$
We set 
$$
\gamma:={\rm pr}_{1}^{*}\gamma_{S}+{\rm pr}_{2}^{*}\gamma_{T}.
$$
Then $\gamma$ is a $\theta\times(-1)_{T}$-invariant Ricci-flat K\"ahler form on $S\times T$ and is regarded as a Ricci-flat K\"ahler form on 
$X_{(S,\theta,T)}=(S\times T)/\theta\times(-1)_{T}$ in the sense of orbifolds. 
Let $\eta$ and $\xi$ be nowhere vanishing canonical forms on $S$ and $T$, respectively.
Then $\eta\wedge\xi$ is a nowhere vanishing canonical form on $S\times T$, which we identify the corresponding nowhere vanishing canonical form
on $X_{(S,\theta,T)}$. By the Ricci-flatness of $\gamma_{S}$ and $\gamma_{T}$, the canonical forms $\eta$, $\xi$, $\eta\wedge\xi$ are parallel
with respect to the Levi-Civita connections associated to the K\"ahler forms $\gamma_{S}$, $\gamma_{T}$, $\gamma$, respectively.
For simplicity, write $X$ for $X_{(S,\theta,T)}$ when there is no possibility of confusion.
\par
Let $\zeta_{p,q}(s)$ be the spectral zeta function of $\square_{p,q}$ acting on $A^{p,q}(X)$, the space of smooth $(p,q)$-forms on $X$ 
in the sense of orbifolds.

\begin{lemma}
\label{lemma:relation:spectral:zeta:1}
The following equality of meromorphic functions on ${\bf C}$ holds
$$
\sum_{p,q\geq0}(-1)^{p+q}pq\,\zeta_{p,q}(s)=9\,\zeta_{0,0}(s)-6\,\zeta_{1,0}(s)+\zeta_{1,1}(s).
$$
\end{lemma}

\begin{pf}
By the ellipticity of the Dolbeault complex $(A^{p,\bullet}(X),\bar{\partial})$, we have the following equality of spectral zeta functions for all $p\geq0$
\begin{equation}
\label{eqn:Dolbeault:relation}
\zeta_{p,0}(s)-\zeta_{p,1}(s)+\zeta_{p,2}(s)-\zeta_{p,3}(s)=0.
\end{equation}
Since $\eta\wedge\xi$ is parallel and nowhere vanishing, the Laplacians $\square_{p,0}$ and $\square_{p,3}$ are isospectral 
by the map $\otimes(\eta\wedge\xi)\colon A^{p,0}(X)\ni\varphi\to\varphi\cdot(\eta\wedge\xi)\in A^{p,3}(X)$. 
In particular, $\zeta_{p,0}(s)=\zeta_{p,3}(s)$, which, together with \eqref{eqn:Dolbeault:relation}, implies $\zeta_{p,1}(s)=\zeta_{p,2}(s)$. Hence
\begin{equation}
\label{eqn:symmetry:spectral:zeta:1}
\zeta_{p,q}(s)=\zeta_{p,3-q}(s)
\end{equation}
for all $0\leq p,q\leq3$. Since $\overline{\square_{p,q}\varphi}=\square_{q,p}\overline{\varphi}$ for all $\varphi\in A^{p,q}(X)$, 
we get for all $p,q\geq0$
\begin{equation}
\label{eqn:symmetry:spectral:zeta:2}
\zeta_{p,q}(s)=\zeta_{q,p}(s).
\end{equation}
By \eqref{eqn:symmetry:spectral:zeta:1}, \eqref{eqn:symmetry:spectral:zeta:2}, every $\zeta_{p,q}(s)$ is equal to one of
$\zeta_{0,0}(s)$, $\zeta_{1,0}(s)$, $\zeta_{1,1}(s)$. Replacing $\zeta_{p,q}(s)$ by the corresponding $\zeta_{0,0}(s)$ or $\zeta_{1,0}(s)$ or $\zeta_{1,1}(s)$,
we get the result.
\end{pf}

\subsubsection
{The spectrum of various Laplacians}
\label{sect:8.2.2}
\par
Let $\square^{S}_{p,q}$ (resp. $\square^{T}_{p,q}$) be the Laplacian acting on $A^{p,q}(S)$ (resp. $A^{p,q}(T)$) 
with respect to $\gamma_{S}$ (resp. $\gamma_{T}$).
Set 
$$
A^{p,q}(S)^{\pm}:=\{\varphi\in A^{p,q}(S);\,\theta^{*}\varphi=\pm\varphi\},
\qquad
\square^{S,\pm}_{p,q}:=\square_{p,q}^{S}|_{A^{p,q}(S)^{\pm}},
$$
$$
A^{p,q}(T)^{\pm}:=\{\varphi\in A^{p,q}(T);\,(-1_{T})^{*}\varphi=\pm\varphi\},
\qquad
\square^{T,\pm}_{p,q}:=\square_{p,q}^{T}|_{A^{p,q}(T)^{\pm}}.
$$
For simplicity, write $\square^{S,\pm}=\square^{S,\pm}_{0,0}$ and $\square^{T,\pm}=\square^{T,\pm}_{0,0}$.
Set $h^{1,1}(S)^{\pm}:=\dim\ker\square_{1,1}^{S,\pm}$. 
Let $\sigma(\square^{S,\pm}_{p,q})$ (resp. $\sigma(\square^{T,\pm}_{p,q})$)
be the set of eigenvalues of $\square^{S,\pm}_{p,q}$ (resp. $\square^{T,\pm}_{p,q}$) counted with multiplicities. 
We can write
$$
\sigma(\square^{S,+})=\sigma(\square_{0,0}^{S,+})=\{0\}\amalg\{\lambda_{i}^{+}(S)\}_{i\in I},
\qquad
\sigma(\square^{S,-})=\sigma(\square_{0,0}^{S,-})=\{\lambda_{j}^{-}(S)\}_{j\in J},
$$
where $\lambda_{i}^{+}(S)>0$ and $\lambda_{j}^{-}(S)>0$ for all $i\in I$ and $j\in J$.
In what follows, for a subset ${\mathcal S}\subset{\bf R}$, the notation $\nu\cdot{\mathcal S}$ implies that every element of ${\mathcal S}$ has multiplicity $\nu$.

\begin{lemma}
\label{lemma:spectrum:Laplacian:2-elementary:K3}
The set $\sigma(\square_{p,q}^{S,\pm})$ is given as follows:
$$
\begin{array}{ll}
(1)&
\sigma(\square_{1,0}^{S,+})=\sigma(\square_{1,0}^{S,-})=\sigma(\square_{0,1}^{S,+})=\sigma(\square_{0,1}^{S,-})
=
\{\lambda_{i}^{+}(S)\}_{i\in I}\amalg\{\lambda_{j}^{-}(S)\}_{j\in J}.
\\
(2)&
\sigma(\square_{1,1}^{S,\pm})
=
h^{1,1}(S)^{\pm}\cdot\{0\}\amalg 2\cdot\{\lambda_{i}^{+}(S)\}_{i\in I}\amalg 2\cdot\{\lambda_{j}^{-}(S)\}_{j\in J}.
\end{array}
$$
\end{lemma}

\begin{pf}
Since the Dolbeault complex $(A^{p,\bullet}(S),\bar{\partial})$ is elliptic and since 
\begin{equation}
\label{eqn:spectrum:(p,0):+=(p,2):-}
\sigma(\square_{p,0}^{S,\pm})=\sigma(\square_{p,2}^{S,\mp})
\end{equation}
via the map $\otimes\eta\colon A^{p,0}(S)^{\pm}\ni\varphi\to\varphi\cdot\overline{\eta}\in A^{p,2}(S)^{\mp}$, we get
$$
\begin{aligned}
\sigma(\square_{0,1}^{S,\pm})
&=
(\sigma(\square_{0,0}^{S,\pm})\setminus\{0\})\amalg(\sigma(\square_{0,2}^{S,\pm})\setminus\{0\})
=
(\sigma(\square_{0,0}^{S,\pm})\setminus\{0\})\amalg(\sigma(\square_{0,0}^{S,\mp})\setminus\{0\})
\\
&=
\{\lambda_{i}^{+}(S)\}_{i\in I}\amalg\{\lambda_{j}^{-}(S)\}_{j\in J}.
\end{aligned}
$$
By taking the complex conjugation, we get
$\sigma(\square_{1,0}^{S,\pm})=\{\lambda_{i}^{+}(S)\}_{i\in I}\amalg\{\lambda_{j}^{-}(S)\}_{j\in J}$.
This proves (1). By the ellipticity of the Dolbeault complex $(A^{1,\bullet}(S),\bar{\partial})$ and \eqref{eqn:spectrum:(p,0):+=(p,2):-},
$$
\begin{aligned}
\sigma(\square_{1,1}^{S,\pm})
&=
h^{1,1}(S)^{\pm}\cdot\{0\}\amalg \sigma(\square_{1,0}^{S,\pm})\amalg \sigma(\square_{1,2}^{S,\pm})
=
h^{1,1}(S)^{\pm}\cdot\{0\}\amalg \sigma(\square_{1,0}^{S,\pm})\amalg \sigma(\square_{1,0}^{S,\mp})
\\
&=
h^{1,1}(S)^{\pm}\cdot\{0\}\amalg 2\cdot\{\lambda_{i}^{+}(S)\}_{i\in I}\amalg 2\cdot\{\lambda_{j}^{-}(S)\}_{j\in J},
\end{aligned}
$$
where the last equality follows from (1). This proves (2).
\end{pf}

For $T={\bf C}/{\bf Z}+\tau{\bf Z}$ and $(m,n)\in{\bf Z}^{2}$, recall that 
$\nu_{m,n}(T)=(2\pi^{2}/\Im\tau)|m\tau+n|^{2}$. Then $\sigma(\square_{0,0}^{T})=\{\nu_{m,n}(T);\,(m,n)\in{\bf Z}^{2}\}$. 
The following is classical (cf. \cite[p.\,166]{RaySinger73}).

\begin{lemma}
\label{lemma:spectrum:Laplacian:elliptic:curve}
Set $({\bf Z}^{2})^{*}:={\bf Z}^{2}\setminus\{(0,0)\}$. Then
$$
\sigma(\square^{T,+})=\sigma(\square_{1,0}^{T,-})=\sigma(\square_{0,1}^{T,-})=\sigma(\square_{1,1}^{T,+})
=
\{0\}\amalg\left\{\nu_{m,n}(T);\,(m,n)\in({\bf Z}^{2})^{*}/\pm1\right\},
$$
$$
\sigma(\square^{T,-})=\sigma(\square_{1,0}^{T,+})=\sigma(\square_{0,1}^{T,+})=\sigma(\square_{1,1}^{T,-})
=
\left\{\nu_{m,n}(T);\,(m,n)\in({\bf Z}^{2})^{*}/\pm1\right\}.
$$
\end{lemma}

\subsubsection
{Various spectral $\zeta$-functions}
\label{sect:8.2.3}
\par
Let $\zeta^{S,\pm}_{p,q}(s)$ and $\zeta^{T,\pm}_{p,q}(s)$ be the spectral zeta functions of $\square^{S,\pm}_{p,q}$ and $\square^{T,\pm}_{p,q}$, respectively.
Then
$$
\zeta^{S,+}(s)=\sum_{i\in I}\lambda_{i}^{+}(S)^{-s},
\qquad
\zeta^{S,-}(s)=\sum_{j\in J}\lambda_{j}^{-}(S)^{-s},
$$
$$
\zeta^{T,+}(s)=\zeta^{T,-}(s)
=
\sum_{(m,n)\in({\bf Z}^{2})^{*}/\pm1}\frac{(\Im\tau)^{s}(2\pi^{2})^{-s}}{|m\tau+n|^{2s}}
=
\sum_{(m,n)\in({\bf Z}^{2})^{*}/\pm1}\nu_{m,n}(T)^{-s}.
$$
Following Harvey-Moore \cite{HarveyMoore98}, we set
$$
\mu^{+}(s)
:=
\sum_{i\in I,\,(m,n)\in({\bf Z}^{2})^{*}/\pm1}
\left(
\lambda_{i}^{+}(S)+\nu_{m,n}(T)
\right)^{-s},
$$
$$
\mu^{-}(s)
:=
\sum_{j\in J,\,(m,n)\in({\bf Z}^{2})^{*}/\pm1}
\left(
\lambda_{j}^{-}(S)+\nu_{m,n}(T)
\right)^{-s}.
$$

\begin{lemma}
\label{lemma:spectral:zeta:Borcea:Voisin:(0,0)}
The following equality of meromorphic functions on ${\bf C}$ holds:
$$
\zeta_{0,0}(s)=\zeta^{T,+}(s)+\zeta^{S,+}(s)+\mu^{+}(s)+\mu^{-}(s).
$$
\end{lemma}

\begin{pf}
Since the decomposition
$$
A^{0,0}(X)^{+}=[A^{0,0}(T)^{+}\hat{\otimes}A^{0,0}(S)^{+}]\oplus[A^{0,0}(T)^{-}\hat{\otimes}A^{0,0}(S)^{-}]
$$
is orthogonal and is preserved by $\square_{0,0}$,
we get 
\begin{equation}
\label{eqn:decomposition:spectrum:Laplacian:(0,0)}
\sigma(\square_{0,0})
=
\sigma(\square_{0,0}|_{A^{0,0}(T)^{+}\hat{\otimes}A^{0,0}(S)^{+}})
\amalg
\sigma(\square_{0,0}|_{A^{0,0}(T)^{-}\hat{\otimes}A^{0,0}(S)^{-}}).
\end{equation}
Since 
\begin{equation}
\label{eqn:action:Laplacian}
\square_{p,q}(\varphi\otimes\psi)=(\square^{S}_{p',q'}\varphi)\otimes\psi+\varphi\otimes(\square^{T}_{p'',q''}\psi)
\end{equation}
for all $\varphi\in A^{p',q'}(S)$ and $\psi\in A^{p'',q''}(T)$ with $p=p'+p'$, $q=q'+q''$,
it follows from \eqref{eqn:action:Laplacian} and 
Lemmas~\ref{lemma:spectrum:Laplacian:2-elementary:K3} and \ref{lemma:spectrum:Laplacian:elliptic:curve} that
\begin{equation}
\label{eqn:spectrum:Laplacian:Borcea:Voisin:(0,0):1}
\begin{aligned}
\sigma(\square_{0,0}|_{A^{0,0}(T)^{+}\hat{\otimes}A^{0,0}(S)^{+}})
&=
\{0\}
\amalg
(\sigma(\square^{T,+})\setminus\{0\})
\amalg
(\sigma(\square^{S,+})\setminus\{0\})
\amalg
\\
&\quad
\left\{
\lambda_{i}^{+}(S)+\nu_{m,n}(T);\,
i\in I,\,(m,n)\in({\bf Z}^{2})^{*}/\pm1
\right\},
\end{aligned}
\end{equation}
\begin{equation}
\label{eqn:spectrum:Laplacian:Borcea:Voisin:(0,0):2}
\sigma(\square_{0,0}|_{A^{0,0}(T)^{-}\hat{\otimes}A^{0,0}(S)^{-}})
=
\left\{
\lambda_{j}^{-}(S)+\nu_{m,n}(T);\,
j\in J,\,(m,n)\in({\bf Z}^{2})^{*}/\pm1
\right\}.
\end{equation}
The result follows from \eqref{eqn:decomposition:spectrum:Laplacian:(0,0)}, \eqref{eqn:spectrum:Laplacian:Borcea:Voisin:(0,0):1},
\eqref{eqn:spectrum:Laplacian:Borcea:Voisin:(0,0):2}.
\end{pf}

\begin{lemma}
\label{lemma:spectral:zeta:Borcea:Voisin:(1,0)}
The following equality of meromorphic functions on ${\bf C}$ holds:
$$
\zeta_{1,0}(s)=\zeta^{T,+}(s)+\zeta^{S,+}(s)+2\zeta^{S,-}(s)+3\mu^{+}(s)+3\mu^{-}(s).
$$
\end{lemma}

\begin{pf}
Since the decomposition
$$
\begin{aligned}
A^{1,0}(X)^{+}
&=
[A^{1,0}(T)^{+}\hat{\otimes}A^{0,0}(S)^{+}]
\oplus
[A^{1,0}(T)^{-}\hat{\otimes}A^{0,0}(S)^{-}]
\\
&\quad
\oplus
[A^{0,0}(T)^{+}\hat{\otimes}A^{1,0}(S)^{+}]
\oplus
[A^{0,0}(T)^{-}\hat{\otimes}A^{1,0}(S)^{-}]
\end{aligned}
$$
is orthogonal and is preserved by $\square_{1,0}$, we get 
\begin{equation}
\label{eqn:decomposition:spectrum:Laplacian:(1,0)}
\begin{aligned}
\sigma(\square_{1,0})
&=
\sigma(\square_{1,0}|_{A^{1,0}(T)^{+}\hat{\otimes}A^{0,0}(S)^{+}})
\amalg
\sigma(\square_{1,0}|_{A^{1,0}(T)^{-}\hat{\otimes}A^{0,0}(S)^{-}})
\\
&\quad
\amalg
\sigma(\square_{1,0}|_{A^{0,0}(T)^{+}\hat{\otimes}A^{1,0}(S)^{+}})
\amalg
\sigma(\square_{1,0}|_{A^{0,0}(T)^{-}\hat{\otimes}A^{1,0}(S)^{-}}).
\end{aligned}
\end{equation}
By \eqref{eqn:action:Laplacian} and Lemmas~\ref{lemma:spectrum:Laplacian:2-elementary:K3} and \ref{lemma:spectrum:Laplacian:elliptic:curve}, we get
\begin{equation}
\label{eqn:spectrum:Laplacian:Borcea:Voisin:(1,0):1}
\begin{aligned}
\,&
\sigma(\square_{1,0}|_{A^{1,0}(T)^{+}\hat{\otimes}A^{0,0}(S)^{+}})
\\
&=
\sigma(\square^{T,+})
\amalg
\left\{
\lambda_{i}^{+}(S)+\nu_{m,n}(T);\,
i\in I,\,(m,n)\in({\bf Z}^{2})^{*}/\pm1
\right\},
\end{aligned}
\end{equation}
\begin{equation}
\label{eqn:spectrum:Laplacian:Borcea:Voisin:(1,0):2}
\begin{aligned}
\,&
\sigma(\square_{1,0}|_{A^{1,0}(T)^{-}\hat{\otimes}A^{0,0}(S)^{-}})
\\
&=
\sigma(\square^{S,-})
\amalg
\left\{
\lambda_{j}^{-}(S)+\nu_{m,n}(T);\,
j\in J,\,(m,n)\in({\bf Z}^{2})^{*}/\pm1
\right\},
\end{aligned}
\end{equation}
\begin{equation}
\label{eqn:spectrum:Laplacian:Borcea:Voisin:(1,0):3}
\begin{aligned}
\,&
\sigma(\square_{1,0}|_{A^{0,0}(T)^{+}\hat{\otimes}A^{1,0}(S)^{+}})
\\
&=
\sigma(\square_{1,0}^{S,+})
\amalg
\left\{
\lambda_{i}^{+}(S)+\nu_{m,n}(T);\,
i\in I,\,(m,n)\in({\bf Z}^{2})^{*}/\pm1
\right\}
\\
&\qquad
\amalg
\left\{
\lambda_{j}^{-}(S)+\nu_{m,n}(T);\,
j\in J,\,(m,n)\in({\bf Z}^{2})^{*}/\pm1
\right\}
\\
&=
\sigma(\square^{S,+})
\amalg
\sigma(\square^{S,-})
\amalg
\left\{
\lambda_{i}^{+}(S)+\nu_{m,n}(T);\,
i\in I,\,(m,n)\in({\bf Z}^{2})^{*}/\pm1
\right\}
\\
&\qquad
\amalg
\left\{
\lambda_{j}^{-}(S)+\nu_{m,n}(T);\,
j\in J,\,(m,n)\in({\bf Z}^{2})^{*}/\pm1
\right\},
\end{aligned}
\end{equation}
\begin{equation}
\label{eqn:spectrum:Laplacian:Borcea:Voisin:(1,0):4}
\begin{aligned}
\sigma(\square_{1,0}|_{A^{0,0}(T)^{-}\hat{\otimes}A^{1,0}(S)^{-}})
&=
\left\{
\lambda_{i}^{+}(S)+\nu_{m,n}(T);\,
i\in I,\,(m,n)\in({\bf Z}^{2})^{*}/\pm1
\right\}
\\
&\amalg
\left\{
\lambda_{j}^{-}(S)+\nu_{m,n}(T);\,
j\in J,\,(m,n)\in({\bf Z}^{2})^{*}/\pm1
\right\}.
\end{aligned}
\end{equation}
The result follows from \eqref{eqn:decomposition:spectrum:Laplacian:(1,0)}, \eqref{eqn:spectrum:Laplacian:Borcea:Voisin:(1,0):1},
\eqref{eqn:spectrum:Laplacian:Borcea:Voisin:(1,0):2}, \eqref{eqn:spectrum:Laplacian:Borcea:Voisin:(1,0):3},
\eqref{eqn:spectrum:Laplacian:Borcea:Voisin:(1,0):4}.
\end{pf}

\begin{lemma}
\label{lemma:spectral:zeta:Borcea:Voisin:(1,1)}
The following equality of meromorphic functions on ${\bf C}$ holds:
$$
\zeta_{1,1}(s)=21\zeta^{T,+}(s)+5\zeta^{S,+}(s)+4\zeta^{S,-}(s)+9\mu^{+}(s)+9\mu^{-}(s).
$$
\end{lemma}

\begin{pf}
Since the decomposition
$$
\begin{aligned}
A^{1,1}(X)^{+}
&=
[A^{1,1}(T)^{+}\hat{\otimes}A^{0,0}(S)^{+}]
\oplus
[A^{1,1}(T)^{-}\hat{\otimes}A^{0,0}(S)^{-}]
\\
&\quad
\oplus
[A^{1,0}(T)^{+}\hat{\otimes}A^{0,1}(S)^{+}]
\oplus
[A^{1,0}(T)^{-}\hat{\otimes}A^{0,1}(S)^{-}]
\\
&\quad
\oplus
[A^{0,1}(T)^{+}\hat{\otimes}A^{1,0}(S)^{+}]
\oplus
[A^{0,1}(T)^{-}\hat{\otimes}A^{1,0}(S)^{-}]
\\
&\quad
\oplus
[A^{0,0}(T)^{+}\hat{\otimes}A^{1,1}(S)^{+}]
\oplus
[A^{0,0}(T)^{-}\hat{\otimes}A^{1,1}(S)^{-}]
\end{aligned}
$$
is orthogonal and is preserved by $\square_{1,1}$, we get 
\begin{equation}
\label{eqn:decomposition:spectrum:Laplacian:(1,1)}
\begin{aligned}
\sigma(\square_{1,1})
&=
\sigma(\square_{1,1}|_{A^{1,1}(T)^{+}\hat{\otimes}A^{0,0}(S)^{+}})
\amalg
\sigma(\square_{1,1}|_{A^{1,1}(T)^{-}\hat{\otimes}A^{0,0}(S)^{-}})
\\
&\quad
\amalg
\sigma(\square_{1,1}|_{A^{1,0}(T)^{+}\hat{\otimes}A^{0,1}(S)^{+}})
\amalg
\sigma(\square_{1,1}|_{A^{1,0}(T)^{-}\hat{\otimes}A^{0,1}(S)^{-}})
\\
&\quad
\amalg
\sigma(\square_{1,1}|_{A^{0,1}(T)^{+}\hat{\otimes}A^{1,0}(S)^{+}})
\amalg
\sigma(\square_{1,1}|_{A^{0,1}(T)^{-}\hat{\otimes}A^{1,0}(S)^{-}})
\\
&\quad
\amalg
\sigma(\square_{1,1}|_{A^{0,0}(T)^{+}\hat{\otimes}A^{1,1}(S)^{+}})
\amalg
\sigma(\square_{1,1}|_{A^{0,0}(T)^{-}\hat{\otimes}A^{1,1}(S)^{-}}).
\end{aligned}
\end{equation}
By \eqref{eqn:action:Laplacian} and Lemmas~\ref{lemma:spectrum:Laplacian:2-elementary:K3} and \ref{lemma:spectrum:Laplacian:elliptic:curve}, we get
\begin{equation}
\label{eqn:spectrum:Laplacian:Borcea:Voisin:(1,1):1}
\begin{aligned}
\sigma(\square_{1,1}|_{A^{1,1}(T)^{+}\hat{\otimes}A^{0,0}(S)^{+}})
&=
\{0\}
\amalg
(\sigma(\square^{T,+})\setminus\{0\})
\amalg
(\sigma(\square^{S,+})\setminus\{0\})
\amalg
\\
&\quad
\left\{
\lambda_{i}^{+}(S)+\nu_{m,n}(T);\,
i\in I,\,(m,n)\in({\bf Z}^{2})^{*}/\pm1
\right\},
\end{aligned}
\end{equation}
\begin{equation}
\label{eqn:spectrum:Laplacian:Borcea:Voisin:(1,1):2}
\sigma(\square_{1,1}|_{A^{1,1}(T)^{-}\hat{\otimes}A^{0,0}(S)^{-}})
=
\left\{
\lambda_{j}^{-}(S)+\nu_{m,n}(T);\,
j\in J,\,(m,n)\in({\bf Z}^{2})^{*}/\pm1
\right\},
\end{equation}
\begin{equation}
\label{eqn:spectrum:Laplacian:Borcea:Voisin:(1,1):3}
\begin{aligned}
\sigma(\square_{1,1}|_{A^{1,0}(T)^{+}\hat{\otimes}A^{0,1}(S)^{+}})
&=
\left\{
\lambda_{i}^{+}(S)+\nu_{m,n}(T);\,
i\in I,\,(m,n)\in({\bf Z}^{2})^{*}/\pm1
\right\}
\\
&\quad
\amalg
\left\{
\lambda_{j}^{-}(S)+\nu_{m,n}(T);\,
j\in J,\,(m,n)\in({\bf Z}^{2})^{*}/\pm1
\right\},
\end{aligned}
\end{equation}
\begin{equation}
\label{eqn:spectrum:Laplacian:Borcea:Voisin:(1,1):4}
\begin{aligned}
\,&
\sigma(\square_{1,1}|_{A^{1,0}(T)^{-}\hat{\otimes}A^{0,1}(S)^{-}})
\\
&=
\sigma(\square^{S,-}_{0,1})
\amalg
\left\{
\lambda_{i}^{+}(S)+\nu_{m,n}(T);\,
i\in I,\,(m,n)\in({\bf Z}^{2})^{*}/\pm1
\right\}
\\
&\qquad
\amalg
\left\{
\lambda_{j}^{-}(S)+\nu_{m,n}(T);\,
j\in J,\,(m,n)\in({\bf Z}^{2})^{*}/\pm1
\right\}
\\
&=
(\sigma(\square^{S,+})\setminus\{0\})
\amalg
\sigma(\square^{S,-})
\amalg
\left\{
\lambda_{i}^{+}(S)+\nu_{m,n}(T);\,
i\in I,\,(m,n)\in({\bf Z}^{2})^{*}/\pm1
\right\}
\\
&\qquad
\amalg
\left\{
\lambda_{j}^{-}(S)+\nu_{m,n}(T);\,
j\in J,\,(m,n)\in({\bf Z}^{2})^{*}/\pm1
\right\},
\end{aligned}
\end{equation}
\begin{equation}
\label{eqn:spectrum:Laplacian:Borcea:Voisin:(1,1):5}
\begin{aligned}
\sigma(\square_{1,1}|_{A^{0,1}(T)^{+}\hat{\otimes}A^{1,0}(S)^{+}})
&=
\left\{
\lambda_{i}^{+}(S)+\nu_{m,n}(T);\,
i\in I,\,(m,n)\in({\bf Z}^{2})^{*}/\pm1
\right\}
\\
&\amalg
\left\{
\lambda_{j}^{-}(S)+\nu_{m,n}(T);\,
j\in J,\,(m,n)\in({\bf Z}^{2})^{*}/\pm1
\right\}.
\end{aligned}
\end{equation}
\begin{equation}
\label{eqn:spectrum:Laplacian:Borcea:Voisin:(1,1):6}
\begin{aligned}
\,&
\sigma(\square_{1,0}|_{A^{0,1}(T)^{-}\hat{\otimes}A^{1,0}(S)^{-}})
\\
&=
\sigma(\square^{S,-}_{1,0})
\amalg
\left\{
\lambda_{i}^{+}(S)+\nu_{m,n}(T);\,
i\in I,\,(m,n)\in({\bf Z}^{2})^{*}/\pm1
\right\}
\\
&\qquad
\amalg
\left\{
\lambda_{j}^{-}(S)+\nu_{m,n}(T);\,
j\in J,\,(m,n)\in({\bf Z}^{2})^{*}/\pm1
\right\}
\\
&=
(\sigma(\square^{S,+})\setminus\{0\})
\amalg
\sigma(\square^{S,-})
\amalg
\left\{
\lambda_{i}^{+}(S)+\nu_{m,n}(T);\,
i\in I,\,(m,n)\in({\bf Z}^{2})^{*}/\pm1
\right\}
\\
&\qquad
\amalg
\left\{
\lambda_{j}^{-}(S)+\nu_{m,n}(T);\,
j\in J,\,(m,n)\in({\bf Z}^{2})^{*}/\pm1
\right\},
\end{aligned}
\end{equation}
\begin{equation}
\label{eqn:spectrum:Laplacian:Borcea:Voisin:(1,1):7}
\begin{aligned}
\sigma(\square_{1,1}|_{A^{0,0}(T)^{+}\hat{\otimes}A^{1,1}(S)^{+}})
&=
h^{1,1}(S)^{+}\cdot\sigma(\square^{T,+}_{0,0})
\amalg
\sigma(\square^{S,+}_{1,1})
\\
&\qquad
\amalg
2\cdot\left\{
\lambda_{i}^{+}(S)+\nu_{m,n}(T);\,
i\in I,\,(m,n)\in({\bf Z}^{2})^{*}/\pm1
\right\}
\\
&\qquad
\amalg
2\cdot\left\{
\lambda_{j}^{-}(S)+\nu_{m,n}(T);\,
j\in J,\,(m,n)\in({\bf Z}^{2})^{*}/\pm1
\right\}
\\
&=
h^{1,1}(S)^{+}\cdot\sigma(\square^{T,+}_{0,0})
\amalg
2\cdot(\sigma(\square^{S,+})\setminus\{0\})
\amalg
2\cdot\sigma(\square^{S,-})
\\
&\qquad
\amalg
2\cdot
\left\{
\lambda_{i}^{+}(S)+\nu_{m,n}(T);\,
i\in I,\,(m,n)\in({\bf Z}^{2})^{*}/\pm1
\right\}
\\
&\qquad
\amalg
2\cdot
\left\{
\lambda_{j}^{-}(S)+\nu_{m,n}(T);\,
j\in J,\,(m,n)\in({\bf Z}^{2})^{*}/\pm1
\right\},
\end{aligned}
\end{equation}
\begin{equation}
\label{eqn:spectrum:Laplacian:Borcea:Voisin:(1,1):8}
\begin{aligned}
\sigma(\square_{1,1}|_{A^{0,0}(T)^{-}\hat{\otimes}A^{1,1}(S)^{-}})
&=
h^{1,1}(S)^{-}\cdot\sigma(\square^{T,-}_{0,0})
\\
&\qquad
\amalg
2\cdot\left\{
\lambda_{i}^{+}(S)+\nu_{m,n}(T);\,
i\in I,\,(m,n)\in({\bf Z}^{2})^{*}/\pm1
\right\}
\\
&\qquad
\amalg
2\cdot\left\{
\lambda_{j}^{-}(S)+\nu_{m,n}(T);\,
j\in J,\,(m,n)\in({\bf Z}^{2})^{*}/\pm1
\right\}.
\end{aligned}
\end{equation}
Since $h^{1,1}(S)^{+}+h^{1,1}(S)^{-}=h^{1,1}(S)=20$,
the result follows from \eqref{eqn:decomposition:spectrum:Laplacian:(1,1)}, \eqref{eqn:spectrum:Laplacian:Borcea:Voisin:(1,1):1},
\eqref{eqn:spectrum:Laplacian:Borcea:Voisin:(1,1):2}, \eqref{eqn:spectrum:Laplacian:Borcea:Voisin:(1,1):3},
\eqref{eqn:spectrum:Laplacian:Borcea:Voisin:(1,1):4}, \eqref{eqn:spectrum:Laplacian:Borcea:Voisin:(1,1):5},
\eqref{eqn:spectrum:Laplacian:Borcea:Voisin:(1,1):6}, \eqref{eqn:spectrum:Laplacian:Borcea:Voisin:(1,1):7},
\eqref{eqn:spectrum:Laplacian:Borcea:Voisin:(1,1):8}.
\end{pf}

\subsection
{Proof of Theorem~\ref{thm:Harvey:Moore}}
\label{sect:8.3}
\par
Set $X:=X_{(S,\theta,T)}$.
Since ${\rm Vol}(S,\gamma_{S})=(2\pi)^{-2}\int_{S}\gamma_{S}^{2}/2!$ and 
${\rm Vol}(T,\gamma_{T})=(2\pi)^{-1}\int_{T}\gamma_{T}=(2\pi)^{-1}$ by definition, we get
\begin{equation}
\label{eqn:volume:Borcea:Voisin:H:0}
{\rm Vol}(X,\gamma)={\rm Vol}(S,\gamma_{S}){\rm Vol}(T,\gamma_{T})/2={\rm Vol}(S,\gamma_{S})/4\pi.
\end{equation}
\par
Following \cite[Lemma 13.4]{FLY08}, we compute the covolume of the lattice 
$$
H^{2}(X,{\bf Z})=H^{2}(S\times T,{\bf Z})^{+}=H^{2}(S,{\bf Z})^{+}\oplus H^{2}(T,{\bf Z})\subset H^{2}(S,{\bf R})^{+}\oplus H^{2}(T,{\bf R})
$$
with respect to $\gamma$. Set $r=r(M)$ and $l=l(M)$. Then $|A_{M}|=|M^{\lor}/M|=2^{l}$. 
Since ${\rm sign}(M)=(1,r-1)$, there is an integral basis $\{{\bf e}_{1},\ldots,{\bf e}_{r}\}$ of $H^{2}(S,{\bf Z})^{+}\cong M$ with
\begin{center}
$\det\left(\langle{\bf e}_{\alpha},{\bf e}_{\beta}\rangle\right)_{1\leq\alpha,\beta\leq r}=(-1)^{r-1}2^{l}=(-1)^{\rho}|A_{M}|$,
\end{center}
where $\langle\cdot,\cdot\rangle$ is the intersection pairing on $H^{2}(S,{\bf Z})$. 
Then $\{{\bf e}_{1},\ldots,{\bf e}_{r},[\gamma_{T}]\}$ is an integral basis of $H^{2}(S\times T,{\bf Z})^{+}$.
Let $\langle\cdot,\cdot\rangle_{L^{2},\gamma}$ be the $L^{2}$-metric on $H^{2}(S\times T,{\bf R})$ with respect to $\gamma$.
By the same computations as in \cite[p.253]{FLY08}, we get
\begin{equation}
\label{eqn:L2:metric:Borcea:Voisin:H:(1,1):No1}
\langle{\bf e}_{\alpha},{\bf e}_{\beta}\rangle_{L^{2},\gamma}
=
(2\pi)^{-3}
\left(
\frac{\langle{\bf e}_{\alpha}[\gamma_{S}]\rangle\langle{\bf e}_{\beta},[\gamma_{S}]\rangle}{\langle[\gamma_{S}],[\gamma_{S}]\rangle}
-
\frac{1}{2}\langle{\bf e}_{\alpha},{\bf e}_{\beta}\rangle
\right),
\end{equation}
\begin{equation}
\label{eqn:L2:metric:Borcea:Voisin:H:(1,1):No2}
\langle{\bf e}_{\alpha},[\gamma_{T}]\rangle_{L^{2},\gamma}=0,
\qquad
\langle[\gamma_{T}],[\gamma_{T}]\rangle_{L^{2},\gamma}
=
\frac{\langle[\gamma_{S}],[\gamma_{S}]\rangle}{(2\pi)^{3}\cdot4}
=
\frac{{\rm Vol}(S,\gamma_{S})}{4\pi}.
\end{equation}
By \eqref{eqn:L2:metric:Borcea:Voisin:H:(1,1):No1}, \eqref{eqn:L2:metric:Borcea:Voisin:H:(1,1):No2} and \cite[p.254, l.1-7]{FLY08}, we get
\begin{equation}
\label{eqn:covolume:Borcea:Voisin:H:(1,1)}
\begin{aligned}
\,&
{\rm Vol}_{L^{2}}(H^{2}(S\times T,{\bf Z})^{+},\gamma)
=
\det
\begin{pmatrix}
\langle{\bf e}_{\alpha},{\bf e}_{\beta}\rangle_{L^{2},\gamma}&\langle{\bf e}_{\alpha},[\gamma_{T}]\rangle_{L^{2},\gamma}
\\
\langle{\bf e}_{\alpha},[\gamma_{T}]\rangle_{L^{2},\gamma}&\langle[\gamma_{T}],[\gamma_{T}]\rangle_{L^{2},\gamma}
\end{pmatrix}
\\
&=
(-2\pi)^{-3r}\frac{{\rm Vol}(S,\gamma_{S})}{4\pi}
\det\left(
\frac{1}{2}\langle{\bf e}_{\alpha},{\bf e}_{\beta}\rangle
-
\frac{\langle{\bf e}_{\alpha}[\gamma_{S}]\rangle\langle{\bf e}_{\beta},[\gamma_{S}]\rangle}{\langle[\gamma_{S}],[\gamma_{S}]\rangle}
\right)
\\
&=
(-1)^{-3r}(2\pi)^{-3r-1}\frac{{\rm Vol}(S,\gamma_{S})}{2}
\cdot(-1)2^{-r}\cdot
\det\left(
\langle{\bf e}_{\alpha},{\bf e}_{\beta}\rangle
\right)
\\
&=
(2\pi)^{-3\rho+2}2^{-\rho}\cdot|A_{M}|{\rm Vol}(S,\gamma_{S}).
\end{aligned}
\end{equation}
\par
Since ${\rm Sing}\,X=(S\times T)^{\theta\times(-1_{T})}=S^{\theta}\times T[2]$ is the $4$ copies of $S^{\theta}$ and 
since $\gamma|_{{\rm Sing}\,X}=\gamma_{S}|_{S^{\theta}}$ for every component of $(S\times T)^{\theta\times(-1_{T})}$, we get
\begin{equation}
\label{eqn:analytic:torsion:singular:locus}
{\rm Vol}\left(
{\rm Sing}\,X,\gamma|_{{\rm Sing}\,X}
\right)\,
\tau\left(
{\rm Sing}\,X,\gamma|_{{\rm Sing}\,X}
\right)
=
\left\{
{\rm Vol}(S^{\theta},\gamma_{S}|_{S^{\theta}})\,\tau(S^{\theta},\gamma_{S}|_{S^{\theta}})
\right\}^{4}.
\end{equation}
\par
By Lemmas~\ref{lemma:spectral:zeta:Borcea:Voisin:(0,0)}, \ref{lemma:spectral:zeta:Borcea:Voisin:(1,0)} and \ref{lemma:spectral:zeta:Borcea:Voisin:(1,1)},
we get 
\begin{equation}
\label{eqn:spectral:zeta:BCOV:Borcea:Voisin}
\sum_{p,q\geq0}(-1)^{p+q}pq\,\zeta_{p,q}(s)
=
24\zeta^{T,+}(s)+8\left\{\zeta^{S,+}(s)-\zeta^{S,-}(s)\right\}.
\end{equation}
Since $\zeta^{T,+}(s)=\zeta^{T,-}(s)=\zeta^{T}_{0,0}(s)/2$, we get by \eqref{eqn:spectral:zeta:BCOV:Borcea:Voisin}
\begin{equation}
\label{eqn:BCOV:torsion:Borcea:Voisin}
\begin{aligned}
T_{\rm BCOV}(X,\gamma)
&=
\exp\left(-12\left.\frac{d}{ds}\right|_{s=0}\zeta^{T}_{0,0}(s)-8\left.\frac{d}{ds}\right|_{s=0}\left\{\zeta^{S,+}(s)-\zeta^{S,-}(s)\right\}\right)
\\
&=
\tau_{{\bf Z}_{2}}(S,\gamma_{S})(\theta)^{-4}\,\tau(T,\gamma_{T})^{-12},
\end{aligned}
\end{equation}
where the second equality follows from \cite[Lemma 4.3]{Yoshikawa08}. 
Since $\gamma$, $\gamma_{S}$, $\gamma_{T}$ are Ricci-flat,
it follows from \eqref{eqn:BCOV:orb:Ricci:flat}, \eqref{eqn:tau:M:Ricci:flat}, \eqref{eqn:torsion:flat:elliptic:curve},
\eqref{eqn:volume:Borcea:Voisin:H:0},
\eqref{eqn:covolume:Borcea:Voisin:H:(1,1)},
\eqref{eqn:analytic:torsion:singular:locus},
\eqref{eqn:BCOV:torsion:Borcea:Voisin}
that
$$
\label{eqn:BCOV:invariant:Borcea:Voisin}
\begin{aligned}
\,&
\tau_{\rm BCOV}^{\rm orb}(X)
=
T_{\rm BCOV}(X,\gamma)
{\rm Vol}(X,\gamma)^{-3+\frac{\chi(S\times T)+3\chi(S^{\theta}\times T[2])}{24}}
{\rm Vol}_{L^{2}}
\left(
H^{2}(S\times T,{\bf Z})^{+},\gamma
\right)^{-1}
\\
&\qquad\qquad\quad\times
\tau
\left(
{\rm Sing}\,X,\gamma|_{{\rm Sing}\,X}
\right)^{-1}
{\rm Vol}
\left(
{\rm Sing}\,X,\gamma|_{{\rm Sing}\,X}
\right)^{-1}
\\
&=
\tau_{{\bf Z}_{2}}(S,\gamma_{S})(\theta)^{-4}\,\tau(T,\gamma_{T})^{-12}
\left\{
{\rm Vol}(S^{\theta},\gamma_{S}|_{S^{\theta}})\,\tau(S^{\theta},\gamma_{S}|_{S^{\theta}})
\right\}^{-4}
\\
&\quad\times
\{{\rm Vol}(S,\gamma_{S})/4\pi\}^{-3+\frac{\chi(S\times T)+3\chi(S^{\theta}\times T[2])}{24}}
\cdot
(2\pi)^{3\rho-2}2^{\rho}\,|A_{M}|^{-1}{\rm Vol}(S,\gamma_{S})^{-1}
\\
&=
2^{14}(2\pi)^{2\rho}|A_{M}|^{-1}\,\tau_{M}(S,\theta)^{-4}\tau_{\rm ell}(T)^{-12}.
\end{aligned}
$$
Here we used $\chi(S\times T)=0$, $\chi(S^{\theta}\times T[2])=8(r-10)$ and $\tau_{\rm ell}(T)=(2\pi)^{-1}\tau(T,\gamma_{T})$
to get the last equality. 
\qed
\newline
\par
After Corollary~\ref{cor:BCOV=orb:BCOV}, we propose the following

\begin{conjecture}
\label{conj:comparison:BCOV:invariants}
For a Calabi-Yau orbifold $Y$ and its crepant resolution $\widetilde{Y}$,
$$
\tau_{\rm BCOV}^{\rm orb}(Y)=C(\widetilde{Y}/Y)\,\tau_{\rm BCOV}(\widetilde{Y}),
$$
where $C(\widetilde{Y}/Y)$ is a constant depending only on the topological type of the crepant resolution $\widetilde{Y}\to Y$.
\end{conjecture}

This conjecture may be viewed an analogue of Roan's theorem \cite{Roan96} on string theoretic Euler characteristic 
for global Calabi-Yau threefolds: $\chi^{\rm orb}(X,G)=\chi(\widetilde{X/G})$.
Conjecture~\ref{conj:comparison:BCOV:invariants} is closely related to another conjecture \cite[Conjecture 4.17]{FLY08}
claiming the birational invariance of BCOV invariants. 
We refer the reader to \cite{MaillotRoessler12} for a current progress on this latter conjecture.


\end{document}